\documentclass[12pt,a4paper,titlepage]{report}
\usepackage[mathscr,mathcal]{eucal}
\usepackage{enumerate}
\usepackage{makeidx}
\makeindex
\usepackage{amsfonts}
\usepackage{amssymb}
\usepackage{amsmath}
\usepackage{xypic}
\newtheorem{theorem}{Theorem}[section]
\newtheorem{definition}[theorem]{Definition}
\newtheorem{ypoth}[theorem]{Assumptions}
\newtheorem{ypothesh}[theorem]{Assumption}
\newtheorem{lemma}[theorem]{Lemma}

\newtheorem{corollary}[theorem]{Corollary}
\newtheorem{proposition}[theorem]{Proposition}
\newtheorem{thedef}[theorem]{Theorem-Definition}

\newcommand{\el}{\lambda}
\newcommand{\eb}{\beta}
\newtheorem{px}[theorem]{Example}

\newcommand{\ie}{\emph{i.e., }}

\begin{document}

$$\textrm{\textbf{Maria Chlouveraki}}$$\\ 
\\

$$\textrm{\textbf{\Large BLOCKS AND FAMILIES}}$$
  $$\textrm{\textbf{\Large FOR CYCLOTOMIC HECKE ALGEBRAS}}$$

$ $\\
$ $\\
$ $\\
$ $\\
$ $\\
$ $\\
$ $\\

  % \vspace{10 cm}
   %$$\textrm{Springer}$$

\chapter*{Preface}

This book contains a thorough study of symmetric algebras, covering topics such as block theory, representation theory and Clifford theory. It can also serve as an introduction to the Hecke algebras of complex reflection groups.
Its aim is the study of the blocks and the determination of the families of characters of the cyclotomic Hecke algebras associated to complex reflection groups.

I would like to thank my thesis' advisor, Michel Brou\'{e}, for his advice. These Springer Lecture Notes were, after all, his idea. I am grateful to  Jean Michel for his help with the realization and presentation of the programming part. I would also like to thank Gunter Malle for his suggestion that I generalize my results on Hecke algebras, which led to the notion of  ``essential algebras''.
Finally, I express my thanks to C\'edric Bonnaf\'e, Meinolf Geck, Nicolas Jacon, Rapha\"el Rouquier and Jacques Th\'evenaz for their useful comments.

\chapter*{Introduction}

The finite groups of matrices with coefficients in $\mathbb{Q}$ generated by reflections, known as \emph{Weyl groups}, are a fundamental building block in the classification of semisimple complex Lie algebras and Lie groups, as well as semisimple algebraic groups over arbitrary algebraically closed fields. They are also a foundation for many other significant mathematical theories, including braid groups and Hecke algebras.

The Weyl groups are particular cases of \emph{complex reflection groups}, finite groups of matrices with coefficients in a finite abelian extension of $\mathbb{Q}$ generated by ``pseudo-reflections'' (elements whose vector space of fixed points is a hyperplane) --- if the coefficients belong to $\mathbb{R}$, then these are the finite Coxeter groups.

The work of Lusztig on the irreducible characters of reductive
groups over finite fields has displayed the important role of the
``families of characters'' of the Weyl groups concerned. However, only
recently was it realized that it would be of great interest to
generalize the notion of families of characters to the complex
reflection groups, or more precisely to the cyclotomic Hecke algebras
associated to complex reflection groups.

On the one hand, the complex reflection groups and their associated
cyclotomic Hecke algebras appear naturally in the classification
of the ``cyclotomic Harish-Chandra series'' of the characters of the
finite reductive groups, generalizing the role of the Weyl group and
its traditional Hecke algebra in the principal series. Since the
families of characters of the Weyl group play an essential role in the
definition of the families of unipotent characters of the
corresponding finite reductive group (cf.~\cite{Lu1}), we can hope that
the families of characters of the cyclotomic Hecke algebras play a key
role in the organization of families of unipotent characters more
generally.

On the other hand, for some complex reflection groups (non-Coxeter)
$W$, some data have been gathered which seem to indicate that behind
the group $W$, there exists another mysterious object - the
\emph{Spets} (cf.~\cite{BMM2}, \cite{Ma3}) - that could play the role
of the ``series of finite reductive groups with Weyl group $W$''. In
some cases, one can define the unipotent characters of the Spets,
which are controlled by the ``spetsial'' Hecke algebra of $W$, a
generalization of the classical Hecke algebra of the Weyl groups.

The main obstacle for this generalization is the lack of
Kazhdan-Lusztig bases for the non-Coxeter complex reflection groups.
However, more 
recent results of Gyoja ($\cite{Gy}$) and Rouquier
($\cite{Rou}$) have made possible the definition of a substitute for
families of characters which can be applied to all complex
reflection groups. Gyoja has shown (case by case) that the partition
into ``$p$-blocks'' of the Iwahori-Hecke algebra of a Weyl group $W$
coincides with the partition into families, when $p$ is the unique
bad prime number for $W$.
Later, Rouquier proved that the families
of characters of a Weyl group $W$ are exactly the blocks of
characters of the Iwahori-Hecke algebra of $W$ over a suitable
coefficient ring, the ``Rouquier ring''. 

Brou\'e, Malle and Rouquier (cf.~\cite{BMR}) have shown that we can define the \emph{generic Hecke algebra} $\mathcal{H}(W)$ associated to a complex reflection group $W$ as a quotient of the group algebra of the braid group of $W$. The algebra $\mathcal{H}(W)$ is  an algebra over a Laurent polynomial ring in a set of indeterminates $\textbf{v}=(v_i)_{0 \leq i \leq m}$ whose cardinality $m$ depends on the group $W$. A \emph{cyclotomic Hecke algebra} is an algebra obtained from $\mathcal{H}(W)$ via a specialization of the form $v_i \mapsto y^{n_i}$, where $y$ is an indeterminate and $n_i \in \mathbb{Z}$ for all $i=0,1,\ldots,m$. The blocks of a cyclotomic Hecke algebra over the Rouquier ring are the \emph{Rouquier blocks} of the cyclotomic Hecke algebra. Thus, the Rouquier blocks generalize the notion of the families of characters   to all complex reflection groups.

In \cite{BK}, Brou{\'e} and Kim presented an algorithm for the determination of the Rouquier blocks of the cyclotomic Hecke algebras of the groups $G(d,1,r)$ and $G(e,e,r)$. Later,
 Kim (cf.~\cite{Kim}) generalized this algorithm to include all the groups of the infinite series $G(de,e,r)$. However,  it was realized recently that their algorithm does not work in general, unless $d$ is a power of a prime number. Moreover, the Rouquier
blocks of the ``spetsial'' cyclotomic Hecke algebra of many
exceptional irreducible complex reflection groups have been calculated by Malle
and Rouquier in \cite{MaRo}.  In this book, we correct and complete the determination of the Rouquier blocks for all cyclotomic Hecke algebras and all complex reflection groups.

The key in our study of the Rouquier blocks has been the proof of the fact that they have the property of ``semi-continuity'' (the name is due to C. Bonnaf\'e). Every complex reflection group $W$ determines some numerical data, which in turn determine the ``essential'' hyperplanes for $W$. To each essential hyperplane $H$, we can associate a partition $\mathcal{B}(H)$ of the set of irreducible characters of $W$ into blocks.
Given a cyclotomic specialization  $v_i \mapsto y^{n_i}$, the Rouquier blocks of the
corresponding cyclotomic Hecke algebra depend only on which essential hyperplanes the integers $n_i$ belong to. In particular, they are unions of the blocks associated with the essential hyperplanes to which the integers $n_i$ belong and they are minimal with respect to that property.

The property of semi-continuity also appears in works on Kazhdan-Lusztig cells  (cf.~\cite{BGIL}, \cite{Bon}, \cite{Jer}) and on Cherednik algebras (cf.~\cite{GoMa}). The common appearance of this as yet unexplained phenomenon implies a connection between these structures and the Rouquier blocks, for which the reason is not yet apparent, but promises to be fruitful when explored thoroughly.
In particular, due to the known relation between Kazhdan-Lusztig cells and families of characters for
Coxeter groups, this could be an indication of the existence of Kazhdan-Lusztig bases for the (non-Coxeter) complex reflection groups.

Another indication of this fact comes from the determination of the Rouquier blocks of the cyclotomic Hecke algebras of all complex reflection groups, obtained in the last chapter of this book with the use of the theory of  ``essential hyperplanes''. In the case of the Weyl groups and their usual Hecke algebra,  Lusztig attaches to every irreducible character two integers, denoted by $a$ and $A$, and shows (cf.~\cite{Lu2}, $3.3$ and $3.4$) that they are constant on the families. 
In an analogous way, we can define integers $a$ and $A$ attached to every irreducible character of a cyclotomic Hecke algebra of a complex reflection group. Using the classification of the Rouquier blocks, it has been proved that the integers $a$ and $A$ are constant on the ``families of characters'' of the cyclotomic Hecke algebras of all complex reflection groups (see end of Chapter $4$).

The first chapter of this book is dedicated to commutative algebra. The need for the results presented in this chapter (some of them are well-known, but others are completely new) arises form the fact that when we are working on Hecke algebras of complex reflection groups, we work over integrally closed rings, which are not necessarily unique factorization domains.

In the second chapter, we present some
classical results of block theory and representation theory of
symmetric algebras. We see that the \emph{Schur elements} associated to the irreducible characters of a symmetric algebra play a crucial role in the determination of its blocks.
Moreover, we generalize the results known as  ``Clifford theory''  (cf.~, for example, \cite{Da}), which determine 
the blocks of certain subalgebras of symmetric algebras, to the case of ``twisted symmetric algebras of finite groups''. Finally, we give a new criterion for a symmetric algebra to be split semisimple.

In the third chapter, we introduce the notion of ``essential algebras''. These are symmetric algebras whose Schur elements have a specific form: they are products of irreducible polynomials evaluated on monomials. We obtain many results on the block theory of these algebras, which we later apply to the Hecke algebras, after we prove that they are essential in Chapter $4$. In particular, we have our first encounter with the phenomenon of semi-continuity (see Theorem $\ref{main theorem}$).

It is in the fourth chapter that we define for the first time the braid group, the generic Hecke algebra and the cyclotomic Hecke algebras associated to a complex reflection group. We show that the generic Hecke algebra of a complex reflection group is essential, by proving that its Schur elements are of the required form. Applying the results of Chapter $3$, we obtain that the Rouquier blocks (\ie the families of characters) of the cyclotomic Hecke algebras have the property of semi-continuity and only depend on some ``essential'' hyperplanes for the group, which are determined by the generic Hecke algebra.

In the fifth and final chapter of this book,  we present the algorithms and the results of the determination of the families of characters for all irreducible complex reflection groups. The use of Clifford theory is essential, since it allows us to restrict ourselves to the study of  only certain cases of complex reflection groups. The computations were made with the use of the GAP package CHEVIE (cf.~\cite{chevie}) for the exceptional irreducible complex reflection groups, whereas combinatorial methods were applied to the groups of the infinite series. In particular, we show that the families of characters for the latter can be obtained from the families of characters of the Weyl groups of type $B$, already determined by Lusztig.

\tableofcontents

\chapter{On Commutative Algebra}

The first chapter contains some known facts and some new results on Commutative Algebra. The former are presented here without their proofs (with the exception of Theorem $\ref{similitude}$) for the convenience of the reader. They are going to be crucial in the proofs of the results of Chapters $3$ and $4$.

In the first section of this chapter, we define the localization of a ring and give some main properties. The second section is dedicated on integrally closed rings. We study particular cases of integrally closed rings, such as valuation rings, discrete valuation rings and Krull rings. We use their properties in order to obtain results on Laurent polynomial rings over integrally closed rings. We state briefly some results on the completions of rings in section $1.3$. In the fourth section, we introduce the notion of \emph{morphisms associated with monomials}. They are morphisms which allow us to pass from a Laurent polynomial ring $A$ in $m+1$ indeterminates to a Laurent polynomial ring $B$ in $m$ indeterminates, while sending a specific monomial to $1$. Moreover, we prove (Proposition $\ref{surjective is associated}$) that every surjective morphism from $A$ to $B$ which sends each indeterminate to a monomial  is associated with a monomial. We call \emph{adapted morphisms} the compositions of morphisms associated with monomials. They play a powerful role in the proof of the main results of Chapters $3$ and $4$. Finally, in the last section of the first chapter, we give a criterion (Theorem $\ref{second irreducible}$) for a polynomial to be irreducible in a Laurent polynomial ring with coefficients in a field.

Throughout this chapter, all rings are assumed to be commutative
with $1$. Moreover, if $R$ is a ring and $x_0,x_1,\ldots,x_m$ is a set
of indeterminates on $R$, then we denote by $R[x_0^{\pm 1},x_1^{\pm
1},\ldots,x_m^{\pm 1}]$ the Laurent polynomial ring in $m+1$
indeterminates
$R[x_0,{x_0}^{-1},x_1,{x_1}^{-1},\ldots,x_m,{x_m}^{-1}]$.

\section{Localizations}

\begin{definition}\label{multiplicatively closed set}
Let $R$ be a commutative ring with $1$. We say that a subset $S$ of
$R$ is a multiplicatively closed set if $0 \notin S$, $1 \in S$ and
every finite product of elements of $S$ belongs to $S$.
\end{definition}

In the set $R \times S$, we introduce an equivalence relation such
that $(r,s)$ is equivalent to $(r',s')$ if and only if there exists
$t \in S$ such that $t(s'r-sr')=0$. We denote the equivalence class
of $(r,s)$ by $r/s$. The set of equivalence classes becomes a ring
under the operations such that the sum and the product of $r/s$ and
$r'/s'$ are given by $(s'r+sr')/ss'$ and $rr'/ss'$ respectively. We
denote this ring by $S^{-1}R$ and we call it the \emph{localization} \index{localization of a ring}
of $R$ at $S$. If $S$ contains no zero divisors of $R$, then any
element $r$ of $R$ can be identified with the element $r/1$ of $S^{-1}R$ and we
can regard the latter as an $R$-algebra.\\
\\
\begin{remarks}\
\emph{\begin{itemize}
  \item If $S$ is the set of all non-zero divisors of $R$, then
  $S^{-1}R$ is called \emph{the total quotient ring} of $R$. If, moreover,
  $R$ is an integral domain, the total quotient ring of $R$ is \emph{the
  field of fractions} of $R$.
  \item If $R$ is Notherian, then $S^{-1}R$ is Noetherian.
  \item If $\mathfrak{p}$ is a prime ideal of $R$, then the set
  $S:=R-\mathfrak{p}$ is a multiplicatively closed subset of $R$. Then the
  ring $S^{-1}R$ is simply denoted by $R_\mathfrak{p}$.
\end{itemize}}
\end{remarks}

The proofs for the following well known results concerning
localizations can be found in \cite{Bou2}.

\begin{proposition}\label{Bourbaki1}
Let $A$ and $B$ be two rings with multiplicative sets $S$ and $T$
respectively and $f$ an homomorphism from $A$ to $B$ such that
$f(S)$ is contained in $T$. There exists a unique homomorphism $f'$
from $S^{-1}A$ to $T^{-1}B$ such that $f'(a/1)=f(a)/1$ for every $a
\in A$. Let us suppose now that $T$ is contained in the
multiplicatively closed set of $B$ generated by $f(S)$. If $f$ is
surjective (resp. injective), then $f'$ is also surjective (resp.
injective).
\end{proposition}

\begin{corollary}\label{inclusion in localizations}
Let $A$ and $B$ be two rings with multiplicative sets $S$ and $T$
respectively such that $A \subseteq B$ and $S \subseteq T$. Then
$S^{-1}A \subseteq T^{-1}B$.
\end{corollary}

\begin{proposition}\label{different multiplicative set}
Let $A$ be a ring and $S,T$ two multiplicative sets of $A$ such that
$S \subseteq T$. We have $S^{-1}A = T^{-1}A$ if and only if every
prime ideal of $R$ that meets $T$ meets $S$.
\end{proposition}

The following proposition and its corollary give us information
about the ideals of the localization of a ring $R$ at a
multiplicatively closed subset $S$ of $R$.

\begin{proposition}\label{prime ideals of a localization}
Let $R$ be a ring and let $S$ be a multiplicatively closed subset of
$R$. Then
\begin{enumerate}
  \item Every ideal $\mathfrak{b}'$ of $S^{-1}R$ is of the form
  $S^{-1}\mathfrak{b}$ for some ideal $\mathfrak{b}$ of $R$.
  \item Let $\mathfrak{b}$ be an ideal of $R$ and let $f$ be the
  canonical surjection $R \twoheadrightarrow R/\mathfrak{b}$. Then $f(S)$
  is a multiplicatively closed subset of $R/\mathfrak{b}$ and the
  homomorphism from $S^{-1}R$ to $(f(S))^{-1}(R/\mathfrak{b})$
  canonically associated with $f$ is surjective with kernel
  $\mathfrak{b}'=S^{-1}\mathfrak{b}$. By passing to quotients, an
  isomorphism between $(S^{-1}R)/\mathfrak{b}'$ and
  $(f(S))^{-1}(R/\mathfrak{b})$ is defined.
  \item The application $\mathfrak{b}' \mapsto \mathfrak{b}$,
  restricted to the set of maximal (resp. prime) ideals of
  $S^{-1}R$, is an isomorphism (for the relation of inclusion)
  between this set and the set of maximal (resp. prime) ideals of
  $R$ that do not meet $S$.
  \item If $\mathfrak{q}'$ is a prime ideal of $S^{-1}R$ and
  $\mathfrak{q}$ is the prime ideal of $R$ such that $\mathfrak{q}'=S^{-1}\mathfrak{q}$
  (we have $\mathfrak{q} \cap S = \emptyset$), then there exists an
  isomorphism from $R_\mathfrak{q}$ to $(S^{-1}R)_{\mathfrak{q}'}$
  which sends $r/s$ to $(r/1)/(s/1)$ for $r \in R$,
  $s \in R-\mathfrak{q}$.
  \end{enumerate}
\end{proposition}

\begin{corollary}\label{an}
Let $R$ be a ring, $\mathfrak{p}$ a prime ideal of $R$ and
$S:=R-\mathfrak{p}$. For every ideal $\mathfrak{b}$ of $R$ which
does not meet $S$, let $\mathfrak{b}':= \mathfrak{b}R_\mathfrak{p}$.
Assume that $\mathfrak{b}' \neq R_\mathfrak{p}$. Then
\begin{enumerate}
  \item Let $f$ be the canonical surjection $R\twoheadrightarrow
  R/\mathfrak{b}$. The ring homomorphism from $R_\mathfrak{p}$ to
  $(R/\mathfrak{b})_{\mathfrak{p}/\mathfrak{b}}$ canonically
  associated with $f$ is surjective and its kernel is
  $\mathfrak{b}'$. Thus it defines, by passing to
  quotients, a canonical isomorphism between
  $R_\mathfrak{p}/\mathfrak{b}'$ and
  $(R/\mathfrak{b})_{\mathfrak{p}/\mathfrak{b}}$.
  \item The application $\mathfrak{b}' \mapsto \mathfrak{b}$,
  restricted to the set of prime ideals of
  $R_\mathfrak{p}$, is an isomorphism (for the relation of inclusion)
  between this set and the set of prime ideals of
  $R$ contained in $\mathfrak{p}$ (thus do not meet $S$). Therefore,
  $\mathfrak{p}R_\mathfrak{p}$ is the only maximal ideal of $R_\mathfrak{p}$.
  \item If now $\mathfrak{b}'$ is
  a prime ideal of $R_\mathfrak{p}$, then there exists an isomorphism
  from $R_\mathfrak{b}$ to $(R_\mathfrak{p})_{\mathfrak{b}'}$ which
  sends $r/s$ to $(r/1)/(s/1)$ for $r \in R$, $s \in
  R-\mathfrak{b}$.
\end{enumerate}
\end{corollary}

The notion of localization can also be extended to the modules over
the ring $R$.

\begin{definition}\label{localization of a module}
Let $R$ be a ring and $S$ a multiplicatively closed set of $R$. If
$M$ is an $R$-module, then we call localization of $M$ \index{localization of a module} at $S$ and
denote by $S^{-1}M$ the $S^{-1}R$-module $M \otimes_R S^{-1}R$.
\end{definition}

\section{Integrally closed rings}

\begin{thedef}\label{integral element}
Let $R$ be a ring, $A$ an $R$-algebra and $a$ an element of $A$. The
following properties are equivalent:
\begin{enumerate}[(i)]
  \item The element $a$ is a root of a monic polynomial with coefficients in
  $R$.
  \item The subalgebra $R[a]$ of $A$ is an $R$-module of
  finite type.
  \item There exists a faithful $R[a]$-module which is an
  $R$-module of finite type.
\end{enumerate}
If $a \in A$ verifies the conditions above, we say that it is
integral \index{integral element} over $R$.
\end{thedef}

\begin{definition}\label{integral closure}
Let $R$ be a ring and $A$ an $R$-algebra. The set of all elements of
$A$ that are integral over $R$ is an $R$-subalgebra of $A$
containing $R$; it is called the integral closure  \index{integral closure} of $R$ in $A$. We
say that $R$ is integrally closed in $A$, if $R$ is an integral
domain and if it coincides with its integral closure in $A$. If now
$R$ is an integral domain and $F$ is its field of fractions, then
the integral closure of $R$ in $F$ is named simply the integral
closure of $R$, and if $R$ is integrally closed in $F$, then $R$ is
said to be integrally closed.  \index{integrally closed}
\end{definition}

The following proposition (\cite{Bou5}, \S 1, Prop.13) implies that
transfer theorem holds for integrally closed rings (corollary
$\ref{integrally closed polynomial ring}$).

\begin{proposition}\label{integral closure of a polynomial ring}
  If $R$ is an integral domain, let us denote by $\bar{R}$ the
  integral closure of $R$. Let $x_0, \ldots,x_m$
  be a set of indeterminates over $R$. Then the integral closure of
  $R[x_0,\ldots,x_m]$ is
  $\bar{R}[x_0,\ldots,x_m]$.
\end{proposition}

\begin{corollary}\label{integrally closed polynomial ring}
Let $R$ be an integral domain. Then $R[x_0,\ldots,x_m]$ is
integrally closed if and only if $R$ is integrally closed.
\end{corollary}

\begin{corollary} If $K$ is a field, then every polynomial ring
  $K[x_0,\ldots,x_m]$ is integrally closed.
\end{corollary}

The next proposition (\cite{Bou5}, \S 1, Prop.16) along with its
corollaries treats the
integral closures of localizations of rings.

\begin{proposition}\label{integral closure of a localization}
Let $R$ be a ring, $A$ an $R$-algebra, $\bar{R}$ the integral
closure of $R$ in $A$ and $S$ a multiplicatively closed subset of
$R$ which contains no zero divisors. Then the integral closure of
$S^{-1}R$ in $S^{-1}A$ is $S^{-1}\bar{R}$.
\end{proposition}

\begin{corollary}Let $R$ be an integral domain, $\bar{R}$ the integral
closure of $R$ and $S$ a multiplicatively closed subset of $R$. Then
the integral closure of $S^{-1}R$ is $S^{-1}\bar{R}$.
\end{corollary}

\begin{corollary}\label{integrally closed localization}
If $R$ is an integrally closed domain and $S$ is a multiplicatively
closed subset of $R$, then $S^{-1}R$ is also integrally closed.
\end{corollary}

\begin{px} \small{\emph{Let $K$ be a finite field extension of $\mathbb{Q}$ and
$\mathbb{Z}_K$ the integral closure of $\mathbb{Z}$ in $K$. Obviously, the ring
$\mathbb{Z}_K$ is integrally
closed. Let $x_0,x_1,\ldots,x_m$ be indeterminates. Then the ring
$\mathbb{Z}_K[x_0^{\pm},x_1^{\pm},\ldots,x_m^{\pm}]$ is also
integrally closed.}}
\end{px}

\subsection{Lifting prime ideals}

\begin{definition}\label{lying over}
Let $R, R'$ be two rings and let $h:R \rightarrow R'$ be a ring
homomorphism. We say that a prime ideal $\mathfrak{a}'$ of $R'$ lies
over a prime ideal $\mathfrak{a}$ of $R$, if
$\mathfrak{a}=h^{-1}(\mathfrak{a}')$.
\end{definition}

The next result is \cite{Bou5}, \S 2, Proposition 2.

\begin{proposition}\label{primes lying over}
Let $h:R \rightarrow R'$ be a ring homomorphism such that $R'$ is
integral over $R$. Let $\mathfrak{p}$ be a prime ideal of $R$,
$S:=R-\mathfrak{p}$ and $(\mathfrak{p}_i')_{i \in I}$ the family of
all the prime ideals of $R'$ lying over $\mathfrak{p}$. If
$S'=\bigcap_{i \in I}(R'-\mathfrak{p}_i')$, then
$S^{-1}R'=S'^{-1}R'$.
\end{proposition}

The following corollary of Proposition $\ref{primes lying over}$ deals with a case we will encounter in a
following chapter, where there exists a unique prime ideal lying
over the prime ideal $\mathfrak{p}$ of $R$. In combination with
Proposition $\ref{integral closure of a localization}$, Proposition
$\ref{primes lying over}$ implies that

\begin{corollary}\label{one prime lying over}
Let $R$ be an integral domain, $A$ an $R$-algebra, $\bar{R}$ the
integral closure of $R$ in $A$. Let $\mathfrak{p}$ be a prime ideal
of $R$ and $S:=R-\mathfrak{p}$. If there exists a unique prime ideal
$\bar{\mathfrak{p}}$ of $\bar{R}$ lying over $\mathfrak{p}$, then
the integral closure of $R_\mathfrak{p}$ in $S^{-1}A$ is
$\bar{R}_{\bar{\mathfrak{p}}}$.
\end{corollary}
\subsection{Valuations}

\begin{definition}\label{valuation}
Let $R$ be a ring and $\Gamma$ a totally ordered abelian group. We
call valuation  \index{valuation} of $R$ with values in $\Gamma$ every application $v:R
\rightarrow \Gamma \cup \{\infty\}$ which satisfies the following
properties:
\begin{description}
\item[(V1)] $v(xy)=v(x)+v(y)$ for $x \in R, y \in R$.
\item[(V2)] $v(x+y) \geq \mathrm{inf}(v(x),v(y))$ for $x \in R, y \in R$.
\item[(V3)] $v(1)=0$ and $v(0)=\infty$.
\end{description}
\end{definition}
In particular, if $v(x) \neq v(y)$, property (V2) gives
$v(x+y)=\textrm{inf}(v(x),v(y))$ for $x \in R, y \in R$. Moreover,
from property (V1), we have that if $z \in R$ with $z^n=1$ for some
integer $n \geq 1$, then $nv(z)=v(z^n)=v(1)=0$ and thus $v(z)=0$.
Consequently, $v(-x)=v(-1)+v(x)=v(x)$ for all $x \in R$.\\

Now let $F$ be a field and let $v:F \rightarrow \Gamma$ be a
valuation of $F$. The set $A$ of $a \in F$ such that $v(a) \geq 0$
is a local subring of $F$. Its maximal ideal $\mathfrak{m}(A)$ is
the set of $a \in A$ such that $v(a) > 0$. For all $a \in F-A$, we
have $a^{-1} \in \mathfrak{m}(A)$. The ring $A$ is called
\emph{the ring of the valuation} $v$ on $F$.\\

We will now introduce the notion of a valuation ring. For more
information about valuation rings and their properties, see
\cite{Bou6}. Some of them will also be discussed in Chapter 2,
Section 2.4.

\begin{definition}\label{valuation ring}
Let $R$ be an integral domain contained in a field $F$. Then $R$ is
a valuation ring  \index{valuation ring} if for all non-zero element $x \in F$, we have $x
\in R$ or $x^{-1} \in R$. Consequently, $F$ is the field of
fractions of $R$.
\end{definition}

If $R$ is a valuation ring, then it has the following properties:
\begin{itemize}
  \item It is an integrally closed local ring.
  \item The set of the principal ideals of $R$ is totally ordered by the relation of inclusion.
  \item The set of the ideals of $R$ is totally ordered by the relation of inclusion.
\end{itemize}

Let $R$ be a valuation ring and $F$ its field of fractions. Let us
denote by $R^\times$ the set of units of $R$. Then the set
$\Gamma_R:=F^\times/R^\times$ is an abelian group, totally ordered
by the relation of inclusion of the corresponding principal ideals.
If we denote by $v_R$ the canonical homomorphism of $F^\times$ onto
$\Gamma_R$ and set $v_R(0)=\infty$, then $v_R$ is a valuation of $F$
whose ring is $R$.

The following proposition gives a characterization of integrally
closed rings in terms of valuation rings (\cite{Bou6}, \S 1, Thm.
3).

\begin{proposition}\label{intersection of valuation rings}
Let $R$ be a subring of a field $F$. The integral closure $\bar{R}$
of $R$ in $F$ is the intersection of all valuation rings in $F$
which contain $R$. Consequently, an integral domain $R$ is
integrally closed if and only if it is the intersection of a family
of valuation rings contained in its field of fractions.
\end{proposition}

This characterization helped us to prove the following result about
integrally closed rings.

\begin{proposition}\label{askhsh 12}
Let $R$ be an integrally closed ring and $f(x)=\sum_ia_ix^i$,
$g(x)=\sum_jb_jx^j$ be two polynomials in $R[x]$. If there exists an
element $c \in R$ such that all the coefficients of $f(x)g(x)$
belong to $cR$, then all the products $a_ib_j$ belong to $cR$.
\end{proposition}
\begin{apod}{Due to Proposition $\ref{intersection of valuation rings}$,
it is enough to prove
the result in the case where $R$ is a valuation ring.

From now on, let $R$ be a valuation ring. Let $v$ be a valuation of
the field of fractions of $R$ such that the ring of valuation of $v$
is $R$. Let $\kappa:=\textrm{inf}_i(v(a_i))$ and
$\lambda:=\textrm{inf}_j(v(b_j))$. Then
$\kappa+\lambda=\textrm{inf}_{i,j}(v(a_ib_j))$. We will show that
$\kappa + \lambda \geq v(c)$ and thus $c$ divides all the products
$a_ib_j$. Argue by contradiction and assume that $\kappa + \lambda <
v(c)$. Let $a_{i_1},a_{i_2},\ldots,a_{i_r}$ with
$i_1<i_2<\ldots<i_r$ be all the elements among the $a_i$ with
valuation equal to $\kappa$. Respectively, let
$b_{j_1},b_{j_2},\ldots,b_{j_s}$ with $j_1<j_2<\ldots<j_s$ be all
the elements among the $b_j$ with valuation equal to $\lambda$. We
have that $i_1+j_1<i_m+j_n$, $\forall (m,n) \neq (1,1)$. Therefore,
the coefficient $c_{i_1+j_1}$ of $x^{i_1+j_1}$ in $f(x)g(x)$ is of
the form $(a_{i_1}b_{j_1}+\sum(\textrm{terms with valuation} >
\kappa+\lambda))$ and since $v(a_{i_1}b_{j_1}) \neq
v(\sum(\textrm{terms with valuation} > \kappa+\lambda))$, we obtain
that $$v(c_{i_1+j_1})=\textrm{inf}(v(a_{i_1}b_{j_1}),
v(\sum(\textrm{terms with valuation} >
\kappa+\lambda)))=\kappa+\lambda.$$ However, since all the
coefficients of $f(x)g(x)$ are divisible by $c$, we have that
$v(c_{i_1+j_1}) \geq v(c) > \kappa + \lambda$, which is a
contradiction.}
\end{apod}

The Propositions $\ref{quotient}$
and $\ref{my second lemma one variable}$ derive from the one above. We will
make use of the results in corollaries $\ref{porisma porismatos}$
and $\ref{my second lemma}$ in Chapter 3.

\begin{proposition}\label{quotient}
Let $R$ be an integrally closed domain and let $F$ be its field of
fractions. Let $\mathfrak{p}$ be a prime ideal of $R$. Then
$$(R[x])_{\mathfrak{p}R[x]} \cap F[x] = R_\mathfrak{p}[x].$$
\end{proposition}
\begin{apod}{The inclusion $R_\mathfrak{p}[x] \subseteq (R[x])_{\mathfrak{p}R[x]} \cap
F[x]$ is obvious. Now, let $f(x)$ be an element of $F[x]$. Then
$f(x)$ can be written in the form $r(x)/\xi$, where $r(x) \in R[x]$
and $\xi \in R$. If, moreover, $f(x)$ belongs to
$(R[x])_{\mathfrak{p}R[x]}$, then there exist $s(x), t(x) \in R[x]$
with $t(x) \notin \mathfrak{p}R[x]$ such that $f(x) = s(x)/t(x)$.
Thus we have $$f(x)=\frac{r(x)}{\xi}= \frac{s(x)}{t(x)}.$$ All the
coefficients of the product $r(x)t(x)$ belong to $\xi R$. Due to
Proposition \ref{askhsh 12}, if $r(x)=\sum_ia_ix^i$ and
$t(x)=\sum_jb_jx^j$, then all the products $a_ib_j$ belong to $\xi
R$. Since $t(x)\notin \mathfrak{p}R[x]$, there exists $j_0$ such
that $b_{j_0} \notin \mathfrak{p}$ and $a_ib_{j_0} \in \xi R,
\forall i$. Consequently, $b_{j_0}f(x)=b_{j_0}(r(x)/\xi) \in R[x]$
and hence all the coefficients of $f(x)$ belong to
$R_\mathfrak{p}$.}
\end{apod}

\begin{corollary}\label{porisma porismatos}
Let $R$ be an integrally closed domain and let $F$ be its field of
fractions. Let $\mathfrak{p}$ be a prime ideal of $R$. Then
\begin{enumerate}
  \item $(R[x,x^{-1}])_{\mathfrak{p}R[x,x^{-1}]} \cap F[x,x^{-1}] =
R_\mathfrak{p}[x,x^{-1}].$
  \item $(R[x_0,\ldots,x_m])_{\mathfrak{p}R[x_0,\ldots,x_m]}
\cap F[x_0,\ldots,x_m] = R_\mathfrak{p}[x_0,\ldots,x_m].$
  \item $(R[x_0^{\pm1},\ldots,x_m^{\pm1}])_{\mathfrak{p}R[x_0^{\pm1},\ldots,x_m^{\pm1}]}
\cap F[x_0^{\pm1},\ldots,x_m^{\pm1}] =
R_\mathfrak{p}[x_0^{\pm1},\ldots,x_m^{\pm1}]$.
\end{enumerate}
\end{corollary}

\begin{proposition}\label{my second lemma one variable}
Let $R$ be an integrally closed domain and let $F$ be its field of
fractions. Let $r(x)$ and $s(x)$ be two elements of $R[x]$ such that
$s(x)$ divides $r(x)$ in $F[x]$. If one of the coefficients of
$s(x)$ is a unit in $R$, then $s(x)$ divides $r(x)$ in $R[x]$.
\end{proposition}
\begin{apod}{Since $s(x)$ divides $r(x)$ in $F[x]$, there exists an
element of the form $t(x)/\xi$ with $t(x) \in R[x]$ and $\xi \in R$
such that
$$r(x) = \frac{s(x)t(x)}{\xi}.$$
All the coefficients of the product $s(x)t(x)$ belong to $\xi R$.
Due to Proposition \ref{askhsh 12}, if $s(x)=\sum_ia_ix^i$ and
$t(x)=\sum_jb_jx^j$, then all the products $a_ib_j$ belong to $\xi
R$. By assumption, there exists $i_0$ such that $a_{i_0}$ is a unit
in $R$ and $a_{i_0}b_j \in \xi R, \forall j$. Consequently, $b_j \in
\xi R, \forall j$ and thus $t(x)/\xi \in R[x]$.}
\end{apod}

\begin{corollary}\label{my second lemma}
Let $R$ be an integrally closed domain and let $F$ be its field of
fractions. Let $r, s$ be two elements of
$R[x_0^{\pm1},\ldots,x_m^{\pm1}]$ such that $s$ divides $r$ in
$F[x_0^{\pm1},\ldots,x_m^{\pm1}]$. If one of the coefficients of $s$
is a unit in $R$, then $s$ divides $r$ in
$R[x_0^{\pm1},\ldots,x_m^{\pm1}]$.
\end{corollary}
\subsection{Discrete valuation rings and Krull rings}

\begin{definition}\label{discrete valuation}
Let $F$ be a field, $\Gamma$ a totally ordered abelian group and $v$
a valuation of $F$ with values in $\Gamma$. We say that the
valuation $v$ is discrete, if $\Gamma$ is isomorphic to
$\mathbb{Z}$.
\end{definition}

\begin{thedef}\label{dvr}
An integral domain $R$ is a discrete valuation ring, \index{discrete valuation ring} if it satisfies
one of the following equivalent conditions:
\begin{enumerate}[(i)]
  \item $R$ is the ring of a discrete valuation.
  \item $R$ is a local Dedekind ring.
  \item $R$ is a local principal ideal domain.
  \item $R$ is a Noetherian valuation ring.
\end{enumerate}
\end{thedef}

By Proposition $\ref{intersection of valuation rings}$, integrally
closed rings are intersections of valuation rings. Krull rings are
essentially intersections of discrete valuation rings.

\begin{definition}\label{Krull ring}
An integral domain $R$ is a Krull ring,  \index{Krull ring} if there exists a family of
valuations $(v_i)_{i \in I}$ of the field of fractions $F$ of $R$
with the following properties:
\begin{description}
  \item[(K1)] The valuations $(v_i)_{i \in I}$ are discrete.
  \item[(K2)] The intersection of the rings of $(v_i)_{i \in I}$ is
  $R$.
  \item[(K3)] For all $x\in F^\times$,there exists a finite number of
  $i \in I$ such that $v_i(x) \neq 0$.
\end{description}
\end{definition}

The proofs of the following results and more information about Krull
rings can be found in \cite{Bou7}, \S 1.

\begin{theorem}\label{Krull-dvr}
Let $R$ be an integral domain and let $\mathrm{Spec}_1(R)$ be the
set of its prime ideals of height $1$. Then $R$ is a Krull ring if and
only if the following properties are satisfied:
\begin{enumerate}
  \item For all $\mathfrak{p} \in \mathrm{Spec}_1(R)$, $R_\mathfrak{p}$ is a
  discrete valuation ring.
  \item $R$ is the intersection of $R_\mathfrak{p}$ for all
  $\mathfrak{p} \in \mathrm{Spec}_1(R)$.
  \item For all $r \neq 0$ in $R$, there exists a finite
  number of ideals $\mathfrak{p} \in \mathrm{Spec}_1(R)$ such that $r \in
  \mathfrak{p}$.
\end{enumerate}
\end{theorem}

Transfer theorem holds also for Krull rings:

\begin{proposition}\label{prime ideals of height 1}
Let $R$ be a Krull ring, $F$ the field of fractions of $R$ and $x$
an indeterminate. Then $R[x]$ is also a Krull ring. Moreover, its
prime ideals of height $1$ are:
\begin{itemize}
  \item the prime ideals of the form $\mathfrak{p}R[x]$, where
   $\mathfrak{p}$ is a prime ideal of height $1$ of $R$,
  \item the prime ideals of the form $\mathfrak{m} \cap R[x]$, where
   $\mathfrak{m}$ is a prime ideal of $F[x]$.
\end{itemize}
\end{proposition}

The following proposition provides us with a simple characterization
of Krull rings, when they are Noetherian.

\begin{proposition}\label{case of Krull}
Let $R$ be a Noetherian ring. Then $R$ is a Krull ring if and only
if it is integrally closed.
\end{proposition}

\begin{px} \small{\emph{Let $K$ be a finite field extension of $\mathbb{Q}$ and
$\mathbb{Z}_K$ the integral closure of $\mathbb{Z}$ in $K$. The ring
$\mathbb{Z}_K$ is a Dedekind ring and thus Noetherian and integrally
closed. Let $x_0,x_1,\ldots,x_m$ be indeterminates. Then the ring
$\mathbb{Z}_K[x_0^{\pm},x_1^{\pm},\ldots,x_m^{\pm}]$ is also
Noetherian and integrally closed and thus a Krull ring.}}
\end{px}

\section{Completions}

For the proofs of all the following results concerning completions, the reader can
refer to \cite{Na}, Chapter II.

Let $I$ be an ideal of a commutative ring $R$ and let $M$ be an
$R$-module. We introduce a topology on $M$ such that the open sets
of $M$ are unions of an arbitrary number of sets of the form $m+I^nM$ 
$(m\in M)$. This topology is called the $I$\emph{-adic topology}  \index{I-adic topology} of
$M$.

\begin{theorem}\label{16.5}
If $M$ is a Noetherian $R$-module, then for any submodule $N$ of
$M$, the $I$-adic topology of $N$ coincides with the topology of $N$
as a subspace of $M$ with the $I$-adic topology.
\end{theorem}

From now on, we will concentrate on semi-local rings and in
particular, on Noetherian semi-local rings.

\begin{definition}\label{semilocal ring}
A ring $R$ is called semi-local,  \index{semi-local ring} if it has only a finite number of
maximal ideals. The Jacobson radical $\mathfrak{m}$ of $R$ is the
intersection of the maximal ideals of $R$.
\end{definition}

\begin{theorem}\label{17.6}
Assume that $R$ is a Noetherian semi-local ring with Jacobson
radical $\mathfrak{m}$ and let $\hat{R}$ be the completion of $R$
with respect to the $\mathfrak{m}$-adic topology. Then $\hat{R}$ is
also a Noetherian semi-local ring and we have $R \subseteq \hat{R}$.
\end{theorem}

\begin{theorem}\label{17.8}
Assume that $R$ is a Noetherian semi-local ring with Jacobson
radical $\mathfrak{m}$ and that $M$ is a finitely generated
$R$-module. Let $\hat{R}$ be the completion of $R$ with respect to
the $\mathfrak{m}$-adic topology. Endow $M$ with the
$\mathfrak{m}$-adic topology. Then $M \otimes_R \hat{R}$ is the
completion of $M$ with respect to that topology.
\end{theorem}

\begin{corollary}\label{17.9}
Let $\mathfrak{a}$ be an ideal of a Noetherian semi-local ring $R$
with Jacobson radical $\mathfrak{m}$. Let $\hat{R}$ be the
completion of $R$ with respect to the $\mathfrak{m}$-adic topology.
Then the completion of $\mathfrak{a}$ is $\mathfrak{a}\hat{R}$ and
$\mathfrak{a}\hat{R}$ is isomorphic to $\mathfrak{a} \otimes_R
\hat{R}$. Furthermore, $\mathfrak{a}\hat{R} \cap R = \mathfrak{a}$
and $\hat{R}/\mathfrak{a}\hat{R}$ is the completion of
$R/\mathfrak{a}$ with respect to the $\mathfrak{m}$-adic topology.
\end{corollary}

\begin{theorem}\label{18.4}
Assume that $R$ is a Noetherian semi-local ring with Jacobson
radical $\mathfrak{m}$ and let $\hat{R}$ be the completion of $R$
with respect to the $\mathfrak{m}$-adic topology. Then
\begin{enumerate}
  \item The total quotient ring $F$ of $R$ (localization of $R$ at
  the set of non-zero divisors) is naturally a subring of the total
  quotient ring $\hat{F}$ of $\hat{R}$.
  \item For any ideal $\mathfrak{a}$ of $R$, $\mathfrak{a}\hat{R} \cap
  F =\mathfrak{a}$.
\end{enumerate}
In particular, $\hat{R} \cap F =R$.
\end{theorem}

\section{Morphisms associated with monomials and adapted morphisms}

We have the following elementary algebra result

\begin{thedef}\label{similitude}
Let $R$ be an integral domain and $M$ a free $R$-module of basis
$(e_i)_{0 \leq i \leq m}$.  Let $x=r_0e_0+r_1e_1+\ldots+r_me_m$ be a
non-zero element of $M$. We set $M^*:=\mathrm{Hom}_R(M,R)$ and
$M^*(x):= \{\varphi(x) \,|\, \varphi \in M^*\}$. Then the following
assertions are equivalent:
\begin{enumerate}[(i)]
  \item  $M^*(x)=R$.
  \item $\sum_{i=0}^m Rr_i =R$.
  \item There exists $\varphi \in M^*$ such that $\varphi(x)=1$.
  \item There exists an $R$-submodule $N$ of $M$ such that $M=Rx\oplus
  N$.
\end{enumerate}
If $x$ satisfies the conditions above, we say that $x$ is a primitive
element of $M$.
\end{thedef}
\begin{apod}{
\begin{description}
  \item[$(i) \Leftrightarrow (ii)$]
  Let $(e_i^*)_{0 \leq i \leq m}$ be the basis of $M^*$ dual to
  $(e_i)_{0 \leq i \leq  m}$.
  Then $M^*(x)$ is generated by $(e_i^*(x))_{0 \leq i \leq m}$
  and $e_i^*(x)=r_i$.
  \item[$(ii) \Rightarrow (iii)$]
  There exist $u_0,u_1,\ldots,u_m \in R$ such that
  $\sum_{i=0}^mu_ir_i=1$. If $\varphi:=\sum_{i=0}^mu_ie_i^*$, then
  $\varphi(x)=1$.
  \item[$(iii) \Rightarrow (iv)$]
  Let $N:=\mathrm{Ker}\varphi$. For all $y \in M$, we have
  $$y = \varphi(y)x + (y -\varphi(y)x).$$
  Since $\varphi(x)=1$, we have $y \in Rx + N$. Obviously, $Rx
  \cap N =\{0\}.$
  \item[$(iv) \Rightarrow (i)$]
  Since $M$ is free and $R$ is an integral domain, $x$ is
  torsion-free (otherwise there exists $r \in R, r\neq 0$ such that
  $\sum_{i=0}^m (rr_i)e_i=0$). Therefore, the map
  $$R \rightarrow M, r \mapsto rx$$
  is an isomorphism of $R$-modules. Its inverse is a linear form on
  $Rx$ which sends $x$ to $1$. Composing it with the map
  $$M \twoheadrightarrow M/N \rightarrow R,$$
  we obtain a linear form $\varphi:M\rightarrow R$ such that
  $\varphi(x)=1$. We have $1 \in M^*(x)$ and thus $M^*(x)=R$.}
\end{description}
\end{apod}

We will apply the above result to the $\mathbb{Z}$-module
$\mathbb{Z}^{m+1}$. Let us consider $(e_i)_{0 \leq i \leq m}$ the
standard basis of $\mathbb{Z}^{m+1}$ and let
$a:=a_0e_0+a_1e_1+\ldots+a_me_m$ be an element of $\mathbb{Z}^{m+1}$
such that $\mathrm{gcd}(a_i)=1$. Then, by Bezout's theorem, there
exist $u_0,u_1,\ldots,u_m \in \mathbb{Z}$ such that
$\sum_{i=0}^mu_ia_i=1$ and hence $\sum_{i=0}^m
\mathbb{Z}a_i=\mathbb{Z}$. By Theorem $\ref{similitude}$, there
exists a $\mathbb{Z}$-submodule $N_a$ of $\mathbb{Z}^{m+1}$ such
that $\mathbb{Z}^{m+1}=\mathbb{Z}a\oplus N_a$. In particular,
$N_a=\mathrm{Ker}\varphi$, where
$\varphi:=\sum_{i=0}^mu_ie_i^*$. We will denote by
$p_a:\mathbb{Z}^{m+1} \twoheadrightarrow N_a$ the projection of
$\mathbb{Z}^{m+1}$ onto $N_a$ such that
$\mathrm{Ker}p_a=\mathbb{Z}a$. We have a $\mathbb{Z}$-module
isomorphism
$i_a:N_a \tilde{\rightarrow} \mathbb{Z}^m$. Then
$f_a:=i_a \circ p_a$ is a surjective $\mathbb{Z}$-module
morphism $\mathbb{Z}^{m+1} \twoheadrightarrow \mathbb{Z}^m$
with $\mathrm{Ker}f_a=\mathbb{Z}a$.

Now let $R$ be an integral domain and let $x_0,x_1,\ldots,x_m$ be
$m+1$ indeterminates over $R$. Let $G$ be the abelian group
generated by all the monomials in $R[x_0^{\pm 1},x_1^{\pm
1},\ldots,x_m^{\pm 1}]$ with group operation the multiplication.
Then $G$ is isomorphic to the additive group $\mathbb{Z}^{m+1}$ by
the isomorphism defined as follows
$$\begin{array}{rccc}
   \theta_G: &G & \tilde{\rightarrow} & \mathbb{Z}^{m+1} \\
    &\prod_{i=0}^m x_i^{l_i} & \mapsto & (l_0,l_1,\ldots,l_m).
  \end{array}$$

\begin{lemma}\label{work on Z}
We have $R[x_0^{\pm 1},x_1^{\pm 1},\ldots,x_m^{\pm 1}]=R[G] \cong
R[\mathbb{Z}^{m+1}]$.
\end{lemma}

Respectively, if $y_1,\ldots,y_m$ are $m$ indeterminates over $R$
and $H$ is the group generated by all the monomials in $R[y_1^{\pm
1},\ldots,y_m^{\pm1}]$, then $H \cong \mathbb{Z}^m$ and $R[y_1^{\pm
1},\ldots,y_m^{\pm1}]=R[H] \cong R[\mathbb{Z}^m].$\\

The morphism $F_a:={\theta_H}^{-1} \circ f_a \circ \theta_G: G
\twoheadrightarrow H$ induces an $R$-algebra morphism
$$\begin{array}{cccc}
    \varphi_a: & R[G] & \rightarrow & R[H] \\
        & \sum_{g \in G}a_gg & \mapsto     & \sum_{g \in G}a_gF_a(g)
  \end{array}$$
Since $F_a$ is surjective, the morphism $\varphi_a$ is also surjective.
Moreover, $\mathrm{Ker}\varphi_a$ is generated (as an $R$-module) by
the set
$$< g-1 \,|\, g \in G \textrm{ such that } \theta_G(g) \in \mathbb{Z}a >.$$

From now on, let $A:=R[x_0^{\pm 1},x_1^{\pm 1},\ldots,x_m^{\pm 1}]$ and
$B:=R[y_1^{\pm 1},y_2^{\pm 1},\ldots,y_m^{\pm 1}]$. Let $M:=\prod_{i=0}^mx_i^{a_i}$
and set $\varphi_M:=\varphi_a$. Then the map $\varphi_M$ is an $R$-algebra morphism from $A$ to $B$ with the following properties:
\begin{enumerate}
\item $\varphi_M$ is surjective.
\item $\varphi_M(x_i)$ is a monomial in $B$ for all $i=0,1,\ldots,m$.
\item $\mathrm{Ker}\varphi_M=(M-1)A$. 
\end{enumerate}

\begin{definition}\label{associated morphism}
Let $M:=\prod_{i=0}^mx_i^{a_i}$ be a monomial in $A$ with
$\mathrm{gcd}(a_i)=1$. An $R$-algebra morphism $\varphi_M:A
\rightarrow B$ defined as above will be called associated with the
monomial $M$.  \index{morphism associated with a monomial}
\end{definition}

\begin{px}\label{x^5y^-3z^-2}
\small{\emph{Let $A:=R[X^{\pm 1},Y^{\pm 1},Z^{\pm 1}]$ and
$M:=X^5Y^{-3}Z^{-2}$. We have $a:=(5,-3,-2)$ and$$ (-1)\cdot
5+(-2)\cdot(-3)+0 \cdot(-2)=1.$$ Hence, with the notations above, we have
$$u_0=-1,u_1=-2,u_3=0.$$
The map $\varphi: \mathbb{Z}^3 \rightarrow \mathbb{Z}$ defined as
$$\varphi:=-e_0^*-2e_1^*.$$
has $\mathrm{Ker}\varphi=\{(x_0,x_1,x_2) \in \mathbb{Z}^3 \,|\, x_0
= -2x_1\} = \{(-2r,r,s) \,|\, r,s \in \mathbb{Z}\}=:N_a$.
\\
\\
By Theorem  $\ref{similitude}$, we have $\mathbb{Z}^3=\mathbb{Z}a
\oplus N_a$
and the projection $p_a:\mathbb{Z}^3\twoheadrightarrow N_a$ is the map
$$(y_0,y_1,y_2)=:y \mapsto
y-\varphi(y)a=(6y_0+10y_1,-3y_0-5y_1,-2y_0-4y_1+y_2).$$ The
$\mathbb{Z}$-module $N_a$ is obviously isomorphic to $\mathbb{Z}^2$
via
$$i_a:(-2r,r,s) \mapsto (r,s).$$
Composing the two previous maps, we obtain a well defined surjection
$$\begin{array}{ccc}
    \mathbb{Z}^3 & \twoheadrightarrow & \mathbb{Z}^2 \\
    (y_0,y_1,y_2) & \mapsto & (-3y_0-5y_1,-2y_0-4y_1+y_2).
  \end{array}$$
The above surjection induces (in the way described before) an
$R$-algebra epimorphism
$$\begin{array}{cccc}
    \varphi_M: & R[X^{\pm 1},Y^{\pm 1},Z^{\pm 1}] & \twoheadrightarrow & R[X^{\pm 1},Y^{\pm 1}] \\
     & X & \mapsto & X^{-3}Y^{-2} \\
     & Y & \mapsto & X^{-5}Y^{-4} \\
     & Z & \mapsto & Y
  \end{array}$$
By straightforward calculations, we can verify that
$\mathrm{Ker}\varphi_M=(M-1)A$.}}
\end{px}

\begin{proposition}\label{primeness of q}
Let $M:=\prod_{i=0}^mx_i^{a_i}$ be a monomial in $A$ such that
$\mathrm{gcd}(a_i)=1$. Then
\begin{enumerate}
  \item The ideal $(M-1)A$ is a prime ideal of $A$.
  \item If $\mathfrak{p}$ is a prime ideal of $R$, then the ideal
$\mathfrak{q}_M:=\mathfrak{p}A+(M-1)A$ is also prime in $A$.
\end{enumerate}
\end{proposition}
\begin{apod}{
\begin{enumerate}
\item Let $\varphi_M:A \rightarrow B$ be a morphism associated with
$M$. Then $\varphi_M$ is surjective and
$\mathrm{Ker}\varphi_M=(M-1)A$. By the $1^{\mathrm{st}}$ isomorphism theorem, we have
$A/(M-1)A \cong B$. Since $B$ is an integral domain, the
ideal generated by $(M-1)$ is prime in $A$.
\item Set $R':=R / \mathfrak{p}$,
$A':=R'[x_0^{\pm 1},x_1^{\pm 1},\ldots,x_m^{\pm 1}]$, $B':=R'[
y_1^{\pm 1},\ldots,y_m^{\pm 1}]$ . The rings $R'$, $A'$ and $B'$ are integral domains
and
$$A/\mathfrak{q}_M \cong A'/(M-1)A' \cong B'.$$
Hence, the ideal $\mathfrak{q}_M$ is
prime in $A$.}
\end{enumerate}
\end{apod}

The following assertions are now straightforward. Nevertheless, they
are stated for further reference.

\begin{proposition}\label{properties of phi}
Let $M:=\prod_{i=0}^mx_i^{a_i}$ be a monomial in $A$ with
$\mathrm{gcd}(a_i)=1$ and let $\varphi_M:A \rightarrow B$ be a
morphism associated with $M$. Let $\mathfrak{p}$ be a prime ideal of
$R$ and set $\mathfrak{q}_M:=\mathfrak{p}A+(M-1)A$. Then the
morphism $\varphi_M$ has the following properties:
\begin{enumerate}
  \item  If $f \in A$, then $\varphi_M(f) \in \mathfrak{p}B$ if and only if $f \in \mathfrak{q}_M$.
  Corollary $\ref{an}$\emph{(1)} implies that $$A_{\mathfrak{q}_M}/(M-1) A_{\mathfrak{q}_M} \cong
   B_{\mathfrak{p}B}.$$
  \item  If $N$ is a monomial in $A$, then $\varphi_M(N)=1$ if and only if
        there exists $k \in \mathbb{Z}$ such that $N=M^k$.
\end{enumerate}
\end{proposition}

\begin{corollary}\label{q intersection zk}
$\mathfrak{q}_M \cap R = \mathfrak{p}$.
\end{corollary}
\begin{apod}{Obviously $\mathfrak{p} \subseteq \mathfrak{q}_M \cap R$.
 Let $\xi \in R$ such that $\xi \in \mathfrak{q}_M$. If
$\varphi_M$ is a morphism associated with $M$, then, by Proposition
\ref{properties of phi}, $\varphi_M(\xi) \in \mathfrak{p}B$. But
$\varphi_M(\xi)=\xi$ and $\mathfrak{p}B \cap R = \mathfrak{p}$. Thus
$\xi \in \mathfrak{p}$.}
\end{apod}\\
\begin{remark}
\emph{If $m=0$ and we set $x:=x_0$, then $A:=R[x,x^{-1}]$ and
$B:=R$. The only monomials that we can associate a morphism $A
\rightarrow B$ with are $x$ and $x^{-1}$. This morphism is unique
and given by $x \mapsto 1$.}
\end{remark}\
\\

The following lemma, whose proof is straightforward when arguing by
contradiction, will be used in the proofs of Propositions
$\ref{surjective is associated}$ and $\ref{surjective is adapted}$.

\begin{lemma}\label{groupsGH}
Let $G$, $H$ be two groups and $p:G \rightarrow H$ a group
homomorphism. If $R$ is an integral domain, let us denote by
$p_R:R[G] \rightarrow R[H]$ the $R$-algebra morphism induced by $p$.
If $p_R$ is surjective, then $p$ is also surjective.
\end{lemma}

\begin{proposition}\label{surjective is associated}
Let $\varphi:A\rightarrow B$ be a surjective $R$-algebra morphism
such that  for all $i\,(0 \leq i \leq m)$, $\varphi(x_i)$ is a monomial
in $B$. Then $\varphi$ is associated with a monomial in $A$.
\end{proposition}
\begin{apod}{Due to the isomorphism of Lemma $\ref{work on Z}$,
$\varphi$ can be considered as a surjective $R$-algebra morphism
$$\varphi: R[\mathbb{Z}^{m+1}] \rightarrow R[\mathbb{Z}^{m}].$$ The
property of $\varphi$ about the $x_i$ implies that the above
morphism is induced by a $\mathbb{Z}$-module morphism $f:
\mathbb{Z}^{m+1} \rightarrow \mathbb{Z}^m$, which is also surjective
by Lemma $\ref{groupsGH}$. Since $\mathbb{Z}^m$ is a free
$\mathbb{Z}$-module the following exact sequence sequence splits
$$0\rightarrow\mathrm{Ker}f\rightarrow\mathbb{Z}^{m+1}
\rightarrow \mathbb{Z}^m\rightarrow0$$ and we obtain that
$\mathbb{Z}^{m+1}\cong \mathrm{Ker}f\,\oplus\,\mathbb{Z}^m$.
Therefore, $\mathrm{Ker}f$ is a $\mathbb{Z}$-module of rank $1$ and
there exists $a:=(a_0,a_1,\ldots,a_m)\in\mathbb{Z}^{m+1}$ such that
$\mathrm{Ker}f=\mathbb{Z}a$. By Theorem $\ref{similitude}$, $a$ is a
primitive element of $\mathbb{Z}^{m+1}$ and we must have
$\sum_{i=0}^m \mathbb{Z}a_i=\mathbb{Z}$, whence
$\mathrm{gcd}(a_i)=1$. By definition, the morphism $\varphi$ is
associated with the monomial $\prod_{i=0}^m x_i^{a_i}$.}
\end{apod}

Now let $r \in \{1,\ldots,m+1\}$ and $C_r:=R[y_r^{\pm1},\ldots,y_m^{\pm 1}]$, 
 where $y_r,\ldots,y_m$ are $m+1-r$ indeterminates over $R$. For $r=m+1$, $C_{r}=R$.

\begin{definition}\label{adapted morphism}
An $R$-algebra morphism $\varphi:A \rightarrow C_r$ is
called adapted, if $\varphi={\varphi_r} \circ {\varphi_{r-1}} \circ
\ldots \circ {\varphi_1}$, where $\varphi_i$ is a morphism
associated with a monomial for all $i=1,\ldots,r$.  \index{adapted morphism} The family
$\mathcal{F}:=\{\varphi_r,\varphi_{r-1},\ldots,\varphi_1\}$ is
called an adapted family  \index{adapted family} for $\varphi$ whose initial morphism is
$\varphi_1$.
\end{definition}

Let us introduce the following notation: If
$M:=\prod_{i=0}^mx_i^{c_i}$ is a monomial such that
$\textrm{gcd}(c_i)=d \in \mathbb{Z}$, then
$M^\circ:=\prod_{i=0}^mx_i^{c_i/d}$.

\begin{proposition}\label{change initial}
Let $\varphi:A \rightarrow C_r$ be an adapted morphism and
$M$ a monomial in $A$ such that $\varphi(M)=1$. Then there exists an
adapted family for $\varphi$ whose initial morphism is associated
with $M^\circ$.
\end{proposition}
\begin{apod}{Let $M=\prod_{i=0}^mx_i^{c_i}$ be a monomial in $A$
such that $\varphi(M)=1$. Note that $\varphi(M)=1$ if and only if
$\varphi(M^\circ)=1.$ Therefore, we can assume that
$\textrm{gcd}(c_i)=1$. We will prove the desired result by induction
on $r$.

If $r=1$, then, due to property $\ref{properties of phi}$(2), $\varphi$ must be a morphism associated with $M$.
If $r=2$, set $B:=R[z_1^{\pm1},\ldots,z_m^{\pm 1}]$.
  Let $\varphi:=\varphi_b \circ \varphi_a$,
  where \begin{itemize}
  \item $\varphi_a: A \rightarrow B$
        is a morphism associated with a monomial
        $\prod_{i=0}^mx_i^{a_i}$ in $A$ such that
        $\textrm{gcd}(a_i)=1$.
  \item $\varphi_b: B \rightarrow C_2$
        is a morphism associated with a monomial
        $\prod_{j=1}^mz_j^{b_j}$ in $B$ such that
        $\textrm{gcd}(b_j)=1$.
\end{itemize}
By Theorem $\ref{similitude}$, the element $a:=(a_0,a_1,\ldots,a_m)$
is a primitive element of $\mathbb{Z}^{m+1}$ and the element
$b:=(b_1,\ldots,b_m)$ is a primitive element of $\mathbb{Z}^{m}$.
Therefore, there exist a $\mathbb{Z}$-submodule $N_a$ of
$\mathbb{Z}^{m+1}$ and a $\mathbb{Z}$-submodule $N_b$ of
$\mathbb{Z}^{m}$ such that $\mathbb{Z}^{m+1}=\mathbb{Z}a \oplus N_a$
and $\mathbb{Z}^{m}=\mathbb{Z}b \oplus N_b$. We will denote by
$p_a:\mathbb{Z}^{m+1}\twoheadrightarrow N_a$ the projection of
$\mathbb{Z}^{m+1}$ onto $N_a$ and by $p_b:\mathbb{Z}^{m}\twoheadrightarrow
N_b$ the projection of $\mathbb{Z}^{m}$ onto $N_b$.
We have isomorphisms $i_a: N_a\tilde{\rightarrow}\mathbb{Z}^{m}$ and
$i_b: N_b \tilde{\rightarrow}\mathbb{Z}^{m-1}$.

By definition of the associated morphism, $\varphi_a$ is induced by
the morphism
$f_a:=i_a \circ p_a:\mathbb{Z}^{m+1}\twoheadrightarrow \mathbb{Z}^{m}$ and
$\varphi_b$ by
 $f_b:=i_b \circ p_b:\mathbb{Z}^{m} \twoheadrightarrow
 \mathbb{Z}^{m-1}.$
Set $f:=f_b \circ f_a$. Then $\varphi$ is the $R$-algebra
morphism induced by $f$.

The morphism $f$ is surjective. Since $\mathbb{Z}^{m-1}$ is
a free $\mathbb{Z}$-module, the following exact sequence sequence
splits
$$0\rightarrow\mathrm{Ker}f\rightarrow\mathbb{Z}^{m+1}
\rightarrow \mathbb{Z}^{m-1}\rightarrow0$$ and we obtain that
$\mathbb{Z}^{m+1}\cong \mathrm{Ker}f\oplus\mathbb{Z}^{m-1}$.

Let $\tilde{b}:=i_a^{-1}(b)$. Then $\mathrm{Ker}f= \mathbb{Z}a
\oplus \mathbb{Z}\tilde{b}$. By assumption, we have that
$c:=(c_0,c_1,\ldots,c_m) \in \mathrm{Ker}f$. Therefore, there exist
unique $\lambda_1, \lambda_2 \in \mathbb{Z}$ such that $c=\lambda_1
a+\lambda_2 \tilde{b}$. Since $\mathrm{gcd}(c_i)=1$, we must also
have $\mathrm{gcd}(\lambda_1,\lambda_2)=1$. Hence $\sum_{i=1}^2
\mathbb{Z}\lambda_i=\mathbb{Z}$. By applying Theorem
$\ref{similitude}$ to the $\mathbb{Z}$-module $\mathrm{Ker} f$, we
obtain that $c$ is a primitive element of $\mathrm{Ker}f$.
Consequently, $f=f\,' \circ f_c$, where $f_c$ is a surjective
$\mathbb{Z}$-module morphism $\mathbb{Z}^{m+1} \twoheadrightarrow
\mathbb{Z}^{m}$ with $\mathrm{Ker}f_c= \mathbb{Z}c$ and
 $f\,'$ is a surjective $\mathbb{Z}$-module
morphism $\mathbb{Z}^{m}
\twoheadrightarrow \mathbb{Z}^{m-1}$. As far as the induced
$R$-algebra morphisms are concerned, we obtain that
$\varphi=\varphi\,' \circ \varphi_c$, where $\varphi_c$ is a morphism
associated with the monomial $M=\prod_{i=0}^mx_i^{c_i}$ and, due to Proposition $\ref{surjective is associated}$,
$\varphi\,'$ is also a morphism associated with a monomial. Thus the assertion is proven for $r=2$.

 Now, let us suppose that  $r>2$ and that the assertion holds for  $1,2,\ldots,r-1$.
  If $\varphi={\varphi_r} \circ {\varphi_{r-1}} \circ \ldots \circ {\varphi_1}$, the induction hypothesis
  implies that there exist morphisms associated with monomials $\varphi_r', \varphi_{r-1}',\ldots,\varphi_2'$
  such that
  \begin{itemize}
  \item $\varphi={\varphi_r'}\circ \ldots\circ {\varphi_2'}\circ{\varphi_1}$.
  \item $\varphi_2'$ is associated with the monomial $(\varphi_1(M))^\circ$.
  \end{itemize}
  We have that $\varphi_2'(\varphi_1(M))=1$. Once more, by induction hypothesis we obtain that
  there exist morphisms associated with monomials $\varphi_2'',\varphi_1''$
  such that
  \begin{itemize}
    \item ${\varphi_2'} \circ{\varphi_1}={\varphi_2''}\circ{\varphi_1''}$.
  \item $\varphi_1''$ is associated with $M^\circ$.
  \end{itemize}
  Thus we have
  $$\varphi={\varphi_r'} \circ \ldots \circ {\varphi_3'} \circ {\varphi_2'}\circ{\varphi_1}=
  {\varphi_r'} \circ \ldots \circ{\varphi_3'} \circ{\varphi_2''} \circ{\varphi_1''}$$
  and $\varphi_1''$ is associated with $M^\circ$.}
\end{apod}

\begin{proposition}\label{surjective is adapted}
Let $\varphi:A\rightarrow C_r$ be a surjective $R$-algebra
morphism such that  for all $i\,(0 \leq i \leq m)$, $\varphi(x_i)$ is a monomial
in $C_r$. Then $\varphi$ is an adapted morphism.
\end{proposition}
\begin{apod}{We will work again by induction on $r$.
For $r=1$, the above result is Proposition $\ref{surjective is
  associated}$. For $r>1$, let us suppose that the result
is true for $1,\ldots,r-1$. Due to the isomorphism of Lemma
$\ref{work on Z}$, $\varphi$ can be considered as a surjective
$R$-algebra morphism
$$\varphi: R[\mathbb{Z}^{m+1}] \rightarrow R[\mathbb{Z}^{m+1-r}].$$ The
property of $\varphi$ about the $x_i$ implies that the above
morphism is induced by a $\mathbb{Z}$-module morphism $f:
\mathbb{Z}^{m+1} \rightarrow \mathbb{Z}^{m+1-r}$, which is also
surjective by Lemma $\ref{groupsGH}$. Since $\mathbb{Z}^{m+1-r}$ is
a free $\mathbb{Z}$-module the following exact sequence 
splits
$$0\rightarrow\mathrm{Ker}f\rightarrow\mathbb{Z}^{m+1}
\rightarrow \mathbb{Z}^{m+1-r}\rightarrow0$$ and we obtain that
$\mathbb{Z}^{m+1}\cong \mathrm{Ker}f\oplus\mathbb{Z}^{m+1-r}$.
Therefore, $\mathrm{Ker}f$ is a $\mathbb{Z}$-module of rank $r$, \ie
$\mathrm{Ker}f \cong \mathbb{Z}^r$. We choose a primitive element
$a$ of $\mathrm{Ker}f$. Then there exists a $\mathbb{Z}$-submodule
$N_a$ of $\mathrm{Ker}f$ such that $\mathrm{Ker}f=\mathbb{Z}a \oplus
N_a$. Since $\mathrm{Ker}f$ is a direct summand of
$\mathbb{Z}^{m+1}$, $a$ is also a primitive element of
$\mathbb{Z}^{m+1}$ and we have
$$\mathbb{Z}^{m+1} \cong \mathbb{Z}a \oplus N_a \oplus \mathbb{Z}^{m+1-r}.$$
Thus, by Theorem $\ref{similitude}$, if $(a_0,a_1,\ldots,a_m)$ are
the coefficients of $a$ with respect to the standard basis of
$\mathbb{Z}^{m+1}$, then $\mathrm{gcd}(a_i)=1$.

Let us denote by $p_a$ the projection $\mathbb{Z}^{m+1}
\twoheadrightarrow N_a \oplus \mathbb{Z}^{m+1-r}$
 and by $p\,'$ the
projection $N_a \oplus \mathbb{Z}^{m+1-r} \twoheadrightarrow
\mathbb{Z}^{m+1-r}$. We have $f=p\,' \circ p_a$. Consequently, if we denote by $i$ the
$\mathbb{Z}$-module isomorphism
$$i: N_a \oplus \mathbb{Z}^{m+1-r} \tilde{\rightarrow}  \mathbb{Z}^{m}$$
and set $f_a:=i \circ p_a$ and $f\,':=p\,' \circ i^{-1}$, we obtain that $f=f\,'
\circ f_a$. If $\varphi_a$ is the $R$-algebra morphism induced by
$f_a$, then, by definition, $\varphi_a$ is a morphism associated
with the monomial $\prod_{i=0}^m x_i^{a_i}$. The $R$-algebra
morphism $\varphi\,'$ induced by $f\,'$ is a surjective morphism
with the same property as $\varphi$ (it sends
all monomials to monomials). By induction hypothesis, $\varphi\,'$ is an
adapted morphism. We have $\varphi=\varphi\,' \circ \varphi_a$ and
so $\varphi$ is also an adapted morphism.}
\end{apod}

\section{Irreducibility}
In this section, we will discuss the irreducibility of polynomials in a Laurent polynomial ring with coefficients in a field.

Let $k$ be a field and $y$ an indeterminate over $k$. We can use the
following theorem in order to determine the irreducibility of a
polynomial of the form $y^n-a$ in $k[y]$ (cf.\cite{La}, Chapter 6,
Thm. 9.1).

\begin{theorem}\label{Lang}
Let $k$ be a field, $a \in k-\{0\}$ and $n \in \mathbb{Z}$ with $n
\geq 2$. The polynomial $y^n-a$ is irreducible in $k[y]$, if for
every prime $p$ dividing $n$, we have $a \notin k^p$ and if $4$
divides $n$, we have $a \notin -4k^4$.
\end{theorem}

Let $x_0,x_1,\ldots,x_m$ be a set of $m+1$ indeterminates over $k$. We will
apply Theorem $\ref{Lang}$ to the field $k(x_1,\ldots,x_m)$.

\begin{proposition}
Let $k$ be a field. The polynomial $x_0^{a_0}-\rho\prod_{i=1}^m
x_i^{a_i}$ with $\rho \in k-\{0\}$, $a_i \in \mathbb{Z}$,
\emph{gcd}$(a_i)=1$ and $a_0 > 0$ is irreducible in
$k(x_1,\ldots,x_m)[x_0]$.
\end{proposition}
\begin{apod}{If $a_0=1$, the polynomial is of degree 1 and thus irreducible in
$k(x_1,\ldots,x_m)[x_0]$. If $a_0 \geq 2$, let $p$ be a prime number dividing $a_0$.
Let us suppose that $\rho\prod_{i=1}^m x_i^{a_i}$ belongs to $k(x_1,\ldots,x_m)^p$, \ie that there exist
 $f(x_1,\ldots,x_m),g(x_1,\ldots,x_m) \in k[x_1,\ldots,x_m]$
prime to each other, with $g(x_1,\ldots,x_m) \neq 0$, such that
$$\rho\prod_{i=1}^m x_i^{a_i}=\left(\frac{f(x_1,\ldots,x_m)}{g(x_1,\ldots,x_m)}\right)^p$$
The above relation can
be rewritten as
$$g(x_1,\ldots,x_m)^p \cdot \left(\rho \prod_{\{i\,|\,a_i\geq 0\}}x_i^{a_i}\right)=f(x_1,\ldots,x_m)^p \cdot\left(\prod_{\{i\,|\,a_i< 0\}}x_i^{-a_i}\right).$$
We have that
$$\textrm{gcd}(f(x_1,\ldots,x_m)^p,g(x_1,\ldots,x_m)^p)=1.$$
Since $k[x_1,\ldots,x_m]$ is a unique factorization domain and the $x_i$
are irreducible in $k[x_1,\ldots,x_m]$, we also have that
$$\textrm{gcd}(\prod_{\{i\,|\,a_i\geq 0\}}x_i^{a_i},\prod_{\{i\,|\,a_i<
0\}}x_i^{-a_i})=1.$$ As a result,
$$f(x_1,\ldots,x_m)^p = \lambda\rho\cdot\prod_{\{i\,|\,a_i\geq 0\}}x_i^{a_i} \textrm{ and }
g(x_1,\ldots,x_m)^p = \lambda\cdot \prod_{\{i\,|\,a_i< 0\}}x_i^{-a_i}.$$ for
some $\lambda \in k-\{0\}$. If $(\lambda\rho)^{1/p},\lambda^{1/p} \notin k$, we have arrived to a contradiction. Suppose that
$(\lambda\rho)^{1/p},\lambda^{1/p} \in k$. Once more, the fact that
$k[x_1,\ldots,x_m]$ is a unique factorization domain and the $x_i$ are
irreducible in $k[x_1,\ldots,x_m]$ implies that
$$f(x_1,\ldots,x_m) = (\lambda\rho)^{1/p}\cdot\prod_{\{i\,|\,a_i\geq 0\}}x_i^{b_i} \textrm{ and }
g(x_1,\ldots,x_m) = \lambda^{1/p}\cdot\prod_{\{i\,|\,a_i< 0\}}x_i^{-b_i},$$
with $b_i \in \mathbb{Z}$ and $b_ip=a_i,\forall i=1,\ldots,m$. Since
$p$ divides $a_0$, this contradicts the fact that $\textrm{gcd}(a_i)=1$. In
the same way, we can show that if $4$ divides $a_0$, then $\rho\prod_{i=1}^m
x_i^{a_i} \notin -4k(x_1,\ldots, x_m)^4$. Thus, by Theorem
$\ref{Lang}$, $x_0^{a_0}-\rho\prod_{i=1}^m x_i^{a_i}$ is irreducible
in $k(x_1,\ldots,x_m)[x_0]$.}
\end{apod}

Thanks to Lemma $\ref{irreducibility in field of fractions}$, the above proposition implies that

\begin{proposition}\label{irreducibility}
Let $k$ be a field. The polynomial $x_0^{a_0}-\rho\prod_{i=1}^m
x_i^{a_i}$ with $\rho \in k-\{0\}$, $a_i \in \mathbb{Z}$,
\emph{gcd}$(a_i)=1$ and $a_0 > 0$ is irreducible in
$k[x_1^{\pm1},\ldots,x_m^{\pm1}][x_0]$.
\end{proposition}

\begin{lemma}\label{irreducibility in field of fractions}
Let $R$ be an integral domain with field of fractions $F$ and $f(x)$
a polynomial in $R[x]$. If $f(x)$ is irreducible in $F[x]$ and at
least one of its coefficients is a unit in $R$, then $f(x)$ is
irreducible in $R[x]$.
\end{lemma}
\begin{apod}{If $f(x)=g(x)h(x)$ for two polynomials
$g(x),h(x) \in R[x]$, then $g(x) \in R$ or $h(x) \in R$. Let us
suppose that $g(x) \in R$. Since one of the coefficients of $f(x)$
is a unit in $R$, $g(x)$ must also be a unit in $R$. Thus, $f(x)$ is
irreducible in $R[x]$.}
\end{apod}

The next result is an immediate consequence of Proposition $\ref{irreducibility}$ and will be used in the proof of Theorem
$\ref{second irreducible}$.

\begin{corollary}\label{M-rho}
Let $M:=\prod_{i=0}^mx_i^{a_i}$ be a monomial in
$k[x_0^{\pm1},x_1^{\pm1},\ldots,x_m^{\pm1}]$ such that
$\mathrm{gcd}(a_i)=1$ and let $\rho \in k-\{0\}$. Then $M-\rho$ is
an irreducible element of
$k[x_0^{\pm1},x_1^{\pm1},\ldots,x_m^{\pm1}]$.
\end{corollary}

\begin{theorem}\label{second irreducible}
Let $M:=\prod_{i=0}^mx_i^{a_i}$ be a monomial in
$k[x_0^{\pm1},x_1^{\pm1},\ldots,x_m^{\pm1}]$ such that
$\mathrm{gcd}(a_i)=1$. If $f(x)$ is an irreducible element of $k[x]$
such that $f(0) \neq 0$, then $f(M)$ is irreducible in
$k[x_0^{\pm1},x_1^{\pm1},\ldots,x_m^{\pm1}]$.
\end{theorem}
\begin{apod}{Suppose that $f(M)=g \cdot h$
with $g,h \in k[x_0^{\pm1},x_1^{\pm1},\ldots,x_m^{\pm1}]$. Let
$\rho_1,\ldots,\rho_n$ be the roots of $f(x)$ in a splitting field
$k'$. Then
$$f(x)=a (x-\rho_1)\ldots (x-\rho_n)$$
for some $a \in k-\{0\}$, whence
$$f(M)=a (M-\rho_1)\ldots (M-\rho_n).$$
By Corollary $\ref{M-rho}$, $M-\rho_j$ is irreducible in
$k'[x_0^{\pm1},x_1^{\pm1},\ldots,x_m^{\pm1}]$ for all $j \in
\{1,\ldots,n\}$. Since $k'[x_0^{\pm1},x_1^{\pm1},\ldots,x_m^{\pm1}]$
is a unique factorization domain, we must have
$$g=r \cdot\left(\prod_{i=0}^m x_i^{b_i}\right)\cdot (M-\rho_{j_1})\ldots(M-\rho_{j_s})$$
for some $r \in k-\{0\}$, $b_i \in \mathbb{Z}$ and $j_1,\ldots,j_s \in
\{1,\ldots,n\}$ with $j_1<j_2<\ldots<j_s$. Thus, there exists $g'(x) \in k[x]$ such that
$$g=\left(\prod_{i=0}^m x_i^{b_i}\right)\cdot g'(M).$$
Respectively, there exists $h'(x) \in k[x]$ such
that$$h=\left(\prod_{i=0}^m x_i^{-b_i}\right)\cdot h'(M).$$ Thus, we obtain that
$$f(M)=g'(M)h'(M).$$
Since $\mathrm{gcd}(a_i)=1$, there exist integers $(u_i)_{0 \leq i
\leq m}$ such that $\sum_{i=0}^m u_ia_i=1$. Let us now consider the
$k$-algebra specialization
$$\begin{array}{cccc}
  \varphi: & k[x_0^{\pm1},x_1^{\pm1},\ldots,x_m^{\pm1}] & \rightarrow & k[x] \\
           & x_i                                        & \mapsto     & x^{u_i}.
\end{array}$$
Then $\varphi(M)=\varphi(\prod_{i=0}^m x_i^{a_i})=x^{\sum_{i=0}^m
u_ia_i}=x$. If we apply $\varphi$ to the equality $f(M)=g'(M)h'(M)$, we
obtain that
$$f(x)=g'(x)h'(x).$$
Since $f(x)$ is irreducible in $k[x]$, we must have that either
$g'(x) \in k$ or $h'(x) \in k$. Respectively, we deduce that either
$g$ or $h$ is a unit in
$k[x_0^{\pm1},x_1^{\pm1},\ldots,x_m^{\pm1}]$.}
\end{apod}

\chapter{On Blocks}

The second chapter focuses on the study of blocks of algebras of finite type and in particular, of symmetric algebras. In the first section, we introduce the notion of blocks and show how we can reduce, in many different ways, the problem of the determination of the blocks of an algebra to easier cases (mostly by changing the ring of definition). In the second section, we give the definition of a symmetric algebra and introduce the \emph{Schur elements} associated to its irreducible characters. We explain why the Schur elements play an essential role in the determination of the blocks of a symmetric algebra.
Now, if we know the blocks of a symmetric algebra, then, given that certain conditions are satisfied, we can obtain the Schur elements and the blocks of its subalgebras with the use of a
 generalization of some classical results, known as ``Clifford theory'', to the case of ``twisted symmetric algebras of finite groups''. This is the object of section 2.3. Finally, in the last section, we present some more results on the block theory of symmetric algebras, which are related to the representation theory of these algebras. We define decomposition maps and generalize results of Brauer theory to the case of symmetric algebras defined over discrete valuation rings.
Moreover, in subsection 2.4.4, we give a new criterion for a symmetric algebra to be split semisimple.

Throughout this chapter, we use the simple example of group algebras in order to illustrate some of the notions. However, when we introduce Hecke algebras in Chapter 4, we will see that they are symmetric too. In particular, the Appendix contains numerous applications of Clifford theory to Hecke algebras. If the reader feels that they are in need of more examples than the ones given in Chapter 2, they should feel free to refer to Chapter 4 or the Appendix.

\section{Generalities}

Let $\mathcal{O}$ be a commutative ring with a unit element and $A$
be an $\mathcal{O}$-algebra. We denote by $ZA$ the center of $A$.

An idempotent in $A$ is an element $e$ such that $e^2=e$. We say
that $e$ is a central idempotent, if it is an idempotent in $ZA$.
Two idempotents $e_1,e_2$ are orthogonal, if $e_1e_2=e_2e_1=0$.
Finally, an idempotent $e$ is primitive, if $e \neq 0$ and $e$ can
not be expressed as the sum of two non-zero orthogonal idempotents.

\begin{definition}\label{blocks}
The block-idempotents \index{block-idempotent} of $A$ are the central primitive idempotents
of $A$.
\end{definition}

Let $e$ be a block-idempotent of $A$. The two sided ideal $Ae$
inherits a structure of algebra, where the composition laws are
those of $A$ and the unit element is $e$. The application
$$\begin{array}{cccc}
  \pi_e: &A &\rightarrow &Ae \\
         &h &\mapsto &he
\end{array}$$
is an epimorphism of algebras. The algebra $Ae$ is called a block of
$A$. From now on, abusing the language, we will also call blocks  \index{block} the
block-idempotents of $A$.

\begin{lemma}\label{orthogonality of blocks}
The blocks of $A$ are mutually orthogonal.
\end{lemma}
\begin{apod}{Let $e$ be a block and $f$ a central idempotent of $A$ with $f \neq e$.
Then $ef$ and $e-ef$ are also central idempotents. We have
$e=ef+(e-ef)$ and due to the primitivity of $e$, we deduce that
either $ef=0$ or $e=ef$. If $f$ is a block too, then either $ef=0$
or $f=ef=e$. Therefore, $f$ is orthogonal to $e$.}
\end{apod}

The above lemma gives rise to the following proposition.

\begin{proposition}\label{1 sum of blocks}
Suppose that the unity element $1$ of $A$ can be expressed as a sum
of blocks: $1=\sum_{e \in E}e$. Then
\begin{enumerate}
  \item The set $E$ is the set of all the blocks of $A$.
  \item The family of morphisms $(\pi_e)_{e \in E}$ defines an
  isomorphism of algebras
  $$A \tilde{\rightarrow} \prod_{e \in E} Ae.$$
\end{enumerate}
\end{proposition}
\begin{apod}{If $f$ is a block, then $f=\sum_{e \in E}ef$. Due to
lemma $\ref{orthogonality of blocks}$, there exists $e \in E$ such
that $f=e$.}
\end{apod}

In the above context (1 is a sum of blocks), let us denote by
$\mathrm{Bl}(A)$ the set of all the blocks of $A$. Proposition
$\ref{1 sum of blocks}$ implies that the category
$_A\mathrm{\textbf{mod}}$ of $A$-modules is a direct sum of the
categories associated with the blocks:
$$_A\mathrm{\textbf{mod}} \tilde{\rightarrow} \bigoplus_{e \in \mathrm{Bl}(A)} \mathrm{}_{Ae}\mathrm{\textbf{mod}}.$$
In particular, every representation of the $\mathcal{O}$-algebra
$Ae$ defines (by composition with $\pi_e$) a representation of $A$
and we say, abusing the language, that it ``belongs to the block
$e$''.

Every indecomposable representation of $A$ belongs to one and only
one block. Thus the following partitions are defined:
$$\mathrm{Ind}(A)=\bigsqcup_{e \in \mathrm{Bl}(A)} \mathrm{Ind}(A,e) \,\,\textrm{ and
}\,\, \mathrm{Irr}(A)=\bigsqcup_{e \in \mathrm{Bl}(A)}
\mathrm{Irr}(A,e),$$ where $\mathrm{Ind}(A)$ (resp.
$\mathrm{Irr}(A)$) denotes the set of indecomposable (resp.
irreducible) representations of $A$ and $\mathrm{Ind}(A,e)$ (resp.
$\mathrm{Irr}(A,e)$) denotes the set of the elements of
$\mathrm{Ind}(A)$ (resp. $\mathrm{Irr}(A)$) which belong to $e$.\\

We will consider two situations where 1 is a sum of blocks.\\
\\
\emph{First case:} Suppose that 1 is a sum of orthogonal primitive
idempotents, \ie $1=\sum_{i \in P}i$, where
\begin{itemize}
  \item every $i \in P$ is a primitive idempotent,
  \item if $i,j \in P$, $i \neq j$, then $ij=ji=0$.
\end{itemize}
Let us consider the equivalence relation $\mathcal{B}$ defined on
$P$ as the symmetric and transitive closure of the relation ``$iAj
\neq \{0\}$''. Thus $(i \mathcal{B} j)$ if and only if there exist
$i_0,i_1,\ldots,i_n \in P$ with $i_0=i$ and $i_n=j$ such that for
all $k \in \{1,\ldots,n\}$, $i_{k-1}Ai_k \neq \{0\}$ or $i_kAi_{k-1}
\neq \{0\}$. To every equivalence class $B$ of $P$ with respect to
$\mathcal{B}$, we associate the idempotent $e_B:=\sum_{i \in B}i$.

\begin{proposition}\label{iAj}
The map $B \mapsto e_B$ is a bijection between the set of
equivalence classes of $\mathcal{B}$ and the set of blocks of $A$.
In particular, we have that $1=\sum_{B \in P/\mathcal{B}}e_B$ and
$1$ is sum of the blocks of $A$.
\end{proposition}
\begin{apod}{It is clear that $1=\sum_{B \in P/\mathcal{B}}e_B$. Let
$a \in A$ and let $B,B'$ be two equivalence classes of $\mathcal{B}$
with $B \neq B'$. Then, by definition of the relation $\mathcal{B}$,
$e_B a e_{B'}=0$. Since $1=\sum_{B \in P/\mathcal{B}}e_B$, we have
that $e_Ba=e_Bae_B=ae_B$. Thus $e_B \in ZA$ for all $B \in
P/\mathcal{B}$.

It remains to show that for all $B \in P/\mathcal{B}$, the central
idempotent $e_B$ is primitive. Suppose that $e_B=e+f$, where $e$ and
$f$ are two orthogonal primitive idempotents in $ZA$. Then we have a
partition $B=B_e \sqcup B_f$, where $B_e:=\{i \in B \,|\, ie=i\}$
and $B_f:=\{j \in B \,|\, jf=j\}$. For all $i \in B_e$ and $j \in
B_f$, we have $iAj=ieAfj=iAefj=\{0\}$ and so no element of $B_e$ can
be $\mathcal{B}$-equivalent to an element of $B_f$. Therefore, we
must have either $B_e=\emptyset$ or $B_f=\emptyset$, which implies
that either $e=0$ or $f=0$.}
\end{apod}
$ $\\
\emph{Second case:} Suppose that $ZA$ is a subalgebra of a
commutative algebra $C$ where 1 is a sum of blocks. For example, if
$A$ is of finite type over $\mathcal{O}$, where $\mathcal{O}$ is an
integral domain with field of fractions $F$, we can choose $C$ to be
the center of the algebra $FA:=F \otimes_\mathcal{O} A$.
We set $1=\sum_{e \in E}e$, where $E$ is the set of blocks of $C$.
For all $S \subseteq E$, set $e_S:=\sum_{e \in S}e$. A subset $S$ of
$E$ is ``on $ZA$'' if $e_S \in ZA$. If $S$ and $T$ are on $ZA$, then
$S \cap T$ is on $ZA$.

\begin{proposition}\label{partition}
Let us denote by $\mathcal{P}_E(ZA)$ the set of non-empty subsets
$B$ of $E$ which are on $ZA$ and are minimal for these two
properties. Then the map $\mathcal{P}_E(ZA) \rightarrow A, B \mapsto
e_B$ induces a bijection between $\mathcal{P}_E(ZA)$ and the set of
blocks of $A$. We have $1=\sum_{B \in \mathcal{P}_E(ZA)}e_B$.
\end{proposition}
\begin{apod}{Since every idempotent in $C$ is of the form $e_S$ for
some $S \subseteq E$, it is clear that $e_B$ is a central primitive
idempotent of $A$, for all $B \in \mathcal{P}_E(ZA)$. It remains to
show that
\begin{enumerate}
  \item If $B$,$B'$ are two distinct elements of
  $\mathcal{P}_E(ZA)$, then $B \cap B' = \emptyset$.
  \item $\mathcal{P}_E(ZA)$ is a partition of $E$.
\end{enumerate}
These two properties, stated in terms of idempotents, mean:
\begin{enumerate}
  \item If $B$,$B'$ are two distinct elements of
  $\mathcal{P}_E(ZA)$, then $e_B$ and $e_{B'}$ are orthogonal.
  \item $1=\sum_{B \in \mathcal{P}_E(ZA)}e_B$.
\end{enumerate}
Let us prove them:
\begin{enumerate}
  \item We have $e_Be_{B'}=e_{B \cap B'}$ and so $B \cap B' = \emptyset$, because $B$ and $B'$ are
  minimal.
  \item Set $F:=\bigcup_{B \in \mathcal{P}_E(ZA)}B$. Then $e_F=\sum_{B \in \mathcal{P}_E(ZA)}e_B \in
  ZA$. Then $1-e_F=e_{E-F} \in ZA$, which means that $E-F$ is on
  $ZA$ . If $E-F \neq \emptyset$, then $E-F$ contains an element
  of $\mathcal{P}_E(ZA)$ in contradiction to the definition of $F$.
  Thus $F=E$ and $\mathcal{P}_E(ZA)$ is a partition of $E$.}
\end{enumerate}
\end{apod}

Now let us assume that
\begin{itemize}
  \item $\mathcal{O}$ is a commutative integral domain with field of fractions
  $F$,
  \item $K$ is a field extension of $F$,
  \item $A$ is an $\mathcal{O}$-algebra, free and finitely generated as an
  $\mathcal{O}$-module.
\end{itemize}

Suppose that the $K$-algebra $KA:=K \otimes_\mathcal{O}A$ is
semisimple.  \index{semisimple algebra} Then $KA$ is isomorphic, by assumption, to a direct
product of simple algebras:
$$KA \cong \prod_{\chi \in \mathrm{Irr}(KA)} M_\chi,$$
where $\mathrm{Irr}(KA)$ denotes the set of irreducible characters
of $KA$ and $M_\chi$ is a simple $K$-algebra.

For all $\chi \in \mathrm{Irr}(KA)$, we denote by $\pi_\chi:KA
\twoheadrightarrow M_\chi$ the projection onto the $\chi$-factor and
by $e_\chi$ the element of $KA$ such that $$\pi_{\chi'}(e_\chi)=
  \left\{
  \begin{array}{ll}
    1_{M_\chi}, & \hbox{if $\chi=\chi'$,} \\
    0, & \hbox{if $\chi \neq \chi'$.}
  \end{array}
\right.$$

The following theorem results directly from Propositions $\ref{1 sum
of blocks}$ and $\ref{partition}$.

\begin{theorem}\label{minimality of blocks}\
\begin{enumerate}
  \item We have $1=\sum_{\chi \in \mathrm{Irr}(KA)}e_\chi$
   and the set $\{e_\chi\}_{\chi \in \mathrm{Irr}(KA)}$ is the set of all the blocks of the algebra $KA$.
  \item There exists a unique partition $\mathrm{Bl}(A)$ of
  $\mathrm{Irr}(KA)$ such that
  \begin{enumerate}[(a)]
    \item For all $B \in \mathrm{Bl}(A)$, the idempotent
    $e_B:=\sum_{\chi \in B}e_\chi$ is a block of $A$.
    \item We have $1=\sum_{B \in \mathrm{Bl}(A)}e_B$ and for
    every central idempotent $e$ of $A$, there exists a subset
    $\mathrm{Bl}(A,e)$ of $\mathrm{Bl}(A)$ such that
    $$e=\sum_{B \in \mathrm{Bl}(A,e)}e_B.$$
    In particular the set $\{e_B\}_{B \in \mathrm{Bl}(A)}$ is the set of all the blocks of $A$.
  \end{enumerate}
\end{enumerate}
\end{theorem}

\begin{remarks}\
\emph{\begin{itemize}
  \item If $\chi \in B$ for some $B \in \mathrm{Bl}(A)$, we say that
  ``$\chi$ belongs to the block $e_B$''.
  \item For all $B \in \mathrm{Bl}(A)$, we have
  $$KAe_B \cong \prod_{\chi \in B}M_\chi.$$
\end{itemize}}
\end{remarks}

From now on, we make the following assumptions

\begin{ypoth}\label{properties of the ring}\
\begin{description}
  \item[(int)] The ring $\mathcal{O}$ is a Noetherian and integrally
  closed domain with field of fractions $F$ and $A$ is an
  $\mathcal{O}$-algebra which is free and finitely generated as an
  $\mathcal{O}$-module.
  \item[(spl)] The field $K$ is a finite Galois extension of $F$ and
  the algebra $KA$ is split  \index{split algebra} (i.e., for every simple $KA$-module $V$, $\mathrm{End}_{KA}(V) \cong K$) semisimple.
\end{description}
\end{ypoth}

We denote by $\mathcal{O}_K$ the integral closure of $\mathcal{O}$
in $K$.

\begin{px}\emph{\small Let $C_6$ be the cyclic group of order $6$ and $A:=\mathbb{Z}[C_6]$. If we take $K:=\mathbb{Q}(\zeta_6)$, where $\zeta_6:=\mathrm{exp}(\frac{2\pi i}{6})$, then the assumptions $\ref{properties of the ring}$ are satisfied.}
\end{px}

\subsection{Blocks and integral closure}

The Galois group $\mathrm{Gal}(K/F)$ acts on $KA=K
\otimes_{\mathcal{O}} A$ (viewed as an $F$-algebra) as follows: if
$\sigma \in \mathrm{Gal}(K/F)$ and $\lambda \otimes a \in KA$, then
$\sigma(\lambda \otimes a):=\sigma(\lambda) \otimes a$.

If $V$ is a $K$-vector space and $\sigma \in \mathrm{Gal}(K/F)$, we
denote by $^\sigma V$ the $K$-vector space defined on the additive
group $V$ with multiplication $\lambda.v:=\sigma^{-1}(\lambda)v$ for
all $\lambda \in K$ and $v \in V$. If $\rho:KA \rightarrow
\mathrm{End}_K(V)$ is a representation of the $K$-algebra $KA$, then
its composition with the action of $\sigma^{-1}$ is also a
representation $^\sigma \rho: KA \rightarrow \mathrm{End}_K(^\sigma
V)$:
$$\diagram KA \rto^{\sigma^{-1}} &KA \rto^{\rho}
&\mathrm{End}_K(V). \enddiagram$$

We denote by $^\sigma \chi$ the character of $^\sigma \rho$ and we
define the action of $\mathrm{Gal}(K/F)$ on $\mathrm{Irr}(KA)$ as
follows: if $\sigma \in \mathrm{Gal}(K/F)$ and $\chi \in
\mathrm{Irr}(KA)$, then $$\sigma(\chi):={}^\sigma\!\chi = \sigma
\circ \chi \circ \sigma^{-1}.$$ This operation induces an action of
$\mathrm{Gal}(K/F)$ on the set of blocks of $KA$:
$$\sigma(e_\chi)=e_{^\sigma \chi} \textrm{ for all } \sigma \in
\mathrm{Gal}(K/F), \chi \in \mathrm{Irr}(KA).$$

Hence, the group $\mathrm{Gal}(K/F)$ acts on the set of idempotents
of $Z\mathcal{O}_KA$ and thus on the set of blocks of
$\mathcal{O}_KA$. Since $F \cap \mathcal{O}_K = \mathcal{O}$, the
idempotents of $ZA$ are the idempotents of $Z\mathcal{O}_KA$ which
are fixed by the action of $\mathrm{Gal}(K/F)$. As a consequence,
the primitive idempotents of $ZA$ are sums of the elements of the
orbits of $\mathrm{Gal}(K/F)$ on the set of primitive idempotents of
$Z\mathcal{O}_KA$. Thus, the blocks of $A$ are in bijection with the
orbits of $\mathrm{Gal}(K/F)$ on the set of blocks of
$\mathcal{O}_KA$. The following proposition is just a reformulation
of this result.

\begin{proposition}\label{Galois action on integral closure}\
\begin{enumerate}
  \item Let $B$ be a block of $A$ and $B'$ a block of
  $\mathcal{O}_KA$ contained in $B$. If $\mathrm{Gal}(K/F)_{B'}$
  denotes the stabilizer of $B'$ in $\mathrm{Gal}(K/F)$, then
  $$B=\bigcup_{\sigma \in
  \mathrm{Gal}(K/F)/\mathrm{Gal}(K/F)_{B'}}\sigma(B')
  \,\,\textrm{ i.e., }\,\,
  e_B=\sum_{\sigma \in
  \mathrm{Gal}(K/F)/\mathrm{Gal}(K/F)_{B'}}\sigma(e_{B'}).$$
  \item Two characters $\chi,\psi \in \mathrm{Irr}(KA)$ are in
  the same block of $A$ if and only if there exists $\sigma \in \mathrm{Gal}(K/F)$
  such that $\sigma(\chi)$ and $\psi$ belong to the same block of
  $\mathcal{O}_KA$.
\end{enumerate}
\end{proposition}\
\begin{remark}\emph{ For all $\chi \in B'$, we have
$\mathrm{Gal}(K/F)_\chi \subseteq \mathrm{Gal}(K/F)_{B'}.$}
\end{remark}\
\\

The assertion (2) of the proposition above allows us to transfer the
problem of the classification of the blocks of $A$ to that of the
classification of the blocks of $\mathcal{O}_KA$.

\subsection{Blocks and prime ideals}

We denote by $\mathrm{Spec}_1(\mathcal{O})$ the set of prime ideals
of height $1$ of $\mathcal{O}$. Since $\mathcal{O}$ is Noetherian and
integrally closed, Proposition $\ref{case of Krull}$ implies that $\mathcal{O}$ is a Krull ring. By Theorem
$\ref{Krull-dvr}$, we have
$$\mathcal{O}=\bigcap_{\mathfrak{p} \in \mathrm{Spec}_1(\mathcal{O})}
\mathcal{O}_\mathfrak{p},$$ where $\mathcal{O}_\mathfrak{p}:=\{x \in
F\,|\,(\exists a \in \mathcal{O}-\mathfrak{p})(ax \in
\mathcal{O})\}$ is the localization of $\mathcal{O}$ at
$\mathfrak{p}$. In particular, $\mathcal{O}_\mathfrak{p}$ is a discrete valuation ring.

Let $\mathfrak{p}$ be a prime ideal of $\mathcal{O}$ and
$\mathcal{O}_\mathfrak{p}A:=\mathcal{O}_\mathfrak{p}
\otimes_{\mathcal{O}}A$. The blocks of $\mathcal{O}_\mathfrak{p}A$
are called the ``$\mathfrak{p}$-blocks of $A$''. If $\chi,\psi \in
\mathrm{Irr}(KA)$ belong to the same block of
$\mathcal{O}_\mathfrak{p}A$, we write $\chi \sim_\mathfrak{p} \psi$.

\begin{proposition}\label{p-blocks}
Two characters $\chi,\psi \in \mathrm{Irr}(KA)$ belong to the same
block of $A$ if and only if there exist a finite sequence
$\chi_0,\chi_1,\ldots,\chi_n \in \mathrm{Irr}(KA)$ and a finite
sequence $\mathfrak{p}_1,\ldots,\mathfrak{p}_n \in
\mathrm{Spec}_1(\mathcal{O})$ such that
\begin{itemize}
  \item $\chi_0=\chi$ and $\chi_n=\psi$,
  \item for all $j$ $(1\leq j \leq n)$, $\chi_{j-1}
  \sim_{\mathfrak{p}_j} \chi_j$.
\end{itemize}
\end{proposition}
\begin{apod}{Let us denote by $\sim$ the equivalence relation on
$\mathrm{Irr}(KA)$ defined as the transitive closure of the relation ``there
exists $\mathfrak{p} \in \mathrm{Spec}_1(\mathcal{O})$ such that $\chi
\sim_\mathfrak{p} \psi$''. Thus, we have to show that $\chi \sim
\psi$ if and only if $\chi$ and $\psi$ belong to the same block of
$A$.

We will first show that the equivalence relation $\sim$ is finer
than the relation ``being in the same block of A''. Let $B$ be a
block of $A$. Then $B$ is a subset of $\mathrm{Irr}(KA)$ such that
$\sum_{\chi \in B}e_\chi \in A$. Since
$\mathcal{O}=\bigcap_{\mathfrak{p} \in \mathrm{Spec}_1(\mathcal{O})}
\mathcal{O}_\mathfrak{p}$, we have that $\sum_{\chi \in B}e_\chi \in
\mathcal{O}_\mathfrak{p}A$ for all $\mathfrak{p} \in
\mathrm{Spec}_1(\mathcal{O})$. Therefore, by Theorem $\ref{minimality
of blocks}$, $C$ is a union of blocks of $\mathcal{O}_\mathfrak{p}A$
for all $\mathfrak{p} \in \mathrm{Spec}_1(\mathcal{O})$ and, hence, a
union of equivalence classes of $\sim$.

Now we will show that the relation ``being in the same block of A''
is finer than the relation $\sim$. Let $C$ be an equivalence class
of $\sim$. Then $\sum_{\chi \in C}e_\chi \in
\mathcal{O}_\mathfrak{p}A$ for all $\mathfrak{p} \in
\mathrm{Spec}_1(\mathcal{O})$. Thus $\sum_{\chi \in C}e_\chi \in
\bigcap_{\mathfrak{p} \in \mathrm{Spec}_1(\mathcal{O})}
\mathcal{O}_\mathfrak{p}A=A$ and $C$ is a union of blocks of $A$.}
\end{apod}\
\\
\emph{Remark:}  More generally, if we denote by
$\mathrm{Spec}(\mathcal{O})$ the set of prime ideals of
$\mathcal{O}$, then
$$\mathcal{O}=\bigcap_{\mathfrak{p} \in \mathrm{Spec}(\mathcal{O})}
\mathcal{O}_\mathfrak{p}.$$
Hence, Proposition $\ref{p-blocks}$ is also true, if we replace $\mathrm{Spec}_1(\mathcal{O})$ by $\mathrm{Spec}(\mathcal{O})$. However, in this case, $\mathcal{O}_\mathfrak{p}$ is not always a discrete valuation ring.\\

Thanks to Proposition $\ref{p-blocks}$, we can transfer the problem of the determination of the blocks of $A$ to that of the determination of the $\mathfrak{p}$-blocks of $A$, where $\mathfrak{p}$ runs over the set or prime ideals (of height $1$) of $\mathcal{O}$. 

\begin{px}\emph{ \small Let $C_6$ be the cyclic group of order $6$ and $A:=\mathbb{Z}[C_6]$.
The group $C_6$ (and thus the algebra $\mathbb{Q}(\zeta_6)[C_6]$) has $6$ irreducible characters $\chi_1,\ldots,\chi_6$, where
\begin{itemize}
\item $\chi_1$ is the trivial character,
\item $\chi_4(g) \in \mathbb{R}$ for all $g \in C_6$,
\item $\chi_2$ is the conjugate of $\chi_6$,
\item $\chi_3$ is the conjugate of $\chi_5$,
\item $\chi_2(g)\chi_3(g) =\chi_4(g)$ for all $g \in C_6$.
\end{itemize}
If $p$ is a prime number, then the localization of $\mathbb{Z}$ at $p\mathbb{Z}$ is a discrete valuation ring. By Brauer theory, the $p$-blocks of $C_6$ are trivial, unless $p$ divides the order of $C_6$.
The $2$-blocks of $C_6$ are $\{\chi_1,\chi_4\}$, $\{\chi_2,\chi_5\}$ and $\{\chi_3,\chi_6\}$, whereas
the  $3$-blocks of $C_6$ are $\{\chi_1,\chi_3,\chi_5\}$ and $\{\chi_2,\chi_4,\chi_6\}$.
Proposition $\ref{p-blocks}$ implies that all the irreducible characters of $C_6$ belong to the same
block of $A$.}
\end{px}

\subsection{Blocks and quotient blocks}

Let $\mathfrak{p}$ be a prime ideal of $\mathcal{O}$. Then $\mathcal{O}_\mathfrak{p}$ is a local ring, whose maximal ideal is $\mathfrak{p}':=\mathfrak{p}\mathcal{O}_\mathfrak{p}$. Let $\mathfrak{q}$ be a prime ideal of $\mathcal{O}$ such that $\mathfrak{q} \subseteq \mathfrak{p}$.
Then $\mathfrak{q'}:=\mathfrak{q}\mathcal{O}_\mathfrak{p}$ is a prime ideal of $\mathcal{O}_\mathfrak{p}$. The natural surjection
$\pi_\mathfrak{q}:\mathcal{O}_\mathfrak{p} \twoheadrightarrow \mathcal{O}_\mathfrak{p}/\mathfrak{q}'$ extends to a morphism 
$\pi_\mathfrak{q}:\mathcal{O}_\mathfrak{p}A \twoheadrightarrow (\mathcal{O}_\mathfrak{p}/\mathfrak{q}')A$, which in turn induces a morphism $$\pi_\mathfrak{q}:Z\mathcal{O}_\mathfrak{p}A \rightarrow Z(\mathcal{O}_\mathfrak{p}/\mathfrak{q}')A. $$ 
The following lemma will serve for the proof of Proposition $\ref{quotient blocks}$.

\begin{lemma}\label{central lifts to central}
Let  $e$ be an idempotent of $\mathcal{O}_\mathfrak{p}A$ whose image $\bar{e}$ in $(\mathcal{O}_\mathfrak{p}/\mathfrak{q}')A$ is central.  Then $e$ is central.
\end{lemma}
\begin{apod}{We have the following equality:
$$\mathcal{O}_\mathfrak{p}A=e\mathcal{O}_\mathfrak{p}Ae \oplus e\mathcal{O}_\mathfrak{p}A(1-e) \oplus (1-e)\mathcal{O}_\mathfrak{p}Ae \oplus (1-e)\mathcal{O}_\mathfrak{p}A(1-e).$$
Since $\bar{e}$ is central, we have
$$\bar{e}(\mathcal{O}_\mathfrak{p}/\mathfrak{q}')A(1-\bar{e})=(1-\bar{e})(\mathcal{O}_\mathfrak{p}/\mathfrak{q}')A\bar{e}=\{0\},$$ \ie
$$e\mathcal{O}_\mathfrak{p}A(1-e) \subseteq \mathfrak{q}\mathcal{O}_\mathfrak{p}A \,\,\textrm{
and }\,\, (1-e)\mathcal{O}_\mathfrak{p}Ae \subseteq \mathfrak{q}\mathcal{O}_\mathfrak{p}A.$$
Since $e$ and $(1-e)$ are idempotents, we get
$$e\mathcal{O}_\mathfrak{p}A(1-e) \subseteq \mathfrak{q}'e\mathcal{O}_\mathfrak{p}A(1-e) \,\,\textrm{ and
}\,\, (1-e)\mathcal{O}_\mathfrak{p}Ae \subseteq
\mathfrak{q}'(1-e)\mathcal{O}_\mathfrak{p}Ae.$$  
However, $\mathfrak{q}' $ is contained in the maximal ideal $\mathfrak{p}'$ of $\mathcal{O}_\mathfrak{p}$.
By Nakayama's lemma, we obtain that
$$e\mathcal{O}_\mathfrak{p}A(1-e)=(1-e)\mathcal{O}_\mathfrak{p}Ae=\{0\}.$$
We deduce that
$$\mathcal{O}_\mathfrak{p}A=e\mathcal{O}_\mathfrak{p}Ae \oplus (1-e)\mathcal{O}_\mathfrak{p}A(1-e)$$
and consequently, $e$ is central.}
\end{apod}
\begin{proposition}\label{quotient blocks}
If $K=F$, then the morphism 
$$\pi_\mathfrak{q}:Z\mathcal{O}_\mathfrak{p}A \rightarrow Z(\mathcal{O}_\mathfrak{p}/\mathfrak{q}')A$$ 
induces a bijection between the set of blocks of $\mathcal{O}_\mathfrak{p}A$ and the set of blocks of
$(\mathcal{O}_\mathfrak{p}/\mathfrak{q}')A$.
\end{proposition}
\begin{apod}{From now on, the symbol $\,\,\,\widehat{}\,\,\,$ will stand for
$\mathfrak{p}$-adic completion. It is immediate that $\pi_\mathfrak{q}$ sends a block of $\mathcal{O}_\mathfrak{p}A$  to a sum of blocks of $(\mathcal{O}_\mathfrak{p}/\mathfrak{q}')A$.
Now let $\bar{e}$ be a block of $(\mathcal{O}_\mathfrak{p}/\mathfrak{q}')A$.
Following Theorem $\ref{17.6}$, all Noetherian local rings are contained in their completions.
By Corollary $\ref{17.9}$, the completion of $\mathcal{O}_\mathfrak{p}/\mathfrak{q}'$ with respect to the $\mathfrak{p}$-adic topology is $\hat{\mathcal{O}_\mathfrak{p}}/\mathfrak{q}\hat{\mathcal{O}_\mathfrak{p}}$. Thus, $\bar{e} \in (\hat{\mathcal{O}_\mathfrak{p}}/\mathfrak{q}\hat{\mathcal{O}_\mathfrak{p}})A \cong \hat{\mathcal{O}_\mathfrak{p}}A/\mathfrak{q}\hat{\mathcal{O}_\mathfrak{p}}A$. In particular, $\bar{e} \in Z(\hat{\mathcal{O}_\mathfrak{p}}/\mathfrak{q}\hat{\mathcal{O}_\mathfrak{p}})A $.
By the theorems of lifting idempotents (cf., for example, \cite{The}, Theorem 3.2) and
the lemma above, $\bar{e}$ is lifted to a central
idempotent $e$ in $\hat{\mathcal{O}}_{\mathfrak{p}}A$.
Following Theorem $\ref{minimality of blocks}$, $e$ is a sum of blocks of  $\hat{\mathcal{O}}_{\mathfrak{p}}A$.
Since the algebra $KA$ is
split semisimple, the blocks of
$\hat{\mathcal{O}}_{\mathfrak{p}}A$ belong to
$KA$. By Theorem $\ref{18.4}$, $K \cap
\hat{\mathcal{O}}_{\mathfrak{p}}=\mathcal{O}_{\mathfrak{p}}$
and thus $e$ is a sum of primitive idempotents of $\mathcal{O}_{\mathfrak{p}}A$ which belong to  $Z \hat{\mathcal{O}}_{\mathfrak{p}}A$. Since
 $\mathcal{O}_{\mathfrak{p}}A \cap
Z \hat{\mathcal{O}}_{\mathfrak{p}}A = Z \mathcal{O}_{\mathfrak{p}}A$, we deduce that $\bar{e}$ is lifted to
a sum of blocks of $\mathcal{O}_{\mathfrak{p}}A$ and this provides
the block bijection.}
\end{apod}

The following well-known result on blocks is the application of our Proposition $\ref{quotient blocks}$ to the case $\mathfrak{q}=\mathfrak{p}$.
\begin{corollary}\label{residue field}
If $K=F$ and $k_{\mathfrak{p}}$ is the residue field of $\mathcal{O}_\mathfrak{p}$, then the morphism
$$\pi_\mathfrak{p}:Z\mathcal{O}_\mathfrak{p}A \rightarrow Zk_\mathfrak{p}A$$
induces a bijection between the set of blocks of $\mathcal{O}_\mathfrak{p}A$ and the set of blocks of
$k_\mathfrak{p}A$.
\end{corollary}

\subsection{Blocks and central characters}

Since $KA$ is a split semisimple $K$-algebra, we have that
$$KA \cong \prod_{\chi \in \mathrm{Irr}(KA)} M_\chi,$$
where $M_\chi$ is a matrix algebra isomorphic to
$\mathrm{Mat}_{\chi(1)}(K)$.

Recall that $A$ is of finite type and thus integral over
$\mathcal{O}$ (\cite{Bou5}, \S1, Def. 2). The map $\pi_\chi: KA
\twoheadrightarrow M_\chi$, restricted to $ZKA$, defines a map
$\omega_\chi:ZKA \twoheadrightarrow K$ (by Schur's lemma), which in
turn, restricted to $ZA$, defines the morphism
$$\omega_\chi: ZA \rightarrow \mathcal{O}_K,$$
where $\mathcal{O}_K$ denotes the integral closure of $\mathcal{O}$
in $K$.  \index{central character}

Now let $\mathfrak{p}$ be a prime ideal of $\mathcal{O}$. In the case where $\mathcal{O}_\mathfrak{p}$ is a discrete valuation ring (for example, when $\mathfrak{p}$ is a prime ideal of height $1$), we
have the following result which is proven later in this chapter,
Proposition $\ref{blocks and central characters}$. For a different
approach to its proof, see \cite{BK}, Proposition 1.18.

\begin{proposition}\label{omega_chi}
Suppose that $\mathcal{O}_\mathfrak{p}$ is a discrete valuation ring with unique
maximal ideal $\mathfrak{p}':=\mathfrak{p}\mathcal{O}_\mathfrak{p}$ and $K=F$. Two characters $\chi,\chi'
\in \mathrm{Irr}(KA)$ belong to the same block of $\mathcal{O}_\mathfrak{p}A$ if and only if
$$\omega_\chi(z) \equiv \omega_{\chi'}(z)\,\, \mathrm{mod} \, \mathfrak{p}'
\textrm{ for all } z \in Z\mathcal{O}_\mathfrak{p}A.$$
\end{proposition}

\section{Symmetric algebras}

Let $\mathcal{O}$ be a ring and let $A$ be an $\mathcal{O}$-algebra.
Suppose again that the assumptions $\ref{properties of the ring}$
are satisfied.

\begin{definition}\label{trace function}
A trace function  \index{trace function} on $A$ is an $\mathcal{O}$-linear map $t:A
\rightarrow \mathcal{O}$ such that $t(ab)=t(ba)$ for all $a,b \in
A$. 
\end{definition}

\begin{definition}\label{symmetric algebra}
We say that a trace function $t:A \rightarrow \mathcal{O}$ is a
symmetrizing form  \index{symmetrizing form} on $A$ or that $A$ is a symmetric algebra  
\index{symmetric algebra} if the
morphism
$$\hat{t}:A \rightarrow \mathrm{Hom}_\mathcal{O}(A,\mathcal{O}),\,\,
  a \mapsto (x \mapsto \hat{t}(a)(x):=t(ax))$$
is an isomorphism of $A$-modules-$A$.
\end{definition}

\begin{px}\label{symmetrizing form of the group algebra}
\small{\emph{In the case where $\mathcal{O}=\mathbb{Z}$ and
$A=\mathbb{Z}[G]$
 ($G$ a finite group), we can define the following symmetrizing form
 (``canonical'')
 on $A$
$$t:\mathbb{Z}[G] \rightarrow \mathbb{Z}, \,\,\, \sum_{g \in G}a_g g \mapsto a_1,$$
where $a_g \in \mathbb{Z}$ for all $g \in G$.}}
\end{px}

Since $A$ is a free $\mathcal{O}$-module of finite rank, we have the
following isomorphism
$$\begin{array}{ccc}
  \mathrm{Hom}_\mathcal{O}(A,\mathcal{O}) \otimes_\mathcal{O} A &\tilde{\rightarrow} &\mathrm{Hom}_\mathcal{O}(A,A)\\
  \varphi \otimes a & \mapsto & (x \mapsto \varphi(x)a).
\end{array}$$
Composing it with the isomorphism
$$\begin{array}{ccc}
  A \otimes_\mathcal{O} A &\tilde{\rightarrow} &\mathrm{Hom}_\mathcal{O}(A,\mathcal{O})\otimes_\mathcal{O} A\\
  a \otimes b & \mapsto & \hat{t}(a) \otimes b,
\end{array}$$
we obtain an isomorphism
$$ A \otimes_\mathcal{O} A \tilde{\rightarrow} \mathrm{Hom}_\mathcal{O}(A,A).$$

\begin{definition}\label{casimir}
We denote by $C_A$ and we call Casimir  \index{Casimir element} of $(A,t)$ the inverse image
of $\mathrm{id}_A$ by the above isomorphism.
\end{definition}

\begin{px}\label{casimir of the group algebra}
\small{\emph{In the case where $\mathcal{O}=\mathbb{Z}$,
$A=\mathbb{Z}[G]$ ($G$ a finite group) and $t$ is the canonical
symmetrizing form, we have $C_{\mathbb{Z}[G]}=\sum_{g \in
G}g^{-1}\otimes g$.}}
\end{px}

More generally, if $(e_i)_{i \in I}$ is a basis of $A$ over
$\mathcal{O}$ and $(e_i')_{i \in I}$ is the dual basis with respect
to $t$ (\ie $t(e_ie_j')=\delta_{ij}$), then
$$C_A = \sum_{i \in I} e_i' \otimes e_i.$$

In this case, let us denote by $c_A$ the image of $C_A$ by the
multiplication $A \otimes A \rightarrow A$, \ie $c_A= \sum_{i \in I}
e_i'e_i$. It is easy to check (see also \cite{BMM2}, 7.9) the
following properties of the Casimir element:

\begin{lemma}\label{properties of the casimir}
For all $a \in A$, we have
\begin{enumerate}
  \item $\sum_i ae_i' \otimes e_i = \sum_i e_i \otimes e_i'a$.
  \item $aC_A=C_Aa$. Consequently, $c_A \in ZA$.
  \item $a=\sum_i t(ae_i')e_i=\sum_i t(ae_i)e_i'=
  \sum_i t(e_i')e_ia=\sum_i t(e_i)e_i'a.$
\end{enumerate}
\end{lemma}

If $\tau:A \rightarrow \mathcal{O}$ is a linear form, we denote by
$\tau^\vee$ its inverse image by the isomorphism $\hat{t}$, \ie
$\tau^\vee$ is the element of $A$ such that
$$t(\tau^\vee a)=\tau(a) \textrm{ for all } a \in A.$$

The element $\tau^\vee$ has the following properties:

\begin{lemma}\label{tau^vee}\
\begin{enumerate}
 \item $\tau$ is a trace function if and only if $\tau^\vee \in ZA$.
 \item We have $\tau^\vee=\sum_i \tau(e_i')e_i=\sum_i \tau(e_i)e_i'$ and
 more generally, for all $a \in A$, we have
 $\tau^\vee a=\sum_i \tau(e_i'a)e_i=\sum_i \tau(e_ia)e_i'$.
\end{enumerate}
\end{lemma}
\begin{apod}{\begin{enumerate}
  \item   Recall that $t$ is a trace function. Let $a \in A$. For all $x \in A$, we have
  $$\hat{t}(\tau^\vee a)(x) = t(\tau^\vee a x)=\tau(ax)$$
  and
  $$\hat{t}(a \tau^\vee )(x) = t(a \tau^\vee x)=t(\tau^\vee
  xa)=\tau(xa)$$
  If $\tau$ is a trace function, then $\tau(ax)=\tau(xa)$ and hence, $\hat{t}(\tau^\vee
  a)= \hat{t}(a \tau^\vee)$. Since $\hat{t}$ is an isomorphism, we
  obtain that $ \tau^\vee a = a \tau^\vee$ and thus $\tau^\vee \in ZA$.

  Now if $\tau^\vee \in ZA$ and $a,b \in A$, then $$\tau(ab)=t(\tau^\vee ab)
  =t(b \tau^\vee a) =t(ba \tau^\vee )=t(\tau^\vee ba)=\tau(ba).$$
  \item It derives from property 3 of Lemma $\ref{properties of the
  casimir}$ and the definition of $\tau^\vee$.}
\end{enumerate}
\end{apod}

Let $\chi_\mathrm{reg}$ be the character of the regular
representation of $A$, \ie the linear form on $A$ defined as
$$\chi_\mathrm{reg}(a):=\mathrm{tr}_{A/\mathcal{O}}(\lambda_a),$$
where $\lambda_a:A\rightarrow A, x \mapsto ax$, is the endomorphism
of left multiplication by $a$.

\begin{proposition}\label{regular^vee}
We have $\chi_\mathrm{reg}^\vee = c_A.$
\end{proposition}
\begin{apod}{Let $a \in A$. The inverse image of $\lambda_a$ by the
isomorphism \\
$A \otimes_\mathcal{O} A \tilde{\rightarrow}
\mathrm{Hom}_\mathcal{O}(A,A)$ is $aC_A$ (by definition of the
Casimir). Hence,
$$\lambda_a = (x \mapsto \sum_i \hat{t}(e_i'a)(x)e_i) = (x \mapsto \sum_i
t(e_i'ax)e_i)$$ and
$$\mathrm{tr}_{A/\mathcal{O}}(\lambda_a)=\sum_i t(e_i'ae_i)=
t(a\sum_i e_i'e_i)=t(ac_A)=t(c_Aa).$$ Therefore, for all $a \in A$,
we have $\chi_\mathrm{reg}(a)=t(c_Aa)$, \ie $\chi_\mathrm{reg}^\vee
= c_A.$}
\end{apod}

If $A$ is a symmetric algebra with a symmetrizing form $t$, we
obtain a symmetrizing form $t^K$ on $KA$ by extension of scalars.
Every irreducible character $\chi \in \mathrm{Irr}(KA)$ is a trace
function on $KA$ and thus we can define $\chi^\vee \in ZKA$.

\begin{definition}\label{Schur element}
For all $\chi \in \mathrm{Irr}(KA)$, we call Schur element of $\chi$  \index{Schur element}
with respect to $t$ and denote by $s_\chi$ the element of $K$
defined by $$s_\chi:=\omega_\chi(\chi^\vee).$$
\end{definition}

\begin{proposition}\label{Schur element belongs to the integral closure}
For all $\chi \in \mathrm{Irr}(KA)$, $s_\chi \in \mathcal{O}_K$.
\end{proposition}

The proof of the above result will be given in Proposition
$\ref{integrality of the Schur elements}$.

\begin{px}\label{Schur elements of the group algebra}
\small{\emph{Let $\mathcal{O}:=\mathbb{Z}$, $A:=\mathbb{Z}[G]$
 ($G$ a finite group) and $t$ the canonical symmetrizing form. If $K$ is an algebraically closed field of
 characteristic 0, then $KA$ is a split semisimple algebra and
 $s_\chi=|G|/\chi(1)$ for all $\chi \in \mathrm{Irr}(KA)$. Because
 of the integrality of the Schur elements, we must have
 $|G|/\chi(1) \in \mathbb{Z}=\mathbb{Z}_K \cap \mathbb{Q}$ for all $\chi \in
 \mathrm{Irr}(KA)$. Thus, we have shown that $\chi(1)$ divides $|G|$.}}
\end{px}

The following properties of the Schur elements can be derived easily
from the above (see also
\cite{Bro},\cite{Ge},\cite{GePf},\cite{GeRo},\cite{BMM2})

\begin{proposition}\label{schur elements and idempotents}\
\begin{enumerate}
  \item We have
  $$t=\sum_{\chi \in \mathrm{Irr}(KA)}\frac{1}{s_\chi}\chi.$$
  \item For all $\chi \in \mathrm{Irr}(KA)$, the central primitive
  idempotent associated with $\chi$ is
  $$e_\chi=\frac{1}{s_\chi}\chi^\vee=\frac{1}{s_\chi}\sum_{i \in I}
  \chi(e_i')e_i.$$
  \item For all $\chi \in \mathrm{Irr}(KA)$, we have
  $$s_\chi \chi(1)=\sum_{i \in I}\chi(e_i')\chi(e_i) \,\textrm{ and  }\,
  s_\chi \chi(1)^2= \chi(\sum_{i \in I}e_i'e_i)=\chi(\chi_{\mathrm{reg}}^\vee).$$
\end{enumerate}
\end{proposition}

\begin{corollary}\label{what we are searching}
The blocks of $A$ are the non-empty subsets $B$ of
$\mathrm{Irr}(KA)$ minimal with respect to the property
$$\sum_{\chi \in B}\frac{1}{s_\chi}\chi(a) \in \mathcal{O} \textrm{ for all } a \in A.$$
\end{corollary}

\section{Twisted symmetric algebras of finite groups}

This part is an adaptation of the section ``Symmetric algebras of
finite groups'' of \cite{BK} to a more general case.\\

Let $A$ be an $\mathcal{O}$-algebra such that the assumptions
$\ref{properties of the ring}$ are satisfied with a symmetrizing
form $t$. Let $\bar{A}$ be a subalgebra of $A$ free and of finite
rank as $\mathcal{O}$-module.

We denote by $\bar{A}^\bot$ the orthogonal of $\bar{A}$ with respect
to $t$, \ie the sub-$\bar{A}$-module-$\bar{A}$ of $A$ defined as
$$\bar{A}^\bot:=\{a \in A \,|\, (\forall \bar{a} \in \bar{A})(t(a\bar{a})=0)\}.$$

\begin{proposition}\label{when symmetric subalgebra}\
\begin{enumerate}
  \item The restriction of $t$ to $\bar{A}$ is a symmetrizing form for $\bar{A}$
  if and only if $\bar{A} \oplus \bar{A}^\bot=A$.
   In this case the projection of $A$ onto $\bar{A}$ parallel to $\bar{A}^\bot$ is
  the map
  $$\mathrm{Br}_{\bar{A}}^A :A \rightarrow \bar{A}
  \textrm{ such that } t(\mathrm{Br}_{\bar{A}}^A(a)\bar{a})=t(a\bar{a})
  \textrm{ for all } a \in A \textrm{ and } \bar{a} \in \bar{A}.$$
  \item If the restriction of $t$ to $\bar{A}$ is a symmetrizing form
  for $\bar{A}$, then $\bar{A}^\bot$ is the sub-$\bar{A}$-module-$\bar{A}$
  of $A$ defined by the following two properties:
  \begin{enumerate}[(a)]
    \item $A=\bar{A} \oplus \bar{A}^\bot$,
    \item $\bar{A}^\bot \subseteq \mathrm{Ker}t$.
  \end{enumerate}
\end{enumerate}
\end{proposition}
\begin{apod}{
\begin{enumerate}
  \item Let us denote by $\bar{t}$ the restriction of $t$ to
  $\bar{A}$. Suppose that $\bar{t}$ is a symmetrizing form for
  $\bar{A}$. Let $a \in A$. Then $\hat{t}(a):=(x \mapsto t(ax)) \in
  \mathrm{Hom}_\mathcal{O}(A,\mathcal{O})$. The restriction of
  $\hat{t}(a)$ to $\bar{A}$ belongs to
  $\mathrm{Hom}_\mathcal{O}(\bar{A},\mathcal{O})$ and therefore,
  there exists $\bar{a} \in \bar{A}$ such that $\bar{t}(\bar{a}\bar{x})=t(\bar{a}\bar{x})=t(a\bar{x})$
  for all $\bar{x} \in \bar{A}$. Thus $a-\bar{a} \in \bar{A}^\bot$
  and since $a=\bar{a}+(a-\bar{a})$, we obtain that $A=\bar{A}+
  \bar{A}^\bot$. If $\bar{a} \in \bar{A} \cap
  \bar{A}^\bot$, then we have $\hat{\bar{t}}(\bar{a})=0 \in
  \mathrm{Hom}_\mathcal{O}(\bar{A},\mathcal{O})$.
  Since $\hat{\bar{t}}$ is an isomorphism, we deduce that
  $\bar{a}=0$. Therefore, $A=\bar{A} \oplus \bar{A}^\bot$ and the
  definition of $\mathrm{Br}_{\bar{A}}^A$ is immediate.

  Now suppose that $A=\bar{A} \oplus \bar{A}^\bot$. We will show that the map
  $$\begin{array}{cccl}
      \hat{\bar{t}}: & \bar{A} & \rightarrow & \mathrm{Hom}_\mathcal{O}(\bar{A},\mathcal{O}) \\
                    & \bar{A} & \mapsto     & (\bar{x} \mapsto \bar{t}(\bar{a}\bar{x})=t(\bar{a}\bar{x}))
    \end{array}$$
  is an isomorphism of $\bar{A}$-modules-$\bar{A}$. The map
  $\hat{\bar{t}}$ is obviously injective, because
  $\hat{\bar{t}}(\bar{a})=0$ implies that $\bar{a} \in \bar{A} \cap \bar{A}^\bot$
  and thus $\bar{a}=0$. Now let $\bar{f}$ be an element of
  $\mathrm{Hom}_\mathcal{O}(\bar{A},\mathcal{O})$. The map $\bar{f}$ can be
  extended to a map $f \in \mathrm{Hom}_\mathcal{O}(A,\mathcal{O})$ such that
  $f(a)=\bar{f}(\mathrm{Br}_{\bar{A}}^A(a))$ for all $a \in A$, where
  $\mathrm{Br}_{\bar{A}}^A$ denotes the projection of $A$ onto $\bar{A}$ parallel to
  $\bar{A}^\bot$. Since $t$ is a symmetrizing form for $A$, there
  exists $a \in A$ such that $\hat{t}(a)=f$, \ie $t(ax)=f(x)$ for
  all $x \in A$. Consequently, if $\bar{x} \in \bar{A}$, we have
  $$t(\mathrm{Br}_{\bar{A}}^A(a)\bar{x})=t(a\bar{x})=f(\bar{x})=\bar{f}(\bar{x})$$
  and thus $\hat{\bar{t}}(\mathrm{Br}_{\bar{A}}^A(a))=\bar{f}$. Hence,
  $\hat{\bar{t}}$ is surjective.
  \item Let $B$ be a sub-$\bar{A}$-module-$\bar{A}$
  of $A$ such that $A=\bar{A} \oplus B$ and $B \subseteq \mathrm{Ker}t$.
  Let $b \in B$. For all $\bar{a} \in
  \bar{A}$, we have $b\bar{a} \in B \subseteq
  \mathrm{Ker}t$ and thus $t(b\bar{a})=0$. Hence $B \subseteq \bar{A}^\bot$.
  Since the restriction of $t$ to $\bar{A}$ is a symmetrizing form
  for $\bar{A}$, we also have $A=\bar{A} \oplus \bar{A}^\bot$. Now
  let $a \in \bar{A}^\bot$. Since $A=\bar{A} \oplus B$, there exist $\bar{a} \in \bar{A}$
  and $b \in B$ such that $a=\bar{a}+b$. Since $b \in \bar{A}^\bot$,
  we must have $a=b \in B$ and therefore, $B=\bar{A}^\bot$.}
\end{enumerate}
\end{apod}

\begin{px}\label{subalgebra of group algebra}
\small{\emph{In the case where $\mathcal{O}=\mathbb{Z}$ and
$A=\mathbb{Z}[G]$
 (G a finite group), let $\bar{A}:=\mathbb{Z}[\bar{G}]$
 be the algebra of a subgroup $\bar{G}$ of $G$.
 Then the morphism $\mathrm{Br}_{\bar{A}}^A$ is the projection given by
 $$\left\{
     \begin{array}{ll}
       g \mapsto g, & \hbox{if $g \in \bar{G}$;} \\
       g \mapsto 0, & \hbox{if $g \notin \bar{G}$.}
     \end{array}
   \right.$$}}
\end{px}\

\begin{definition}\label{symmetric subalgebra}
Let $A$ be a symmetric $\mathcal{O}$-algebra with symmetrizing form
$t$. Let $\bar{A}$ be a subalgebra of $A$. We say that $\bar{A}$ is
a symmetric subalgebra of $A$,  \index{symmetric subalgebra} if it satisfies the following two
conditions:
\begin{enumerate}
  \item $\bar{A}$ is free (of finite rank) as an $\mathcal{O}$-module and the
  restriction $\mathrm{Res}_{\bar{A}}^A(t)$ of the form $t$ to $\bar{A}$ is a symmetrizing form
  on $\bar{A}$,
  \item $A$ is free (of finite rank) as an $\bar{A}$-module for the action
  of left multiplication by the elements of $\bar{A}$.
\end{enumerate}
\end{definition}

From now on, let us suppose that $\bar{A}$ is a symmetric subalgebra
of $A$ and set $\bar{t}:=\mathrm{Res}_{\bar{A}}^A(t)$. We denote by
$$\mathrm{Ind}_{\bar{A}}^A: _{\bar{A}}\mathrm{\textbf{mod}} \rightarrow _A\mathrm{\textbf{mod}}
\,\textrm{ and }\, \mathrm{Res}_{\bar{A}}^A: _A\mathrm{\textbf{mod}}
\rightarrow _{\bar{A}}\mathrm{\textbf{mod}} $$ the functors defined
as usual by
$$\mathrm{Ind}_{\bar{A}}^A:=A \otimes_{\bar{A}}- \textrm{ where $A$ is viewed as an $A$-module-$\bar{A}$}$$
and
$$\mathrm{Res}_{\bar{A}}^A:=A \otimes_A - \textrm{ where $A$ is viewed as an $\bar{A}$-module-$A$.}$$
Since $A$ is free as $\bar{A}$-module and as module-$\bar{A}$, the
functors $\mathrm{Res}_{\bar{A}}^A$ and $\mathrm{Ind}_{\bar{A}}^A$
are adjoint from both sides.

Moreover, let $K$ be a finite Galois extension of the field of fractions of $\mathcal{O}$ such that the algebras $KA$ and
$K\bar{A}$ are both split semisimple.

We denote by $\langle-,-\rangle_{KA}$ the scalar product on the
$K$-vector space of trace functions for which the family
$(\chi)_{\chi \in \mathrm{Irr}(KA)}$ is orthonormal and
$\langle-,-\rangle_{K\bar{A}}$ the scalar product on the $K$-vector
space of trace functions for which the family
$(\bar{\chi})_{\bar{\chi} \in \mathrm{Irr}(K\bar{A})}$ is
orthonormal.

Since the functors $\mathrm{Res}_{\bar{A}}^A$ and
$\mathrm{Ind}_{\bar{A}}^A$ are adjoint from both sides, we obtain
the \emph{Frobenius reciprocity} formula:
$$\langle \chi, \mathrm{Ind}_{K\bar{A}}^{KA}(\bar{\chi}) \rangle_{KA} =
\langle \mathrm{Res}_{K\bar{A}}^{KA}(\chi), \bar{\chi}
\rangle_{K\bar{A}} .$$

For every element $\chi \in \mathrm{Irr}(KA)$, let
$$\mathrm{Res}_{K\bar{A}}^{KA}(\chi)=\sum_{\bar{\chi} \in
\mathrm{Irr}(K\bar{A})}m_{\chi,\bar{\chi}}\bar{\chi} \textrm{ (where
} m_{\chi,\bar{\chi}} \in \mathbb{N}).$$

Frobenius reciprocity implies that, for all $\bar{\chi} \in
\mathrm{Irr}(K\bar{A})$,

$$\mathrm{Ind}_{K\bar{A}}^{KA}(\bar{\chi})=\sum_{\chi \in
\mathrm{Irr}(KA)}m_{\chi,\bar{\chi}}\chi .$$

The following property is immediate.

\begin{lemma}\label{mxx'}
For $\chi \in \mathrm{Irr}(KA)$ and $\bar{\chi} \in
\mathrm{Irr}(K\bar{A})$, let $e(\chi)$ and $\bar{e}(\bar{\chi})$ be
respectively the block-idempotents of $KA$ and $K\bar{A}$ associated
with $\chi$ and $\bar{\chi}$. The following conditions are
equivalent:
\begin{enumerate}[(i)]
  \item  $m_{\chi,\bar{\chi}} \neq 0$,
  \item  $e(\chi)\bar{e}(\bar{\chi}) \neq 0$.
\end{enumerate}
\end{lemma}

For all $\bar{\chi} \in \mathrm{Irr}(K\bar{A})$, we set
$$\mathrm{Irr}(KA,\bar{\chi}):=\{\chi \in \mathrm{Irr}(KA) \,|\, m_{\chi,\bar{\chi}} \neq 0\},$$

and for all $\chi \in \mathrm{Irr}(KA)$,
$$\mathrm{Irr}(K\bar{A},\chi):=\{\bar{\chi} \in \mathrm{Irr}(K\bar{A}) \,|\, m_{\chi,\bar{\chi}} \neq 0\}.$$

We denote respectively by $s_\chi$ and $s_{\bar{\chi}}$ the Schur
elements of $\chi$ and $\bar{\chi}$ ( with respect to the
symmetrizing forms $t$ for $A$ and $\bar{t}$ for $\bar{A}$).

\begin{lemma}\label{induction and Schur elements}
For all $\bar{\chi} \in \mathrm{Irr}(K\bar{A})$ we have
$$\frac{1}{s_{\bar{\chi}}}=\sum_{\chi \in
\mathrm{Irr}(KA,\bar{\chi})}\frac{m_{\chi,\bar{\chi}}}{s_\chi}  .$$
\end{lemma}
\begin{apod}
{It derives from the relations
$$t=\sum_{\chi \in \mathrm{Irr}(KA)}\frac{1}{s_\chi}\chi,\,\,
\bar{t}=\sum_{\bar{\chi} \in
\mathrm{Irr}(K\bar{A})}\frac{1}{s_{\bar{\chi}}}\bar{\chi},\,\,
\bar{t}= \mathrm{Res}_{\bar{A}}^{A}(t).$$}
\end{apod}

In the next chapters, we will work on the Hecke algebras of complex
reflection groups, which, under certain assumptions, are symmetric.
Sometimes the Hecke algebra of a group $W$ appears as a symmetric
subalgebra of the Hecke algebra of another group $W'$, which
contains $W$. Since we will be mostly interested in the
determination of the blocks of these algebras, it would be helpful,
if we could obtain the blocks of the former from the blocks of the
latter or vice versa. This is possible with the use of a generalization of some
classical results known as ``Clifford theory''.

``Clifford theory'' was originally developed by A. H. Clifford in 1937 for finite group algebras over a field
(cf.~\cite{Clif}): if $G$ is a finite group, $F$ is a field and $N$ is a normal subgroup of $G$, then important information concerning simple and indecomposable $KG$-modules can be obtained by applying (perhaps repeatedly) three basic operations: (a) restriction to $FN$, (b) extension from $FN$, (c) induction from $FN$. In the past twenty years, the theory has enjoyed vigorous development. The foundations have been strengthened and reorganized from new points of view and Clifford's results have been generalized to various types of algebras (cf., for example, \cite{Da}). Here we will generalize these results to the case of, what we are about to define as, the twisted symmetric algebras of finite groups and in particular, of finite cyclic groups.

\begin{definition}\label{symmetric algebra of a finite group}
Let $A$ be a symmetric $\mathcal{O}$-algebra with symmetrizing form
$t$. We say that $A$ is the twisted
symmetric algebra of a finite group $G$  \index{twisted symmetric algebra} over the subalgebra
$\bar{A}$, if the following conditions are satisfied:
\begin{enumerate}
  \item $\bar{A}$ is a symmetric subalgebra of $A$.
  \item There exists a family $\{A_g \,|\, g \in G\}$ of
  $\mathcal{O}$-submodules of $A$ such that
  \begin{enumerate}
    \item $A= \bigoplus_{g \in G}A_g$,
    \item $A_gA_h=A_{gh}$ for all $g,h \in G$,
    \item $A_1=\bar{A}$,
    \item $t(A_g)=0$ for all $g \in G\setminus \{1\}$.
    \end{enumerate}
\end{enumerate}
\end{definition}
If that is the case, then Proposition $\ref{when symmetric
subalgebra}$ implies that
$$\bigoplus_{g \in G\setminus\{1\}} A_g =
\bar{A}^\bot.$$

\begin{proposition}\label{factor group}
Let $A$ be a symmetric $\mathcal{O}$-algebra with symmetrizing form
$t$. The algebra $A$ is the twisted
symmetric algebra of a finite group $G$ over the subalgebra
$\bar{A}$ if and only if the following conditions are satisfied:
\begin{enumerate}
  \item $\bar{A}$ is a symmetric subalgebra of $A$.
  \item There exists a family $\{a_g \,|\, g \in G\}$ of
  invertible elements of $A$ such that
  \begin{enumerate}
    \item $A= \bigoplus_{g \in G}a_g\bar{A}$,
    \item $a_g\bar{A}=\bar{A}a_g$ for all $g \in G$,
    \item $a_ga_h\bar{A}=a_{gh}\bar{A}$ for all $g,h \in G$,
    \item $a_1=1$,
    \item $a_g \in \bar{A}^\bot$ for all $g \in G\setminus \{1\}$.
 \end{enumerate}
\end{enumerate}
\end{proposition}
\begin{apod}{If $A$ is the twisted
symmetric algebra of a finite group $G$ over the subalgebra
$\bar{A}$, then $\bar{A}$ is a symmetric subalgebra of $A$ and
there exists a family $\{A_g \,|\, g \in G\}$ of
  $\mathcal{O}$-submodules of $A$ which satisfy the conditions (a)--(d) of Definition
  $\ref{symmetric algebra of a finite group}$.
  Let $g \in G$. Then, by property (b), $A_gA_{g^{-1}}=A_1=\bar{A}$. Since $1 \in \bar{A}$, there exists 
  $a_g \in A_g$ and $a_{g^{-1}} \in A_{g^{-1}}$ such that $a_ga_{g^{-1}}=1$.   Now, if $x \in A_g$, then
  $x=a_ga_g^{-1}x \in a_gA_{g^{-1}}A_g=a_g\bar{A}$ and thus, $A_g\subseteq a_g\bar{A}$. Property (b) implies the inverse inclusion. In the same way, we show that $A_g=\bar{A}a_g$. Hence, we obtain that there exists a family $\{a_g \,|\, g \in G\}$ of invertible elements of $A$ such that
 \begin{center}
 $A= \bigoplus_{g \in G}a_g\bar{A}$ \,and\, $a_g\bar{A}=\bar{A}a_g$\, for all $g \in G$.
 \end{center}
 We can choose $a_1:=1$. Moreover, for all $g,h \in G$, we have $a_ga_h \in A_gA_h=A_{gh}$. Since $a_ga_h$ is a unit in $A$, we obtain that $A_{gh}=a_ga_h\bar{A}$, \ie $a_{gh}\bar{A}=a_ga_h\bar{A}$.
 Finally, for all $g \in G\setminus\{1\}$ and $\bar{a} \in \bar{A}$, we have  $a_g\bar{a} \in a_g\bar{A}=A_g$ and property (d) implies that $t(a_g\bar{a})=0$. Thus, $a_g \in \bar{A}^\bot$ for all $g \in G\setminus \{1\}$.

Now, let us suppose that the conditions $1$ and $2$ are satisfied. For all $g \in G$, we set $A_g:=a_g\bar{A}=\bar{A}a_g$. We only need to show that the family $\{A_g \,|\, g \in G\}$ satisfies the conditions (b), (c) and (d) of Definition  $\ref{symmetric algebra of a finite group}$. Obviously, $A_1=a_1\bar{A}=1\bar{A}=\bar{A}$. Moreover, if $g \in G\setminus \{1\}$ and $x \in A_g$, then there exists $\bar{a} \in \bar{A}$ such that $x=a_g\bar{a}$. Since $a_g \in \bar{A}^\bot$, we obtain that $t(x)=t(a_g\bar{a})=0$. Hence, $t(A_g)=0$ for all $g \in G\setminus \{1\}$. It remains to show that  $A_gA_h=A_{gh}$ for all $g,h \in G$. We have $A_gA_h=(a_g\bar{A})(a_h\bar{A}) = a_g(\bar{A}a_h)\bar{A}=
a_g(a_h\bar{A})\bar{A}=a_ga_h\bar{A}=a_{gh}\bar{A}=A_{gh}$.}
\end{apod}

From now on, let $(A,t)$ be the twisted symmetric algebra of a
finite group $G$ over the subalgebra $\bar{A}$. Due to the proposition above, 
for all $g \in G$, there exists $a_g \in A_g
\cap A^\times$ such that $A_g=a_g\bar{A}=\bar{A}a_g$. We fix a system of representatives
$\mathrm{Rep}(A/\bar{A}):=\{a_g \,|\, g \in G\}.$

\begin{proposition}\label{dual basis of a symmetric algebra of a finite group}
Let $(\bar{e}_i)_{i \in I}$ be a basis of $\bar{A}$ over
$\mathcal{O}$ and $(\bar{e}'_i)_{i \in I}$ its dual with respect to
the symmetrizing form $\bar{t}$.  Then the families
$$ (\bar{e}_ia_g)_{i \in I, a_g \in \mathrm{Rep}(A/\bar{A})}
 \textrm{ and } (a_g^{-1}\bar{e}'_i)_{i \in I,a_g \in \mathrm{Rep}(A/\bar{A})}$$
are two $\mathcal{O}$-bases of $A$ dual to each other.
\end{proposition}

\subsection{Action of $G$ on $Z\bar{A}$}

\begin{lemma}\label{definition of ga}
Let $\bar{a} \in Z\bar{A}$ and $g \in G$. There exists a unique
element $g(\bar{a})$ of $\bar{A}$ satisfying
$$g(\bar{a})x_g=x_g\bar{a} \textrm{ for all } x_g \in
A_g.\,\,\,\,\,\,\, (\dag)$$ In particular, $g(\bar{a})=a_g\bar{a}a_g^{-1}.$
\end{lemma}
\begin{apod}{
 For all $x_g \in A_g$, we have
$a_g^{-1}x_g \in
  \bar{A}$. Set $g(\bar{a}):=a_g\bar{a}a_g^{-1}.$ 
  Since $\bar{a} \in Z\bar{A}$, we obtain that  $g(\bar{a})x_g=a_g\bar{a}a_g^{-1}x_g=a_ga_g^{-1}x_g\bar{a}=x_g\bar{a}$.
Now, let $y$ be another element of $A$ such that $y
x_g=x_g\bar{a} \textrm{ for all } x_g \in
A_g.$ Then $ya_g=a_g\bar{a}$, whence $y=g(\bar{a})$. Therefore,
$g(\bar{a})$ is the unique element of $A$ which satisfies $(\dag)$.}
\end{apod}$ $\\
\begin{remark} \emph{Note that $g(\bar{a})$ does not depend on the choice of $a_g$.}
\end{remark}
\begin{proposition}\label{action of G on ZA}
The map $\bar{a} \mapsto g(\bar{a})$ defines an action of $G$ as
ring automorphism of $Z\bar{A}$.
\end{proposition}
\begin{apod}{Let $\bar{a} \in Z\bar{A}$ and $g \in G$. We will show that $g(\bar{a}) \in Z\bar{A}$. If
$\bar{x} \in \bar{A}$, then $\bar{x}a_g \in A_g$ and we have
$$\bar{x}(g(\bar{a})a_g)=\bar{x}(a_g\bar{a})= (\bar{x}a_g)\bar{a}=g(\bar{a})(\bar{x}a_g).$$
Multiplying both sides by $a_g^{-1}$, we obtain that
$$\bar{x}g(\bar{a})=g(\bar{a})\bar{x}$$
and hence, $g(\bar{a}) \in Z\bar{A}$.

Since $a_1=1$, we have
$1_G(\bar{a})=\bar{a}$. If $g_1,g_2 \in G$, then equation ($\dag$)
gives
$$g_1(g_2(\bar{a}))a_{g_1}a_{g_2}=a_{g_1}g_2(\bar{a})a_{g_2}=
a_{g_1}a_{g_2}\bar{a}.$$ Due to property (b) of  Definition
$\ref{symmetric algebra of a finite group}$, the product
$a_{g_1}a_{g_2}$ generates the submodule $A_{g_1g_2}$. Therefore,
$g_1(g_2(\bar{a}))u=u\bar{a}$ for all $u \in A_{g_1g_2}$. By Lemma
$\ref{definition of ga}$, we obtain that
$(g_1g_2)(\bar{a})=g_1(g_2(\bar{a})).$

Finally, let us fix $g \in G$. By definition, the map $\bar{a}
\mapsto g(\bar{a})$ is an additive automorphism of $Z\bar{A}$. If
$\bar{a}_1,\bar{a}_2 \in Z\bar{A}$, then
$$x_g\bar{a}_1\bar{a}_2=g(\bar{a}_2)x_g\bar{a}_1=g(\bar{a}_1)g(\bar{a}_2)x_g
\textrm{ for all } x_g \in A_g.$$ By Lemma $\ref{definition
of ga}$, we obtain that
$g(\bar{a}_1\bar{a}_2)=g(\bar{a}_1)g(\bar{a}_2).$}
\end{apod}

Now let $\bar{b}$ be a block(-idempotent) of $\bar{A}$. If $g \in
G$, then $g(\bar{b})$ is also a block of $\bar{A}$. So we must have
either $g(\bar{b})=\bar{b}$ or $g(\bar{b})$ orthogonal to $\bar{b}$.
Set
$$\mathrm{Tr}(G,\bar{b}):=\sum_{g \in G/G_{\bar{b}}}g(\bar{b}),$$
where $G_{\bar{b}}:=\{g \in G \,|\, g(\bar{b})=\bar{b}\}$. It is
clear that
\begin{itemize}
  \item $\bar{b}$ is a central idempotent of $\bigoplus_{g \in G_{\bar{b}}}A_g=:A_{G_{\bar{b}}}$,
  \item $\mathrm{Tr}(G,\bar{b})$ is a central idempotent of $A$.
\end{itemize}
From now on, let $b:=\mathrm{Tr}(G,\bar{b})$ and $x_g \in
A_g$. We have
\begin{itemize}
\item $\bar{b}x_g\bar{b}= \left\{
                           \begin{array}{ll}
                             x_g\bar{b}=\bar{b}x_g, & \hbox{if $g \in G_{\bar{b}}$;} \\
                             0, & \hbox{if $g \notin G_{\bar{b}}$,}
                           \end{array}
                         \right.$
\item $\bar{b}x_gb=\bar{b}x_g$ and $bx_g\bar{b}=x_g\bar{b}.$
\end{itemize}

\begin{proposition}\label{morita}
The applications
$$\left\{
    \begin{array}{ll}
      bA\bar{b} \otimes_{A_{G_{\bar{b}}}\bar{b}} \bar{b}Ab \rightarrow Ab \\
      ba_g\bar{a}\bar{b} \otimes \bar{b}\bar{a}'a_{g'}b \mapsto a_g\bar{a}\bar{b}\bar{a}'a_{g'}\\
      ab \mapsto \sum_{g \in G/G_{\bar{b}}} baa_g\bar{b} \otimes \bar{b}a_g^{-1},
    \end{array}
  \right.
 $$
and
$$\left\{
    \begin{array}{ll}
      \bar{b}Ab \otimes_{Ab} bA\bar{b} \rightarrow A_{G_{\bar{b}}}\bar{b}\\
      \bar{b}\bar{a}a_gb \otimes ba_{g'}\bar{a}'\bar{b} \mapsto \left\{
                                                            \begin{array}{ll}
                                                              \bar{b}\bar{a}a_ga_{g'}\bar{a}'\bar{b}, & \hbox{if $gg' \in G_{\bar{b}}$;} \\
                                                              0, & \hbox{if not.}
                                                            \end{array}
                                                          \right. \\
      \bar{a}a_g\bar{b} \mapsto \bar{a}a_g\bar{b} \otimes \bar{b} \textrm{ } (\textrm{ where } g \in G_{\bar{b}}),
    \end{array}
  \right.$$ define isomorphisms inverse to each other
$$ bA\bar{b} \otimes_{A_{G_{\bar{b}}}\bar{b}} \bar{b}Ab\, \tilde{\leftrightarrow}\, Ab  \,\,\,\textrm{ and
}\,\,\, \bar{b}Ab \otimes_{Ab} bA\bar{b} \,\tilde{\leftrightarrow}\,
A_{G_{\bar{b}}}\bar{b}.$$ Therefore, $bA\bar{b}$ and $\bar{b}Ab$ are
Morita equivalent. In particular, the functors
$$\mathrm{Ind}_{\bar{A}}^A=(bA\bar{b} \otimes_{A_{G_{\bar{b}}}\bar{b}}-) \,\,\,\textrm{ and
}\,\,\, \bar{b} \cdot \mathrm{Res}_{\bar{A}}^A= (\bar{b}Ab
\otimes_{Ab} -)$$ define category equivalences inverse to each other
between $_{A_{G_{\bar{b}}}\bar{b}}\emph{\textbf{mod}}$ and\\
$_{Ab}\emph{\textbf{mod}}$.
\end{proposition}

\subsection{Multiplication of an $A$-module by an $\mathcal{O}G$-module}

Let $X$ be an $A$-module and $\rho:A \rightarrow
\mathrm{End}_\mathcal{O}(X)$ be the structural morphism. We define
an additive functor
$$ X \cdot - : _{\mathcal{O}G}\mathrm{\textbf{mod}} \rightarrow _A
\mathrm{\textbf{mod}}, Y \mapsto X \cdot Y $$ as follows:
If $Y$ is an $\mathcal{O}G$-module and $\sigma: \mathcal{O}G
\rightarrow \mathrm{End}_\mathcal{O}(Y)$ is the structural morphism,
we denote by $X \cdot Y$ the $\mathcal{O}$-module $X
\otimes_\mathcal{O} Y$. The action of $A$ on the latter is given by
the morphism
$$\rho\cdot\sigma:A\rightarrow\mathrm{End}_\mathcal{O}(X
\otimes Y), \bar{a}a_g \mapsto \rho(\bar{a}a_g) \otimes \sigma(g).$$

\begin{proposition}\label{multiplication of A and G modules}
Let $X$ be an $\bar{A}$-module. The application $$A
\otimes_{\bar{A}}X \rightarrow X \cdot \mathcal{O}G$$ defined by
$$a_g\otimes_{\bar{A}}x \mapsto \rho(a_g)(x)\otimes_\mathcal{O}g
\textrm{ (for all $x \in X$ and $g \in G$)}$$ is an isomorphism of
$A$-modules
$$\mathrm{Ind}_{\bar{A}}^A(X) \tilde{\rightarrow} X \cdot \mathcal{O}G.$$
\end{proposition}

\subsection{Induction and restriction of $KA$-modules and
$K\bar{A}$-modules}

Let $X$ be a $KA$-module of character $\chi$ and $Y$ a $KG$-module
of character $\xi$. We denote by $\chi \cdot \xi$ the character of
the $KA$-module $X \cdot Y$. From now on, all group algebras over
$K$ will be considered split semisimple.

\begin{proposition}\label{1.39}
Let $\chi$ be an irreducible character of $KA$ which restricts to an
irreducible character $\bar{\chi}$ of $K\bar{A}$. Then
\begin{enumerate}
  \item The characters $(\chi \cdot \xi)_{\xi \in \mathrm{Irr}(KG)}$
  are distinct irreducible characters of $KA$.
  \item We have $$\mathrm{Ind}_{K\bar{A}}^{KA}(\bar{\chi})=
  \sum_{\xi \in \mathrm{Irr}(KG)} \xi(1)(\chi \cdot \xi).$$
\end{enumerate}
\end{proposition}
\begin{apod}{The second relation results from Proposition $\ref{multiplication of A and G
modules}$. We have
$$\langle \mathrm{Ind}_{K\bar{A}}^{KA}(\bar{\chi}),
\mathrm{Ind}_{K\bar{A}}^{KA}(\bar{\chi}) \rangle_{KA} =
\sum_{\xi,\xi' \in \mathrm{Irr}(KG)} \xi(1)\xi'(1)\langle \chi
\cdot \xi,\chi \cdot \xi'\rangle_{KA}.$$
Frobenius reciprocity now gives
$$\begin{array}{ccc}
\langle \mathrm{Ind}_{K\bar{A}}^{KA}(\bar{\chi}),
\mathrm{Ind}_{K\bar{A}}^{KA}(\bar{\chi}) \rangle_{KA} & = & \langle
\mathrm{Res}_{K\bar{A}}^{KA}(\sum_{\xi \in \mathrm{Irr}(KG)}
\xi(1)\chi \cdot \xi), \bar{\chi} \rangle_{K\bar{A}} \\
     & = & \langle
\sum_{\xi \in \mathrm{Irr}(KG)} \xi(1)^2 \bar{\chi},\bar{\chi}
\rangle_{K\bar{A}} \\
     & = & \sum_{\xi \in \mathrm{Irr}(KG)} \xi(1)^2 = |G|,
  \end{array}$$
whence we obtain
$$|G|=\sum_{\xi,\xi' \in \mathrm{Irr}(KG)} \xi(1)\xi'(1)\langle \chi
\cdot \xi,\chi \cdot \xi'\rangle_{KA}.$$ Since $|G|=\sum_{\xi \in
\mathrm{Irr}(KG)} \xi(1)^2$, we must have $\langle \chi \cdot
\xi,\chi \cdot \xi'\rangle_{KA}=\delta_{\xi,\xi'}$ and the proof is
complete.}
\end{apod}

For all $\bar{\chi} \in \mathrm{Irr}(K\bar{A})$, we denote by
$\bar{e}(\bar{\chi})$ the block of $K\bar{A}$ associated with
$\bar{\chi}$. We have seen that if $g \in G$, then
$g(\bar{e}(\bar{\chi}))$ is also a block of $K\bar{A}$. Since
$K\bar{A}$ is split semisimple, it must be associated with an
irreducible character $g(\bar{\chi})$ of $K\bar{A}$. Thus, we can
define an action of $G$ on $\mathrm{Irr}(K\bar{A})$ such that for
all $g \in G$, $\bar{e}(g(\bar{\chi}))=g(\bar{e}(\bar{\chi}))$. We
denote by $G_{\bar{\chi}}$ the stabilizer of $\bar{\chi}$ in $G$. Obviously, $G_{\bar{\chi}}=
G_{\bar{e}(\bar{\chi})}$.

\begin{proposition}\label{1.40}
Let $\bar{\chi} \in \mathrm{Irr}(K\bar{A})$ and suppose that
$\bar{\chi}$ is extended to a character
$\tilde{\chi} \in  \mathrm{Irr}(KA_{G_{\bar{\chi}}})$ (i.e.,
$\mathrm{Res}_{K\bar{A}}^{KA_{G_{\bar{\chi}}}}(\tilde{\chi})=\bar{\chi}$).
We set
$$\chi:=\mathrm{Ind}_{KA_{G_{\bar{\chi}}}}^{KA}(\tilde{\chi}) \,\textrm{
and }\,
\chi_\xi:=\mathrm{Ind}_{KA_{G_{\bar{\chi}}}}^{KA}(\tilde{\chi} \cdot
\xi) \textrm{ for all } \xi \in \mathrm{Irr}(KG_{\bar{\chi}}).$$
Then
\begin{enumerate}
  \item The characters $(\chi_\xi)_{\xi \in \mathrm{Irr}(KG_{\bar{\chi}})}$
  are distinct irreducible characters of $KA$.
  \item We have $$\mathrm{Ind}_{K\bar{A}}^{KA}(\bar{\chi})=
  \sum_{\xi \in \mathrm{Irr}(KG_{\bar{\chi}})} \xi(1)\chi_\xi.$$
  In particular,
  $$m_{\chi_\xi,\bar{\chi}}=\xi(1) \,\textrm{ and }\,
  \chi_\xi(1)=|G:G_{\bar{\chi}}|\bar{\chi}(1)\xi(1).$$
  \item For all $\xi \in \mathrm{Irr}(KG_{\bar{\chi}})$, we have
  $$s_{\chi_\xi}=\frac{|G_{\bar{\chi}}|}{\xi(1)}s_{\bar{\chi}}.$$
\end{enumerate}
\end{proposition}
\begin{apod}{
\begin{enumerate}
  \item By Proposition $\ref{1.39}$, we obtain that the characters
   $(\tilde{\chi} \cdot \xi)_{\xi \in \mathrm{Irr}(KG_{\bar{\chi}})}$
   are distinct irreducible characters of
   $\mathrm{Irr}(KA_{G_{\bar{\chi}}})$. Now let $\bar{e}(\bar{\chi})$ be the
   block of $K\bar{A}$ associated with the irreducible character
   $\bar{\chi}$. We have seen that $\bar{e}(\bar{\chi})$ is a central
   idempotent of $KA_{G_{\bar{\chi}}}$. Proposition $\ref{morita}$
   implies that the functor
   $\mathrm{Ind}_{K\bar{A}}^{KA}$ defines a Morita equivalence
   between the category
   $_{KA_{G_{\bar{\chi}}}\bar{e}(\bar{\chi})}\mathrm{\textbf{mod}}$
   and its image. Therefore, the characters $(\mathrm{Ind}_{KA_{G_{\bar{\chi}}}}^{KA}(\tilde{\chi}
   \cdot \xi))_{\xi \in \mathrm{Irr}(KG_{\bar{\chi}})}$ are distinct
   irreducible characters of $KA$.
  \item By Proposition $\ref{1.39}$, we obtain that
   $$\mathrm{Ind}_{K\bar{A}}^{KA_{G_{\bar{\chi}}}}(\bar{\chi})=
   \sum_{\xi \in \mathrm{Irr}(KG_{\bar{\chi}})} \xi(1)\tilde{\chi} \cdot \xi.$$
   Applying $\mathrm{Ind}_{KA_{G_{\bar{\chi}}}}^{KA}$ to both sides gives us the
   required relation. Obviously, $m_{\chi_\xi,\bar{\chi}}=\xi(1)$.

   Now let us calculate the value of $\chi_\xi(\bar{a})$ for any
   $\bar{a} \in \bar{A}$.
   Let $Y$ be an irreducible $KA_{G_{\bar{\chi}}}$-module of character $\psi$.Then
   $\mathrm{Ind}_{KA_{G_{\bar{\chi}}}}^{KA}(Y)=
   KA \otimes_{KA_{G_{\bar{\chi}}}} Y$ has character
   $\mathrm{Ind}_{KA_{G_{\bar{\chi}}}}^{KA}(\psi)$.
   We have $KA = \bigoplus_{g \in G/G_{\bar{\chi}}}a_g K\bar{A}$.
   Let $\bar{a} \in \bar{A}$. Then
   $$\begin{array}{ccl}
       \bar{a}\mathrm{Ind}_{KA_{G_{\bar{\chi}}}}^{KA}(Y) & = & \bigoplus_{g \in G/G_{\bar{\chi}}}\bar{a}a_g K\bar{A} \otimes_{KA_{G_{\bar{\chi}}}} Y \\
        & = & \bigoplus_{g \in G/G_{\bar{\chi}}}a_g(a_g^{-1}\bar{a}a_g) K\bar{A} \otimes_{KA_{G_{\bar{\chi}}}} Y\\
        & = & \bigoplus_{g \in G/G_{\bar{\chi}}}a_g K\bar{A} \otimes_{KA_{G_{\bar{\chi}}}} (a_g^{-1}\bar{a}a_g)Y.
     \end{array}$$
   Thus,
   $\mathrm{Ind}_{KA_{G_{\bar{\chi}}}}^{KA}(\psi)(\bar{a})=
   \sum_{g \in G/G_{\bar{\chi}}}\psi(a_g^{-1}\bar{a}a_g)$
   and
   $$\chi_\xi(\bar{a})=\sum_{g \in G/G_{\bar{\chi}}}(\tilde{\chi} \cdot
   \xi)(a_g^{-1}\bar{a}a_g)= \sum_{g \in
   G/G_{\bar{\chi}}}\bar{\chi}(a_g^{-1}\bar{a}a_g)\xi(1).$$
   Therefore,
   $$\chi_\xi(1)=\sum_{g \in G/G_{\bar{\chi}}}\bar{\chi}(1)\xi(1)=
   |G:G_{\bar{\chi}}|\bar{\chi}(1)\xi(1).$$
  \item Let $(\bar{e}_i)_{i \in I}$ be a basis of $\bar{A}$ as
   $\mathcal{O}$-module and let $(\bar{e}_i')_{i \in I}$ be its
   dual with respect to the symmetrizing form $\bar{t}$.
   Proposition $\ref{schur elements and idempotents}$(3), in combination with Proposition
   $\ref{dual basis of a symmetric algebra of a finite group}$, gives
   $$s_{\chi_\xi}\chi_\xi(1)^2=\chi_\xi(\sum_{i \in I,g \in G}\bar{e}_i'a_ga_g^{-1}\bar{e}_i)=
     \chi_\xi(|G|\sum_{i \in I}\bar{e}_i'\bar{e}_i).$$
   However, $\sum_{i \in I,g \in G}\bar{e}_i'a_ga_g^{-1}\bar{e}_i$
   belongs to to center of $A$ (by Lemma $\ref{properties of the casimir}$) and thus, for all $h \in G$,
   $${a_h}^{-1}(\sum_{i \in I}\bar{e}_i'a_ga_g^{-1}\bar{e}_i)a_h=
     \sum_{i \in I}\bar{e}_i'a_ga_g^{-1}\bar{e}_i=
     |G|\sum_{i \in I}\bar{e}_i'\bar{e}_i.$$
   Since $\sum_{i \in I}\bar{e}_i'\bar{e}_i \in \bar{A}$, we must have
      $$\begin{array}{ccl}
       \chi_\xi(|G|\sum_{i \in I}\bar{e}_i'\bar{e}_i) & = & \sum_{h \in G/G_{\bar{\chi}}}\bar{\chi}({a_h}^{-1}(|G|\sum_{i \in
   I}\bar{e}_i'\bar{e}_i)a_h)\xi(1) \\
    & & \\
        & = & \sum_{h \in G/G_{\bar{\chi}}}\bar{\chi}(|G|\sum_{i \in
   I}\bar{e}_i'\bar{e}_i)\xi(1) \\
   & & \\
        & = & |G:G_{\bar{\chi}}| |G| \xi(1) \bar{\chi}(\sum_{i \in
   I}\bar{e}_i'\bar{e}_i) \\
   & & \\
        & = & |G:G_{\bar{\chi}}|^2 |G_{\bar{\chi}}| \xi(1) s_{\bar{\chi}} \bar{\chi}(1)^2.
     \end{array}$$
   So we obtain
   $$s_{\chi_\xi}\chi_\xi(1)^2=|G:G_{\bar{\chi}}|^2 |G_{\bar{\chi}}| \xi(1) \bar{\chi}(1)^2 s_{\bar{\chi}}.$$
   Replacing $\chi_\xi(1)=|G:G_{\bar{\chi}}|\bar{\chi}(1)\xi(1)$
   gives
   $$s_{\chi_\xi} \xi(1) = |G_{\bar{\chi}}| s_{\bar{\chi}} .$$}
\end{enumerate}
\end{apod}

Now let $\bar{\Omega}$ be the orbit of the character $\bar{\chi} \in
\mathrm{Irr}(K\bar{A})$ under the action of $G$. We have
$|\bar{\Omega}|=|G|/|G_{\bar{\chi}}|$. Define
$$\bar{e}(\bar{\Omega})=\sum_{g \in
G/G_{\bar{\chi}}}\bar{e}(g(\bar{\chi}))= \sum_{g \in
G/G_{\bar{\chi}}}g(\bar{e}(\bar{\chi})).$$ 

If $\bar{\chi} \in \bar{\Omega}$, the set
$\mathrm{Irr}(KA,\bar{\chi})$ depends only on $\bar{\Omega}$ and we
set $\mathrm{Irr}(KA,\bar{\Omega}):=\mathrm{Irr}(KA,\bar{\chi})$.
The idempotent $\bar{e}(\bar{\Omega})$ belongs to the algebra
$(ZK\bar{A})^G$ of the elements in the center of $K\bar{A}$ fixed by
$G$ and thus to the center of $KA$ (since it commutes with all
elements of $\bar{A}$ and all $a_g$, $g \in G$). Therefore, it must
be a sum of blocks of $KA$, \ie
$$\bar{e}(\bar{\Omega}) = \sum_{\chi \in
\mathrm{Irr}(KA,\bar{\Omega})}e(\chi).$$

Let $X$ be an irreducible $KA$-module of character $\chi$ and
$\bar{X}$ an irreducible $K\bar{A}$-submodule of
$\mathrm{Res}_{K\bar{A}}^{KA}(X)$ of character $\bar{\chi}$. For $g
\in G$, the $K\bar{A}$-submodule $g(\bar{X})$ of
$\mathrm{Res}_{K\bar{A}}^{KA}(X)$ has character $g(\bar{\chi})$.
Then $\sum_{g \in G}g(\bar{X})$ is a $KA$-submodule of $X$. We
deduce that
$$\mathrm{Res}_{K\bar{A}}^{KA}(X)=( \bigoplus_{g \in
G/G_{\bar{\chi}}}g(\bar{X}))^{m_{\chi,\bar{\chi}}},$$ \ie
$$\mathrm{Res}_{K\bar{A}}^{KA}(\chi)= m_{\chi,\bar{\chi}}\sum_{g \in
G/G_{\bar{\chi}}}g(\bar{\chi}).$$ In particular, we see that
$\mathrm{Irr}(K\bar{A},\chi)$ is an orbit of $G$ on
$\mathrm{Irr}(K\bar{A})$. Notice that
$\chi(1)=m_{\chi,\bar{\chi}}|\bar{\Omega}|\bar{\chi}(1)$.\\
\\
\emph{Case where $G$ is cyclic}\\
\\
Suppose that $G$ is a cyclic group of order $d$ and let $g \in G$ be a generator of
$G$ (we can choose
$\mathrm{Rep}(A/\bar{A})=\{1,a_g,a_g^2,\ldots,a_g^{d-1}\}$). We
will show that the assumptions of Proposition $\ref{1.40}$ are
satisfied for all irreducible characters of $K\bar{A}$.

Let $\bar{X}$ be an irreducible $K\bar{A}$-module and let
$\bar{\rho}:K\bar{A}\rightarrow \mathrm{End}_K(\bar{X})$ be the
structural morphism. Since the representation of $\bar{X}$ is
invariant by the action of $G_{\bar{\chi}}$, there exists an
automorphism $\alpha$ of the $K$-vector space $\bar{X}$ such that
$$\alpha \bar{\rho}(\bar{a}) \alpha^{-1}= g(\bar{\rho})(\bar{a}),$$
for all $g \in G_{\bar{\chi}}$.

The subgroup $G_{\bar{\chi}}$ is also cyclic. Let
$d(\bar{\chi}):=|G_{\bar{\chi}}|$. Then
$$\bar{\rho}(\bar{a})=\alpha^{d(\bar{\chi})} \bar{\rho}(\bar{a}) \alpha^{-d(\bar{\chi})}.$$
Since $\bar{X}$ is irreducible and $K\bar{A}$ is split semisimple,
$\alpha^{d(\bar{\chi})}$ must be a scalar. Instead of enlarging the
field $K$, we can assume that $K$ contains a $d(\bar{\chi})$-th root
of that scalar. By dividing $\alpha$ by that root, we reduce to the
case where $\alpha^{d(\bar{\chi})}=1$.

This allows us to extend the structural morphism
$\bar{\rho}:K\bar{A}\rightarrow \mathrm{End}_K(\bar{X})$ to a
morphism
$$\tilde{\rho}:KA_{G_{\bar{\chi}}}\rightarrow
\mathrm{End}_K(\bar{X})$$ such that
$$\tilde{\rho}(\bar{a}a_h^j):=\bar{\rho}(\bar{a})\alpha^j
 \textrm{ for } 0\leq j<d(\bar{\chi}),$$
 where $h:=g^{d/d(\bar{\chi})}$ generates $G_{\bar{\chi}}$.
The morphism $\tilde{\rho}$ defines a $KA_{G_{\bar{\chi}}}$-module
$\tilde{X}$ of character $\tilde{\chi}$. By definition of $\tilde{\rho}$,
$\mathrm{Res}_{K\bar{A}}^{KA_{G_{\bar{\chi}}}}(\tilde{\chi})=\bar{\chi}$.

Since the group $G$ is abelian, the set $\mathrm{Irr}(KG)$ forms a
group, which we denote by $G^\vee$. The application $\psi \mapsto
\psi \cdot \xi$, where $\psi \in \mathrm{Irr}(KA)$ and $\xi \in
G^\vee$, defines an action of $G^\vee$ on $\mathrm{Irr}(KA)$.

Let $\Omega$ be the orbit of $\tilde{\chi}$ under the action of
$(G_{\bar{\chi}})^\vee$. By Proposition $\ref{1.39}$, we obtain that
$\Omega$ is a regular orbit (\ie $|\Omega|=|G_{\bar{\chi}}|$) and
that $\Omega=\mathrm{Irr}(KA_{G_{\bar{\chi}}},\bar{\chi})$.
Like in Proposition $\ref{1.40}$, we introduce the notations
$$\chi:=\mathrm{Ind}_{KA_{G_{\bar{\chi}}}}^{KA}(\tilde{\chi})\textrm{ and }
\chi_\xi:=\mathrm{Ind}_{KA_{G_{\bar{\chi}}}}^{KA}(\tilde{\chi} \cdot
\xi) \textrm{ for all } \xi \in (G_{\bar{\chi}})^\vee.$$
Then
$$\mathrm{Irr}(KA,\bar{\chi})=\{\chi_\xi\,|\,\xi\in(G_{\bar{\chi}})^\vee\}
\textrm{ and } m_{\chi_\xi,\bar{\chi}}=\xi(1)=1 \textrm{ for all }
\xi \in (G_{\bar{\chi}})^\vee.$$

Recall that $|G_{\bar{\chi}}|=d(\bar{\chi})$. There exists a
surjective morphism $G \twoheadrightarrow G_{\bar{\chi}}$ defined by
$g \mapsto g^{d/d(\bar{\chi})}$, which induces an inclusion
$(G_{\bar{\chi}})^\vee \hookrightarrow G^\vee$. If $\xi \in
(G_{\bar{\chi}})^\vee$, we denote (abusing notation) by $\xi$ its
image in $G^\vee$ by the above injection. It is easy to check that
$\chi_\xi=\chi \cdot \xi$.

Hence, Proposition $\ref{1.40}$ implies the following result.

\begin{proposition}\label{1.42}
If the group $G$ is cyclic, then there exists a bijection
$$\begin{array}{ccc}
    \mathrm{Irr}(K\bar{A})/G & \tilde{\leftrightarrow} & \mathrm{Irr}(KA)/G^\vee \\
    \bar{\Omega} & \leftrightarrow & \Omega
  \end{array}$$
such that
$$\bar{e}(\bar{\Omega})=e(\Omega),\, |\bar{\Omega}||\Omega|=|G| \textrm{ and }
\left\{
  \begin{array}{ll}
    \forall \chi \in \Omega, & \mathrm{Res}_{K\bar{A}}^{KA}(\chi)=\sum_{\bar{\chi} \in \bar{\Omega}}\bar{\chi}\\
    &\\
    \forall \bar{\chi} \in \bar{\Omega}, &\mathrm{Ind}_{K\bar{A}}^{KA}(\bar{\chi})=\sum_{\chi \in \Omega}\chi
  \end{array}
\right.
$$
Moreover, for all $\chi \in \Omega$ and $\bar{\chi} \in
\bar{\Omega}$, we have
$$s_\chi = |\Omega| s_{\bar{\chi}}.$$
\end{proposition}

\subsection{Blocks of $A$ and blocks of $\bar{A}$}

Let us denote by $\mathrm{Bl}(A)$ the set of blocks of $A$ and by
$\mathrm{Bl}(\bar{A})$ the set of blocks of $\bar{A}$. For $\bar{b}
\in \mathrm{Bl}(\bar{A})$, we have set
$$\mathrm{Tr}(G,\bar{b}):=\sum_{g \in G/G_{\bar{b}}}g(\bar{b}).$$

The algebra $(Z\bar{A})^G$ is contained in both $Z\bar{A}$ and $ZA$
and the set of its blocks is
$$\mathrm{Bl}((Z\bar{A})^G)=\{\mathrm{Tr}(G,\bar{b}) \,|\, \bar{b} \in
\mathrm{Bl}(\bar{A})/G\}.$$ Moreover, $\mathrm{Tr}(G,\bar{b})$ is a
sum of blocks of $A$ and we define the subset
$\mathrm{Bl}(A,\bar{b})$ of $\mathrm{Bl}(A)$ as follows:
$$\mathrm{Tr}(G,\bar{b}):=\sum_{b \in \mathrm{Bl}(A,\bar{b})}b.$$

\begin{lemma}\label{1.43}
Let $\bar{b}$ be a block of $\bar{A}$ and
$\bar{B}:=\mathrm{Irr}(K\bar{A}\bar{b})$. Then
\begin{enumerate}
  \item For all $\bar{\chi} \in \bar{B}$, we have $G_{\bar{\chi}}
  \subseteq G_{\bar{b}}$.
  \item We have
  $$\mathrm{Tr}(G,\bar{b})=\sum_{\bar{\chi} \in \bar{B}/G}
  \mathrm{Tr}(G,\bar{e}(\bar{\chi}))=
  \sum_{\{\bar{\Omega}\,|\,\bar{\Omega} \cap \bar{B} \neq \emptyset \}}\bar{e}(\bar{\Omega}).$$
\end{enumerate}
\end{lemma}
\begin{apod}{
\begin{enumerate}
  \item We have $\bar{b}=\sum_{\bar{\chi} \in \bar{B}}\bar{e}(\bar{\chi})$.
  If $g \notin G_{\bar{b}}$, then the blocks $\bar{b}$ and $g(\bar{b})$
  are orthogonal. Hence, $g \notin G_{\bar{\chi}}$ for all $\bar{\chi}\in \bar{B}$.
  \item Note that $\bar{b}=\sum_{\bar{\chi} \in \bar{B}}\bar{e}(\bar{\chi})=
  \sum_{\bar{\chi} \in
  \bar{B}/G_{\bar{b}}}\mathrm{Tr}(G_{\bar{b}},\bar{e}(\bar{\chi}))$.
  Thus,
  $$\mathrm{Tr}(G,\bar{b})=\sum_{\bar{\chi} \in \bar{B}/G}
  \mathrm{Tr}(G,\bar{e}(\bar{\chi}))=
  \sum_{\{\bar{\Omega}\,|\,\bar{\Omega} \cap \bar{B} \neq \emptyset \}}\bar{e}(\bar{\Omega}),$$
  by the definition of $\bar{e}(\bar{\Omega})$.}
\end{enumerate}
\end{apod}

Now let $G^\vee:=\mathrm{Hom}(G,K^\times)$. We suppose that $K=F$.
The multiplication of the characters of $KA$ by the characters of
$KG$ defines an action of the group $G^\vee$ on $\mathrm{Irr}(KA)$.
This action is induced by the operation of $G^\vee$ on the algebra
$A$, which is defined in the following way:
$$ \xi \cdot (\bar{a}a_g) := \xi(g)\bar{a}a_g \,\,\textrm{ for all }
\xi \in G^\vee, \bar{a} \in \bar{A}, g \in G.$$ In particular,
$G^\vee$ acts on the set of blocks of $A$. If $b$ is a block of
$A$, we denote by $\xi \cdot b$ the product of $\xi$ and $b$ and by
$(G^\vee)_b$ the stabilizer of $b$ in $G^\vee$. We set
$$\mathrm{Tr}(G^\vee,b):=\sum_{\xi\in G^\vee/(G^\vee)_b}\xi\cdot b.$$

The set of blocks of the algebra $(ZA)^{G^\vee}$ is given by
$$\mathrm{Bl}((ZA)^{G^\vee})=\{\mathrm{Tr}(G^\vee,b) \,|\, b \in
\mathrm{Bl}(A)/G^\vee\}.$$

The following lemma is the analogue of Lemma $\ref{1.43}$

\begin{lemma}\label{1.44}
Let $b$ be a block of $A$ and $B:=\mathrm{Irr}(KAb)$. Then
\begin{enumerate}
  \item For all $\chi \in B$, we have $(G^\vee)_\chi \subseteq (G^\vee)_b$.
  \item We have
  $$\mathrm{Tr}(G^\vee,b)=\sum_{\chi\in B/G^\vee} \mathrm{Tr}(G^\vee,e(\chi))=
  \sum_{\{\Omega\,|\,\Omega\cap B \neq \emptyset \}}e(\Omega).$$
\end{enumerate}
\end{lemma}
\emph{Case where $G$ is cyclic}\\
\\
For every orbit $\bar{\mathcal{Y}}$ of $G$ on
$\mathrm{Bl}(\bar{A})$, we denote by $\bar{b}(\bar{\mathcal{Y}})$
the block of $(Z\bar{A})^G$ defined as
$$\bar{b}(\bar{\mathcal{Y}}):=\sum_{\bar{b} \in \bar{\mathcal{Y}}}\bar{b}.$$
For every orbit $\mathcal{Y}$ of $G^\vee$ on $\mathrm{Bl}(A)$, we
denote by $b(\mathcal{Y})$ the block of $(ZA)^{G^\vee}$ defined as
$$b(\mathcal{Y}):=\sum_{b \in \mathcal{Y}}b.$$

The following proposition results from Proposition $\ref{1.42}$ and
Lemmas $\ref{1.43}$ and $\ref{1.44}$.

\begin{proposition}\label{1.45}
If the group $G$ is cyclic, then there exists a bijection
$$\begin{array}{ccc}
    \mathrm{Bl}(\bar{A})/G & \tilde{\leftrightarrow} & \mathrm{Bl}(A)/G^\vee \\
    \bar{\mathcal{Y}} & \leftrightarrow & \mathcal{Y}
  \end{array}$$
such that
$$\bar{b}(\bar{\mathcal{Y}})=b(\mathcal{Y}),$$
i.e.,
$$\mathrm{Tr}(G,\bar{b})=\mathrm{Tr}(G^\vee,b) \textrm{ for all }
\bar{b} \in \bar{\mathcal{Y}} \textrm{ and } b \in \mathcal{Y}.$$ In
particular, the algebras $(Z\bar{A})^G$ and $(ZA)^{G^\vee}$ have the
same blocks.
\end{proposition}

\begin{corollary}\label{clifford}
If the blocks of $A$ are stable under the action of $G^\vee$, then the
blocks of $A$ coincide with the blocks of $(Z\bar{A})^G$.
\end{corollary}

\section{Representation theory of symmetric algebras}

In the last section of Chapter $2$, we present some results concerning the representation theory of symmetric algebras. If a symmetric algebra satisfies certain conditions, we can define a decomposition map (and consequently, a decomposition matrix) and obtain the blocks with the use of a Brauer graph. In order to check whether the required conditions are satisfied, we have to know when a symmetric algebra is split or semisimple. In subsection $2.4.4$, we prove a theorem which gives us a new criterion for a symmetric algebra to be split and semisimple. With only this exception, all the results in this section are well-known and mostly taken from the seventh chapter  of
\cite{GePf}.

\subsection{Grothendieck groups}

Let $\mathcal{O}$ be an integral domain and $K$ a field containing
$\mathcal{O}$. Let $A$ be an $\mathcal{O}$-algebra free and finitely
generated as an $\mathcal{O}$-module.

Let $R_0(KA)$ be the Grothendieck group  \index{Grothendieck group} of finite-dimensional
$KA$-modules. Thus, $R_0(KA)$ is generated by expressions $[V]$, one
for each $KA$-module $V$ (up to isomorphism), with relations
$[V]=[V']+[V'']$ for each exact sequence $0 \rightarrow V'
\rightarrow V \rightarrow V'' \rightarrow 0$ of $KA$-modules. Two
$KA$-modules $V,V'$ give rise to the same element in $R_0(KA)$, if
$V$ and $V'$ have the same composition factors, counting
multiplicities. It follows that $R_0(KA)$ is free abelian with basis
given by the isomorphism classes of simple modules. Finally, let
$R_0^+(KA)$ be the subset of $R_0(KA)$ consisting of elements $[V]$,
where $V$ is a finite-dimensional $KA$-module.

\begin{definition}\label{pk}
Let $x$ be an indeterminate over $K$ and $\mathrm{Maps}(A,K[x])$ the
$K$-algebra of maps from $A$ to $K[x]$ (with pointwise
multiplication of maps as algebra multiplication). If $V$ is a
$KA$-module, let $\rho_V:KA \rightarrow \mathrm{End}_K(V)$ denote
its structural morphism. We define the map
$$\begin{array}{cccl}
    \mathfrak{p}_K: & R_0^+(KA) & \rightarrow & \emph{Maps}(A,K[x]) \\
                    & [V] & \mapsto & (a \mapsto \textrm{\emph{characteristic
                    polynomial of} } \rho_V(a)).
  \end{array}$$
Considering $\emph{Maps}(A,K[x])$ as a semigroup with respect to
multiplication, the map $\mathfrak{p}_K$ is a semigroup
homomorphism.
\end{definition}

Let $\mathrm{Irr}(KA)$ be the set of all irreducible characters of the algebra $KA$ (\ie the set of all characters $\chi_V$, where $V$ is a simple $KA$-module). The following result is known as the
``Brauer-Nesbitt lemma''\index{Brauer-Nesbitt lemma} (cf.~\cite{BrNe}, Lemma $2$).

\begin{lemma}\label{Brauer-Nesbitt} Assume that
$\mathrm{Irr}(KA)$ is a linearly independent subset of
 $\mathrm{Hom}_K(KA,K)$. Then the map $\mathfrak{p}_K$ is injective.
\end{lemma}
\begin{apod}{Let $V,V'$ be two $KA$-modules such that
$\mathfrak{p}_K([V])=\mathfrak{p}_K([V'])$. Since $[V], [V']$ only
depend on the composition factors of $V,V'$, we may assume that
$V,V'$ are semisimple modules. Let
$$V=\bigoplus_{i=1}^n a_iV_i \,\,\textrm{ and }\,\, V'=\bigoplus_{i=1}^n b_iV_i,$$
where the $V_i$ are pairwise non-isomorphic simple $KA$-modules and
$a_i,b_i \geq 0$ for all $i$. We have to show that $a_i=b_i$ for all
$i$.

If, for some $i$, we have both $a_i>0$ and $b_i>0$, then we can
write $V=V_i \oplus \tilde{V}$ and $V'=V_i \oplus \tilde{V}'$. Since
$\mathfrak{p}_K$ is a semigroup homomorphism, we obtain
$$\mathfrak{p}_K([V_i]) \cdot \mathfrak{p}_K([\tilde{V}]) = \mathfrak{p}_K([V])
= \mathfrak{p}_K([V'])= \mathfrak{p}_K([V_i]) \cdot
\mathfrak{p}_K([\tilde{V}']),$$ and, dividing by
$\mathfrak{p}_K([V_i])$, we conclude that
$\mathfrak{p}_K([\tilde{V}])=\mathfrak{p}_K([\tilde{V}'])$. Thus, we
can suppose that, for all $i$, we have $a_i=0$ or $b_i=0$. Taking
characters yields that
$$\chi_V = \sum_i a_i \chi_{V_i} \,\,\textrm{ and }\,\, \chi_{V'} = \sum_i b_i \chi_{V_i}.$$
For each $a \in A$, the character values $\chi_V(a)$
and $\chi_{V'}(a)$ appear as coefficients in the polynomials
$\mathfrak{p}_K([V])(a)$ and $\mathfrak{p}_K([V'])(a)$ respectively.
Since $\mathfrak{p}_K([V])=\mathfrak{p}_K([V'])$, we deduce that $\sum_i (a_i-b_i)\chi_{V_i}=0 $. By assumption,
the characters $\chi_{V_i}$ are linearly independent. So we must
have $(a_i-b_i) 1_K=0$ for all $i$. Since for all $i$, $a_i=0$ or
$b_i=0$, this means that $a_i 1_K=0$ and $b_i 1_K=0$ for all $i$. If
the field $K$ has characteristic 0, we conclude that $a_i=b_i=0$ for
all $i$ and we are done. If $K$ has characteristic $p>0$, we
conclude that $p$ divides all $a_i$ and all $b_i$ and so
$\frac{1}{p}[V]$ and $\frac{1}{p}[V']$ exist in $R_0^+(KA)$.
We also have
$\mathfrak{p}_K(\frac{1}{p}[V])=\mathfrak{p}_K(\frac{1}{p}[V'])$.
Repeating the above argument for $\frac{1}{p}[V]$ and
$\frac{1}{p}[V']$ yields that the multiplicity of $V_i$ in each of
these modules is still divisible by $p$. If we repeat this again and
again, we deduce that $a_i$ and $b_i$ should be divisible by
arbitrary powers of $p$. This forces $a_i=b_i=0$ for all $i$, as
desired.}
\end{apod}\\
\begin{remark} \emph{The assumption of the Brauer-Nesbitt lemma is satisfied when
(but not only when) 
 $KA$ is split or
 $K$ is a perfect field.}
\end{remark}\\

The following lemma (cf.~\cite{GePf},
Lemma 7.3.4) implies the compatibility of the map
$\mathfrak{p}_K$ with the field extensions of $K$ .

\begin{lemma}\label{compatibility with field extensions}
Let $K \subseteq K'$ be a field extension. Then there is a canonical
map $d_K^{K'}:R_0(KA) \rightarrow R_0(K'A)$ given by $[V] \mapsto [K'
\otimes_K V]$. Furthermore, we have a commutative diagram
$$\diagram R_0^+(KA) \dto^{d_K^{K'}} \rto^{\mathfrak{p}_K} &\emph{Maps}(A,K[x]) \dto^{\tau_K^{K'}} \\
           R_0^+(K'A) \rto^{\mathfrak{p}_{K'}} &\emph{Maps}(A,K'[x])
\enddiagram$$
where $\tau_K^{K'}$ is the canonical embedding. If, moreover, $KA$
is split, then $d_K^{K'}$ is an isomorphism which preserves
isomorphism classes of simple modules.
\end{lemma}
\subsection{Integrality}

We have seen in Chapter $1$ that a subring $\mathcal{R} \subseteq K$
is a valuation ring if, for each non-zero element $x \in K$, we have
$x \in \mathcal{R}$ or $x^{-1} \in \mathcal{R}$. Consequently, $K$
is the field of fractions of $\mathcal{R}$.

Such a valuation ring is a local ring whose maximal ideal we will
denote by $\mathfrak{m}_\mathcal{R}$. Valuation rings have interesting
properties, some of which are:
\begin{description}
  \item[(V1)] If $I$ is a prime ideal of $\mathcal{O}$, then there
  exists a valuation ring $\mathcal{R} \subseteq K$ such that
  $\mathcal{O} \subseteq \mathcal{R}$ and $\mathfrak{m}_\mathcal{R} \cap
  \mathcal{O}=I$.
  \item[(V2)] Every finitely generated torsion-free module over a
  valuation ring in $K$ is free.
  \item[(V3)] The intersection of all valuation rings $\mathcal{R}
  \subseteq K$ with $\mathcal{O} \subseteq \mathcal{R}$ is the
  integral closure of $\mathcal{O}$ in $K$; each valuation ring
  itself is integrally closed in $K$ (Proposition $\ref{intersection of valuation rings}$).
\end{description}

\begin{lemma}\label{realizing modules over O}
Let $V$ be a $KA$-module. Choosing a $K$-basis of $V$, we obtain a
corresponding matrix representation $\rho:KA \rightarrow
\emph{M}_n(K)$, where $n=\emph{dim}_K(V)$. If $\mathcal{R} \subseteq
K$ is a valuation ring with $\mathcal{O} \subseteq \mathcal{R}$,
then a basis of $V$ can be chosen so that $\rho(a) \in
\emph{M}_n(\mathcal{R})$ for all $a \in A$. In that case, we say
that $V$ is realized over $\mathcal{R}$. \index{realized over}
\end{lemma}
\begin{apod}{Let $(v_1,\ldots,v_n)$ be a $K$-basis of $V$ and $\mathcal{B}$
an $\mathcal{O}$-basis for $A$. Let $\tilde{V}$ be the
$\mathcal{R}$-submodule of $V$ spanned by the finite set $\{v_ib \,
|\, 1 \leq i \leq n, b \in \mathcal{B}\}$. Then $\tilde{V}$ is
invariant under the action of $\mathcal{R}A$ and hence a finitely
generated $\mathcal{R}A$-module. Since it is contained in a
$K$-vector space, it is also torsion-free. So (V2) implies that
$\tilde{V}$ is an $\mathcal{R}A$-lattice (a finitely generated
$\mathcal{R}A$-module which is free as an $\mathcal{R}$-module) such
that $K \otimes_\mathcal{R}  \tilde{V} \cong V$. Thus, any
$\mathcal{R}$-basis of $\tilde{V}$ is also a $K$-basis of $V$ with
the required property.}
\end{apod}\\
\begin{remark}
\emph{Note that the above argument only requires that $\mathcal{R}$
is a subring of $K$ such that $K$ is the field of fractions of
$\mathcal{R}$ and $\mathcal{R}$ satisfies (V2). These conditions
also hold, for example, when $\mathcal{R}$ is a principal ideal
domain with $K$ as field of fractions.}
\end{remark}\\

The following two important results (\cite{GePf}, Propositions 7.3.8 and 7.3.9)  derive from the above lemma.

\begin{proposition}\label{integrality of pk}
Let $V$ be a $KA$-module and $\mathcal{O}_K$ be the integral closure
of $\mathcal{O}$ in $K$. Then we have $\mathfrak{p}_K([V])(a) \in
\mathcal{O}_K[x]$ for all $a \in A$. Thus the map $\mathfrak{p}_K$
of Definition $\ref{pk}$ is in fact a map $ R_0^+(KA)  \rightarrow
\emph{Maps}(A,\mathcal{O}_K[x])$.
\end{proposition}
\begin{apod}{Fix $a \in A$. Let $\mathcal{R} \subseteq K$ be a
valuation ring with $\mathcal{O} \subseteq \mathcal{R}$. By Lemma
$\ref{realizing modules over O}$, there exists a basis of $V$ such
that the action of $a$ on $V$ with respect to that basis is given by
a matrix with coefficients in $\mathcal{R}$. Therefore, we have that
$\mathfrak{p}_K([V])(a) \in \mathcal{R}[x]$. Since this holds for
all valuation rings $\mathcal{R}$ in $K$ containing $\mathcal{O}$,
property (V3) implies that $\mathfrak{p}_K([V])(a) \in
\mathcal{O}_K[x]$.}
\end{apod}

Note that, in particular, Proposition $\ref{integrality of pk}$
implies that $\chi_V(a)\in \mathcal{O}_K$ for all $a \in A$, where
$\chi_V$ is the character of the representation $\rho_V$.

The next proposition is a result on symmetric algebras already mentioned in section 2.2 
(Proposition $\ref{Schur element belongs to the integral closure}$): the integrality of the Schur elements.

\begin{proposition}\label{integrality of the Schur elements}
Assume that we have a
symmetrizing form $t$ on $A$. Let $V$ be a split simple $KA$-module
(i.e., $\mathrm{End}_{KA}(V) \cong K$) and let $s_V$ be its Schur
element with respect to the induced form $t^K$ on $KA$. Then $s_V
\in \mathcal{O}_K$.
\end{proposition}
\begin{apod}{Let $\mathcal{R} \subseteq K$ be a
valuation ring with $\mathcal{O} \subseteq \mathcal{R}$. By Lemma
$\ref{realizing modules over O}$, we can assume that $V$ affords a
representation $\rho:KA \rightarrow \mathrm{M}_n(K)$ such that
$\rho(a) \in \mathrm{M}_n(\mathcal{R})$ for all $a \in A$. Let
$\mathcal{B}$ be an $\mathcal{O}$-basis of $A$ and let
$\mathcal{B}'$ be its dual with respect to $t$. Then $s_V=\sum_{b
\in \mathcal{B}}\rho(b)_{ij}\rho(b')_{ji}$ for all $1 \leq i,j \leq
n$ (\cite{GePf}, Cor. 7.2.2). All terms in the sum lie in
$\mathcal{R}$ and so $s_V \in \mathcal{R}$. Since this holds for all
valuation rings $\mathcal{R}$ in $K$ containing $\mathcal{O}$,
property (V3) implies that $s_V \in \mathcal{O}_K$.}
\end{apod}
\subsection{The decomposition map}

Now, we moreover assume that the ring $\mathcal{O}$ is integrally
closed in $K$, \ie $\mathcal{O}_K=\mathcal{O}$. Throughout we will fix a ring homomorphism $\theta:
\mathcal{O} \rightarrow L$ into a field $L$ such that $L$ is the
field of fractions of $\theta(\mathcal{O})$. We call such a ring
homomorphism a \emph{specialization} of $\mathcal{O}$.\index{specialization}

Let $\mathcal{R} \subseteq K$ be a valuation ring with $\mathcal{O}
\subseteq \mathcal{R}$ and $\mathfrak{m}_\mathcal{R} \cap \mathcal{O} =
\mathrm{Ker}\theta$ (note that $\mathrm{Ker}\theta$ is a prime
ideal, since $\theta(\mathcal{O})$ is contained in a field). Let $k$
be the residue field of $\mathcal{R}$. Then the restriction of the
canonical map $\pi:\mathcal{R} \rightarrow k$ to $\mathcal{O}$ has
kernel $\mathfrak{m}_\mathcal{R} \cap \mathcal{O} = \mathrm{Ker}\theta$. Since
$L$ is the field of fractions of $\theta(\mathcal{O})$, we may
regard $L$ as a subfield of $k$. Thus, we obtain a commutative diagram
$$\diagram\mathcal{O}\dto^{\theta}&\subseteq&
          \mathcal{R}\dto^{\pi}&\subseteq &K\\
           L & \subseteq & k & & \enddiagram $$

From now on, we make the following assumption:

\begin{ypothesh}\label{split assumption}
\emph{(a)} $LA \textrm{ is split}$ \,\,\,\,or\,\,\,\,  \emph{(b)} $L=k \textrm{
and } k \textrm{ is perfect}$.
\end{ypothesh}

The map $\theta:\mathcal{O} \rightarrow L$ induces a map $A
\rightarrow LA, a \mapsto 1 \otimes a$. One consequence of the
assumption $\ref{split assumption}$ is that, due to Lemma
$\ref{compatibility with field extensions}$, the map $d_L^k:R_0(LA)
\rightarrow R_0(kA)$ is an isomorphism which preserves isomorphism
classes of simple modules. Thus we can identify $R_0(LA)$ and
$R_0(kA)$. Moreover, the Brauer-Nesbitt lemma holds for $LA$, \ie
the map $\mathfrak{p}_L:R_0^+(LA) \rightarrow \mathrm{Maps}(A,L[x])$
is injective.

Let $V$ be a $KA$-module and $\mathcal{R} \subseteq K$ be a
valuation ring with $\mathcal{O} \subseteq \mathcal{R}$. By Lemma
$\ref{realizing modules over O}$, there exists a $K$-basis of $V$
such that the corresponding matrix representation $\rho:KA
\rightarrow \mathrm{M}_n(K)$ ($n=\mathrm{dim}_K(V)$) has the
property that $\rho(a) \in \mathrm{M}_n(\mathcal{R})$ for all $a\in A$. Following the proof of
Lemma $\ref{realizing modules over O}$, that
basis generates an $\mathcal{R}A$-lattice $\tilde{V}$ such that
$K \otimes_\mathcal{R} \tilde{V}=V$. The $k$-vector space
$k \otimes_\mathcal{R} \tilde{V}$ is a $kA$-module via
$(1\otimes v)(1\otimes a) = 1\otimes va$ $(v \in \tilde{V},a \in A)$,
which we call the \emph{modular reduction} of $\tilde{V}$.
To simplify notation, we shall write
$$K\tilde{V}:=K \otimes_\mathcal{R} \tilde{V} \textrm{ and }
k\tilde{V}:=k \otimes_\mathcal{R} \tilde{V}  .$$
The matrix representation $\rho^k:kA \rightarrow \mathrm{M}_n(k)$
afforded by $k\tilde{V}$ is given by
$$\rho^k(1\otimes a)=(\pi(a_{ij})) \textrm{ where } a \in A \textrm{ and } \rho(a)=(a_{ij}).$$
Note that if
$\tilde{V}'$ is another $\mathcal{R}A$-lattice such that $K
\otimes_\mathcal{R} \tilde{V}' \cong V$, then $\tilde{V}$ and $\tilde{V}'$
need not be isomorphic. The same hold for the $kA$-modules
$k \otimes_\mathcal{R} \tilde{V} $ and
$k \otimes_\mathcal{R} \tilde{V}'$.\\

Now we are ready to state and prove the following result (cf.~\cite{GePf}, Theorem 7.4.3), which associates to $A$ a decomposition map, in the case where $\mathcal{O}$ is integrally closed in $K$ .

\begin{thedef}\label{existence of decomposition maps}
Let $\theta:\mathcal{O} \rightarrow L$ be a ring homomorphism into a
field $L$ such that $L$ is the field of fractions of
$\theta(\mathcal{O})$ and $\mathcal{O}$ is integrally closed in $K$.
Assume that we have chosen a valuation ring $\mathcal{R}$ with
$\mathcal{O} \subseteq \mathcal{R} \subseteq K$ and $\mathfrak{m}_\mathcal{R}
\cap \mathcal{O} = \mathrm{Ker}\theta$ and that the assumption
$\ref{split assumption}$ is satisfied. Then
\begin{enumerate}[(a)]
  \item The modular reduction induces an additive map $d_\theta:R_0^+(KA) \rightarrow R_0^+(LA)$
  such that $d_\theta([K\tilde{V}])=[k\tilde{V}]$, where $\tilde{V}$
  is an $\mathcal{R}A$-lattice and $[k\tilde{V}]$ is regarded as an element of
  $R_0^+(LA)$ via the identification of $R_0(kA)$ and $R_0(LA)$.
  \item By Proposition $\ref{integrality of pk}$, the image of $\mathfrak{p}_K$ is
  contained in $\emph{Maps}(A,\mathcal{O}[x])$ and we have the following commutative
  diagram
  $$\diagram
  R_0^+(KA) \dto^{d_\theta}  \rto^{\mathfrak{p}_K} &\emph{Maps}(A,\mathcal{O}[x]) \dto^{\tau_\theta}\\
  R_0^+(LA) \rto^{\mathfrak{p}_L} & \mathrm{Maps}(A,L[x]) \enddiagram $$
  where $\tau_\theta:\emph{Maps}(A,\mathcal{O}[x]) \rightarrow \emph{Maps}(A,L[x])$
  is the map induced by $\theta$.
  \item The map $d_\theta$ is uniquely determined by the commutativity of
  the above diagram. In particular, the map $d_\theta$ depends only on
  $\theta$ and not on the choice of $\mathcal{R}$.
\end{enumerate}
The map $d_\theta$ will be called the decomposition map  \index{decomposition map} associated
with the specialization $\theta:\mathcal{O} \rightarrow L$. The
matrix of that map with respect to the bases of $R_0(KA)$ and
$R_0(LA)$ consisting of the classes of the simple modules is called
the decomposition matrix  \index{decomposition matrix} associated with $\theta$.
\end{thedef}
\begin{apod}{Let $\tilde{V}$ be an $\mathcal{R}A$-lattice and $a \in A$.
Let $(a_{ij}) \in \mathrm{M}_n(\mathcal{R})$
 be the matrix describing the action of $a$
on $\tilde{V}$ with respect to a chosen $\mathcal{R}$-basis of
$\tilde{V}$. Due to the properties of modular reduction, the action
of $1 \otimes a \in kA$ on $k\tilde{V}$ is given by the matrix
$(\pi(a_{ij}))$. Then, by definition,
$\mathfrak{p}_L([k\tilde{V}])(a)$ is the characteristic polynomial
of $(\pi(a_{ij}))$. On the other hand, applying $\theta$ (which is
the restriction of $\pi$ to $\mathcal{O}$) to the coefficients of
the characteristic polynomial of $(a_{ij})$ returns $(\tau_\theta
\circ \mathfrak{p}_K)([K\tilde{V}])(a)$. Since the two actions
just described commute, the two polynomials obtained are equal. Thus
the following relation is established:
$$\mathfrak{p}_L([k\tilde{V}])=\tau_\theta \circ \mathfrak{p}_K([K\tilde{V}]) \textrm{ for all }
\mathcal{R}A\textrm{-lattices } \tilde{V} \textrm{    } (\dag) $$

Now let us prove \emph{(a)}. We have to show that the map $d_\theta$ is
well-defined \ie if $\tilde{V},\tilde{V}'$ are two
$\mathcal{R}A$-lattices such that $K\tilde{V}$ and $K\tilde{V}'$
have the same composition factors (counting multiplicities), then
the classes of $k\tilde{V}$ and $k\tilde{V}'$ in $R_0(LA)$ are the
same. Since $[K\tilde{V}]=[K\tilde{V}']$, the endomorphisms
$\rho_{K\tilde{V}}(a)$ and $\rho_{K\tilde{V}'}(a)$ are conjugate for all $a \in A$. 
The equality $(\dag)$ implies that
$$\mathfrak{p}_L([k\tilde{V}])(a)=\mathfrak{p}_L([k\tilde{V}'])(a) \textrm{ for all } a\in A.$$
We have already remarked that, since the assumption $\ref{split
assumption}$ is satisfied, the Brauer-Nesbitt lemma holds for $LA$.
So we conclude that $[k\tilde{V}]=[k\tilde{V}']$, as desired.

Having established the existence of $d_\theta$, we have
$[k\tilde{V}]=d_\theta([K\tilde{V}])$ for any $\mathcal{R}A$-lattice
$\tilde{V}$. Hence $(\dag)$ yields the commutativity of the diagram
in \emph{(b)}.

Finally, by the Brauer-Nesbitt lemma, the map $\mathfrak{p}_L$ is
injective. Hence there exists at most one map which makes the
diagram in \emph{(b)} commutative. This proves \emph{(c)}.}
\end{apod}\
\\
\begin{remark} \emph{ Note that if $\mathcal{O}$ is a discrete valuation ring and $L$ its
residue field, we do not need the assumption $\ref{split
assumption}$ in order to define a decomposition map from $R_0^+(KA)$
to $R_0^+(LA)$ associated with the canonical map $\theta:\mathcal{O}
\rightarrow L$. For a given $KA$-module $V$, there exists an
$A$-lattice $\tilde{V}$ such that $V=K \otimes_\mathcal{O}
\tilde{V}$. The map $d_\theta:R_0^+(KA) \rightarrow
R_0^+(LA),\,[K\tilde{V}]\mapsto[L\tilde{V}]$ is well and
uniquely defined.
For the details of this construction, see \cite{CuRe}, \S16C}.
\end{remark}\\

Recall from Proposition $\ref{integrality of pk}$ that if $V$ is a
$KA$-module, then its character $\chi_V$ restricts to a trace
function $\dot{\chi}_V: A \rightarrow \mathcal{O}$. Now, any linear
map $\lambda:A \rightarrow \mathcal{O}$ induces an $L$-linear map
$$\lambda^L:LA \rightarrow L, 1\otimes a\mapsto \theta(\lambda(a))
(a\in A).$$ It is clear that if $\lambda$ is a trace function, so is
$\lambda^L$. Applying this to $\dot{\chi}_V$ shows that
$\dot{\chi}_V^L$ is a trace function on $LA$. Since character values
occur as coefficients in characteristic polynomials, Theorem
$\ref{existence of decomposition maps}$ implies that
$\dot{\chi}_V^L$ is the character of $d_\theta([V])$. Moreover, for
any simple $KA$-module $V$, we have
$$\dot{\chi}_V^L = \sum_{V'} d_{VV'} \chi_{V'},$$
where the sum is over all simple $LA$-modules $V'$ (up to
isomorphism) and $D=(d_{VV'})$ is the decomposition matrix
associated with $\theta$.

The following result gives a criterion for $d_\theta$ to be trivial.
It is known as ``Tits' deformation theorem''.
For its proof, the reader may refer, for example, to \cite{GePf},
Theorem 7.4.6.

\begin{theorem}\label{Tits} Assume that $KA$ and $LA$ are split. If $LA$ is semisimple, then
$KA$ is also semisimple and the decomposition map $d_\theta$ is an
isomorphism which preserves isomorphism classes of simple modules.
In particular, the map $\emph{Irr}(KA) \rightarrow \emph{Irr}(LA),
\chi \mapsto \dot{\chi}^L$ is a bijection.  \index{Tits' deformation theorem}
\end{theorem}

Finally, if $A$ is symmetric, we can check whether the assumption of
Tits' deformation theorem is satisfied, using the following semisimplicity criterion
(cf.~\cite{GePf}, Thm. 7.4.7).

\begin{theorem}\label{semisimplicity}
Assume that $KA$ and $LA$ are split and that $A$ is
symmetric with symmetrizing form $t$. For any simple $KA$-module
$V$, let $s_V \in \mathcal{O}$ be the Schur element with respect to
the induced symmetrizing form $t^K$ on $KA$. Then $LA$ is semisimple
if and only if $\theta(s_V) \neq 0$ for all $V$.
\end{theorem}

\begin{corollary}\label{injective preserves splitness}
Let $K$ be the field of fractions of $\mathcal{O}$.
Assume that $KA$ is split semisimple and that $A$ is
symmetric with symmetrizing form $t$. If the map $\theta$ is
injective, then $LA$ is split semisimple.
\end{corollary}

\subsection{A variation for Tits' deformation
theorem}

Let us suppose that $\mathcal{O}$ is a Krull ring and
$\theta:\mathcal{O} \rightarrow L$ is a ring homomorphism into a field $L$ such that $L$ is the
field of fractions of $\theta(\mathcal{O})$. We will give a new criterion for the algebra $LA$ to be split semisimple.

\begin{theorem}\label{Lehrer}
Let $K$ be the field of fractions of $\mathcal{O}$. Assume that $KA$ is split semisimple and that $A$ is symmetric with
symmetrizing form $t$. For any simple $KA$-module $V$, let $s_V \in
\mathcal{O}$ be the Schur element with respect to the induced
symmetrizing form $t^K$ on $KA$. If $\mathrm{Ker}\theta$ is a prime
ideal of $\mathcal{O}$ of height $1$, then $LA$ is split semisimple if
and only if $\theta(s_V) \neq 0$ for all $V$.
\end{theorem}
\begin{apod}{If $LA$ is split semisimple, then Theorem
$\ref{semisimplicity}$ implies that $\theta(s_V) \neq 0$ for all
$V$. Now let us denote by $\mathrm{Irr}(KA)$ the set of irreducible
characters of $KA$. If $\chi$ is the character afforded by a simple
$KA$-module $V_\chi$, then $s_\chi:=s_{V_\chi}$. We set
$\mathfrak{q}:=\mathrm{Ker}\theta$ and suppose that $s_\chi \notin
\mathfrak{q}$ for all $\chi \in \mathrm{Irr}(KA)$. Since $KA$ is
split semisimple, it is isomorphic to a product of matrix algebras
over $K$:
$$KA\,\, \cong \prod_{\chi \in \mathrm{Irr}(KA)}\mathrm{End}_{K}(V_\chi)$$
Let us denote by $\pi_\chi:KA \twoheadrightarrow
\mathrm{End}_{K}(V_\chi)$ the projection onto the $\chi$-factor,
such that $\pi:=\prod_{\chi \in \mathrm{Irr}(KA)}\pi_\chi$ is the
above isomorphism. Then $\chi=\mathrm{tr}_{V_\chi} \circ \pi_\chi$,
where $\mathrm{tr}_{V_\chi}$ denotes the standard trace on
$\mathrm{End}_{K}(V_\chi)$.

Let $\mathcal{B},\mathcal{B}'$ be two dual bases of $A$ with respect
to the symmetrizing form $t$. By Lemma $\ref{tau^vee}$, for all $a
\in KA$ and $\chi \in \mathrm{Irr}(KA)$, we have
$$\chi^\vee a= \sum_{b \in \mathcal{B}} \chi(b'a)b.$$
Applying $\pi$ to both sides yields
$$\pi(\chi^\vee) \pi(a)= \sum_{b \in \mathcal{B}} \chi(b'a)\pi(b).$$
By definition of the Schur element, if $\omega_\chi$ denotes the central character associated with $\chi$, then $s_\chi=\omega_\chi(\chi^\vee)=\pi_\chi(\chi^\vee)=
\pi(\chi^\vee)$.
Thus, if $\alpha \in \mathrm{End}_{K}(V_\chi)$, then
$$\pi^{-1}(\alpha)=\frac{1}{s_\chi}\sum_{b \in \mathcal{B}}\mathrm{tr}_{V_\chi}(\pi_\chi(b')\alpha)b. \,\,\,\,\,\,\,\,\,\,(\dag)$$

Since $\mathcal{O}$ is a Krull ring and $\mathfrak{q}$ is a prime
ideal of height $1$ of $\mathcal{O}$, the ring
$\mathcal{O}_{\mathfrak{q}}$ is, by Theorem $\ref{Krull-dvr}$, a
discrete valuation ring. Thanks to Lemma $\ref{realizing modules
over O}$, there exists an $\mathcal{O}_{\mathfrak{q}}A$-lattice
$\tilde{V}_\chi$ such that $K \otimes_{\mathcal{O}_{\mathfrak{q}}}
\tilde{V}_\chi \cong V_\chi$.

Moreover, $1/s_\chi \in \mathcal{O}_\mathfrak{q}$ for all ${\chi \in
\mathrm{Irr}(KA)}$. Due to the relation $(\dag)$, the map $\pi$
induces an isomorphism
$$\mathcal{O}_{\mathfrak{q}}A \,\,\cong
\prod_{\chi \in
\mathrm{Irr}(KA)}\mathrm{End}_{\mathcal{O}_{\mathfrak{q}}}(\tilde{V}_\chi),$$
\ie $\mathcal{O}_{\mathfrak{q}}A$ is isomorphic to a product of matrix algebras
over $\mathcal{O}_{\mathfrak{q}}$. Since
$\mathrm{Ker}\theta=\mathfrak{q}$, the above isomorphism remains
after applying $\theta$ to both sides. Therefore, we obtain that $LA$ is isomorphic to a product
of matrix algebras over $L$ and thus split semisimple.}
\end{apod}

If that is the case, then the assumption of Tits' deformation
theorem is satisfied and there exists a bijection $\mathrm{Irr}(KA)
\leftrightarrow \mathrm{Irr}(LA)$.

\subsection{Symmetric algebras over discrete valuation rings}

From now on, we assume that the following conditions are satisfied:
\begin{itemize}
  \item $\mathcal{O}$ is a discrete valuation ring in $K$ and $K$ is
  perfect; let $v:K\rightarrow \mathbb{Z} \cup \{\infty\}$ be the
  corresponding valuation and $\mathfrak{p}$ the maximal ideal of $\mathcal{O}$.
  \item $KA$ is split semisimple.
  \item $\theta:\mathcal{O} \rightarrow L$ is the canonical map onto
  the residue field $L$ of $\mathcal{O}$.
  \item $A$ is a symmetric algebra with symmetrizing form $t$.
\end{itemize}

We have already seen that we have a well-defined decomposition map
$d_\theta:R_0^+(KA) \rightarrow R_0^+(LA)$. The decomposition matrix
associated with $d_\theta$ is the $|\mathrm{Irr}(KA)| \times
|\mathrm{Irr}(LA)|$ matrix $D=(d_{\chi\phi})$ with non-negative
integer entries such that
$$d_\theta([V_\chi])=\sum_{\phi \in
\mathrm{Irr}(LA)}d_{\chi\phi}[V_\phi'] \,\,\textrm{ for all } \chi \in
\mathrm{Irr}(KA),$$ where $V_\chi$ is a simple $KA$-module with
character $\chi$ and $V_\phi'$ is a simple $LA$-module with
character $\phi$. We sometimes call the characters of $KA$
``ordinary'' and the characters of $LA$ ``modular''. We say that
$\phi \in \mathrm{Irr}(LA)$ is a \emph{modular constituent} of $\chi
\in \mathrm{Irr}(KA)$, if $d_{\chi\phi} \neq 0$.

The rows of $D$ describe the decomposition of $d_\theta([V_\chi])$
in the standard basis of $R_0(LA)$. An interpretation of the columns
is given by the following result (cf.~\cite{GePf}, Theorem 7.5.2),
which is part of Brauer's classical theory of modular
representations (``Brauer reciprocity'').\index{Brauer reciprocity}

\begin{theorem}\label{Brauer reciprocity}
For each $\phi \in \mathrm{Irr}(LA)$, there exists some primitive
idempotent $e_\phi \in A$ such that
$$[e_\phi KA]=\sum_{\chi \in \mathrm{Irr}(KA)}d_{\chi\phi}[V_\chi]
\in R_0^+(KA).$$
\end{theorem}

Let $\phi \in \mathrm{Irr}(LA)$. Consider the map $\psi(\phi):ZKA
\rightarrow K$ defined by
$$\psi(\phi):=\sum_{\chi \in \mathrm{Irr}(KA)}
\frac{d_{\chi\phi}}{s_\chi}\omega_\chi,$$ where $\omega_\chi: ZKA
\twoheadrightarrow K$ is the central character associated with $\chi
\in \mathrm{Irr}(KA)$, as defined in subsection 2.1.4. The next result is due to
Geck and Rouquier (cf.~\cite{GeRo}, Proposition 4.4).

\begin{theorem}\label{Geck-Rouquier}
The map $\psi(\phi)$ restricts to a map $ZA \rightarrow
\mathcal{O}$. In particular,
$$\psi(\phi)(1)=\sum_{\chi \in \mathrm{Irr}(KA)}
\frac{d_{\chi\phi}}{s_\chi} \in \mathcal{O}.$$
\end{theorem}
\begin{apod}{Let us denote by $t^K$ the induced symmetrizing
form on $KA$. If $e_\phi$ is an idempotent as in Theorem
$\ref{Brauer reciprocity}$, then we can define a $K$-linear map
$\lambda_\phi:ZKA \rightarrow K,\,z \mapsto t^K(ze_\phi)$. We claim
that $\lambda_\phi=\psi(\phi)$. Since $KA$ is split semisimple, the
elements $\{\chi^\vee\,|\, \chi \in \mathrm{Irr}(KA)\}$ form a basis
of $ZKA$ (recall that $\chi^\vee$ is the element of $ZKA$ such that
$t^K(\chi^\vee x)=\chi(x)$ for all $x \in KA$). It is, therefore, sufficient to show that
$$\lambda_\phi(\chi^\vee)=\psi(\phi)(\chi^\vee) \textrm{ for all } \chi
\in \mathrm{Irr}(KA).$$ We have
$\psi(\phi)(\chi^\vee)=d_{\chi\phi}1_K$. Now consider the left-hand
side.
$$\begin{array}{rl}
    \lambda_\phi(\chi^\vee) & =t^K(\chi^\vee e_\phi)=\chi(e_\phi) =\mathrm{dim}_K(V_\chi e_\phi)1_K \\
     & \\
     & =\mathrm{dim}_K(\mathrm{Hom}_K(e_\phi
KA,V_\chi))1_K=d_{\chi\phi}1_K.
  \end{array}$$
Hence the above claim is established.

Finally, it remains to observe that since $e_\phi \in A$, the
function $\lambda_\phi$ takes values in $\mathcal{O}$ on all
elements of $A$.}
\end{apod}

Finally, we will treat the block distribution of characters. For
this purpose, we introduce the following notions.

\begin{definition}\label{brauer graph}\
\begin{enumerate}
  \item The Brauer graph  \index{Brauer graph} associated with $A$ has vertices labeled by
  the irreducible characters of $KA$ and an edge joining
  $\chi,\chi' \in \emph{Irr}(KA)$ if $\chi \neq \chi'$ and there
  exists some $\phi \in \emph{Irr}(LA)$ such that $d_{\chi\phi} \neq 0 \neq
  d_{\chi'\phi}$, i.e., there. A connected component of the Brauer graph is called
  a block.  \index{block}
  \item Let $\chi \in \mathrm{Irr}(KA)$. Recall that $0 \neq s_\chi
  \in \mathcal{O}$. Let $\delta_\chi:=v(s_\chi)$, where $v$
  is the given valuation. Then $\delta_\chi$ is called the defect of
  $\chi$  \index{defect of a character} and we have $\delta_\chi \geq 0$ for all $\chi \in \mathrm{Irr}(KA)$.  If $B$ is a block, then
  $\delta_B:=\mathrm{max}\{\delta_\chi \,|\, \chi \in B\}$ is called the
  defect of $B$.  \index{defect of a block}
\end{enumerate}
\end{definition}

Following \cite{Fe}, I.17.9, each block $B$ of $A$ corresponds to a central
primitive idempotent (\ie block-idempotent) $e_B$ of $A$. If $\chi \in B$ and $e_\chi$ is its
corresponding central primitive idempotent in $KA$, then $e_Be_\chi
\neq 0$.

Every $\chi \in \mathrm{Irr}(KA)$ determines a central character
$\omega_\chi:ZKA \rightarrow K$. Since $\mathcal{O}$ is integrally
closed, we have $\omega_\chi(z) \in \mathcal{O}$ for all $z \in ZA$.
We have the following standard result (already presented in subsection 2.1.4)
 relating blocks with central characters.

\begin{proposition}\label{blocks and central characters}
Let $\chi,\chi' \in \emph{Irr}(KA)$. Then $\chi$ and $\chi'$ belong
to the same block of $A$ if and only if
$$\theta(\omega_\chi(z))=\theta(\omega_{\chi'}(z)) \textrm{ for all
} z \in ZA.$$ i.e.,
$$\omega_\chi(z) \equiv \omega_{\chi'}(z) \,\,\emph{mod}\,\mathfrak{p} \textrm{ for all
} z \in ZA.$$
\end{proposition}
\begin{apod}{First assume that $\chi,\chi'$ belong to the same
block of $A$, \ie they belong to a connected component of the Brauer
graph. It is sufficient to consider the case where $\chi,\chi'$ are
directly linked on the Brauer graph, \ie there exists some $\phi \in
\mathrm{Irr}(LA)$ such that $d_{\chi\phi} \neq 0 \neq
d_{\chi'\phi}$. Let $\tilde{V}_\chi$ (resp. $\tilde{V}_{\chi'}$) be an $A$-lattice such that
$K\tilde{V}_\chi$  (resp. $K\tilde{V}_{\chi'}$) 
affords $\chi$ (resp. $\chi'$). Let $z \in ZA$. Then $1 \otimes z$
acts by the scalar $\theta(\omega_\chi(z))$ on every modular
constituent of $k\tilde{V}_\chi$. Similarly, $1 \otimes z$ acts by
the scalar $\theta(\omega_{\chi'}(z))$ on every modular constituent
of $k\tilde{V}_{\chi'}$. Since, by
assumption, $K\tilde{V}_\chi$ and $K\tilde{V}_{\chi'}$ have a
modular constituent in common, we have
$\theta(\omega_\chi(z))=\theta(\omega_{\chi'}(z))$, as desired.

Now assume that $\chi$ belongs to the block $B$ and $\chi'$ to the
block $B'$, with $B \neq B'$. Let $e_B$, $e_{B'}$ be the
corresponding central primitive idempotents. Then
$\omega_\chi(e_B)=1$ and $\omega_{\chi'}(e_B)=0$. Consequently,
$\theta(\omega_\chi(e_B)) \neq \theta(\omega_{\chi'}(e_B))$.}
\end{apod}

Using the above characterization of blocks, we can prove the following result about the characters
of defect $0$, \ie the characters whose Schur elements do not belong to the maximal ideal
$\mathfrak{p}$ of $\mathcal{O}$.

\begin{theorem}\label{blocks of defect 0}
Let $\chi \in \emph{Irr}(KA)$ with $\theta(s_\chi) \neq 0$. Then
$\chi$ is an isolated vertex in the Brauer graph and the
corresponding decomposition matrix is just \emph{(1)}.
\end{theorem}
\begin{apod}{Let $t^K$ be the induced symmetrizing
form on $KA$ and $\hat{t}^K$ the isomorphism from $KA$ to
$\mathrm{Hom}_K(KA,K)$ induced by $t^K$. The irreducible character
$\chi \in \mathrm{Irr}(KA)$ is a trace function on $KA$ and thus we
can define $\chi^\vee:=(\hat{t}^K)^{-1}(\chi) \in ZKA$. Since $\chi$
restricts to a trace function $A \rightarrow \mathcal{O}$, we have
in fact $\chi^\vee \in ZA$. By definition, we have that
$\omega_\chi(\chi^\vee)=s_\chi$ and $\omega_{\chi'}(\chi^\vee)=0$
for any $\chi' \in \mathrm{Irr}(KA), \chi' \neq \chi$. Now assume
that there exists some character $\chi'$ which is linked to $\chi$
in the Brauer graph. Proposition $\ref{blocks and central
characters}$ implies that $0 \neq
\theta(s_\chi)=\theta(\omega_\chi(\chi^\vee))=\theta(\omega_{\chi'}(\chi^\vee))=0$,
which is absurd.

It remains to show that $d_\theta([V_\chi])$ is the class of a
simple module in $R_0^+(LA)$.  By Lemma $\ref{realizing modules over
O}$, there exists a basis of $V_\chi$ and a corresponding
representation $\rho:KA \rightarrow \mathrm{M}_n(K)$ afforded by
$V_\chi$ such that $\rho(a) \in \mathrm{M}_n(\mathcal{O})$ for all
$a \in A$. Let $\mathcal{B}$ be an $\mathcal{O}$-basis of $A$ and
let $\mathcal{B}'$ be its dual with respect to $t$. We have seen in
the proof of Proposition $\ref{integrality of the Schur elements}$ that
$s_\chi=\sum_{b \in \mathcal{B}}\rho(b)_{ij}\rho(b')_{ji}$ for all
$1 \leq i,j \leq n$. All terms in this expression lie in
$\mathcal{O}$. So we can apply the map $\theta$ and obtain a similar
relation for $\theta(s_\chi)$ with respect to the module
$L\tilde{V}_\chi$, where $\tilde{V}_\chi \subseteq V_\chi$ is the
$A$-lattice spanned by the above basis of $V_\chi$. Since
$\theta(s_\chi) \neq 0$, the module $L\tilde{V}_\chi$ is simple
(\cite{GePf}, Lemma 7.2.3).}
\end{apod}

The next result (\cite{MaRo}, Lemma 2.6(b)) is a consequence of Theorems
$\ref{Geck-Rouquier}$ and $\ref{blocks of defect 0}$.

\begin{proposition}\label{Malle-Rouquier}
Assume that the canonical map $ZA \rightarrow ZLA$ is surjective and
let $\chi \in \mathrm{Irr}(KA)$. Then $\chi$ is a block by itself if
and only if $\theta(s_\chi) \neq 0.$
\end{proposition}

\chapter{On Essential algebras}

In this chapter we introduce the notion of  ``essential algebras''. These are symmetric algebras defined over a Laurent polynomial ring whose Schur elements are polynomials of a specific form (described by Definition $\ref{essential algebra}$).  This form gives rise to the definition of the ``essential monomials'' for the algebra. As we have seen in the previous chapter, the Schur elements play an important role in the determination of the blocks of a symmetric algebra. In the following sections, we see how the form of the Schur elements affects the behavior of the blocks of an essential algebra when specialized via different types of morphisms (a morphism associated with a monomial in $3.2$, an adapted morphism in $3.3$, the morphism $I^n$ defined in $3.4$). In particular, in the first two cases, we show that the blocks depend only on the essential monomials for the algebra. 

In the next chapter, we will see that the generic Hecke algebras of complex reflection groups are a particular case of essential algebras. 

\section{Generalities}

Let $R$ be a Noetherian integrally closed domain with field of fractions $K$. Let $\textbf{x}=(x_i)_{0 \leq i \leq m}$ be a set of $m+1$ indeterminates over $R$. We set $A:=R[\textbf{x},\textbf{x}^{-1}]$ the Laurent polynomial ring in these indeterminates. The ring $A$ is also a Noetherian integrally closed domain and thus a Krull ring, by Proposition $\ref{case of Krull}$. The field of fractions of $A$ is $K(\textbf{x})$.
 Let $H$ be an $A$-algebra such that
\begin{itemize}
\item $H$ is free and finitely generated as an $A$-module.
\item There exists a linear form
$t: H \rightarrow A $ which is symmetrizing on $H$.
\item The algebra $K(\textbf{x})H:=K(\textbf{x})\otimes_AH$ is split semisimple.
\end{itemize}

Due to Proposition $\ref{schur elements and idempotents}$, we have
that the symmetrizing form $t$ is of the form
  $$t=\sum_{\chi \in \mathrm{Irr}(K(\textbf{x})H)}
  \frac{1}{s_{\chi}}\chi,$$
  where $s_\chi$ denotes the Schur element of $\chi$ with respect to $t$.
  We know that $s_\chi \in A$ by Proposition $\ref{Schur element belongs to the integral closure}$.
 Moreover, for all $\chi \in \mathrm{Irr}(K(\textbf{x})H)$, the
  block-idempotent of $K(\textbf{x})H$ associated to $\chi$ is
  $e_\chi=\chi^\vee/s_\chi$ (for the definition and properties of $\chi^\vee$, see Lemma $\ref{tau^vee}$).

\begin{definition}\label{essential algebra}
We say that the algebra $H$ is essential  \index{essential algebra} if, for each irreducible character $\chi \in \textrm{Irr}(K(\textbf{\emph{x}})H)$, 
the Schur element $s_\chi$ associated to $\chi$ is an element of $A$ of the form
$$s_\chi = \xi_\chi N_\chi \prod_{i \in I_\chi} \Psi_{\chi,i}(M_{\chi,i})^{n_{\chi,i}}$$
where
\begin{enumerate}[(a)]
    \item $\xi_\chi$ is an element of $R$,
    \item $N_\chi$ is a monomial in $A$, 
    \item $I_\chi$ is an index set,
    \item $(\Psi_{\chi,i})_{i \in I_\chi}$ is a family of  monic polynomials  in $R[x]$ of degree at least $1$, which are irreducible over $K$, prime to $x$ and $x-1$ and whose constant term is a unit in $R$,
    \item $(M_{\chi,i})_{i \in I_\chi}$ is a family of primitive monomials \index{primitive monomial} in $A$, i.e., if $M_{\chi,i} = \prod_{i=0}^m x_i^{a_i}$, then $\textrm{\emph{gcd}}(a_i)=1$,
     \item  ($n_{\chi,i})_{i \in I_\chi}$ is a family of positive integers.
    \end{enumerate}
\end{definition}

Following Theorem $\ref{second irreducible}$, Definition $\ref{essential algebra}$ describes the factorization of $s_\chi$ into irreducible factors in
$K[\textbf{x},\textbf{x}^{-1}]$. This factorization is unique. However, this does not mean that
the monomials $M_{\chi,i}$ appearing in it are unique. Suppose that
$$\Psi_{\chi}(M_{\chi})= u \,\Phi_{\chi}(N_{\chi}),$$
where 
\begin{itemize}
\item $\Psi_{\chi}, \Phi_{\chi}$ are two
$K$-irreducible polynomials as in Definition $\ref{essential algebra}$(d), 
\item $M_{\chi},N_{\chi}$ are two primitive
monomials in $A$,
\item
 $u$ is a unit of $K[\textbf{x},\textbf{x}^{-1}]$.
 \end{itemize}
Since the coefficients of $\Psi_\chi$ belong to $R$ and the constant term of $\Phi_\chi$ is a unit in $R$, we deduce that $u \in R[ \textbf{x},\textbf{x}^{-1}]$.
  Let $\varphi_M$ be an $R$-algebra morphism from
$A$ to a Laurent polynomial
ring $R[\textbf{y},\textbf{y}^{-1}]$ in $m$ indeterminates
associated with the monomial $M_{\chi}$ (recall Definition
$\ref{associated morphism}$). If we apply $\varphi_M$ to the
above equality, we obtain
$$\Psi_{\chi}(1)=\varphi_M(u)\Phi_{\chi}(\varphi_M(N_{\chi})).$$
The morphism $\varphi_M$ sends $u$ to a unit  and $N_\chi$ to a monomial in $K[\textbf{y},\textbf{y}^{-1}]$.
Since $\Psi_{\chi}(1) \in R$ and 
$\Phi_\chi$ is prime to $x$, we deduce that
$\varphi_M(N_{\chi})=1$. By Proposition $\ref{properties of
phi}$(2) and the fact that $N_{\chi}$ is primitive, we obtain that
$$N_{\chi} = M_{\chi}^{\pm 1}.$$
Now, if $M_{\chi} = N_{\chi}$, then $\Psi_{\chi} = \Phi_{\chi}$
and $u=1$. If $N_{\chi} = M_{\chi}^{- 1}$, then $\rm{deg}(\Psi_{\chi})=\rm{deg}(\Phi_{\chi})$ 
and $u = \Psi_{\chi}(0)
M_{\chi}^{\mathrm{deg}(\Psi_{\chi})}$. We summarize the above results into the following proposition.

\begin{proposition}\label{uniqueness up to inversion}
Let $\chi \in \textrm{Irr}(K(\textbf{\emph{x}})H)$
and assume that the Schur element $s_\chi$ associated to $\chi$ has the factorization described in Definition $\ref{essential algebra}$. Then
\begin{enumerate}
\item The monomials $(M_{\chi,i})_{i \in I_\chi}$ are unique up to inversion.
\item The coefficient $\xi_\chi$ is unique up to multiplication by a unit of $R$.
\end{enumerate}
\end{proposition}

Now let $\mathfrak{p}$ be a prime ideal of $R$.
If $\Psi_{\chi,i}(M_{\chi,i})$ is a factor of $s_\chi$ and
$\Psi_{\chi,i}(1) \in \mathfrak{p}$, then the monomial $M_{\chi,i}$
is called $\mathfrak{p}$\emph{-essential for}  \index{p-essential monomial}  $\chi$ in $A$. By
Proposition $\ref{properties of phi}$(1), we have
$$\Psi_{\chi,i}(1) \in \mathfrak{p} \Leftrightarrow
  \Psi_{\chi,i}(M_{\chi,i}) \in \mathfrak{q}_{\chi,i},\,\,\,\,\,(\dag)$$
where $\mathfrak{q}_{\chi,i}:=(M_{\chi,i}-1)A +\mathfrak{p}A.$ 
We deduce that $M_{\chi,i}$ is $\mathfrak{p}$-essential for $\chi$ if and only if $M_{\chi,i}^{-1}$ is  $\mathfrak{p}$-essential for $\chi$.

Recall that,  by Proposition $\ref{primeness of
q}$,
$\mathfrak{q}_{\chi,i}$ is a prime ideal of $A$. Due to the primeness of $\mathfrak{q}_{\chi,i}$
and the proposition above,
 the following
result is an immediate consequence of $(\dag)$.

\begin{proposition}\label{p-essential}
Let $M$ be a primitive monomial in $A$
and $\mathfrak{q}_M:=(M-1)A +\mathfrak{p}A$. Then  $M$ is
$\mathfrak{p}$-essential for $\chi$ in $A$ if and only if
$s_\chi/\xi_\chi \in \mathfrak{q}_M$, where $\xi_\chi$ denotes the
coefficient of $s_\chi$.
\end{proposition}
\begin{apod}{If $M$ is a $\mathfrak{p}$-essential monomial for $\chi$ in $A$, then there exists a polynomial $\Psi$ as in Definition $\ref{essential algebra}$(d) such that $\Psi(M)$ is a factor of $s_\chi$ and $\Psi(1) \in \mathfrak{p}$. By $(\dag)$, we have that $\Psi(M) \in \mathfrak{q}_M$, whence $s_\chi/\xi_\chi \in \mathfrak{q}_M$.
Now, if $s_\chi/\xi_\chi \in \mathfrak{q}_M$, then, since $\mathfrak{q}_M$ is a prime ideal of $A$, there exists
a polynomial $\Phi$ as in Definition $\ref{essential algebra}$(d) and a primitive monomial $N$ such that
$\Phi(N)$ is an irreducible factor of $s_\chi$ and $\Phi(N) \in \mathfrak{q}_M$.
Let $\varphi_M$ be an $R$-algebra morphism from $A$
to a Laurent polynomial ring $B$ in $m$ indeterminates
associated with the monomial $M$. By Proposition  $\ref{properties of phi}$(1), we obtain that
$\varphi_M(\Phi(N))=\Phi(\varphi_M(N)) \in \mathfrak{p}B$. Since $\phi_M(N)$ is a monomial in $B$ and the constant term of $\Phi$ is a unit in $R$, we must have $\varphi_M(N)=1$ and $\Phi(1) \in \mathfrak{p}$. Since $N$ is primitive, Proposition  $\ref{properties of phi}$(2) implies that $N=M^{\pm 1}$. By definition, $M$ is a $\mathfrak{p}$-essential monomial for $\chi$.}
\end{apod}

A primitive monomial $M$ in $A$ is called $\mathfrak{p}$-\emph{essential for} $H$, if there exists an irreducible character $\chi \in \textrm{Irr}(K(\textbf{x})H)$ such that $M$ is $\mathfrak{p}$-essential for $\chi$.

\section{Specialization via morphisms associated with monomials}

From now on, we assume that the algebra $H$ is essential. 
Let $M:=\prod_{i=0}^m x_i^{a_i}$ be a monomial in $A$ such that
$a_i \in \mathbb{Z}$ and gcd$(a_i)=1$, \ie $M$ is primitive. Let $\textbf{y}=(y_j)_{1 \leq j \leq m}$ be a set of $m$ indeterminates over $R$. We set 
 $B:=R[y_1^{\pm 1},\ldots,y_m^{\pm 1}]$ and consider $\varphi_M: A \rightarrow B$ an
$R$-algebra morphism associated with $M$. Let us denote
by $H_{\varphi_M}$ the algebra obtained from $H$
via the specialization $\varphi_M$. The algebra $H_{\varphi_M}$ has also a symmetrizing form defined as the specialization of $t$ via $\varphi_M$.

\begin{proposition}\label{associated morphism preserves splitness}
The algebra $K(\textbf{\emph{y}})H_{\varphi_M}$ is split
semisimple.
\end{proposition}
\begin{apod}{By assumption, the algebra
$K(\textbf{x})H$ is split semisimple. The ring $A$
is a Krull ring and $\mathrm{Ker}\varphi_M=(M-1)A$ is a prime ideal
of height 1 of $A$. Let $\chi \in \textrm{Irr}(K(\textbf{x})H)$. Using the same description for the Schur element $s_\chi$ associated to $\chi$ as in Definition $\ref{essential algebra}$, we obtain that
$$\varphi_M(s_\chi)=\xi_\chi \varphi_M(N_\chi) \prod_{i \in I_\chi} \Psi_{\chi,i}(\varphi_M(M_{\chi,i}))^{n_{\chi,i}}.$$
Since $\varphi_M(N_\chi)$  and $\varphi_M(M_{\chi,i})$ are monomials in $B$ and $\Psi_{\chi,i}(1)\neq 0$ for all $i \in I_\chi$, we deduce that $\varphi_M(s_\chi) \neq 0$.
 Thus, we can apply
Theorem $\ref{Lehrer}$ and obtain that the algebra
$K(\textbf{y})H_{\varphi_M}$ is split semisimple.}
\end{apod}

By ``Tits' deformation theorem'' (Theorem $\ref{Tits}$), the map $\varphi_M$ induces a bijection between the set of irreducible characters of the algebra $K(\textbf{x})H$ and the set of  irreducible characters of the algebra  $K(\textbf{y})H_{\varphi_M}$. The Schur elements of the latter are the
specializations of the Schur elements of
$K(\textbf{x})H$ via $\varphi_M$ and hence, $H_{\varphi_M}$ is also essential.

From now on, whenever we refer to irreducible characters, we mean
irreducible characters of the algebra $K(\textbf{x})H$. Due to the bijection induced by ``Tits' deformation theorem'', we can compare the blocks of $H$ and
$H_{\varphi_M}$ (in terms of partitions of $\mathrm{Irr}(K(\textbf{x})H)$) over suitable rings.

Let $\mathfrak{p}$ be a prime ideal of $R$ and
$\mathfrak{q}_M:=(M-1)A + \mathfrak{p}A$.

\begin{proposition}\label{Aq B}
The  blocks of $B_{\mathfrak{p}B}H_{\varphi_M}$ coincide
with the blocks of $A_{\mathfrak{q}_M}H$.
\end{proposition}
\begin{apod}{Let us denote by $\mathfrak{n}_M$ the kernel of $\varphi_M$, \ie
$\mathfrak{n}_M:=(M-1)A$. By Proposition $\ref{properties of
phi}$(1), we have that $A_{\mathfrak{q}_M}/\mathfrak{n}_M
A_{\mathfrak{q}_M} \cong B_{\mathfrak{p}B}$. Proposition $\ref{quotient blocks}$ implies that the canonical surjection
$A_{\mathfrak{q}_M}H \twoheadrightarrow
(A_{\mathfrak{q}_M}/\mathfrak{n}_M A_{\mathfrak{q}_M})H$
induces a block bijection between these two algebras, whence the desired result.}
\end{apod}$ $\\
\begin{remark}\emph{ Proposition $\ref{Aq B}$ implies that the $\mathfrak{p}$-blocks of
an algebra obtained as the specialization of $H$ via a morphism associated with a monomial $M$
do not depend on the actual choice of the morphism.}
\end{remark}

\begin{proposition}\label{simple inclusion}
If two irreducible characters $\chi$ and $\psi$ are in the same
block of $A_{\mathfrak{p}A}H$, then they are in the same
block of $A_{\mathfrak{q}_M}H$.
\end{proposition}
\begin{apod}{Let $C$ be a block of $A_{\mathfrak{q}_M}H$. Then $\sum_{\chi \in C} e_\chi
\in A_{\mathfrak{q}_M}H \subset
A_{\mathfrak{p}A}H$. Thus $C$ is a union of blocks of
$A_{\mathfrak{p}A}H$.}
\end{apod}
\begin{corollary}\label{union of blocks}
If two irreducible characters $\chi$ and $\psi$ are in the same
block of $A_{\mathfrak{p}A}H$, then they are in the same
block of $B_{\mathfrak{p}B}H_{\varphi_M}$.
\end{corollary}

The corollary above implies that the size of $\mathfrak{p}$-blocks
grows larger as the number of indeterminates becomes smaller.
However, we will now see that the size of blocks remains the same,
if our specialization is not associated with a
$\mathfrak{p}$-essential monomial.

\begin{proposition}\label{not essential for block}
Let $C$ be a block of $A_{\mathfrak{p}A}H$. If $M$ is not
a $\mathfrak{p}$-essential monomial for any $\chi \in C$, then $C$
is a block of $A_{\mathfrak{q}_M}H$ (and thus of
$B_{\mathfrak{p}B}H_{\varphi_M}$).
\end{proposition}
\begin{apod}{ Using the notations of Definition $\ref{essential algebra}$, 
Proposition $\ref{p-essential}$ implies that, for all $\chi \in C$, we have $s_\chi/\xi_\chi \notin
\mathfrak{q}_M$.  Since $C$ is a block of
$A_{\mathfrak{p}A}H$, we have  $$\sum_{\chi \in C} e_\chi
=\sum_{\chi \in C} \frac{\chi^\vee}{s_\chi} \in
A_{\mathfrak{p}A}H.$$ If $\mathcal{B},\mathcal{B}'$ are
two $A$-bases of $H$ dual to each other, then
$\chi^\vee=\sum_{b \in \mathcal{B}}\chi(b)b'$ and the above relation
implies that $$\sum_{\chi \in C}\frac{\chi(b)}{s_\chi} \in
A_{\mathfrak{p}A}, \forall b \in \mathcal{B}.$$ Set $f_b:=\sum_{\chi
\in C}(\chi(b)/s_\chi) \in A_{\mathfrak{p}A}$. Then $f_b$ is of the
form $r_b/(\xi s)$, where 
$$\xi:=\prod_{\chi \in C}\xi_\chi \in
R \textrm{ and }s:=\prod_{\chi \in C} s_\chi/\xi_\chi \in A.$$
Since $\mathfrak{q}_M$ is a prime ideal of $A$, the element $s$, by
assumption, doesn't belong to $\mathfrak{q}_M$. Moreover, we have that
$r_b/\xi \in A_{\mathfrak{p}A}$. By Corollary $\ref{porisma
porismatos}$, there exists $\xi' \in R-\mathfrak{p}$ such
that $r_b/\xi=r_b'/\xi'$ for some $r_b' \in A$. Since
$\mathfrak{q}_M \cap R= \mathfrak{p}$ (cf.~Corollary $\ref{q
intersection zk}$), the element $\xi'$ doesn't belong to the ideal
$\mathfrak{q}_M$ either. Therefore, $f_b=r_b'/(\xi' s) \in
A_{\mathfrak{q}_M} \,\forall b \in \mathcal{B}$, whence
$\sum_{\chi \in C} e_\chi \in A_{\mathfrak{q}_M}H$. Thus,
$C$ is a union of blocks of $A_{\mathfrak{q}_M}H$. Since
the blocks of $A_{\mathfrak{q}_M}H$ are unions of blocks
of $A_{\mathfrak{p}A}H$, by Proposition $\ref{simple inclusion}$, we eventually obtain that $C$ is
a block of $A_{\mathfrak{q}_M}H$.}
\end{apod}
\begin{corollary}\label{not essential for all}
If $M$ is not a $\mathfrak{p}$-essential monomial for any $\chi \in
\mathrm{Irr}(K(\textbf{\emph{x}})H)$, then the blocks of $A_{\mathfrak{q}_M}H$
coincide with the blocks of $A_{\mathfrak{p}A}H$.
\end{corollary}

Of course, all the above results hold for $B$ in the place of $A$, if
we further specialize $B$ (and $H_{\varphi_M}$) via a
morphism associated with a monomial
in $B$.

\section{Specialization via adapted morphisms}

For $r \in \{1,\ldots,m+1\}$, we set $C_r:=R[\textbf{y},\textbf{y}^{-1}]$, 
 where $\textbf{y}:=(y_j)_{r \leq j \leq m}$ is a set of $m+1-r$ indeterminates over $R$. For $r=m+1$, $C_{r}=R$.
 
 From now on, we fix $r\in \{1,\ldots,m+1\}$ and set $\mathcal{R}:=C_r$.
 We  recall that
an $R$-algebra morphism $\varphi:A \rightarrow \mathcal{R}$ is
called \emph{adapted}, if $\varphi={\varphi_r} \circ {\varphi_{r-1}} \circ
\ldots \circ {\varphi_1}$, where $\varphi_i$ is a morphism
associated with a monomial for all $i=1,\ldots,r$. The family
$\mathcal{F}:=\{\varphi_r,\varphi_{r-1},\ldots,\varphi_1\}$ is
called an \emph{adapted family} for $\varphi$ whose \emph{initial morphism} is
$\varphi_1$.

Let $\varphi:A \rightarrow \mathcal{R}$ be an adapted morphism and let us denote by
$H_\varphi$ the algebra obtained as the specialization of
$H$ via $\varphi$. Applying Proposition $\ref{associated
morphism preserves splitness}$ $r$ times, we obtain that the algebra
$K(\textbf{y})H_\varphi$ is split semisimple. By
``Tits' deformation theorem'', the morphism $\varphi$
induces a bijection from the set
$\mathrm{Irr}(K(\textbf{x})H)$ to the set  $\mathrm{Irr}(K(\textbf{y)}H_\varphi)$  of irreducible
characters of $K(\textbf{y})H_{\varphi}$.
Therefore, whenever we refer to irreducible characters, we mean
irreducible characters of the algebra  $K(\textbf{x})H$.

If $M:=\prod_{i=0}^mx_i^{b_i}$ is a monomial
such that $\textrm{gcd}(b_i)=d \in \mathbb{Z}$, then we denote
$M^\circ:=\prod_{i=0}^mx_i^{b_i/d}$. Once more, let $\mathfrak{p}$ be a prime ideal of $R$.

\begin{proposition}\label{blocks of initial monomial}
Let $\varphi:A \rightarrow \mathcal{R}$ be an adapted morphism and
$H_\varphi$ the algebra obtained as the specialization of
$H$ via $\varphi$. If $M$ is a monomial in $A$ such that
$\varphi(M)=1$ and
$\mathfrak{q}_{M^\circ}:=(M^\circ-1)A+\mathfrak{p}A$, then the
blocks of $\mathcal{R}_{\mathfrak{p}\mathcal{R}}H_\varphi$ are unions of
blocks of $A_{\mathfrak{q}_{M^\circ}}H$.
\end{proposition}
\begin{apod}{Let $M$ be a monomial
in $A$ such that $\varphi(M)=1$. Due to proposition $\ref{change
initial}$, there exists an adapted family for $\varphi$ whose initial
morphism $\varphi_1$ is associated with $M^\circ$. Let us denote by
$B$ the image of $\varphi_1$ and by $H_{\varphi_1}$ the
algebra obtained as the specialization of $H$ via
$\varphi_1$. Due to Proposition $\ref{Aq B}$, the blocks of
$B_{\mathfrak{p}B}H_{\varphi_1}$ coincide with the blocks
of $A_{\mathfrak{q}_{M^\circ}}H$. Now, by corollary
$\ref{union of blocks}$, if two irreducible characters belong to the
same $\mathfrak{p}$-block of an essential algebra, then they belong to
the same $\mathfrak{p}$-block of its specialization via a morphism
associated with a monomial. Inductively, we obtain that the blocks
of $\mathcal{R}_{\mathfrak{p}\mathcal{R}}H_\varphi$ are unions of blocks of
$B_{\mathfrak{p}B}H_{\varphi_1}$ and thus of
$A_{\mathfrak{q}_{M^\circ}}H$.}
\end{apod}

We will now state and prove our main result concerning the
$\mathfrak{p}$-blocks of essential algebras.  
Let $\varphi:A \rightarrow \mathcal{R}$ be an adapted morphism and
$H_\varphi$ the algebra obtained as the specialization of
$H$ via $\varphi$. Let $M_1,\ldots,M_k$ be all the
$\mathfrak{p}$-essential monomials for $H$ such that $\phi(M_j)=1$
for all $j=1,\ldots,k$. Note that we can easily find all $\mathfrak{p}$-essential monomials for
$H$ by looking at the unique factorization of its Schur
elements in $K[\textbf{x},\textbf{x}^{-1}]$.
We have $M_j^\circ=M_j$ for all
$j=1,\ldots,k$. Set $\mathfrak{q}_0:=\mathfrak{p}A$,
$\mathfrak{q}_j:=\mathfrak{p}A+(M_j-1)A$ for $j=1,\ldots,k$ and
$\mathcal{Q}:=\{\mathfrak{q}_0,\mathfrak{q}_1,\ldots,\mathfrak{q}_k\}$.

Now let $\mathfrak{q} \in \mathcal{Q}$. If two irreducible
characters $\chi,\psi$ belong to the same block of
$A_\mathfrak{q}H$, we write $\chi \sim_{\mathfrak{q}}
\psi$.

\begin{theorem}\label{main theorem}
Two irreducible characters $\chi,\psi \in \textrm{\emph{Irr}}(K(\textbf{\emph{x}})H)$
are in the same block of $\mathcal{R}_{\mathfrak{p}\mathcal{R}}H_\varphi$ if
and only if there exist a finite sequence
$\chi_0,\chi_1,\ldots,\chi_n \in \textrm{\emph{Irr}}(K(\textbf{\emph{x}})H)$ and a
finite sequence $\mathfrak{q}_{j_1},\ldots,\mathfrak{q}_{j_n} \in
\mathcal{Q}$ such that
\begin{itemize}
  \item $\chi_0=\chi$ and $\chi_n=\psi$,
  \item for all $i$ $(1\leq i \leq n)$, $\chi_{i-1} \sim_{\mathfrak{q}_{j_i}}\chi_i$.
\end{itemize}
\end{theorem}
\begin{apod}{Let us denote by $\sim$ the equivalence relation on
$\textrm{Irr}(K(\textbf{x})H)$ defined as the transitive closure of the relation ``there
exists $\mathfrak{q} \in \mathcal{Q}$ such that $\chi
\sim_\mathfrak{q} \psi$''. We have to show that $\chi$
and $\psi$ are in the same block of
$\mathcal{R}_{\mathfrak{p}\mathcal{R}}H_\varphi$ if and only if $\chi \sim
\psi$.

If $\chi \sim \psi$, then Proposition $\ref{blocks of initial
monomial}$ implies that $\chi$ and $\psi$ are in the same block of
$\mathcal{R}_{\mathfrak{p}\mathcal{R}}H_\varphi$. Now let $C$ be an
equivalence class of $\sim$. We have that $C$ is a union of blocks
of $A_\mathfrak{q}H$, for all $\mathfrak{q} \in \mathcal{Q}$.
Therefore,
$$\sum_{\theta \in C}\frac{\theta^\vee}{s_\theta} \in
A_\mathfrak{q}H, \forall \mathfrak{q} \in \mathcal{Q}.$$
If $\mathcal{B},\mathcal{B}'$ are two dual bases of $H$
with respect to the symmetrizing form $t$, then
$\theta^\vee=\sum_{b\in \mathcal{B}}\theta(b)b'$ and hence,
$$\sum_{\theta \in C}\frac{\theta(b)}{s_\theta} \in
A_\mathfrak{q}, \forall \mathfrak{q} \in \mathcal{Q}, \forall b \in
\mathcal{B}.$$ Using the notations for
$s_\theta$ of Definition $\ref{essential algebra}$, set
$\xi_C:=\prod_{\theta \in C}\xi_\theta$ and  $s_C:=\prod_{\theta \in
C}(s_\theta/\xi_\theta)$. Then, for all $b \in \mathcal{B}$, there
exists an element $r_{C,b} \in A$ such that
$$\sum_{\theta \in C}\frac{\theta(b)}{s_\theta} = \frac{r_{C,b}}{\xi_C s_C}.$$

The element $s_C \in A$ is product of terms (monic $K$-irreducible
polynomials taking values on monomials) which are irreducible in
$K[\textbf{x},\textbf{x}^{-1}]$, due to Theorem
$\ref{second irreducible}$. We also have $s_C \notin \mathfrak{p}A$.

Fix $b \in \mathcal{B}$. The ring $K[\textbf{x},\textbf{x}^{-1}]$ is a unique factorization domain and thus
the quotient $r_{C,b}/s_C$ can be written uniquely in the form
$r/\alpha s$, where
\begin{itemize}
\item $r,s \in A$,
\item $\alpha \in R$,
\item $s$ divides $s_C$ in $A$,
\item $\mathrm{gcd}(r,s)=1$ in $K[\textbf{x},\textbf{x}^{-1}]$.
\end{itemize}
Setting  $\xi:= \alpha \xi_C $, we obtain
$$\frac{r_{C,b}}{\xi_C s_C} = \frac{r}{\xi s}\in
A_\mathfrak{q}, \forall \mathfrak{q} \in \mathcal{Q}.$$ Thus, for
all $\mathfrak{q} \in \mathcal{Q}$, there exist
$r_\mathfrak{q},s_\mathfrak{q} \in A$ with $s_\mathfrak{q} \notin
\mathfrak{q}$ such that
$$\frac{r}{\xi s}=\frac{r_\mathfrak{q}}{s_\mathfrak{q}}.$$

Since $\mathrm{gcd}(r,s)=1$, we obtain that $s$ divides $s_\mathfrak{q}$ in
$K[\textbf{x},\textbf{x}^{-1}]$. However, $s$ divides $s_C$ in
$A$ and hence, $s$ is a product of monic $K$-irreducible polynomials taking
values on monomials. Consequently, at least one of the coefficients
of $s$ is a unit in $A$. Corollary $\ref{my second lemma}$ implies
that $s$ divides $s_\mathfrak{q}$ in $A$. Therefore, $s \notin \mathfrak{q}$
for all $\mathfrak{q} \in \mathcal{Q}$.

Moreover, we have that $r/\xi \in A_{\mathfrak{q}_0} =
A_{\mathfrak{p}A}$. By Corollary $\ref{porisma porismatos}$, there
exist $r' \in A$ and $\xi' \in R-\mathfrak{p}$ such that
$r/\xi=r'/\xi'$. 
Since $\mathfrak{q} \cap R=\mathfrak{p}$, we deduce that 
$\xi' \notin \mathfrak{q}$
for all $\mathfrak{q} \in \mathcal{Q}$.
Thus,
$$\frac{r}{\xi s}=\frac{r'}{\xi' s}\in
A_\mathfrak{q}, \forall \mathfrak{q} \in \mathcal{Q}.$$

Now let us suppose that $\varphi(\xi' s)=\xi' \varphi(s)$ belongs to
$\mathfrak{p}\mathcal{R}$. Since $\xi' \notin \mathfrak{p}$, we must have
$\varphi(s) \in \mathfrak{p}\mathcal{R}$. However, the morphism $\varphi$
sends every monomial of $A$ to a monomial in $\mathcal{R}$. Since $s \notin \mathfrak{p}A$
and $s$ divides $s_C$, $s$ must have a factor of the form $\Psi(M)$, where
\begin{itemize}
  \item $M$ is a primitive monomial in $A$ such that $\varphi(M)=1$,
  \item $\Psi$ is a monic $K$-irreducible polynomial as in Definition $\ref{essential algebra}$(d)
  such that $\Psi(1) \in \mathfrak{p}$.
\end{itemize}
By Proposition $\ref{properties of phi}$(1), we obtain that $s \in \mathfrak{q}_M:=(M-1)A+\mathfrak{p}A$. 

Since $s$ divides $s_C$, $\Psi(M)$ also divides
$s_C$. By definition, $M$ is a $\mathfrak{p}$-essential monomial for some
irreducible character $\theta \in C$. Consequently, $M \in
\{M_1,\ldots,M_k\}$. This contradicts the fact that $s \notin
\mathfrak{q}$ for all $\mathfrak{q} \in \mathcal{Q}$. Therefore,
$\varphi(\xi' s) \notin \mathfrak{p}\mathcal{R}$ and
so we have
$$\frac{\varphi(r_{C,b})}{\varphi(\xi_C s_C)}=
\frac{\varphi(r')}{\varphi(\xi' s)} \in \mathcal{R}_{\mathfrak{p}\mathcal{R}}.$$ 
This equality 
holds for all $b \in \mathcal{B}$ and hence,
$$\sum_{\theta \in C}\frac{\varphi(\theta^\vee)}{\varphi(s_\theta)} \in
\mathcal{R}_{\mathfrak{p}\mathcal{R}}H_\varphi.$$ Thus, $C$ is a union of
blocks of $\mathcal{R}_{\mathfrak{p}\mathcal{R}}H_\varphi$.}
\end{apod}\\
\begin{remark}
\emph{We can obtain Corollary $\ref{not essential for all}$ as an
application of the above theorem for
$\mathcal{Q}=\{\mathfrak{q}_0\}$.}
\end{remark}
\\

Theorem $\ref{main theorem}$ allows us to calculate
the blocks of $\mathcal{R}_{\mathfrak{p}\mathcal{R}}H_\varphi$ for any adapted
morphism $\varphi:A \rightarrow \mathcal{R}$, if we know the blocks of
$A_{\mathfrak{p}A}H$ and the blocks of
$A_{\mathfrak{q}_M}H$ for all $\mathfrak{p}$-essential
monomials $M$. Thus, the study of the blocks of the 
algebra $H$  in a finite number of cases gives us the $\mathfrak{p}$-blocks of all essential algebras obtained via such specializations. 

\section{The map $I^n$}

Let $n$ be an integer, $n\neq 0$, and let $\textbf{y}:=(y_i)_{0 \leq i \leq m}$ be a set of $m+1$ indeterminates over $R$. Define $I^n:A \rightarrow
A':=R[\textbf{y},\textbf{y}^{-1}]$ to be
the $R$-algebra morphism $x_i \mapsto y_i^n$. Obviously,
$I^n$ is injective. Therefore, if we denote by $H'$ the
algebra obtained as the specialization of $H$ via $I^n$,  Corollary $\ref{injective preserves splitness}$ implies that the algebra $K(\textbf{y})H'$ is split semisimple. Again, by 
``Tits' deformation theorem'', the morphism $I^n$
induces a bijection from the set
$\mathrm{Irr}(K(\textbf{x})H)$ to the set  $\mathrm{Irr}(K(\textbf{y})H')$.
Again, let $\mathfrak{p}$ be a prime ideal of $R$.
 
\begin{lemma}\label{power}
The blocks of $A'_{\mathfrak{p}A'}H'$ coincide with
the blocks of $A_{\mathfrak{p}A}H$.
\end{lemma}
\begin{apod}{Since the map $I^n$ is injective, we can consider $A$ as a subring of $A'$ via the
identification $x_i = y_i^n$ for all $i=0,1,\ldots,m$. By Corollary
$\ref{inclusion in localizations}$, we obtain that
$A_{\mathfrak{p}A}$ is contained in $A'_{\mathfrak{p}A'}$ and hence,
the blocks of $A_{\mathfrak{p}A}H$ are unions of blocks of
$A'_{\mathfrak{p}A'}H'$.

Now let $C$ be a block of $A'_{\mathfrak{p}A'}H'$. Since
the field of fractions of $A$ is a splitting field for $H$
(and thus for $H'$), we obtain that
$$\sum_{\chi \in C} e_\chi \in (A'_{\mathfrak{p}A'} \cap
K(\textbf{x}))H'.$$ If $A'_{\mathfrak{p}A'} \cap
K(\textbf{x})=A_{\mathfrak{p}A}$, then $C$ is also a union
of blocks of $A_{\mathfrak{p}A}H$ and we obtain the
desired result.

In order to prove that $A'_{\mathfrak{p}A'} \cap
K(\textbf{x})=A_{\mathfrak{p}A}$, it suffices to show that:
\begin{enumerate}[(a)]
  \item The ring $A_{\mathfrak{p}A}$ is integrally closed.
  \item The ring $A'_{\mathfrak{p}A'}$ is integral over
  $A_{\mathfrak{p}A}$.
\end{enumerate}

Since the ring $A$ is integrally closed, (a) is immediate by
Corollary $\ref{integrally closed localization}$. For (b), we have that $A'$ is integral over $A$, since
$y_i^n-x_i=0$, for all $i=0,1,\ldots,m$. Moreover, $A'$ is
integrally closed and thus the integral closure of $A$ in
$K(\textbf{y})$. The only prime ideal of $A'$ lying over
$\mathfrak{p}A$ is $\mathfrak{p}A'$. Following Corollary
$\ref{one prime lying over}$, we obtain that the integral closure of
$A_{\mathfrak{p}A}$ in $K(\textbf{y})$ is
$A'_{\mathfrak{p}A'}$. Thus, $A'_{\mathfrak{p}A'}$ is integral over
$A_{\mathfrak{p}A}$.}
\end{apod}

We can consider $I^n$ as an endomorphism of $A$ and denote it by
$I^n_A$. If $k$ is another non-zero integer, then ${I^k_A} \circ
I^n_A = {I^n_A}\circ I^k_A = I^{kn}_A$. If now $\varphi:A
\rightarrow \mathcal{R}$ is an adapted morphism, we can easily check that
$\varphi \circ I^n_A = {I^n_\mathcal{R}} \circ \varphi$. Abusing notation, we
write $\varphi \circ I^n = {I^n} \circ \varphi$.

\begin{corollary}\label{i gives same blocks}
Let $\varphi:A \rightarrow \mathcal{R}$ be an adapted morphism and
$H_\varphi$ the algebra obtained as the specialization of
$H$ via $\varphi$. Let $\phi: A \rightarrow \mathcal{R}$ be an
$R$-algebra morphism such that ${I^\alpha}\circ
\varphi={I^\beta}\circ \phi$ for some $\alpha,\beta \in
\mathbb{Z}\setminus\{0\}$. If $H_\phi$ is the algebra obtained as
the specialization of $H$ via $\phi$, then the blocks of
$\mathcal{R}_{\mathfrak{p}\mathcal{R}}H_\phi$ coincide with the blocks of
$\mathcal{R}_{\mathfrak{p}\mathcal{R}}H_\varphi$ and we can use Theorem
$\ref{main theorem}$ to calculate them.
\end{corollary}

\chapter{On Hecke algebras}

We will start this chapter by  giving the definition  and the classification of \emph{complex reflection groups}. We will also define the \emph{braid group} and the \emph{pure braid group} associated to a complex reflection group. We will then introduce the \emph{generic Hecke algebra} of a complex reflection group, which is a quotient of the group algebra of the associated braid group defined over a Laurent polynomial in a finite number of indeterminates. Under certain assumptions, which have been verified for all but a finite number of cases, we prove (Theorem $\ref{Schur element generic}$) that the generic Hecke algebras of complex reflection groups are essential. Therefore, all results obtained in Chapter $3$ apply to the case of the generic Hecke algebras. 

A \emph{cyclotomic Hecke algebra} is obtained from the generic Hecke algebra via  a \emph{cyclotomic specialization} (Definition $\ref{specialization}$).  We prove (Theorem $\ref{cyclotomic}$) that any cyclotomic specialization is essentially an adapted morphism. Thus, we can use Theorem $\ref{main theorem}$ in order to obtain the \emph{Rouquier blocks} of a cyclotomic Hecke algebra (\ie its blocks over the \emph{Rouquier ring}, defined in section  4.4), which are a substitute for the families of characters that can be applied to all complex reflection groups. We will see that the Rouquier blocks have the property of \emph{semi-continuity}, thus depending only on some ``essential'' hyperplanes for the group, which are determined by the generic Hecke algebra. 

The theory developed in this chapter will allow us to determine the Rouquier blocks of the cyclotomic Hecke algebras of all (irreducible) complex reflection groups in the next and final chapter.

\section{Complex reflection groups and associated braid groups}

Let $\mu_\infty$ be the group of all the roots of unity in
$\mathbb{C}$ and $K$ a number field contained in
$\mathbb{Q}(\mu_\infty)$. We denote by $\mu(K)$ the group of all the
roots of unity of $K$. For every integer $d>1$, we set
$\zeta_d:=\mathrm{exp}(2\pi i/d)$ and denote by $\mu_d$ the group of
all the $d$-th roots of unity. Let $V$ be a $K$-vector space of
finite dimension $r$.

\subsection{Complex reflection groups}

\begin{definition}\label{pseudo-reflection}
A pseudo-reflection \index{reflection} is a non-trivial element $s$
of $\mathrm{GL}(V)$ which acts trivially on a hyperplane, called the
reflecting hyperplane of $s$.

If $W$ is a finite subgroup of $\mathrm{GL}(V)$ generated by
pseudo-reflections, then $(V,W)$ is called a $K$-reflection group \index{reflection group} of
rank $r$.
\end{definition}

We have the following classification of complex reflection groups,
also known as the ``Shephard-Todd classification''. For more details
about the classification, one may refer to \cite{ShTo}.\index{Shephard-Todd classification}

\begin{theorem} Let $(V,W)$ be an irreducible complex
reflection group (i.e., $W$ acts irreducibly on $V$). Then one of
the following assertions is true:
\begin{itemize}
  \item There exist non-zero integers $d,e,r$ such
  that $(V,W) \cong G(de,e,r)$, where $G(de,e,r)$ is the group of all 
  $r \times r$ monomial matrices with non-zero entries in $\mu_{de}$ such that the product of all non-zero
  entries lies in $\mu_d$.
  \item $(V,W)$ is isomorphic to one of the $34$ exceptional groups
  $G_n$ $(n=4,\ldots,37)$.
\end{itemize}
\end{theorem}
\begin{remark} \emph{Among the irreducible complex reflection groups, we encounter the irreducible real reflection groups. In particular, we have:
\begin{itemize}
\item $G(1,1,r) \cong A_{r-1}$ for $r \geq 2$, 
\item $G(2,1,r) \cong B_r$ (or $C_r$) for $r \geq 2$,
\item $G(2,2,r) \cong D_r$ for $r \geq 4$,
\item $G(e,e,2)\cong I_2(e)$, where $I_2(e)$ denotes the dihedral group of order $2e$,
\item $G_{23}=H_3$, $G_{28}=F_4$, $G_{30}=H_4$, $G_{35}=E_6$, $G_{36}=E_7$, $G_{37}=E_8$.
\end{itemize}}
\end{remark}\

The following theorem has been proved (using a case by case
analysis) by Benard (\cite{Ben}) and Bessis (\cite{Bes1}) and
generalizes a well known result for Weyl groups.

\begin{thedef}\label{field of definition}
Let $(V,W)$ be a reflection group. Let $K$ be the field generated by
the traces on $V$ of all the elements of $W$. Then all irreducible
$KW$-representations are absolutely irreducible i.e., $K$ is a
splitting field for $W$. The field $K$ is called the field of
definition of the group $W$. \index{field of definition}
\end{thedef}

\begin{itemize}
  \item If $K \subseteq \mathbb{R}$, then $W$ is a (finite) Coxeter
  group.
  \item If $K=\mathbb{Q}$, then $W$ is a Weyl group.
\end{itemize}

\subsection{Braid groups associated to complex reflection groups}
 
For all definitions and results about braid groups we
follow \cite{BMR}. Note that for a given topological space $X$ and a point $x_0 \in X$, we denote
by $\Pi_1(X,x_0)$ the fundamental group with base point $x_0$.

Let $V$ be a $K$-vector space of
finite dimension $r$.
Let $W$ be a finite
subgroup of $\mathrm{GL}(V)$ generated by pseudo-reflections and
acting irreducibly on $V$. We denote by $\mathcal{A}$ the set of its
reflecting hyperplanes. We define the \emph{regular variety} \index{regular variety}
$V^{\textrm{reg}}:= \mathbb{C} \otimes V-\bigcup_{H \in
\mathcal{A}}\mathbb{C} \otimes H$. For $x_0 \in V^{\textrm{reg}}$,
we define $P:=\Pi_1(V^{\textrm{reg}},x_0)$ the \emph{pure braid
group} \index{pure braid group} (at $x_0$) associated with $W$. If $p:V^{\textrm{reg}}
\rightarrow V^{\textrm{reg}}/W$ denotes the canonical surjection, we
define $B:=\Pi_1(V^{\textrm{reg}}/W,p(x_0))$ the \emph{braid group} \index{braid group}
(at $x_0$) associated with $W$.

The projection $p$ induces a surjective map $B\twoheadrightarrow W,
\sigma \mapsto \bar{\sigma}$ as follows: Let
$\tilde{\sigma}:[0,1]\rightarrow V^{\textrm{reg}}$ be a path in
$V^{\textrm{reg}}$ such that $\tilde{\sigma}(0)=x_0$, which lifts
$\sigma$. Then $\bar{\sigma}$ is defined by the equality
$\bar{\sigma}(x_0)=\tilde{\sigma}(1)$. Note that the map $\sigma
\mapsto \bar{\sigma}$ is an anti-morphism.

Denoting by $W^\mathrm{op}$ the group opposite to $W$, we have the
following short exact sequence
$$1\rightarrow P\rightarrow B\rightarrow
W^\mathrm{op}\rightarrow1,$$ where the map $B\rightarrow
W^\mathrm{op}$ is defined by $\sigma \mapsto \bar{\sigma}.$

Now, for every hyperplane $H \in \mathcal{A}$, we set $e_H$ the
order of the group $W_H$, where $W_H$ is the subgroup of $W$ formed
by $\mathrm{id}_V$ and all the reflections fixing the hyperplane $H$. The group
$W_H$ is cyclic: if $s_H$ denotes an element of $W_H$ with
determinant $\zeta_H:=\zeta_{e_H}$, then $W_H=<s_H>$ and $s_H$ is
called a \emph{distinguished reflection} \index{distinguished reflection} in $W$.

Let $L_H:=\mathrm{Im}(s-\mathrm{id}_V)$. Then, for all $x \in V$, we
have $x=\mathrm{pr}_H(x)+\mathrm{pr}_{L_H}(x)$ with
$\mathrm{pr}_H(x) \in H$ and $\mathrm{pr}_{L_H}(x) \in L_H$. Thus,
$s_H(x)=\mathrm{pr}_H(x)+\zeta_H\mathrm{pr}_{L_H}(x)$.

If $t \in \mathbb{R}$, we set $\zeta_H^t:=\mathrm{exp}(2\pi it/e_H)$
and we denote by $s_H^t$ the element of $\mathrm{GL}(V)$ (a
pseudo-reflection if $t\neq 0$) defined by
$$s_H^t(x):=\mathrm{pr}_H(x)+\zeta_H^t\mathrm{pr}_{L_H}(x).$$

For $x \in V$, we denote by $\sigma_{H,x}$ the path in $V$ from $x$
to $s_H(x)$ defined by
$$\sigma_{H,x}:[0,1] \rightarrow V,\,\, t \mapsto s_H^t(x).$$

Let $\gamma$ be a path in $V^\mathrm{reg}$ with initial point $x_0$
and terminal point $x_H$. Then $\gamma^{-1}$ is the path in
$V^\mathrm{reg}$ with initial point $x_H$ and terminal point $x_0$
such that
$$\gamma^{-1}(t)=\gamma(1-t) \textrm{ for all } t \in [0,1].$$
Thus, we can define the path $s_H(\gamma^{-1}):t \mapsto
s_H(\gamma^{-1}(t))$, which goes from $s_H(x_H)$ to $s_H(x_0)$ and
lies also in $V^\mathrm{reg}$, since for all $x \in V^\mathrm{reg}$,
$s_H(x) \in V^\mathrm{reg}$ (If $s_H(x) \notin V^\mathrm{reg}$, then
$s_H(x)$ must belong to a hyperplane $H'$. If $s_{H'}$ is a
distinguished pseudo-reflection with reflecting hyperplane $H'$,
then $s_{H'}(s_H(x))=s_H(x)$ and ${s_H}^{-1}(s_{H'}(s_H(x)))=x$.
However, ${s_H}^{-1}s_{H'}s_H$ is a reflection and $x$ belongs to
its reflecting hyperplane, $s_H^{-1}(H')$. This contradicts the fact
that $x$ belongs to $V^\mathrm{reg}$.). Now we define a path from
$x_0$ to $s_H(x_0)$ as follows:
$$\sigma_{H,\gamma}:=s_H(\gamma^{-1}(t)) \cdot \sigma_{H,x_H} \cdot
\gamma$$

If $x_H$ is chosen ``close to $H$ and far from the other reflecting
hyperplanes'', the path $\sigma_{H,\gamma}$ lies in $V^\mathrm{reg}$
and its homotopy class does not depend on the choice of $x_H$. The
element it induces in the braid group $B$, $\textbf{s}_{H,\gamma}$,
is a distinguished braid reflection around the image of $H$ in
$V^\mathrm{reg}/W$.

\begin{proposition}\label{braid reflections}\
\begin{enumerate}
  \item The braid group $B$ is generated by the distinguished braid
   reflections around the images of the hyperplanes $H \in
   \mathcal{A}$ in $V^\mathrm{reg}/W$.
  \item The image of $\emph{\textbf{s}}_{H,\gamma}$ in $W$ is $s_H$.
  \item Whenever $\gamma'$ is a path in $V^\mathrm{reg}$ from $x_0$
   to $x_H$, if $\lambda$ denotes the loop in $V^\mathrm{reg}$ defined
   by $\lambda:=\gamma'^{-1}\gamma$, then
   $$\sigma_{H,\gamma'}=s_H(\lambda) \cdot \sigma_{H,\gamma} \cdot
   \lambda^{-1}.$$
   In particular, $\emph{\textbf{s}}_{H,\gamma}$ and
  $\emph{\textbf{s}}_{H,\gamma}$ are conjugate in $P$.
  \item The path $\prod_{j=e_H-1}^{j=0}\sigma_{H,s_H^j(\gamma)}$,
  a loop in $V^\mathrm{reg}$, induces the element
  $\emph{\textbf{s}}_{H,\gamma}^{e_H}$ in the braid group $B$ and belongs to the pure braid
  group $P$. It is a distinguished braid reflection around $H$ in
  $P$.
\end{enumerate}
\end{proposition}

\begin{definition}\label{s-distinguished braid reflection}
Let $s$ be a distinguished pseudo-reflection in $W$ with reflecting
hyperplane $H$. An $s$-distinguished braid reflection or monodromy
generator is a distinguished braid reflection $\emph{\textbf{s}}$
around the image of $H$ in $V^\mathrm{reg}/W$ such that
$\bar{\emph{\textbf{s}}}=s$.
\end{definition}

\begin{definition}\label{pi}
Let $x_0 \in V^\mathrm{reg}$ as before. We denote by \textbf{$\tau$}
the element of $P$ defined by the loop  $t \mapsto
x_0\mathrm{exp}(2\pi it).$
\end{definition}

\begin{lemma}\label{pi in ZP}
We have $\tau \in ZP$.
\end{lemma}

\begin{thedef}\label{length function}
Given $\mathcal{C} \in \mathcal{A}/W$, there exists a unique length
function $l_\mathcal{C}:B \rightarrow \mathbb{Z}$ defined as
follows: if $b=\emph{\textbf{s}}_1^{n_1} \cdot
\emph{\textbf{s}}_2^{n_2} \cdot \cdot \cdot
\emph{\textbf{s}}_m^{n_m}$ where (for all $j$) $n_j \in \mathbb{Z}$
and $\emph{\textbf{s}}_j$ is a distinguished braid reflection around
an element of $\mathcal{C}_j$, then
$$l_\mathcal{C}(b)=\sum_{\{j\,|\,\mathcal{C}_j=\mathcal{C}\}}n_j.$$
The length function $l:B \rightarrow \mathbb{Z}$ is defined,
for all $b \in B$, as
$$l(b)=\sum_{\mathcal{C} \in \mathcal{A}/W}l_\mathcal{C}(b).$$
\end{thedef}

We say that $B$ has an \emph{Artin-like} presentation \index{Artin-like presentation} (cf.~\cite{Op2},
5.2), if it has a presentation of the form
$$<\textbf{s} \in \textbf{S} \,|\, \{\textbf{v}_i=\textbf{w}_i\}_{i \in
I}>,$$ where $\textbf{S}$ is a finite set of distinguished braid
reflections and $I$ is a finite set of relations which are
multi-homogeneous, \ie such that, for each $i$, $\textbf{v}_i$ and
$\textbf{w}_i$ are positive words in elements of $\textbf{S}$ (and
hence, for each $\mathcal{C} \in \mathcal{A}/W$, we have
$l_\mathcal{C}(\textbf{v}_i)=l_\mathcal{C}(\textbf{w}_i)$).

The following result by Bessis (\cite{Bes3}, Theorem 0.1) shows that any
braid group has an Artin-like presentation.

\begin{theorem}\label{bessis}
Let $W$ be a complex reflection group with associated braid group
$B$. Then there exists a subset
$\emph{\textbf{S}}=\{\emph{\textbf{s}}_1,\ldots,\emph{\textbf{s}}_n\}$
of $B$ such that
\begin{enumerate}
  \item The elements $\emph{\textbf{s}}_1,\ldots,\emph{\textbf{s}}_n$ are
  distinguished braid reflection and therefore, their images $s_1,\ldots,s_n$ in
  $W$ are distinguished reflections.
  \item The set $\emph{\textbf{S}}$ generates $B$ and therefore,
  $S:=\{s_1,\ldots,s_n\}$ generates $W$.
  \item There exists a set $\mathcal{R}$ of relations of the form
  $\emph{\textbf{w}}_1=\emph{\textbf{w}}_2$, where $\emph{\textbf{w}}_1$ and
  $\emph{\textbf{w}}_2$ are positive words of equal length in the elements
  of $\emph{\textbf{S}}$, such that $<\emph{\textbf{S} }\,|\, \mathcal{R}>$ is a
  presentation of $B$.
  \item Viewing now $\mathcal{R}$ as a set of relations in $S$, the
  group $W$ is presented by
  $$<S \,|\, \mathcal{R}; (\forall s \in S)(s^{e_s}=1)>,$$
  where $e_s$ denotes the order of $s$ in $W$.
\end{enumerate}
\end{theorem}

\section{Generic Hecke algebras}

Let $K,V,W,\mathcal{A},P,B$ be defined as in the previous section.
For every orbit $\mathcal{C}$ of $W$ on $\mathcal{A}$, we set
$e_{\mathcal{C}}$ the common order of the subgroups $W_H$, where $H$
is any element of $\mathcal{C}$ and $W_H$ the subgroup formed by $\mathrm{id}_V$
and all the reflections fixing the hyperplane $H$.

We choose a set of indeterminates
$\textbf{u}=(u_{\mathcal{C},j})_{(\mathcal{C} \in
\mathcal{A}/W)(0\leq j \leq e_{\mathcal{C}}-1)}$ and we denote by
$\mathbb{Z}[\textbf{u},\textbf{u}^{-1}]$ the Laurent polynomial ring
in all the indeterminates $\textbf{u}$. We define the \emph{generic
Hecke algebra} \index{generic Hecke algebra} $\mathcal{H}$ of $W$ to be the quotient of the group
algebra $\mathbb{Z}[\textbf{u},\textbf{u}^{-1}]B$ by the ideal
generated by the elements of the form
$$(\textbf{s}-u_{\mathcal{C},0})(\textbf{s}-u_{\mathcal{C},1}) \ldots (\textbf{s}-u_{\mathcal{C},e_{\mathcal{C}}-1}),$$
where $\mathcal{C}$ runs over the set $\mathcal{A}/W$ and
$\textbf{s}$ runs over the set of monodromy generators around the
images in $V^{\textrm{reg}}/W$ of the elements of the hyperplane
orbit $\mathcal{C}$.

\begin{px}
\small{\emph{Let $W:=G_2=<s,t \,|\, \,ststst=tststs, s^2=t^2=1>$ be the dihedral group of order $12$.  Then the generic Hecke algebra of
$W$ is defined over the Laurent polynomial ring in four indeterminates
$\mathbb{Z}[u_0,u_0^{-1},u_1,u_1^{-1},w_0,w_0^{-1},w_1,w_1^{-1}]$ and
can be presented as follows:
$$\begin{array}{rll}
   \mathcal{H}(G_2)=<S,T \,\,|&STSTST=TSTSTS, &(S-u_0)(S-u_1)=0, \\
                                &         &(T-w_0)(T-w_1)=0>.
  \end{array}$$}}
\end{px}

\begin{px}
\small{\emph{Let $W:=G_4=<s,t \,|\, sts=tst, s^3=t^3=1>$. Then $s$
and $t$ are conjugate in $W$ and their reflecting hyperplanes belong
to the same orbit of $W$ on $\mathcal{A}$. The generic Hecke algebra of
$W$ can be presented as follows:
$$\begin{array}{rll}
   \mathcal{H}(G_4)=<S,T \,\,|&STS=TST, &(S-u_0)(S-u_1)(S-u_2)=0, \\
                                &         &(T-u_0)(T-u_1)(T-u_2)=0>.
  \end{array}$$}}
\end{px}

We make some assumptions for the generic Hecke algebra $\mathcal{H}$. Note that
they have been verified for all but a finite number of irreducible
complex reflection groups (\cite{BMM2}, remarks before 1.17, $\S$ 2;
\cite{GIM}).

\begin{ypoth}\label{ypo}
The algebra $\mathcal{H}$ is a free
$\mathbb{Z}[\textbf{\emph{u}},\textbf{\emph{u}}^{-1}]$-module of
rank $|W|$. Moreover, there exists a linear form
$t:\mathcal{H}\rightarrow
\mathbb{Z}[\textbf{\emph{u}},\textbf{\emph{u}}^{-1}]$ with the
following properties:
\begin{enumerate}
    \item $t$ is a symmetrizing form on $\mathcal{H}$, i.e.,
     $t(hh')=t(h'h)$ for all $h,h' \in \mathcal{H}$ and the map
     $$\begin{array}{cccc}
     \hat{t}: & \mathcal{H} & \rightarrow & \textrm{\emph{Hom}}(\mathcal{H},\mathbb{Z}[\textbf{\emph{u}},\textbf{\emph{u}}^{-1}]) \\
              & h & \mapsto & (h' \mapsto t(hh'))
     \end{array}$$
     is an isomorphism.
    \item Via the specialization $u_{\mathcal{C},j} \mapsto
     \zeta_{e_\mathcal{C}}^j$, the form $t$ becomes the canonical
     symmetrizing form on the group algebra $\mathbb{Z}_K[W]$.
    \item If we denote by $\alpha \mapsto \alpha^*$ the automorphism of
     $\mathbb{Z}[\emph{\textbf{u}},\emph{\textbf{u}}^{-1}]$ consisting of the
     simultaneous inversion of the indeterminates, then for all $b \in B$, we
     have
          $$t(b^{-1})^*=\frac{t(b\tau)}{t(\tau)},$$
     where $\tau$ is the (central) element of $P$ defined by the loop  $t \mapsto
x_0\mathrm{exp}(2\pi it).$\end{enumerate}
\end{ypoth}

We know that the form $t$ is unique (\cite{BMM2}, 2.1). From now on,
we suppose that the assumptions $\ref{ypo}$ are satisfied. Then
we have the following result by G. Malle (\cite{Ma4}, 5.2).

\begin{theorem}\label{Semisimplicity Malle}
Let $\textbf{\emph{v}}=(v_{\mathcal{C},j})_{(\mathcal{C} \in
\mathcal{A}/W)(0\leq j \leq e_{\mathcal{C}}-1)}$ be a set of
$\sum_{\mathcal{C} \in \mathcal{A}/W}e_{\mathcal{C}}$ indeterminates
such that, for every $\mathcal{C},j$, we have
$v_{\mathcal{C},j}^{|\mu(K)|}=\zeta_{e_\mathcal{C}}^{-j}u_{\mathcal{C},j}$.
Then the $K(\textbf{\emph{v}})$-algebra
$K(\textbf{\emph{v}})\mathcal{H}$ is split semisimple.
\end{theorem}

By ``Tits' deformation theorem'' (Theorem $\ref{Tits}$), it follows
that the specialization $v_{\mathcal{C},j}\mapsto 1$ induces a
bijection $\chi_{\textbf{v}} \mapsto \chi$ from the set
$\mathrm{Irr}(K(\textbf{v})\mathcal{H})$ of absolutely irreducible characters of $K(\textbf{v})\mathcal{H}$ to the
set $\mathrm{Irr}(W)$ of absolutely
irreducible characters of $W$, such that the
following diagram is commutative $$\begin{array}{rccc}
  \chi_\textbf{v} : & \mathcal{H} & \rightarrow & \mathbb{Z}_K[\textbf{v},\textbf{v}^{-1}] \\
  & \downarrow &  & \downarrow \\
  \chi: & \mathbb{Z}_K[W] &\rightarrow &\mathbb{Z}_K .
\end{array}$$

Since the assumptions $\ref{ypo}$ are satisfied and the algebra
$K(\textbf{v})\mathcal{H}$ is split semisimple, we can define the
Schur element $s_{\chi}({\textbf{v}})$ for every irreducible
character $\chi_{\textbf{v}}$ of $K(\textbf{v})\mathcal{H}$ with
respect to the symmetrizing form $t$. 
The following result describes the form of the Schur elements associated
to the irreducible characters of $K(\textbf{v})\mathcal{H}$.

\begin{theorem}\label{Schur element generic}
The Schur element $s_\chi(\textbf{\emph{v}})$ associated to the irreducible
character $\chi_{\textbf{\emph{v}}}$ of
$K(\textbf{\emph{v}})\mathcal{H}$ is an element of
$\mathbb{Z}_K[\textbf{\emph{v}},\textbf{\emph{v}}^{-1}]$ of the form
$$s_\chi({\textbf{\emph{v}}})=\xi_\chi N_\chi \prod_{i \in I_\chi} \Psi_{\chi,i}(M_{\chi,i})^{n_{\chi,i}}$$
where
\begin{enumerate}[(a)]
    \item $\xi_\chi$ is an element of $\mathbb{Z}_K$,
    \item $N_\chi= \prod_{\mathcal{C},j} v_{\mathcal{C},j}^{b_{\mathcal{C},j}}$ is a monomial in $\mathbb{Z}_K[\textbf{\emph{v}},\textbf{\emph{v}}^{-1}]$
          with $\sum_{j=0}^{e_\mathcal{C}-1}b_{\mathcal{C},j}=0$
          for all $\mathcal{C} \in \mathcal{A}/W$,
    \item $I_\chi$ is an index set,
    \item $(\Psi_{\chi,i})_{i \in I_\chi}$ is a family of $K$-cyclotomic polynomials in one variable
           (i.e., minimal polynomials of the roots of unity over $K$),
    \item $(M_{\chi,i})_{i \in I_\chi}$ is a family of monomials in $\mathbb{Z}_K[\textbf{\emph{v}},\textbf{\emph{v}}^{-1}]$ such that
          if $M_{\chi,i} = \prod_{\mathcal{C},j} v_{\mathcal{C},j}^{a_{\mathcal{C},j}}$,
          then $\textrm{\emph{gcd}}(a_{\mathcal{C},j})=1$
          and $\sum_{j=0}^{e_\mathcal{C}-1}a_{\mathcal{C},j}=0$
          for all $\mathcal{C} \in \mathcal{A}/W$,
    \item ($n_{\chi,i})_{i \in I_\chi}$ is a family of positive integers.
\end{enumerate}
\end{theorem}
\begin{apod}{By Proposition $\ref{Schur element belongs to the integral closure}$,
we have that $s_\chi(\textbf{v}) \in
\mathbb{Z}_K[\textbf{v},\textbf{v}^{-1}]$. The rest is a case by
case analysis: Let us first consider the group $G(d,1,r)$.
The Schur elements of $\mathcal{H}(G(d,1,r))$ have been calculated independently by Geck, Iancu and Malle (\cite{GIM}) and by Mathas (\cite{Mat}). Following Theorem $\ref{schur elements of ArikiKoike}$, they are obviously of the desired form. Moreover, in the Appendix we give the generic Schur elements for the groups $G(2d,2,2)$, $G_7$, $G_{11}$, $G_{19}$, $G_{26}$, $G_{32}$ (calculated by Malle in \cite{Ma2} and \cite{Ma5}) and $F_4$ (calculated by Lusztig in \cite{Lu79b}) and show that they are of the form described above.
In the Appendix, we also give the specializations of the parameters which make
\begin{itemize}
\item $\mathcal{H}(G(de,1,r))$ the twisted symmetric algebra of the cyclic group $C_e$ over
$\mathcal{H}(G(de,e,r))$ in the case where $r>2$ or $r=2$ and $e$ is odd.
\item  $\mathcal{H}(G(de,2,2))$ the twisted symmetric algebra of the cyclic group $C_{e/2}$ over
$\mathcal{H}(G(de,e,2))$ in the case where $e$ is even.
\item $\mathcal{H}(G_7)$ the twisted symmetric algebra of some finite cyclic group over $\mathcal{H}(G_4)$, $\mathcal{H}(G_5)$ and $\mathcal{H}(G_6)$.
\item $\mathcal{H}(G_{11})$ the twisted symmetric algebra of some finite cyclic group over $\mathcal{H}(G_8)$, $\mathcal{H}(G_9)$, $\mathcal{H}(G_{10})$,
$\mathcal{H}(G_{12})$, $\mathcal{H}(G_{13})$, $\mathcal{H}(G_{14})$ and
$\mathcal{H}(G_{15})$.
\item $\mathcal{H}(G_{19})$ the twisted symmetric algebra of some finite cyclic group over
$\mathcal{H}(G_{16})$, $\mathcal{H}(G_{17})$, $\mathcal{H}(G_{18})$,
$\mathcal{H}(G_{20})$, $\mathcal{H}(G_{21})$ and $\mathcal{H}(G_{22})$.
\item $\mathcal{H}(G_{26})$ the twisted symmetric algebra of the cyclic group $C_2$ over
$\mathcal{H}(G_{25})$.
\end{itemize}
In all these cases, Proposition $\ref{1.42}$ implies that the Schur elements of the twisted symmetric algebra are scalar multiples of the Schur elements of the subalgebra. Due to the nature of the specializations (each indeterminate is sent to an indeterminate or a root of unity or a product of the two), the Schur elements of the subalgebra are also of the desired form.

Finally, if $W$ is one of the remaining exceptional irreducible complex reflection groups, then $W$ has
one hyperplane orbit $\mathcal{C}$ with $e_{\mathcal{C}}=2$. The generic Hecke algebra of $W$ is defined over a Laurent polynomial ring in two indeterminates $v_{\mathcal{C},0}$ and $v_{\mathcal{C},1}$. Its Schur elements should be products of $K$-cyclotomic polynomials in one variable
$v:=v_{\mathcal{C},0}v_{\mathcal{C},1}^{-1}$.
The generic Schur elements have been calculated
\begin{itemize}
  \item for $E_6$ and $E_7$ by Surowski (\cite{Sur78}),
  \item for $E_8$ by Benson (\cite{Ben79}),
  \item for $H_3$ by Lusztig  (\cite{Lu82}),
  \item for $H_4$ by Alvis and Lusztig (\cite{AlLu82}),
  \item for $G_{24}$, $G_{27}$, $G_{29}$, $G_{31}$, $G_{33}$ and $G_{34}$ by Malle (\cite{Ma5}),
\end{itemize}
and they are, indeed, products of $K$-cyclotomic polynomials in ``one'' variable.

Note that in order to write  the Schur elements in the desired form, we have used the GAP Package CHEVIE (where some mistakes in the 
articles cited above have been corrected).}
\end{apod}$ $\\
\begin{remark}\emph{ It is a consequence of \cite{Raph}, Theorem 3.5, that the
irreducible factors of the generic Schur elements over
$\mathbb{C}[\textbf{v},\textbf{v}^{-1}]$ are divisors of Laurent
polynomials of the form $M(\textbf{v})^n-1$, where
\begin{itemize}
\item $M(\textbf{v})$ is a monomial in $\mathbb{C}[\textbf{v},\textbf{v}^{-1}]$,
\item $n$ is a positive integer.
\end{itemize}}
\end{remark}\

We have seen that
the specialization $v_{\mathcal{C},j}\mapsto 1$ induces a
bijection  $\chi_{\textbf{v}} \mapsto \chi$ from
$\mathrm{Irr}(K(\textbf{v})\mathcal{H})$ to $\mathrm{Irr}(W)$.
Due to the assumptions $\ref{ypo}$, it sends $s_\chi(\textbf{v})$ to $|W|/\chi(1)$, which is the Schur element of $\chi$ with respect to the canonical symmetrizing form. Therefore, 
the first cyclotomic polynomial does not appear in the
factorization of  $s_\chi(\textbf{v})$ (otherwise the
specialization $v_{\mathcal{C},j} \mapsto 1$ would send
$s_\chi(\textbf{v})$ to $0$). 

The following result is an immediate application of Definition
$\ref{essential algebra}$.

\begin{theorem}\label{Hecke is essential} 
 The algebra $\mathcal{H}$, defined over the ring 
$\mathbb{Z}_K[\textbf{\emph{v}},\textbf{\emph{v}}^{-1}]$, is an essential algebra.
\end{theorem}

Thanks to Theorem $\ref{Hecke is essential}$, all the results of Chapter 3 can be applied to the generic Hecke algebra of an irreducible complex reflection group.

\begin{definition}\label{essential for group}
Let $\mathfrak{p}$ be a prime ideal of $\mathbb{Z}_K$. We say that a (primitive) monomial $M$ in $\mathbb{Z}_K[\textbf{\emph{v}},\textbf{\emph{v}}^{-1}]$ is 
$\mathfrak{p}$-essential for $W$, \index{p-essential monomial} if  $M$ is $\mathfrak{p}$-essential for $\mathcal{H}$.
\end{definition}

\begin{px}\label{Schur elements for G2}
\emph{\small Let $W:=G_2$.  The group $G_2$ is a Weyl group. We have seen that
$$\begin{array}{rll}
   \mathcal{H}(G_2)=<S,T \,\,|&STSTST=TSTSTS, &(S-u_0)(S-u_1)=0, \\
                                &         &(T-w_0)(T-w_1)=0>.
  \end{array}$$
Set $x_0^2:=u_0$, $x_1^2:=-u_1$, $y_0^2:=w_0$, $y_1^2:=-w_1$.
 By Theorem $\ref{Semisimplicity Malle}$, the algebra
 $\mathbb{Q}(x_0,x_1,y_0,y_1)\mathcal{H}(G_2)$ is split semisimple and hence, there exists a bijection between its irreducible characters and the irreducible characters of $G_2$. The group
 $G_2$ has $4$ irreducible characters of degree $1$ and $2$ irreducible characters of degree $2$.
Set
\begin{description}
\item[$s_1(x_0,x_1,y_0,y_1):=$]${\Phi_4(x_0x_1^{-1})} \cdot  
               { \Phi_4(y_0y_1^{-1})} \cdot  
               {\Phi_3(x_0x_1^{-1}y_0y_1^{-1})} \cdot 
                \Phi_6(x_0x_1^{-1}y_0y_1^{-1}),$
\item[$s_2(x_0,x_1,y_0,y_1):=$]$2  x_0^{-2}x_1^2\cdot {\Phi_3(x_0x_1^{-1}y_0y_1^{-1})} \cdot  \Phi_6(x_0x_1^{-1}y_0^{-1}y_1),$
\end{description}
where
$  \Phi_3(x)=x^2+x+1,\, \Phi_4(x)=x^2+1,\, \Phi_6(x)=x^2-x+1.$\\ \\
The Schur elements of $\mathcal{H}(G_2)$ are
$$s_1(x_0,x_1,y_0,y_1), s_1(x_0,x_1,y_1,y_0), s_1(x_1,x_0,y_0,y_1), s_1(x_1,x_0,y_1,y_0),$$                         
$$s_2(x_0,x_1,y_0,y_1), s_2(x_0,x_1,y_1,y_0).$$
Since $\Phi_3(1)=3$, $\Phi_4(1)=2$ and $\Phi_6(1)=1$, we obtain that
\begin{itemize}
\item the $(2)$-essential monomials for $G_2$ are $x_0x_1^{-1}$ and $y_0y_1^{-1}$ (and their inverses),
\item the $(3)$-essential monomials for $G_2$ are $x_0x_1^{-1}y_0y_1^{-1}$ and $x_0x_1^{-1}y_0^{-1}y_1$ (and their inverses).
 \end{itemize}      }
\end{px}

\begin{px}\label{Schur elements for G4}
\emph{\small Let $W:=G_4$.  The field of definition of $G_4$ is $\mathbb{Q}(\zeta_3)$. We have seen that
$$\begin{array}{rll}
   \mathcal{H}(G_4)=<S,T \,\,|&STS=TST, &(S-u_0)(S-u_1)(S-u_2)=0, \\
                                &         &(T-u_0)(T-u_1)(T-u_2)=0>.
  \end{array}$$
Set $v_0^6:=u_0$, $v_1^6:=\zeta_3^2u_1$, $v_2^6:=\zeta_3u_2$.
 By Theorem $\ref{Semisimplicity Malle}$, the algebra
 $\mathbb{Q}(\zeta_3)(v_0,v_1,v_2)\mathcal{H}(G_4)$ is split semisimple and hence, there exists a bijection between its irreducible characters and the irreducible characters of $G_4$. The group
 $G_4$ has $3$ irreducible characters of degree $1$,  $3$ irreducible characters of degree $2$ and $1$ irreducible character of degree $3$.
Set}
\begin{description}{\small
\item[$s_1(v_0,v_1,v_2)=$]
                $\Phi_9''(v_0v_1^{-1}) \cdot \Phi_{18}'(v_0v_1^{-1}) \cdot
                \Phi_4(v_0v_1^{-1})\cdot  \Phi_{12}'(v_0v_1^{-1}) \cdot
                \Phi_{12}''(v_0v_1^{-1}) \cdot \Phi_{36}'(v_0v_1^{-1}) \cdot
                {\Phi_9'(v_0v_2^{-1})} \cdot \Phi_{18}''(v_0v_2^{-1}) \cdot
               {\Phi_4(v_0v_2^{-1})} \cdot  \Phi_{12}'(v_0v_2^{-1}) \cdot
                \Phi_{12}''(v_0v_2^{-1}) \cdot \Phi_{36}''(v_0v_2^{-1})
                \cdot
                {\Phi_4(v_0^2v_1^{-1}v_2^{-1})} \cdot
                \Phi_{12}'(v_0^2v_1^{-1}v_2^{-1}) \cdot
                \Phi_{12}''(v_0^2v_1^{-1}v_2^{-1})$,
\item[$s_2(v_0,v_1,v_2)=$]$-\zeta_3^2v_2^6v_1^{-6}\cdot
                 {\Phi_9'(v_1v_0^{-1})} \cdot \Phi_{18}''(v_1v_0^{-1}) \cdot
               {\Phi_9''(v_2v_0^{-1})} \cdot \Phi_{18}'(v_2v_0^{-1}) \cdot
                { \Phi_4(v_1v_2^{-1})} \cdot  \Phi_{12}'(v_1v_2^{-1}) \cdot
                 \Phi_{12}''(v_1v_2^{-1}) \cdot \Phi_{36}'(v_1v_2^{-1}) \cdot
                {\Phi_4(v_0^{-2}v_1v_2)} \cdot
                 \Phi_{12}'(v_0^{-2}v_1v_2) \cdot
                 \Phi_{12}''(v_0^{-2}v_1v_2)$,
\item[$s_3(v_0,v_1,v_2)=$]${\Phi_4(v_0^2v_1^{-1}v_2^{-1})} \cdot
                \Phi_{12}'(v_0^2v_1^{-1}v_2^{-1}) \cdot
                \Phi_{12}''(v_0^2v_1^{-1}v_2^{-1}) \cdot
                {\Phi_4(v_1^2v_2^{-1}v_0^{-1})} \cdot
                \Phi_{12}'(v_1^2v_2^{-1}v_0^{-1}) \cdot
                \Phi_{12}''(v_1^2v_2^{-1}v_0^{-1}) \cdot
                {\Phi_4(v_2^2v_0^{-1}v_1^{-1})} \cdot
                \Phi_{12}'(v_2^2v_0^{-1}v_1^{-1}) \cdot
                \Phi_{12}''(v_2^2v_0^{-1}v_1^{-1})$},
\end{description}
\emph{\small where $\Phi_4(x)=x^2+1$, $\Phi_9'(x)=x^3-\zeta_3$,
$\Phi_9''(x)=x^3-\zeta_3^2$, $\Phi_{12}''(x)=x^2+\zeta_3$,
$\Phi_{12}'(x)=x^2+\zeta_3^2$, $\Phi_{18}''(x)=x^3+\zeta_3$,
$\Phi_{18}'(x)=x^3+\zeta_3^2$, $\Phi_{36}''(x)=x^6+\zeta_3$,
$\Phi_{36}'(x)=x^6+\zeta_3^2$.\\ \\The Schur elements of $\mathcal{H}(G_4)$ are
$$s_1(v_0,v_1,v_2), s_1(v_1,v_2,v_0), s_1(v_2,v_0,v_1),$$                         
$$s_2(v_0,v_1,v_2), s_2(v_1,v_2,v_0), s_2(v_2,v_0,v_1),
s_3(v_0,v_1,v_2).$$  
We deduce that the $(2)$-essential monomials for $G_4$ are
$$v_0v_1^{-1}, 
v_0v_2^{-1}, v_1v_2^{-1},
v_0^2v_1^{-1}v_2^{-1},v_1^2v_2^{-1}v_0^{-1},
v_2^2v_0^{-1}v_1^{-1}.$$
The first three are also the $(1-\zeta_3)$-essential monomials for $G_4$.}
\end{px}

\section{Cyclotomic Hecke algebras}

Let $y$ be an indeterminate. We set $x:=y^{|\mu(K)|}.$

\begin{definition}\label{specialization}
A cyclotomic specialization \index{cyclotomic specialization} of $\mathcal{H}$ is a
$\mathbb{Z}_K$-algebra morphism $\phi:
\mathbb{Z}_K[\textbf{\emph{v}},\textbf{\emph{v}}^{-1}]\rightarrow
\mathbb{Z}_K[y,y^{-1}]$ with the following properties:
\begin{itemize}
  \item $\phi: v_{\mathcal{C},j} \mapsto y^{n_{\mathcal{C},j}}$ where
  $n_{\mathcal{C},j} \in \mathbb{Z}$ for all $\mathcal{C}$ and $j$.
  \item For all $\mathcal{C} \in \mathcal{A}/W$, if $z$ is another
  indeterminate, the element of $\mathbb{Z}_K[y,y^{-1},z]$ defined by
  $$\Gamma_\mathcal{C}(y,z):=\prod_{j=0}^{e_\mathcal{C}-1}(z-\zeta_{e_\mathcal{C}}^jy^{n_{\mathcal{C},j}})$$
  is invariant by the action of $\textrm{\emph{Gal}}(K(y)/K(x))$.
\end{itemize}
\end{definition}

If $\phi$ is a cyclotomic specialization of $\mathcal{H}$,
the corresponding \emph{cyclotomic Hecke algebra} \index{cyclotomic Hecke algebra} is the
$\mathbb{Z}_K[y,y^{-1}]$-algebra, denoted by $\mathcal{H}_\phi$,
which is obtained as the specialization of the
$\mathbb{Z}_K[\textbf{v},\textbf{v}^{-1}]$-algebra $\mathcal{H}$ via
the morphism $\phi$. It also has a symmetrizing form $t_\phi$
defined as the specialization of the canonical form $t$.\\
\\
\begin{remark} \emph{Sometimes we describe the morphism $\phi$ by the
formula}
$$u_{\mathcal{C},j} \mapsto \zeta_{e_\mathcal{C}}^j x^{n_{\mathcal{C},j}}.$$
\emph{If now we set $q:=\zeta x$ for some root of unity $\zeta \in
\mu(K)$, then the cyclotomic specialization $\phi$ becomes a
$\zeta$-\emph{cyclotomic specialization} and $\mathcal{H}_\phi$ can
be also considered over $\mathbb{Z}_K[q,q^{-1}]$.}
\end{remark}
\begin{px}\label{spetsial}
\small{\emph{The ``spetsial'' cyclotomic Hecke algebra $\mathcal{H}_q^s(W)$ is the
1-cyclotomic algebra obtained by the specialization\index{spetsial Hecke algebra}
$$u_{\mathcal{C},0} \mapsto q,\,\, u_{\mathcal{C},j} \mapsto \zeta_{e_\mathcal{C}}^j \textrm{
for } 1 \leq j \leq e_\mathcal{C}-1, \textrm{ for all } \mathcal{C}
\in \mathcal{A}/W.$$ For example, 
$$\mathcal{H}_q^s(G_2)=<S,T \textrm{ }|\, STSTST=TSTSTS,
(S-q)(S+1)=(T-q)(T+1)=0>.$$
and
$$\mathcal{H}_q^s(G_4)=<S,T \textrm{ }|\, STS=TST,
(S-q)(S^2+S+1)=(T-q)(T^2+T+1)=0>.$$}}
\end{px}\

Set $A:=\mathbb{Z}_K[\textbf{v},\textbf{v}^{-1}]$ and
$\Omega:=\mathbb{Z}_K[y,y^{-1}]$. Let $\phi: A \rightarrow \Omega$
be a cyclotomic specialization such that $\phi(v_{\mathcal{C},j})=
y^{n_{\mathcal{C},j}}$ for all $\mathcal{C},j$. Recall that, for $\alpha \in
\mathbb{Z}\setminus\{0\}$, we denote by $I^\alpha:\Omega \rightarrow \Omega$
the monomorphism $y \mapsto y^\alpha$.

\begin{theorem}\label{cyclotomic}
Let $\phi:A \rightarrow \Omega$ be a cyclotomic specialization like
above. Then there exist an adapted $\mathbb{Z}_K$-algebra morphism
$\varphi:A \rightarrow \Omega$ and $\alpha \in \mathbb{Z}\setminus\{0\}$
such that
$$\phi={I^\alpha}\circ\varphi.$$
\end{theorem}
\begin{apod}{We set $d:=\mathrm{gcd}(n_{\mathcal{C},j})$ and consider the
cyclotomic specialization $\varphi: v_{\mathcal{C},j} \mapsto
y^{n_{\mathcal{C},j}/d}$. We have $\phi=I^d \circ\varphi$. Since
$\mathrm{gcd}(n_{\mathcal{C},j}/d)=1$,
there exist $a_{\mathcal{C},j} \in \mathbb{Z}$ such that
$$\sum_{\mathcal{C},j}a_{\mathcal{C},j}(n_{\mathcal{C},j}/d)=1.$$
We have
$y=\varphi(\prod_{\mathcal{C},j}v_{\mathcal{C},j}^{a_{\mathcal{C},j}})$
an hence, $\varphi$ is surjective.  Then, by Proposition
$\ref{surjective is adapted}$, $\varphi$ is adapted.}
\end{apod}

Let $\varphi$ be defined as in Theorem $\ref{cyclotomic}$ and
$\mathcal{H}_\varphi$ the corresponding cyclotomic Hecke algebra.
Proposition $\ref{associated morphism preserves splitness}$ implies
that the algebra $K(y)\mathcal{H}_\varphi$ is split semisimple. Due
to Corollary $\ref{injective preserves splitness}$ and the theorem
above, we deduce that 

\begin{proposition}\label{cyclotomic Hecke is split ss}
The algebra $K(y)\mathcal{H}_\phi$ is 
split semisimple. 
\end{proposition}

For $y=1$, the algebra $K(y)\mathcal{H}_\phi$ specializes to the group
algebra $KW$ (the form $t_\phi$ becoming the canonical form on the
group algebra). Thus, by ``Tits' deformation theorem'', the
specialization $v_{\mathcal{C},j} \mapsto 1$ defines the following bijections
$$\begin{array}{ccccc}
    \textrm{Irr}(K(\textbf{v})\mathcal{H}) & \leftrightarrow & \textrm{Irr}(K(y)\mathcal{H}_\phi) & \leftrightarrow & \textrm{Irr}(W) \\
    \chi_{\textbf{v}} & \mapsto & \chi_{\phi} & \mapsto & \chi.
  \end{array}$$

The following result is an immediate consequence of Theorem
$\ref{Schur element generic}$.

\begin{proposition}\label{Schur element cyclotomic}
The Schur element $s_{\chi_\phi}(y)$ associated to the irreducible
character $\chi_\phi$ of $K(y)\mathcal{H}_\phi$ is a Laurent
polynomial in $y$ of the form
$$s_{\chi_\phi}(y)=\psi_{\chi,\phi} y^{a_{\chi,\phi}} \prod_{\Phi \in
C_K}\Phi(y)^{n_{\chi,\phi,\Phi}}$$ where $\psi_{\chi,\phi} \in
\mathbb{Z}_K$, $a_{\chi,\phi} \in \mathbb{Z}$, $n_{\chi,\phi,\Phi} \in
\mathbb{N}$ and $C_K$ is a set of $K$-cyclotomic polynomials.
\end{proposition}

\subsection{Essential hyperplanes}

Let $\mathfrak{p}$ be a prime ideal of $\mathbb{Z}_K$. Let
$\phi: v_{\mathcal{C},j} \mapsto y^{n_{\mathcal{C},j}}$ be a cyclotomic specialization of $\mathcal{H}$
and let $\varphi$ be an adapted morphism as in
Theorem
$\ref{cyclotomic}$. By Corollary $\ref{i gives same blocks}$, the blocks of
$\Omega_{\mathfrak{p}\Omega}\mathcal{H}_\phi$ coincide
with the blocks of $\Omega_{\mathfrak{p}\Omega}\mathcal{H}_\varphi$
and the latter can be calculated with the use of 
Theorem $\ref{main theorem}$. Therefore, we need to know
which $\mathfrak{p}$-essential monomials are sent to 1 by $\varphi$.

Let $M:=\prod_{\mathcal{C},j}v_{\mathcal{C},j}^{a_{\mathcal{C},j}}$
be a $\mathfrak{p}$-essential monomial for $W$. Then
$$\varphi(M)=1 \Leftrightarrow \phi(M)=1 \Leftrightarrow \sum_{\mathcal{C},j}a_{\mathcal{C},j}n_{\mathcal{C},j}=0.$$
Set $m:=\sum_{\mathcal{C}\in \mathcal{A}/W}e_\mathcal{C}$. The
hyperplane defined in $\mathbb{C}^m$ by the relation
$$\sum_{\mathcal{C},j}a_{\mathcal{C},j}t_{\mathcal{C},j}=0,$$ where
$(t_ {\mathcal{C},j})_{ \mathcal{C},j}$ is a set of $m$
indeterminates, is called \emph{$\mathfrak{p}$-essential hyperplane} \index{p-essential hyperplane}
for $W$. A hyperplane in $\mathbb{C}^m$ is called \emph{essential} \index{essential hyperplane}
for $W$, if it is $\mathfrak{p}$-essential for some prime ideal
$\mathfrak{p}$ of $\mathbb{Z}_K$.

\begin{px}\label{Essential hyperplanes for G2}
\emph{\small Let $W:=G_2$.  Following Example $\ref{Schur elements for G2}$, let
$$\phi: x_0 \mapsto y^{n_0}, x_1 \mapsto y^{n_1}, y_0 \mapsto y^{m_0}, y_1 \mapsto y^{m_1} $$ be a cyclotomic specialization. Then
\begin{itemize}
\item the $(2)$-essential hyperplanes for $G_2$ are $N_0-N_1=0$ and $M_0-M_1=0$,
\item the $(3)$-essential hyperplanes for $G_2$ are $N_0-N_1+M_0-M_1=0$ and $N_0-N_1-M_0+M_1=0$.
 \end{itemize}      }
\end{px}

\begin{px}\label{Essential hyperplanes for G4}
\emph{\small Let $W:=G_4$.  Following Example $\ref{Schur elements for G4}$, let
$\phi: v_i \mapsto y^{n_i}$ for $i=0,1,2$ be a cyclotomic specialization. Then the hyperplanes}
\emph{\small\begin{itemize}
\item $N_0-N_1=0$, $N_0-N_2=0$ and $N_1-N_2=0$ are $(2)$-essential and $(1-\zeta_3)$-essential for $G_4$, 
\item $2N_0-N_1-N_2=0$, $2N_1-N_2-N_0=0$ and $2N_2-N_0-N_1=0$ are just $(2)$-essential for $G_4$.
\end{itemize}}
\end{px}

In order to calculate the blocks of
$\Omega_{\mathfrak{p}\Omega}\mathcal{H}_\phi$, we check to which
$\mathfrak{p}$-essential hyperplanes the $n_{\mathcal{C},j}$ belong and we apply Theorem $\ref{main theorem}$:
\begin{itemize}
  \item If the $n_{\mathcal{C},j}$ belong to no $\mathfrak{p}$-essential hyperplane, then the
blocks of $\Omega_{\mathfrak{p}\Omega}\mathcal{H}_\phi$ coincide
with the blocks of $A_{\mathfrak{p}A}\mathcal{H}$. We call these blocks \index{p-blocks associated with no essential hyperplane}$\mathfrak{p}$-\emph{blocks associated with no essential hyperplane}. 
  \item If the
 $n_{\mathcal{C},j}$ belong to exactly one $\mathfrak{p}$-essential
 hyperplane $H_M$, corresponding to the $\mathfrak{p}$-essential monomial
 $M$, then the blocks of
 $\Omega_{\mathfrak{p}\Omega}\mathcal{H}_\phi$ coincide
 with the blocks of $A_{\mathfrak{q}_M}\mathcal{H}$, where $\mathfrak{q}_M:=\frak{p}A+(M-1)A$.
 We call these blocks \index{p-blocks associated with an essential hyperplane}$\mathfrak{p}$-\emph{blocks associated with the essential hyperplane} $H_M$. 
\item If the $n_{\mathcal{C},j}$ belong to more than one $\mathfrak{p}$-essential
hyperplane, then, following Theorem $\ref{main theorem}$, 
the blocks of $\Omega_{\mathfrak{p}\Omega}\mathcal{H}_\phi$
are unions of the $\mathfrak{p}$-blocks associated with the $\mathfrak{p}$-essential hyperplanes to which the $n_{\mathcal{C},j}$ belong and they are minimal with respect to that property.
\end{itemize}

This last property of the  $\mathfrak{p}$-blocks is called \index{semi-continuity} ``property of \emph{semi-continuity}'' (the name is due to C. Bonnaf\'e). The property of semi-continuity also appears in works on Kazhdan-Lusztig cells  (cf.~\cite{BGIL}, \cite{Bon}, \cite{Jer}) and on Cherednik algebras (cf.~\cite{GoMa}). In the next section, we will see that the Rouquier blocks of the cyclotomic Hecke algebras also have this property.

\subsection{Group algebra}

Let $\mathfrak{p}$ be a prime ideal of $\mathbb{Z}_K$ lying over a prime number $p$ and let
$\phi: v_{\mathcal{C},j} \mapsto y^{n_{\mathcal{C},j}}$ be a cyclotomic specialization of $\mathcal{H}$.
If $n_{\mathcal{C},j}=n \in \mathbb{Z}$ for all $\mathcal{C}$ and $j$,
then the $n_{\mathcal{C},j}$ belong to
all $\mathfrak{p}$-essential hyperplanes for $W$ and we have $\Omega_{\mathfrak{p}\Omega}\mathcal{H}_\phi \cong
\Omega_{\mathfrak{p}\Omega}W$.
Note that, since the ring
$\Omega_{\mathfrak{p}\Omega}$ is a discrete valuation ring (by
Theorem $\ref{Krull-dvr}$), the blocks of
$\Omega_{\mathfrak{p}\Omega}W$ are the $p$-blocks of $W$ as
determined by Brauer theory. 
 Due to Theorem $\ref{main
theorem}$, we obtain the following result which relates the $\mathfrak{p}$-blocks of any cyclotomic Hecke algebra to the $p$-blocks of $W$.

\begin{proposition}\label{group algebra}
 Let
$\phi: v_{\mathcal{C},j} \mapsto y^{n_{\mathcal{C},j}}$ be a cyclotomic specialization of $\mathcal{H}$. If two irreducible characters $\chi,\psi \in \mathrm{Irr}(W)$
are in the same block of
$\Omega_{\mathfrak{p}\Omega}\mathcal{H}_\phi$, then they are in the
same $p$-block of $W$.
\end{proposition}
\begin{apod}{The blocks of $\Omega_{\mathfrak{p}\Omega}{H}_\phi$ are unions of the blocks of
$A_{\mathfrak{q}_M}\mathcal{H}$ for all $\mathfrak{p}$-essential
monomials $M$ such that $\phi(M)=1$, whereas the $p$-blocks of $W$
are unions of the blocks of $A_{\mathfrak{q}_M}\mathcal{H}$ for all
$\mathfrak{p}$-essential monomials $M$.}
\end{apod}

However, we know from Brauer theory that if the order of the group
$W$ is prime to $p$, then every character of $W$ is a $p$-block by
itself (see, for example, \cite{Se}, 15.5, Proposition 43). It is an
immediate consequence of Proposition $\ref{group algebra}$ that

\begin{corollary}\label{p prime to the order of the group}
If $\mathfrak{p}$ is a prime ideal of $\mathbb{Z}_K$ lying over a
prime number $p$ which does not divide the order of the group $W$,
then the blocks of $\Omega_{\mathfrak{p}\Omega}\mathcal{H}_\phi$ are
singletons.
\end{corollary}

\section{Rouquier blocks of the cyclotomic Hecke algebras}

\begin{definition}\label{Rouquier ring}
We call Rouquier ring \index{Rouquier ring} of $K$ and denote by $\mathcal{R}_K(y)$ the
$\mathbb{Z}_K$-subalgebra of $K(y)$
$$\mathcal{R}_K(y):=\mathbb{Z}_K[y,y^{-1},(y^n-1)^{-1}_{n\geq 1}]$$
\end{definition}

Let $\phi: v_{\mathcal{C},j} \mapsto y^{n_{\mathcal{C},j}}$ be a
cyclotomic specialization and $\mathcal{H}_\phi$ the corresponding
cyclotomic Hecke algebra. The \emph{Rouquier blocks}  \index{Rouquier block} of
$\mathcal{H}_\phi$ are the blocks of the algebra
$\mathcal{R}_K(y)\mathcal{H}_\phi$.

 It has been shown by Rouquier
(cf.~\cite{Rou}), that if $W$ is a Weyl group and $\mathcal{H}_\phi$ is
obtained via the ``spetsial'' cyclotomic specialization (see Example
$\ref{spetsial}$), then the Rouquier blocks of  $\mathcal{H}_\phi$ coincide with the
families of characters defined by Lusztig. Thus, the Rouquier
blocks generalize the notion of ``families of characters'' to all complex reflection groups. 
\\ \\
\begin{remark}\emph{ We have seen that if we set $q:=\zeta y^{|\mu(K)|}$ for some root of unity $\zeta \in \mu(K)$, then the cyclotomic Hecke algebra $\mathcal{H}_\phi$ can be also considered over the  ring $\mathbb{Z}_{K}[q,q^{-1}]$. We define the Rouquier blocks of  $\mathcal{H}_\phi$ to be the blocks of $\mathcal{R}_{K}(y)\mathcal{H}_\phi$. However, in other texts, as, for example, in \cite{BK}, the Rouquier blocks are defined to be the blocks of  $\mathcal{R}_{K}(q)\mathcal{H}_\phi$. Since $\mathcal{R}_{K}(y)$ is the integral closure of $\mathcal{R}_{K}(q)$ in the splitting field  $K(y)$ for $\mathcal{H}_\phi$, Proposition $\ref{Galois action on integral closure}$ establishes the (determining) relation between the blocks of
$\mathcal{R}_{K}(y)\mathcal{H}_\phi$ and the blocks of $\mathcal{R}_{K}(q)\mathcal{H}_\phi$.}
\end{remark}$ $\\

The Rouquier ring $\mathcal{R}_K(y)$ has many interesting properties. The next result describes some of them.

\begin{proposition}\label{Some properties of the Rouquier
ring}\
\begin{enumerate}
  \item The group of units $\mathcal{R}_K(y)^\times$ of the Rouquier ring $\mathcal{R}_K(y)$
  consists of the elements of the form
  $$u y^n \prod_{\Phi \in \mathrm{Cycl}(K)} \Phi(y)^{n_\Phi},$$
  where $u \in \mathbb{Z}_K^\times$, $n, n_\Phi \in \mathbb{Z}$, $\mathrm{Cycl}(K)$ is the set
  of $K$-cyclotomic polynomials and
  $n_\Phi=0$ for all but a finite number of $\Phi$.
  \item The prime ideals of $\mathcal{R}_K(y)$ are
  \begin{itemize}
    \item the zero ideal $\{0\}$,
    \item the ideals of the form $\mathfrak{p}\mathcal{R}_K(y)$,
    where $\mathfrak{p}$ is a prime ideal of $\mathbb{Z}_K$,
    \item the ideals of the form $P(y)\mathcal{R}_K(y)$, where
    $P(y)$ is an irreducible element of $\mathbb{Z}_K[y]$ of degree
    at least $1$, prime to $y$ and to $\Phi(y)$ for all $\Phi \in
    \mathrm{Cycl}(K)$.
  \end{itemize}
  \item The Rouquier ring $\mathcal{R}_K(y)$ is a Dedekind ring.
\end{enumerate}
\end{proposition}
\begin{apod}{
\begin{enumerate}
   \item This part is immediate from the definition of
   $K$-cyclotomic polynomials.
   \item Since $\mathcal{R}_K(y)$ is an integral domain, the zero
   ideal is prime.

   The ring $\mathbb{Z}_K$ is a Dedekind ring
   and thus a Krull ring, by Proposition $\ref{case of
   Krull}$. Proposition $\ref{prime ideals of height 1}$ implies that the ring
   $\mathbb{Z}_K[y]$ is also a Krull ring whose prime ideals of height 1
    are of the form $\mathfrak{p}\mathbb{Z}_K[y]$
   ($\mathfrak{p}$ prime in $\mathbb{Z}_K$) and
   $P(y)\mathbb{Z}_K[y]$ ($P(y)$ irreducible in $\mathbb{Z}_K[y]$ of degree at least 1).
   Moreover, $\mathbb{Z}_K$ has an infinite
   number of non-zero prime ideals whose intersection is the zero ideal. Since
   all non-zero prime ideals of $\mathbb{Z}_K$ are maximal, we
   obtain that every prime ideal of $\mathbb{Z}_K$ is the intersection of
   maximal ideals. Thus $\mathbb{Z}_K$ is, by definition, a Jacobson
   ring (cf.~\cite{Eis}, \S 4.5). The general form of the
   Nullstellensatz (\cite{Eis}, Theorem 4.19) implies that
   for every maximal ideal $\mathfrak{m}$ of $\mathbb{Z}_K[y]$, the
   ideal $\mathfrak{m} \cap \mathbb{Z}_K$ is a maximal ideal of
   $\mathbb{Z}_K$. We deduce that the maximal ideals of
   $\mathbb{Z}_K[y]$ are of the form $\mathfrak{p}\mathbb{Z}_K[y]+P(y)\mathbb{Z}_K[y]$
   ($\mathfrak{p}$ prime in $\mathbb{Z}_K$ and $P(y)$ of degree at least 1 irreducible
   modulo $\mathfrak{p}$). Since
   $\mathbb{Z}_K[y]$ has Krull dimension 2, we have now described
   all its prime ideals.

   The Rouquier ring $\mathcal{R}_K(y)$ is a localization of $\mathbb{Z}_K[y]$. Therefore,
   in order to prove that the non-zero prime ideals of
   $\mathcal{R}_K(y)$ are the ones described above, it is enough to
   show that $\mathfrak{m}\mathcal{R}_K(y)=\mathcal{R}_K(y)$ for all
   maximal ideals $\mathfrak{m}$ of $\mathbb{Z}_K[y]$. For this, it
   suffices to show that $\mathfrak{p}\mathcal{R}_K(y)$ is a maximal
   ideal of $\mathcal{R}_K(y)$ for all prime ideals $\mathfrak{p}$
   of $\mathbb{Z}_K$.

   Let $\mathfrak{p}$ be a prime ideal of $\mathbb{Z}_K$. Then
   $$\mathcal{R}_K(y)/\mathfrak{p}\mathcal{R}_K(y) \cong \mathbb{F}_\mathfrak{p}[y,y^{-1},(y^n-1)^{-1}_{n\geq 1}],$$
   where $\mathbb{F}_\mathfrak{p}$ denotes the finite field
   $\mathbb{Z}_K/\mathfrak{p}$. Since $\mathbb{F}_\mathfrak{p}$ is finite,
   every non-zero polynomial in $\mathbb{F}_\mathfrak{p}[y]$ is a product of
   elements which divide $y$ or $y^n-1$ for some $n \in \mathbb{N}$.
   Thus every non-zero element of $\mathbb{F}_\mathfrak{p}[y]$ is invertible in
   $\mathcal{R}_K(y)/\mathfrak{p}\mathcal{R}_K(y)$. Consequently, we
   obtain that
   $$\mathcal{R}_K(y)/\mathfrak{p}\mathcal{R}_K(y) \cong \mathbb{F}_\mathfrak{p}(y),$$
   whence $\mathfrak{p}$ generates a maximal ideal in
   $\mathcal{R}_K(y)$.
   \item  The ring $\mathcal{R}_K(y)$ is the localization of a
   Noetherian integrally closed ring and thus Noetherian and
   integrally closed itself. Moreover, following the description of
   its prime ideals in part 2, it has Krull dimension 1.}
\end{enumerate}
\end{apod}
\begin{remark}
 \emph{If $P(y)$ is an irreducible element of $\mathbb{Z}_K[y]$ of degree
 at least 1, prime to $y$ and to $\Phi(y)$ for all $\Phi \in
 \mathrm{Cycl}(K)$, then the field $\mathcal{R}_K(y)/P(y)\mathcal{R}_K(y)$
 is isomorphic to the field of fractions of the ring
 $\mathbb{Z}_K[y]/P(y)\mathbb{Z}_K[y]$}.
\end{remark}
\\

Now let us recall the form of the Schur elements of the cyclotomic
Hecke algebra $\mathcal{H}_\phi$ given in Proposition $\ref{Schur
element cyclotomic}$. If $\chi_\phi$ is an irreducible character of
$K(y)\mathcal{H}_\phi$, then its Schur element $s_{\chi_\phi}(y)$ is
of the form
$$s_{\chi_\phi}(y)=\psi_{\chi,\phi} y^{a_{\chi,\phi}} \prod_{\Phi \in
C_K}\Phi(y)^{n_{\chi,\phi,\Phi}}$$ where $\psi_{\chi,\phi} \in
\mathbb{Z}_K$, $a_{\chi,\phi} \in \mathbb{Z}$, $n_{\chi,\phi,\Phi} \in
\mathbb{N}$ and $C_K$ is a set of $K$-cyclotomic polynomials.

\begin{definition}\label{bad}
A prime ideal $\mathfrak{p}$ of $\mathbb{Z}_K$ lying over a prime
number $p$ is $\phi$-bad \index{$\phi$-bad prime ideal} for $W$, if there exists $\chi_\phi \in
\textrm{\emph{Irr}}(K(y)\mathcal{H}_\phi)$ with $\psi_{\chi,\phi}
\in \mathfrak{p}$. If $\mathfrak{p}$ is $\phi$-bad for $W$, we say
that $p$ is a $\phi$-bad prime number \index{$\phi$-bad prime number} for $W$.
\end{definition}
\begin{remark}\emph{ If $W$ is a Weyl group and $\phi$ is the
``spetsial'' cyclotomic specialization, then the $\phi$-bad prime
ideals are the ideals generated by the bad prime numbers (in the
``usual'' sense) for $W$ (see \cite{GeRo}, 5.2).}
\end{remark}$ $\\

Note that if $p$ is a $\phi$-bad prime number for $W$, then $p$ must
divide the order of the group (since
$s_{\chi_\phi}(1)=|W|/\chi(1)$).\\

Let us denote by $\mathcal{O}$ the Rouquier ring. By Proposition
$\ref{p-blocks}$, the Rouquier blocks of $\mathcal{H}_\phi$ are
unions of the blocks of $\mathcal{O}_\mathcal{P}\mathcal{H}_\phi$, where $\mathcal{P}$ runs over the set of prime ideals of $\mathcal{O}$. However, in all
of the following cases, due to the form of the Schur elements, the
blocks of $\mathcal{O}_\mathcal{P}\mathcal{H}_\phi$ are singletons
(\ie $e_{\chi_\phi}=\chi_\phi^\vee /s_{\chi_\phi} \in
\mathcal{O}_\mathcal{P}\mathcal{H}_\phi$ for all $\chi_\phi \in
\mathrm{Irr}(K(y)\mathcal{H}_\phi)$):
\begin{itemize}
  \item $\mathcal{P}$ is the zero ideal $\{0\}$.
  \item $\mathcal{P}$ is of the form $P(y)\mathcal{O}$, where
$P(y)$ is an irreducible element of $\mathbb{Z}_K[y]$ of degree at
least $1$, prime to $y$ and to $\Phi(y)$ for all $\Phi \in
\mathrm{Cycl}(K)$.
  \item $\mathcal{P}$ is of the form $\mathfrak{p}\mathcal{O}$, where
$\mathfrak{p}$ is a prime ideal of $\mathbb{Z}_K$ which is not
$\phi$-bad for $W$.
\end{itemize}
Therefore, the blocks of $\mathcal{O}\mathcal{H}_\phi$ are, simply,
unions of the blocks of
$\mathcal{O}_{\mathfrak{p}\mathcal{O}}\mathcal{H}_\phi$, where $\mathfrak{p}$ runs over the set of
$\phi$-bad prime ideals $\mathfrak{p}$ of $\mathbb{Z}_K$. More precisely, we have that

\begin{proposition}\label{semicontinuity of Rouquier blocks}
Let $\chi,\psi \in \emph{Irr}(W)$. The characters $\chi_\phi$ and
$\psi_\phi$ are in the same Rouquier block of $\mathcal{H}_\phi$ if
and only if there exist a finite sequence
$\chi_0,\chi_1,\ldots,\chi_n \in \emph{Irr}(W)$ and a finite
sequence $\mathfrak{p}_1,\ldots,\mathfrak{p}_n$ of $\phi$-bad prime
ideals for $W$ such that
\begin{itemize}
  \item $(\chi_0)_\phi=\chi_\phi$ and $(\chi_n)_\phi=\psi_\phi$,
  \item for all $j$ $(1\leq j \leq n)$,\,\,
        $(\chi_{j-1})_\phi$ and $(\chi_j)_\phi$ are in the same block of
       $\mathcal{O}_{\mathfrak{p}_j\mathcal{O}}\mathcal{H}_\phi.$
\end{itemize}
\end{proposition}

By
Proposition $\ref{prime ideals of a localization}$(4), we obtain
that $\mathcal{O}_{\mathfrak{p}\mathcal{O}} \cong
\Omega_{\mathfrak{p}\Omega}$, where
$\Omega:=\mathbb{Z}_K[y,y^{-1}]$. In the previous section we saw how
we can use Theorem $\ref{main theorem}$ to calculate the blocks of
$\Omega_{\mathfrak{p}\Omega}\mathcal{H}_\phi$ and thus obtain
the Rouquier blocks of $\mathcal{H}_\phi$. We deduce that the Rouquier blocks of the
cyclotomic Hecke algebras also have the property of \emph{semi-continuity}:
\begin{itemize}
  \item If the $n_{\mathcal{C},j}$ belong to no essential hyperplane for $W$, then the
  Rouquier blocks of $\mathcal{H}_\phi$ are the \index{Rouquier blocks associated with no essential hyperplane}\emph{Rouquier blocks associated with no essential hyperplane}. 
  \item If the
 $n_{\mathcal{C},j}$ belong to exactly one essential
 hyperplane $H$ for $W$, then the
  Rouquier blocks of $\mathcal{H}_\phi$ are the \index{Rouquier blocks associated with an essential hyperplane}\emph{Rouquier blocks associated with the essential hyperplane} $H$. 
\item If the $n_{\mathcal{C},j}$ belong to more than one essential
hyperplane, then the Rouquier blocks of $\mathcal{H}_\phi$
are unions of the  Rouquier blocks associated with the essential hyperplanes to which the $n_{\mathcal{C},j}$ belong and they are minimal with respect to that property.
\end{itemize}

\subsection{Rouquier blocks and central morphisms}

The following description of the Rouquier blocks results from
Proposition $\ref{omega_chi}$ and the description of $\phi$-bad
prime ideals for $W$.

\begin{proposition}\label{Rouquier blocks and central characters}
Let $\chi,\psi \in \emph{Irr}(W)$. The characters $\chi_\phi$ and
$\psi_\phi$ are in the same Rouquier block of $\mathcal{H}_\phi$ if
and only if there exists a finite sequence
$\chi_0,\chi_1,\ldots,\chi_n \in \emph{Irr}(W)$ and a finite
sequence $\mathfrak{p}_1,\ldots,\mathfrak{p}_n$ of $\phi$-bad prime
ideals for $W$ such that
\begin{itemize}
  \item $(\chi_0)_\phi=\chi_\phi$ and $(\chi_n)_\phi=\psi_\phi$,
  \item for all $j$ $(1\leq j \leq n)$,\,\,
        $\omega_{(\chi_{j-1})_\phi} \equiv \omega_{(\chi_j)_\phi}
        \emph{ mod } \mathfrak{p}_j\mathcal{O}_{\mathfrak{p}_j\mathcal{O}}.$
\end{itemize}
\end{proposition}

\subsection{Rouquier blocks and functions $a$ and $A$}
\index{functions $a$ and $A$}
Following the notations of \cite{BMM2}, 6B, for every element $P(y)
\in \mathbb{C}(y)$, we call
\begin{itemize}
  \item \emph{valuation of $P(y)$ at $y$} and denote by $\mathrm{val}_y(P)$ the order of $P(y)$
  at 0 (we have $\mathrm{val}_y(P)<0$ if 0 is a pole of $P(y)$ and $\mathrm{val}_y(P)>0$ if 0 is a zero of $P(y)$),
  \item \emph{degree of $P(y)$ at $y$} and denote by $\mathrm{deg}_y(P)$ the opposite of the
  valuation of $P(1/y)$.
\end{itemize}
Moreover, if $x:=y^{|\mu(K)|}$, then
$$\mathrm{val}_x(P(y)):=\frac{\mathrm{val}_y(P)}{|\mu(K)|} \,\textrm{ and
}\, \mathrm{deg}_x(P(y)):=\frac{\mathrm{deg}_y(P)}{|\mu(K)|}.$$ For $\chi
\in \mathrm{Irr}(W)$, we define
$$a_{\chi_\phi}:=\mathrm{val}_x(s_{\chi_\phi}(y)) \,\textrm{ and }\,
A_{\chi_\phi}:=\mathrm{deg}_x(s_{\chi_\phi}(y)).$$ The following
result is proven in \cite{BK}, Proposition 2.9.

\begin{proposition}\label{aA}\
\begin{enumerate}
  \item For all $\chi \in \emph{Irr}(W)$, we have
        $$\omega_{\chi_\phi}(\tau)=t_\phi(\tau)x^{a_{\chi_\phi}+A_{\chi_\phi}},$$
        where $\tau$ is the central element of the pure braid group defined in
        $\ref{pi}$.
  \item Let $\chi,\psi \in \mathrm{Irr}(W)$. If $\chi_\phi$ and
        $\psi_\phi$ belong to the same Rouquier block, then
        $$a_{\chi_\phi}+A_{\chi_\phi}=a_{\psi_\phi}+A_{\psi_\phi}.$$
\end{enumerate}
\end{proposition}
\begin{apod}
{\begin{enumerate}
      \item If $P(y) \in \mathbb{C}[y,y^{-1}]$, we denote by
      $P(y)^*$ the polynomial whose coefficients are the complex
      conjugates of those of $P(y)$. By \cite{BMM2}, 2.8, we know
      that the Schur element $s_{\chi_\phi}(y)$ is semi-palindromic and satisfies
      $$s_{\chi_\phi}(y^{-1})^*=\frac{t_\phi(\tau)}{\omega_{\chi_\phi}(\tau)}s_{\chi_\phi}(y).$$
      We deduce (\cite{BMM2}, 6.5, 6.6) that
      $$\frac{t_\phi(\tau)}{\omega_{\chi_\phi}(\tau)}=\xi
      x^{-(a_{\chi_\phi}+A_{\chi_\phi})},$$
      for some $\xi \in \mathbb{C}$.
      For $y=x=1$, the first equation gives
      $t_\phi(\tau)=\omega_{\chi_\phi}(\tau)$ and the second
      one $\xi=1$. Thus we obtain
      $$\omega_{\chi_\phi}(\tau)=t_\phi(\tau)x^{a_{\chi_\phi}+A_{\chi_\phi}}.$$
      \item Suppose that $\chi_\phi$ and $\psi_\phi$ belong to the same Rouquier
      block. Due to Proposition $\ref{Rouquier blocks and central
      characters}$, it is enough to show that if there exists a
      $\phi$-bad prime ideal $\mathfrak{p}$ of $\mathbb{Z}_K$ such that
      $\omega_{\chi_\phi} \equiv \omega_{\psi_\phi}
      \textrm{ mod } \mathfrak{p}\mathcal{O}_{\mathfrak{p}\mathcal{O}}$,
      then $a_{\chi_\phi}+A_{\chi_\phi}=a_{\psi_\phi}+A_{\psi_\phi}.$
      If $\omega_{\chi_\phi} \equiv \omega_{\psi_\phi}
      \textrm{ mod } \mathfrak{p}\mathcal{O}_{\mathfrak{p}\mathcal{O}}$, then, in
      particular, $\omega_{\chi_\phi}(\tau) \equiv
      \omega_{\psi_\phi}(\tau)
      \textrm{ mod } \mathfrak{p}\mathcal{O}_{\mathfrak{p}\mathcal{O}}$. Part 1 implies
      that
      $$t_\phi(\tau)x^{a_{\chi_\phi}+A_{\chi_\phi}} \equiv
      t_\phi(\tau)x^{a_{\psi_\phi}+A_{\psi_\phi}} \textrm{ mod } \mathfrak{p}\mathcal{O}_{\mathfrak{p}\mathcal{O}}.$$
      We know by \cite{BMM2}, 2.1 that $t_\phi(\tau)$ is of the
      form $\xi x^M$, where $\xi$ is a root of unity and $M \in
      \mathbb{Z}$. Thus $t_\phi(\tau) \notin
      \mathfrak{p}\mathcal{O}_{\mathfrak{p}\mathcal{O}}$ and the above congruence gives
      $$x^{a_{\chi_\phi}+A_{\chi_\phi}} \equiv
        x^{a_{\psi_\phi}+A_{\psi_\phi}} \textrm{ mod } \mathfrak{p}\mathcal{O}_{\mathfrak{p}\mathcal{O}},$$
      whence
      $$a_{\chi_\phi}+A_{\chi_\phi}=a_{\psi_\phi}+A_{\psi_\phi}.$$}
\end{enumerate}
\end{apod}
\begin{remark}\emph{ For all Coxeter groups, Lusztig has proved (cf., for example, \cite{Lu2}, 3.3 and 3.4) that if $\chi_\phi$ and
$\psi_\phi$ belong to the same Rouquier block of the Iwahori-Hecke algebra, then
$a_{\chi_\phi}=a_{\psi_\phi}$ and $A_{\chi_\phi}=A_{\psi_\phi}$. This assertion has also been proved
\begin{itemize}
  \item for almost all cyclotomic Hecke algebras of the groups $G(d,1,r)$ and $G(e,e,r)$ in \cite{BK},
  \item for the ``spetsial'' cyclotomic Hecke algebra of the ``spetsial'' exceptional complex reflection groups in
  \cite{MaRo}.
\end{itemize}
Using the results of the next chapter, we have been able to obtain the same result for all cyclotomic Hecke algebras 
\begin{itemize}
  \item of the groups $G(d,1,r)$ in \cite{Chlou2},
  \item of the groups $G(de,e,r)$ in \cite{Chlou3}, and
  \item of all exceptional irreducible complex reflection groups in \cite{DegVal},
\end{itemize}
thus completing its proof for all complex reflection groups.}
\end{remark}

\chapter{On the determination of the Rouquier blocks}

The aim of this chapter is the determination of the Rouquier blocks of the cyclotomic Hecke algebras of all irreducible complex reflection groups. In the last chapter, we saw that the Rouquier blocks have the property of  ``semi-continuity''. This property allows us to obtain the Rouquier blocks for any cyclotomic Hecke algebra by actually calculating them in a small number of cases. Following the theory developed in the two previous chapters, we only need to determine the Rouquier blocks ``associated with no and each essential hyperplane'' for all irreducible complex reflection groups.

For  the exceptional irreducible complex reflection groups, the computations were made 
with the use of the GAP package CHEVIE (cf.~\cite{chevie}). In section 5.2, we give the algorithm which has been used. 
This algorithm is heuristic and was applied only to the groups $G_7$, $G_{11}$, $G_{19}$, $G_{26}$, $G_{28}$ and $G_{32}$. The results presented in the Appendix allow us to use Clifford Theory in order to obtain the Rouquier blocks for the groups $G_4,\ldots,G_{22}$ and $G_{25}$. The remaining groups have already been studied by Malle and Rouquier in \cite{MaRo}. We have stored all the calculated data in a computer file and created GAP functions to display them. These functions are presented in this chapter and can be found on my webpage (\cite{chlouveraki}).

As far as the groups of the infinite series are concerned, Clifford theory again allows us to obtain the Rouquier blocks of the cyclotomic Hecke algebras of $G(de,e,r)$ (when $r>2$ or $r=2$ and $e$ is odd) and $G(2fd,2f,2)$ from those of $G(de,1,r)$ and $G(2fd,2,2)$ respectively. Therefore, only the last two cases need to be studied thoroughly.

In section 5.3, we determine the Rouquier blocks associated with the essential hyperplanes for the group $G(d,1,r)$. The algorithm of Lyle and Mathas (cf.~\cite{LyMa}) for the determination of the blocks of an Ariki-Koike algebra over a field has played a key role in the achievement of this goal.
The description of the Rouquier blocks for $G(d,1,r)$ is combinatorial and demonstrates an unexpected relation between them and the families of characters of the Weyl groups of type $B_n$, $n \leq r$. 

In section 5.4, we calculate the Rouquier blocks associated with no and each essential hyperplane for the group $G(2d,2,2)$. The method used follows the same principles as the algorithm for the exceptional irreducible complex reflection groups.

Finally, in section 5.5, we explain how exactly we apply the results of Clifford theory (Propositions $\ref{1.42}$ and $\ref{1.45}$) to obtain the Rouquier blocks of the cyclotomic Hecke algebras of the groups $G(de,e,r)$.

\section{General principles}

Let $W$ be an irreducible complex reflection group with field of definition $K$
and let $\mathcal{H}$ be its generic Hecke
algebra. We suppose that the assumptions
$\ref{ypo}$ are satisfied. Following Theorem $\ref{Semisimplicity Malle}$, we can find a set of indeterminates $\textbf{v}$ such that the algebra $K(\textbf{v})\mathcal{H}$ is split semisimple.
Set
$A:=\mathbb{Z}_K[\textbf{v},\textbf{v}^{-1}]$ and let us denote by $\mathcal{O}$ the Rouquier ring $\mathcal{R}_K(y)$ of $K$. Let $\mathfrak{p}$ be
a prime ideal of $\mathbb{Z}_K$ lying over a prime number $p$ which
divides the order of the group $W$.  
We can determine the
$\mathfrak{p}$-essential hyperplanes $W$ from the factorization of the
Schur elements of $\mathcal{H}$ over $K[\textbf{v},\textbf{v}^{-1}]$.

Let $\phi_{\emptyset} :v_{\mathcal{C},j} \mapsto y^{n_{\mathcal{C},j}}$ be a
cyclotomic specialization such that the integers $n_{\mathcal{C},j}$ belong to no essential hyperplane for $W$.
Such a cyclotomic specialization will be called  \index{cyclotomic specialization associated with no essential hyperplane}\emph{associated with no essential hyperplane}. The blocks of $\mathcal{O}_{\mathfrak{p}\mathcal{O}}\mathcal{H}_{\phi_\emptyset}$ are called \index{p-blocks associated with no essential hyperplane}$\mathfrak{p}$-\emph{blocks associated with no essential hyperplane} and coincide with the blocks of $A_{\mathfrak{p}A}\mathcal{H}$.

Let $\phi_{H} :v_{\mathcal{C},j} \mapsto y^{n_{\mathcal{C},j}}$ be a
cyclotomic specialization such that the integers $n_{\mathcal{C},j}$ belong to exactly one essential hyperplane $H$, corresponding to the essential monomial $M$. Such a cyclotomic specialization will be called  \index{cyclotomic specialization associated with an essential hyperplane}\emph{associated with the essential hyperplane $H$}.  The blocks of $\mathcal{O}_{\mathfrak{p}\mathcal{O}}\mathcal{H}_{\phi_H}$ are called \index{p-blocks associated with an essential hyperplane}$\mathfrak{p}$-\emph{blocks associated with the essential hyperplane} $H$. 
If $H$ is not $\mathfrak{p}$-essential for $W$, then the blocks of $\mathcal{O}_{\mathfrak{p}\mathcal{O}}\mathcal{H}_{\phi_H}$
coincide with the $\mathfrak{p}$-blocks associated with no essential hyperplane. 
If $H$ is $\mathfrak{p}$-essential for $W$, then the blocks of $\mathcal{O}_{\mathfrak{p}\mathcal{O}}\mathcal{H}_{\phi_H}$ coincide with the blocks of $A_{\mathfrak{q}_M}\mathcal{H}$, where
 $\mathfrak{q}_M=(M-1)A+\mathfrak{p}A$. By Proposition $\ref{simple inclusion}$, the  $\mathfrak{p}$-blocks associated with the essential hyperplane $H$ are unions of $\mathfrak{p}$-blocks associated with no essential hyperplane.
 
Following Proposition $\ref{semicontinuity of Rouquier blocks}$, the Rouquier blocks of $\mathcal{H}_{\phi_\emptyset}$
can be obtained as unions of $\mathfrak{p}$-blocks associated with no essential hyperplane, where
$\mathfrak{p}$ runs over the set of prime ideals of $\mathbb{Z}_K$ lying over the prime divisors of $|W|$ (if $\mathfrak{p}$ is not $\phi_\emptyset$-bad, then the corresponding $\mathfrak{p}$-blocks are trivial). The Rouquier blocks of
$\mathcal{H}_{\phi_\emptyset}$ are the \index{Rouquier blocks associated with no essential hyperplane}\emph{Rouquier blocks associated with no essential hyperplane}. Respectively, the Rouquier blocks of $\mathcal{H}_{\phi_H}$ are the \index{Rouquier blocks associated with an essential hyperplane}\emph{Rouquier blocks associated with the essential hyperplane} $H$. Like the $\mathfrak{p}$-blocks, the  Rouquier blocks associated with the essential hyperplane $H$ are unions of Rouquier blocks associated with no essential hyperplane.

The following result is a consequence of Theorem $\ref{main theorem}$ and summarizes the results of Chapter $4$.

\begin{theorem}\label{summary}
Let $\phi :v_{\mathcal{C},j} \mapsto y^{n_{\mathcal{C},j}}$ be a cyclotomic specialization which is not associated with no essential hyperplane. Let $\mathcal{E}$ be the set of all essential hyperplanes to which the integers $n_{\mathcal{C},j}$ belong. Let $\chi, \psi \in \mathrm{Irr}(W)$. 
The characters
$\chi_\phi$ and $\psi_\phi$ belong to the same block of $\mathcal{O}_{\mathfrak{p}\mathcal{O}}\mathcal{H}_{\phi}$ if and only if there exist a finite sequence
$\chi_0,\chi_1,\ldots,\chi_n \in \emph{Irr}(W)$ and a finite
sequence $H_1,\ldots,H_n \in \mathcal{E}$  such that
\begin{itemize}
  \item $(\chi_0)_\phi=\chi_\phi$ and $(\chi_n)_\phi=\psi_\phi$,
  \item for all $j$ $(1\leq j \leq n)$,\,\,
        $(\chi_{j-1})_\phi$ and $(\chi_j)_\phi$ are in the same $\mathfrak{p}$-block associated with the essential hyperplane $H_j$.
\end{itemize}
Moreover, the characters
$\chi_\phi$ and $\psi_\phi$ belong to the same Rouquier block of $\mathcal{H}_{\phi}$ if and only if there exist a finite sequence
$\chi_0,\chi_1,\ldots,\chi_n \in \emph{Irr}(W)$ and a finite
sequence $H_1,\ldots,H_n \in \mathcal{E}$  such that
\begin{itemize}
  \item $(\chi_0)_\phi=\chi_\phi$ and $(\chi_n)_\phi=\psi_\phi$,
  \item for all $j$ $(1\leq j \leq n)$,\,\,
        $(\chi_{j-1})_\phi$ and $(\chi_j)_\phi$ are in the same  Rouquier block associated with the essential hyperplane $H_j$.
\end{itemize}
\end{theorem}

Thanks to the above theorem, in order to determine the Rouquier blocks of any cyclotomic Hecke algebra associated to $W$, we only need to consider a cyclotomic specialization associated with no and each essential hyperplane and 
\begin{itemize}
\item either calculate their $\mathfrak{p}$-blocks, for all prime ideals $\mathfrak{p}$ lying over the prime divisors of $|W|$, and use Proposition $\ref{semicontinuity of Rouquier blocks}$ in order to obtain their Rouquier blocks,
\item or calculate directly their Rouquier blocks.
\end{itemize}

In the case of the exceptional groups, we will use the first method, whereas in the case of the groups of the infinite series, we will mostly use the second one.
In both cases, we will need some criteria in order to determine the corresponding partitions of $\mathrm{Irr}(W)$ into blocks. These are results which have already been presented in previous chapters, but we are going to repeat here for the convenience of the reader. Once more, let $\phi :v_{\mathcal{C},j} \mapsto y^{n_{\mathcal{C},j}}$ be a cyclotomic specialization and let $\mathfrak{p}$ be a prime ideal of $\mathbb{Z}_K$ lying over the prime number $p$.
\\ \\
\textbf{Proposition $\ref{Malle-Rouquier}$}. An irreducible character $\chi \in \mathrm{Irr}(W)$ is a block of $\mathcal{O}_{\mathfrak{p}\mathcal{O}}\mathcal{H}_\phi$ by
  itself 
  if and only if $s_{\chi_\phi} \notin \mathfrak{p}\mathbb{Z}_K[y,y^{-1}]$.\\
  \\
\textbf{Proposition $\ref{not essential for block}$}. Let $C$ be a block of $A_{\mathfrak{p}A}\mathcal{H}$. If $M$ is an essential monomial for $W$ which is not $\mathfrak{p}$-essential for any $\chi \in C$ , then $C$ is a block of $A_{\mathfrak{q}_M}\mathcal{H}$, where $\mathfrak{q}_M=(M-1)A+\mathfrak{p}A$.\\
\\
\textbf{Proposition $\ref{group algebra}$}. If $\chi,\psi \in \mathrm{Irr}(W)$ belong to the same block of
  $\mathcal{O}_{\mathfrak{p}\mathcal{O}}\mathcal{H}_\phi$, then they are in the
  same $p$-block of $W$.\\ 
  \\
 \textbf{Proposition $\ref{aA}$}. If $\chi,\psi \in \mathrm{Irr}(W)$ are in the same block of
  $\mathcal{O}_{\mathfrak{p}\mathcal{O}}\mathcal{H}_\phi$, then they are in the same Rouquier
  block of $\mathcal{H}_\phi$ and we have
  $$a_{\chi_\phi}+A_{\chi_\phi}=a_{\psi_\phi}+A_{\psi_\phi}.$$

\section{The exceptional irreducible complex reflection groups}

Let $W:=G_n$ $(4 \leq n \leq 37)$ be an irreducible exceptional complex reflection group with field of definition $K$.

If $n \in \{23,24,27,29,30,31,33,34,35,36,37\}$, then $W$ has
only one hyperplane orbit $\mathcal{C}$ with $e_{\mathcal{C}}=2$. The generic Hecke algebra of $W$ is defined over a Laurent polynomial ring in two indeterminates $v_{\mathcal{C},0}$ and $v_{\mathcal{C},1}$ and the only essential monomial for $W$ is
$v_{\mathcal{C},0}v_{\mathcal{C},1}^{-1}$.

If $\phi$ is the ``spetsial'' cyclotomic specialization (see Example $\ref{spetsial}$), then $\phi$ is associated with no essential hyperplane for $W$. The $\mathfrak{p}$-blocks, for all $\phi$-bad prime ideals $\mathfrak{p}$, and the Rouquier blocks of the spetsial cyclotomic Hecke algebra of these groups have been calculated by Malle and Rouquier in \cite{MaRo}.

If $\phi$ is a cyclotomic specialization associated with the unique essential hyperplane for $W$, then $\mathcal{H}_\phi$ is isomorphic to the group algebra $\mathbb{Z}_KW$. Its $p$-blocks are known from Brauer theory, whereas there exists a single Rouquier block (see also \cite{Rou}, \S3, Rem.1).

Therefore, we will only study in detail the remaining cases. 
 
\subsection{Essential hyperplanes}

Let $\mathfrak{p}_1$, $\mathfrak{p}_2$ be two prime ideal of $\mathbb{Z}_K$ lying over the same prime number $p$. If $\Psi$ is a $K$-cyclotomic polynomial, then $\Psi(1) \in \mathfrak{p}_1$ if and only if $\Psi_1 \in \mathfrak{p}_2$. We deduce that an essential hyperplane is $\mathfrak{p}_1$-essential for $W$ if and only if it is $\mathfrak{p}_2$-essential for $W$. Therefore, we can talk about determining the $p$-essential hyperplanes for $W$, where $p$ runs over the set of prime divisors of $|W|$.

Together with Jean Michel, we have programmed into the GAP package CHEVIE the Schur elements of the generic Hecke algebras of all exceptional irreducible complex
reflection groups in factorized form
(function \emph{SchurModels} and \emph{SchurData}). Given a prime ideal $\mathfrak{p}$ of $\mathbb{Z}_K$, GAP
provides us with a way to determine whether an element of
$\mathbb{Z}_K$ belongs to $\mathfrak{p}$. Therefore, we can easily determine the $p$-essential monomials and thus, the $p$-essential hyperplanes for $W$.

In particular, we only need to follow this procedure for the groups $G_7$, $G_{11}$, $G_{19}$, $G_{26}$, $G_{28}$ and $G_{32}$.  In the Appendix, we give the specializations of the parameters which make
\begin{itemize}
\item $\mathcal{H}(G_7)$ the twisted symmetric algebra of some finite cyclic group over $\mathcal{H}(G_4)$, $\mathcal{H}(G_5)$ and $\mathcal{H}(G_6)$.
\item $\mathcal{H}(G_{11})$ the twisted symmetric algebra of some finite cyclic group over $\mathcal{H}(G_8)$, $\mathcal{H}(G_9)$, $\mathcal{H}(G_{10})$,
$\mathcal{H}(G_{12})$, $\mathcal{H}(G_{13})$, $\mathcal{H}(G_{14})$ and
$\mathcal{H}(G_{15})$.
\item $\mathcal{H}(G_{19})$ the twisted symmetric algebra of some finite cyclic group over
$\mathcal{H}(G_{16})$, $\mathcal{H}(G_{17})$, $\mathcal{H}(G_{18})$,
$\mathcal{H}(G_{20})$, $\mathcal{H}(G_{21})$ and $\mathcal{H}(G_{22})$.
\item $\mathcal{H}(G_{26})$ the twisted symmetric algebra of the cyclic group $C_2$ over
$\mathcal{H}(G_{25})$.
\end{itemize}
In all these cases, Proposition $\ref{1.42}$ implies that the Schur elements of the twisted symmetric algebra are scalar multiples of the Schur elements of the subalgebra. Due to the nature of the specializations, we can obtain the essential hyperplanes for the smaller group from the ones of the larger.

\begin{px} \small{\emph{ The essential hyperplanes for $G_7$ are 
given in Example $\ref{example g7}$ (note that different letters
represent different hyperplane orbits).
The only 3-essential hyperplanes for $G_7$ are:
$$\begin{array}{ccc}
    c_1-c_2=0, & c_0-c_1=0, & c_0-c_2=0  \\
    b_1-b_2=0, & b_0-b_1=0, & b_0-b_2=0
  \end{array}$$
All its remaining essential hyperplanes are strictly 2-essential.
From these, we can obtain the $p$-essential hyperplanes (where
$p=2,3$)
\begin{itemize}
  \item for $G_6$ by setting $b_0=b_1=b_2=0$,
  \item for $G_5$ by setting $a_0=a_1=0$,
  \item for $G_4$ by setting $a_0=a_1=b_0=b_1=b_2=0$.
\end{itemize}
}}
\end{px}

We have created the GAP function \emph{EssentialHyperplanes}
which is applied as follows:
\begin{verbatim}
gap> EssentialHyperplanes(W,p);
\end{verbatim}
and returns
\begin{itemize}
  \item the essential hyperplanes for $W$, if $p=0$.
  \item the $p$-essential hyperplanes for $W$, if $p$ divides the order of
  $W$.
  \item error, if $p$ does not divide the order of $W$.
\end{itemize}

\begin{px}\
\begin{verbatim}
gap> W:=ComplexReflectionGroup(4);
gap> EssentialHyperplanes(W,0);
c_1-c_2=0
c_0-c_1=0
c_0-c_2=0
2c_0-c_1-c_2=0
c_0-2c_1+c_2=0
c_0+c_1-2c_2=0
gap> EssentialHyperplanes(W,2);
2c_0-c_1-c_2=0
c_0-2c_1+c_2=0
c_0+c_1-2c_2=0
c_0-c_1=0
c_1-c_2=0
c_0-c_2=0
gap> EssentialHyperplanes(W,3);
c_1-c_2=0
c_0-c_1=0
c_0-c_2=0
gap> EssentialHyperplanes(W,5);
Error, The number p should divide the order of the group
\end{verbatim}
\end{px}

\subsection{Algorithm}

Let $\mathfrak{p}$ be a prime ideal of $\mathbb{Z}_K$ lying over a
prime number $p$ which divides the order of the group $W$. 
In this section, we will present an algorithm for the determination of the $\mathfrak{p}$-blocks associated with no and each essential hyperplane for $W$.  We retake here the notations of section 5.1.\\

If we are interested in calculating the blocks of
$A_{\mathfrak{p}A}\mathcal{H}$, we follow the steps below:
\begin{enumerate}
  \item We select the characters $\chi \in \mathrm{Irr}(W)$ whose generic Schur elements
  belong to $\mathfrak{p}A$.
  The remaining ones will be blocks of $A_{\mathfrak{p}A}\mathcal{H}$ by themselves, due to Proposition
  $\ref{Malle-Rouquier}$. Thus we form a first partition $\lambda_1$
  of $\mathrm{Irr}(W)$; one part formed by the selected characters,
  each remaining character forming a part by itself.
  \item We calculate the $p$-blocks of $W$.
  By Proposition $\ref{group algebra}$, if two irreducible characters are not
  in the same $p$-block of $W$, then they can not be in the same block of $A_{\mathfrak{p}A}\mathcal{H}$.
  We intersect the partition $\lambda_1$ with the partition obtained
  by the $p$-blocks of $W$ and we obtain a finer partition, named
  $\lambda_2$.
  \item We find a cyclotomic specialization $\phi:v_{\mathcal{C},j} \mapsto y^{n_{\mathcal{C},j}}$
  associated with no essential hyperplane by trying and
  checking random values for the $n_{\mathcal{C},j}$.
  Following Proposition $\ref{aA}$, we take the intersection of the
  partition we already have with the subsets of $\mathrm{Irr}(W)$,
  where the sum $a_{\chi_\phi}+A_{\chi_\phi}$ remains constant. This
  procedure is repeated several times, because sometimes the
  partition becomes finer after some repetitions. Finally, we obtain
  the partition $\lambda_3$, which is the finest of all.
\end{enumerate}

If we are interested in calculating the blocks of
$A_{\mathfrak{q}_M}\mathcal{H}$ for some $\mathfrak{p}$-essential
monomial $M$, the procedure is more or less the same:
\begin{enumerate}
  \item We select the characters $\chi \in \mathrm{Irr}(W)$ for which $M$ is a
  $\mathfrak{p}$-essential monomial.
  We form a first partition $\lambda_1$ of $\mathrm{Irr}(W)$; one part formed by the selected characters,
  each remaining character forming a part by itself. The idea is
  that, by Proposition $\ref{not essential for block}$, if $M$ is not
  $\mathfrak{p}$-essential for any character in a block $C$ of
  $A_{\mathfrak{p}A}\mathcal{H}$, then $C$ is a block of
  $A_{\mathfrak{q}_M}\mathcal{H}$. This explains step 4.
  \item We calculate the $p$-blocks of $W$.
  By Proposition $\ref{group algebra}$, if two irreducible characters are not
  in the same $p$-block of $W$, then they can not be in the same block of $A_{\mathfrak{q}_M}\mathcal{H}$.
  We intersect the partition $\lambda_1$ with the partition obtained
  by the $p$-blocks of $W$ and we obtain a finer partition, named
  $\lambda_2$.
  \item We find a cyclotomic specialization $\phi:v_{\mathcal{C},j} \mapsto y^{n_{\mathcal{C},j}}$
  associated with the  $\mathfrak{p}$-essential hyperplane defined by $M$ (again by trying and
  checking random values for the $n_{\mathcal{C},j}$). We repeat the
  third step as described for $A_{\mathfrak{p}A}\mathcal{H}$ to obtain partition $\lambda_3$.
  \item We take the union of $\lambda_3$ and the partition defined by the blocks
  of $A_{\mathfrak{p}A}\mathcal{H}$.
\end{enumerate}

The above algorithm is, due to step 3, heuristic.
However, we will see in the next section that
we only need to apply this algorithm 
to the groups $G_7$, $G_{11}$, $G_{19}$, $G_{26}$, $G_{28}$ and $G_{32}$.
In these cases, we have
been able to determine (using again the criteria presented in section 5.1) that the partition obtained at the end is minimal
and corresponds to the blocks we are looking for.\\
\\
\begin{remark} \emph{Eventually, the above algorithm provides us with the correct Rouquier
blocks for all exceptional irreducible complex reflection groups,
except for $G_{34}$.}
\end{remark}\
\\ \\
\begin{remark}
\emph{If $\mathfrak{p}_1$, $\mathfrak{p}_2$ are two prime ideals of $\mathbb{Z}_K$ lying over the same prime number $p$, we have observed that, for all exceptional irreducible complex reflection groups, the $\mathfrak{p}_1$-blocks always coincide with the $\mathfrak{p}_2$-blocks.  Therefore, we can talk about determining the $p$-blocks associated with no and each essential hyperplane.}
\end{remark}

\subsection{Results}

With the help of the GAP package CHEVIE, we created a program which implements the algorithm of the previous section. Using Proposition $\ref{semicontinuity of Rouquier blocks}$, we have been able to determine the Rouquier blocks associated with no and each essential hyperplane for the
groups $G_7$, $G_{11}$, $G_{19}$, $G_{26}$, $G_{28}$ and $G_{32}$.

Now, Clifford theory allows us to calculate the Rouquier blocks associated with no and each essential hyperplane for the remaining exceptional irreducible complex reflection groups. In all the cases presented in the Appendix, the explicit calculation of the blocks of the twisted symmetric algebras with the use of the algorithm of the previous section has shown that
the assumptions of Corollary $\ref{clifford}$ are satisfied. Moreover, in all these cases,
if $H$ is the twisted symmetric algebra of the finite cyclic group $G$ over $\bar{H}$, then
each irreducible character of $H$ restricts to an irreducible character of $\bar{H}$.
Using the notations of Proposition $\ref{1.42}$, this means that $|\bar{\Omega}|=1$, whence
the blocks of $\bar{H}$ are stable under the action of $G$. We deduce that the
block-idempotents of $H$ and $\bar{H}$ over the Rouquier ring coincide.
In particular, if $C$ is a block (of characters) of $H$, then $\{\mathrm{Res}^H_{\bar{H}}(\chi)\,|\,\chi \in C\}$ is a block of $\bar{H}$. \\

We will give here the example of $G_7$ and show how we obtain the
blocks of $G_6$ from those of $G_7$. Nevertheless, let us first
explain the notations of characters used by the CHEVIE package.

Let $W$ be an exceptional irreducible complex reflection group. For
$\chi \in \mathrm{Irr}(W)$, we set $d(\chi):=\chi(1)$ and we denote
by $b(\chi)$ the valuation of the fake degree of $\chi$ (for the
definition of the fake degree see \cite{Brou}, 1.20). The
irreducible characters $\chi$ of $W$ are determined by the
corresponding pairs $(d(\chi),b(\chi))$ and we write
$\chi=\phi_{d,b}$, where $d:=d(\chi)$ and $b:=b(\chi)$. If two
irreducible characters $\chi$ and $\chi'$ have $d(\chi)=d(\chi')$
and $b(\chi)=b(\chi')$, we use primes `` $'$ '' to distinguish them
(following \cite{Ma3}, \cite{MaRo}).

\begin{px}\label{example g7}
\emph{\small The generic Hecke algebra of $G_7$ is
$$\begin{array}{rccl}
    \mathcal{H}(G_7) & = & <S,T,U \,\,| &  STU=TUS=UST  \\
     &  &  & (S-x_0)(S-x_1)=0 \\
     &  &  & (T-y_0)(T-y_1)(T-y_2)=0 \\
     &  &  & (U-z_0)(U-z_1)(U-z_2)=0>
  \end{array}$$
Let $$\phi : \left\{ 
\begin{array}{ll} 
x_i \mapsto (\zeta_2)^i q^{a_i} &(0 \leq i <2),\\ 
y_j \mapsto (\zeta_3)^j q^{b_j} &(0 \leq j <3),\\
z_k \mapsto (\zeta_3)^k q^{c_k}  &(0 \leq k <3)
\end{array} \right. 
$$
 be a cyclotomic specialization of $\mathcal{H}(G_7)$. The only prime numbers which divide the order of $G_7$ are 2 and 3. Using the algorithm of the previous section, we have determined  the  Rouquier blocks
associated with no and each essential hyperplane for $G_7$. We present here only the non-trivial ones:}
\begin{description}\scriptsize
\item[No essential hyperplane]\hfil\break
$\{\phi_{2,9'},\phi_{2,15}\}$, $\{\phi_{2,7'},\phi_{2,13'}\}$,
$\{\phi_{2, 11'},\phi_{2,5'}\}$, $\{\phi_{2,7''},\phi_{2,13''}\}$,
$\{\phi_{2,11''}, \phi_{2,5''}\}$, $\{\phi_{2,9''},\phi_{2,3'}\}$,
$\{\phi_{2,11'''},\phi_{2, 5'''}\}$,
$\{\phi_{2,9'''},\phi_{2,3''}\}$, $\{\phi_{2,7'''},\phi_{2,1}\}$,
 $\{\phi_{3,6},\phi_{3,10},\phi_{3,2}\}$, $\{\phi_{3,4},\phi_{3,8},\phi_{3,
12}\}$\item[$c_1-c_2=0$]\hfil\break $\{\phi_{1,4'},\phi_{1,8'}\}$,
$\{\phi_{1,8''},\phi_{1,12'}\}$, $\{\phi_{1, 12''},\phi_{1,16}\}$,
$\{\phi_{1,10'},\phi_{1,14'}\}$, $\{\phi_{1,14''}, \phi_{1,18'}\}$,
$\{\phi_{1,18''},\phi_{1,22}\}$, $\{\phi_{2,9'},\phi_{2, 15}\}$,
$\{\phi_{2,7'},\phi_{2,11'},\phi_{2,13'},\phi_{2,5'}\}$, $\{\phi_{2,
7''},\phi_{2,13''}\}$,
$\{\phi_{2,11''},\phi_{2,9''},\phi_{2,5''},\phi_{2, 3'}\}$,
$\{\phi_{2,11'''},\phi_{2,5'''}\}$, $\{\phi_{2,9'''},\phi_{2,7'''},
\phi_{2,3''},\phi_{2,1}\}$, $\{\phi_{3,6},\phi_{3,10},\phi_{3,2}\}$,
 $\{\phi_{3,4},\phi_{3,8},\phi_{3,12}\}$\item[$c_0-c_1=0$]\hfil\break
$\{\phi_{1,0},\phi_{1,4'}\}$, $\{\phi_{1,4''},\phi_{1,8''}\}$,
$\{\phi_{1, 8'''},\phi_{1,12''}\}$, $\{\phi_{1,6},\phi_{1,10'}\}$,
$\{\phi_{1,10''}, \phi_{1,14''}\}$,
$\{\phi_{1,14'''},\phi_{1,18''}\}$,\\ $\{\phi_{2,9'},\phi_{2,
7'},\phi_{2,15},\phi_{2,13'}\}$, $\{\phi_{2,11'},\phi_{2,5'}\}$,
$\{\phi_{2, 7''},\phi_{2,11''},\phi_{2,13''},\phi_{2,5''}\}$,
$\{\phi_{2,9''},\phi_{2, 3'}\}$,\\
$\{\phi_{2,11'''},\phi_{2,9'''},\phi_{2,5'''},\phi_{2,3''}\}$,
 $\{\phi_{2,7'''},\phi_{2,1}\}$, $\{\phi_{3,6},\phi_{3,10},\phi_{3,2}\}$,
 $\{\phi_{3,4},\phi_{3,8},\phi_{3,12}\}$\item[$c_0-c_2=0$]\hfil\break
$\{\phi_{1,0},\phi_{1,8'}\}$, $\{\phi_{1,4''},\phi_{1,12'}\}$,
$\{\phi_{1, 8'''},\phi_{1,16}\}$, $\{\phi_{1,6},\phi_{1,14'}\}$,
$\{\phi_{1,10''},\phi_{1, 18'}\}$,
$\{\phi_{1,14'''},\phi_{1,22}\}$,\\ $\{\phi_{2,9'},\phi_{2,11'},
\phi_{2,15},\phi_{2,5'}\}$, $\{\phi_{2,7'},\phi_{2,13'}\}$,
$\{\phi_{2,7''}, \phi_{2,9''},\phi_{2,13''},\phi_{2,3'}\}$,
$\{\phi_{2,11''},\phi_{2,5''}\}$,\\
 $\{\phi_{2,11'''},\phi_{2,7'''},\phi_{2,5'''},\phi_{2,1}\}$, $\{\phi_{2,
9'''},\phi_{2,3''}\}$, $\{\phi_{3,6},\phi_{3,10},\phi_{3,2}\}$,
$\{\phi_{3,4},
\phi_{3,8},\phi_{3,12}\}$\item[$b_1-b_2=0$]\hfil\break
$\{\phi_{1,4''},\phi_{1,8'''}\}$, $\{\phi_{1,8''},\phi_{1,12''}\}$,
 $\{\phi_{1,12'},\phi_{1,16}\}$, $\{\phi_{1,10''},\phi_{1,14'''}\}$,
 $\{\phi_{1,14''},\phi_{1,18''}\}$,\\ $\{\phi_{1,18'},\phi_{1,22}\}$,
 $\{\phi_{2,9'},\phi_{2,15}\}$, $\{\phi_{2,7'},\phi_{2,13'}\}$, $\{\phi_{2,
11'},\phi_{2,5'}\}$,
$\{\phi_{2,7''},\phi_{2,11'''},\phi_{2,13''},\phi_{2, 5'''}\}$,\\
$\{\phi_{2,11''},\phi_{2,9'''},\phi_{2,5''},\phi_{2,3''}\}$,
 $\{\phi_{2,9''},\phi_{2,7'''},\phi_{2,3'},\phi_{2,1}\}$, $\{\phi_{3,6},
\phi_{3,10},\phi_{3,2}\}$,
$\{\phi_{3,4},\phi_{3,8},\phi_{3,12}\}$\item[$b_0-b_1=0$]\hfil\break
$\{\phi_{1,0},\phi_{1,4''}\}$, $\{\phi_{1,4'},\phi_{1,8''}\}$,
$\{\phi_{1,8'}, \phi_{1,12'}\}$, $\{\phi_{1,6},\phi_{1,10''}\}$,
$\{\phi_{1,10'},\phi_{1, 14''}\}$,
$\{\phi_{1,14'},\phi_{1,18'}\}$,\\
$\{\phi_{2,9'},\phi_{2,7''}, \phi_{2,15},\phi_{2,13''}\}$,
$\{\phi_{2,7'},\phi_{2,11''},\phi_{2,13'}, \phi_{2,5''}\}$,
$\{\phi_{2,11'},\phi_{2,9''},\phi_{2,5'},\phi_{2,3'}\}$,\\
 $\{\phi_{2,11'''},\phi_{2,5'''}\}$, $\{\phi_{2,9'''},\phi_{2,3''}\}$,
 $\{\phi_{2,7'''},\phi_{2,1}\}$, $\{\phi_{3,6},\phi_{3,10},\phi_{3,2}\}$,
 $\{\phi_{3,4},\phi_{3,8},\phi_{3,12}\}$\item[$b_0-b_2=0$]\hfil\break
$\{\phi_{1,0},\phi_{1,8'''}\}$, $\{\phi_{1,4'},\phi_{1,12''}\}$,
$\{\phi_{1, 8'},\phi_{1,16}\}$, $\{\phi_{1,6},\phi_{1,14'''}\}$,
$\{\phi_{1,10'},\phi_{1, 18''}\}$, $\{\phi_{1,14'},\phi_{1,22}\}$,\\
$\{\phi_{2,9'},\phi_{2,11'''}, \phi_{2,15},\phi_{2,5'''}\}$,
$\{\phi_{2,7'},\phi_{2,9'''},\phi_{2,13'}, \phi_{2,3''}\}$,
$\{\phi_{2,11'},\phi_{2,7'''},\phi_{2,5'},\phi_{2,1}\}$,
 $\{\phi_{2,7''},\phi_{2,13''}\}$,\\ $\{\phi_{2,11''},\phi_{2,5''}\}$,
 $\{\phi_{2,9''},\phi_{2,3'}\}$, $\{\phi_{3,6},\phi_{3,10},\phi_{3,2}\}$,
$\{\phi_{3,4},\phi_{3,8},\phi_{3,12}\}$\item[$a_0-a_1-2b_0+b_1+b_2-2c_0+c_1+c_2=0$]\hfil\break
$\{\phi_{1,6},\phi_{2,9'},\phi_{2,15},\phi_{3,4},\phi_{3,8},\phi_{3,12}\}$,
 $\{\phi_{2,7'},\phi_{2,13'}\}$, $\{\phi_{2,11'},\phi_{2,5'}\}$, $\{\phi_{2,
7''},\phi_{2,13''}\}$, $\{\phi_{2,11''},\phi_{2,5''}\}$,\\
$\{\phi_{2,9''}, \phi_{2,3'}\}$, $\{\phi_{2,11'''},\phi_{2,5'''}\}$,
$\{\phi_{2,9'''},\phi_{2, 3''}\}$, $\{\phi_{2,7'''},\phi_{2,1}\}$,
$\{\phi_{3,6},\phi_{3,10},\phi_{3,
2}\}$\item[$a_0-a_1-2b_0+b_1+b_2+c_0-2c_1+c_2=0$]\hfil\break
$\{\phi_{1,10'},\phi_{2,7'},\phi_{2,13'},\phi_{3,4},\phi_{3,8},\phi_{3,12}\}$,
 $\{\phi_{2,9'},\phi_{2,15}\}$, $\{\phi_{2,11'},\phi_{2,5'}\}$, $\{\phi_{2,
7''},\phi_{2,13''}\}$, $\{\phi_{2,11''},\phi_{2,5''}\}$,\\
$\{\phi_{2,9''}, \phi_{2,3'}\}$, $\{\phi_{2,11'''},\phi_{2,5'''}\}$,
$\{\phi_{2,9'''},\phi_{2, 3''}\}$, $\{\phi_{2,7'''},\phi_{2,1}\}$,
$\{\phi_{3,6},\phi_{3,10},\phi_{3,
2}\}$\item[$a_0-a_1-2b_0+b_1+b_2+c_0+c_1-2c_2=0$]\hfil\break
$\{\phi_{1,14'},\phi_{2,11'},\phi_{2,5'},\phi_{3,4},\phi_{3,8},\phi_{3,12}\}$,
 $\{\phi_{2,9'},\phi_{2,15}\}$, $\{\phi_{2,7'},\phi_{2,13'}\}$, $\{\phi_{2,
7''},\phi_{2,13''}\}$, $\{\phi_{2,11''},\phi_{2,5''}\}$,\\
$\{\phi_{2,9''}, \phi_{2,3'}\}$, $\{\phi_{2,11'''},\phi_{2,5'''}\}$,
$\{\phi_{2,9'''},\phi_{2, 3''}\}$, $\{\phi_{2,7'''},\phi_{2,1}\}$,
$\{\phi_{3,6},\phi_{3,10},\phi_{3,
2}\}$\item[$a_0-a_1-b_0-b_1+2b_2-c_0-c_1+2c_2=0$]\hfil\break
$\{\phi_{1,16},\phi_{2,7'''},\phi_{2,1},\phi_{3,6},\phi_{3,10},\phi_{3,2}\}$,
 $\{\phi_{2,9'},\phi_{2,15}\}$, $\{\phi_{2,7'},\phi_{2,13'}\}$, $\{\phi_{2,
11'},\phi_{2,5'}\}$, $\{\phi_{2,7''},\phi_{2,13''}\}$,\\
$\{\phi_{2,11''}, \phi_{2,5''}\}$, $\{\phi_{2,9''},\phi_{2,3'}\}$,
$\{\phi_{2,11'''},\phi_{2, 5'''}\}$,
$\{\phi_{2,9'''},\phi_{2,3''}\}$, $\{\phi_{3,4},\phi_{3,8},\phi_{3,
12}\}$\item[$a_0-a_1-b_0-b_1+2b_2-c_0+2c_1-c_2=0$]\hfil\break
$\{\phi_{1,12''},\phi_{2,9'''},\phi_{2,3''},\phi_{3,6},\phi_{3,10},\phi_{3,
2}\}$, $\{\phi_{2,9'},\phi_{2,15}\}$,
$\{\phi_{2,7'},\phi_{2,13'}\}$,
 $\{\phi_{2,11'},\phi_{2,5'}\}$, $\{\phi_{2,7''},\phi_{2,13''}\}$,\\ $\{\phi_{2,
11''},\phi_{2,5''}\}$, $\{\phi_{2,9''},\phi_{2,3'}\}$,
$\{\phi_{2,11'''}, \phi_{2,5'''}\}$, $\{\phi_{2,7'''},\phi_{2,1}\}$,
$\{\phi_{3,4},\phi_{3,8},
\phi_{3,12}\}$\item[$a_0-a_1-b_0-b_1+2b_2+2c_0-c_1-c_2=0$]\hfil\break
$\{\phi_{1,8'''},\phi_{2,11'''},\phi_{2,5'''},\phi_{3,6},\phi_{3,10},\phi_{3,
2}\}$, $\{\phi_{2,9'},\phi_{2,15}\}$,
$\{\phi_{2,7'},\phi_{2,13'}\}$,
 $\{\phi_{2,11'},\phi_{2,5'}\}$, $\{\phi_{2,7''},\phi_{2,13''}\}$,\\ $\{\phi_{2,
11''},\phi_{2,5''}\}$, $\{\phi_{2,9''},\phi_{2,3'}\}$,
$\{\phi_{2,9'''}, \phi_{2,3''}\}$, $\{\phi_{2,7'''},\phi_{2,1}\}$,
$\{\phi_{3,4},\phi_{3,8},
\phi_{3,12}\}$\item[$a_0-a_1-b_0+b_2-c_0+c_1=0$]\hfil\break
$\{\phi_{1,12''},\phi_{1,6},\phi_{2,9''},\phi_{2,3'}\}$,
$\{\phi_{2,9'}, \phi_{2,15}\}$, $\{\phi_{2,7'},\phi_{2,13'}\}$,
$\{\phi_{2,11'},\phi_{2, 5'}\}$, $\{\phi_{2,7''},\phi_{2,13''}\}$,
$\{\phi_{2,11''},\phi_{2,5''}\}$,\\
 $\{\phi_{2,11'''},\phi_{2,5'''}\}$, $\{\phi_{2,9'''},\phi_{2,3''}\}$,
 $\{\phi_{2,7'''},\phi_{2,1}\}$, $\{\phi_{3,6},\phi_{3,10},\phi_{3,2}\}$,
$\{\phi_{3,4},\phi_{3,8},\phi_{3,12}\}$\item[$a_0-a_1-b_0+b_2-c_1+c_2=0$]\hfil\break
$\{\phi_{1,16},\phi_{1,10'},\phi_{2,7''},\phi_{2,13''}\}$,
$\{\phi_{2,9'}, \phi_{2,15}\}$, $\{\phi_{2,7'},\phi_{2,13'}\}$,
$\{\phi_{2,11'},\phi_{2, 5'}\}$, $\{\phi_{2,11''},\phi_{2,5''}\}$,
$\{\phi_{2,9''},\phi_{2,3'}\}$,\\
 $\{\phi_{2,11'''},\phi_{2,5'''}\}$, $\{\phi_{2,9'''},\phi_{2,3''}\}$,
 $\{\phi_{2,7'''},\phi_{2,1}\}$, $\{\phi_{3,6},\phi_{3,10},\phi_{3,2}\}$,
$\{\phi_{3,4},\phi_{3,8},\phi_{3,12}\}$\item[$a_0-a_1-b_0+b_2+c_0-c_2=0$]\hfil\break
$\{\phi_{1,8'''},\phi_{1,14'},\phi_{2,11''},\phi_{2,5''}\}$,
$\{\phi_{2,9'}, \phi_{2,15}\}$, $\{\phi_{2,7'},\phi_{2,13'}\}$,
$\{\phi_{2,11'},\phi_{2, 5'}\}$, $\{\phi_{2,7''},\phi_{2,13''}\}$,
$\{\phi_{2,9''},\phi_{2,3'}\}$,\\
 $\{\phi_{2,11'''},\phi_{2,5'''}\}$, $\{\phi_{2,9'''},\phi_{2,3''}\}$,
 $\{\phi_{2,7'''},\phi_{2,1}\}$, $\{\phi_{3,6},\phi_{3,10},\phi_{3,2}\}$,
$\{\phi_{3,4},\phi_{3,8},\phi_{3,12}\}$\item[$a_0-a_1-b_0+b_1-c_0+c_2=0$]\hfil\break
$\{\phi_{1,12'},\phi_{1,6},\phi_{2,9'''},\phi_{2,3''}\}$,
$\{\phi_{2,9'}, \phi_{2,15}\}$, $\{\phi_{2,7'},\phi_{2,13'}\}$,
$\{\phi_{2,11'},\phi_{2, 5'}\}$, $\{\phi_{2,7''},\phi_{2,13''}\}$,
$\{\phi_{2,11''},\phi_{2,5''}\}$,\\
 $\{\phi_{2,9''},\phi_{2,3'}\}$, $\{\phi_{2,11'''},\phi_{2,5'''}\}$,
 $\{\phi_{2,7'''},\phi_{2,1}\}$, $\{\phi_{3,6},\phi_{3,10},\phi_{3,2}\}$,
$\{\phi_{3,4},\phi_{3,8},\phi_{3,12}\}$\item[$a_0-a_1-b_0+b_1+c_1-c_2=0$]\hfil\break
$\{\phi_{1,8''},\phi_{1,14'},\phi_{2,11'''},\phi_{2,5'''}\}$,
$\{\phi_{2,9'}, \phi_{2,15}\}$, $\{\phi_{2,7'},\phi_{2,13'}\}$,
$\{\phi_{2,11'},\phi_{2, 5'}\}$, $\{\phi_{2,7''},\phi_{2,13''}\}$,
$\{\phi_{2,11''},\phi_{2,5''}\}$,\\
 $\{\phi_{2,9''},\phi_{2,3'}\}$, $\{\phi_{2,9'''},\phi_{2,3''}\}$, $\{\phi_{2,
7'''},\phi_{2,1}\}$, $\{\phi_{3,6},\phi_{3,10},\phi_{3,2}\}$,
$\{\phi_{3,4},
\phi_{3,8},\phi_{3,12}\}$\item[$a_0-a_1-b_0+b_1+c_0-c_1=0$]\hfil\break
$\{\phi_{1,4''},\phi_{1,10'},\phi_{2,7'''},\phi_{2,1}\}$,
$\{\phi_{2,9'}, \phi_{2,15}\}$, $\{\phi_{2,7'},\phi_{2,13'}\}$,
$\{\phi_{2,11'},\phi_{2, 5'}\}$, $\{\phi_{2,7''},\phi_{2,13''}\}$,
$\{\phi_{2,11''},\phi_{2,5''}\}$,\\
 $\{\phi_{2,9''},\phi_{2,3'}\}$, $\{\phi_{2,11'''},\phi_{2,5'''}\}$,
 $\{\phi_{2,9'''},\phi_{2,3''}\}$, $\{\phi_{3,6},\phi_{3,10},\phi_{3,2}\}$,
$\{\phi_{3,4},\phi_{3,8},\phi_{3,12}\}$\item[$a_0-a_1-b_0+2b_1-b_2-c_0-c_1+2c_2=0$]\hfil\break
$\{\phi_{1,12'},\phi_{2,9''},\phi_{2,3'},\phi_{3,6},\phi_{3,10},\phi_{3,2}\}$,
 $\{\phi_{2,9'},\phi_{2,15}\}$, $\{\phi_{2,7'},\phi_{2,13'}\}$, $\{\phi_{2,
11'},\phi_{2,5'}\}$, $\{\phi_{2,7''},\phi_{2,13''}\}$,\\
$\{\phi_{2,11''}, \phi_{2,5''}\}$,
$\{\phi_{2,11'''},\phi_{2,5'''}\}$, $\{\phi_{2,9'''},\phi_{2,
3''}\}$, $\{\phi_{2,7'''},\phi_{2,1}\}$,
$\{\phi_{3,4},\phi_{3,8},\phi_{3,
12}\}$\item[$a_0-a_1-b_0+2b_1-b_2-c_0+2c_1-c_2=0$]\hfil\break
$\{\phi_{1,8''},\phi_{2,11''},\phi_{2,5''},\phi_{3,6},\phi_{3,10},\phi_{3,
2}\}$, $\{\phi_{2,9'},\phi_{2,15}\}$,
$\{\phi_{2,7'},\phi_{2,13'}\}$,
 $\{\phi_{2,11'},\phi_{2,5'}\}$, $\{\phi_{2,7''},\phi_{2,13''}\}$,\\ $\{\phi_{2,
9''},\phi_{2,3'}\}$, $\{\phi_{2,11'''},\phi_{2,5'''}\}$,
$\{\phi_{2,9'''}, \phi_{2,3''}\}$, $\{\phi_{2,7'''},\phi_{2,1}\}$,
$\{\phi_{3,4},\phi_{3,8},
\phi_{3,12}\}$\item[$a_0-a_1-b_0+2b_1-b_2+2c_0-c_1-c_2=0$]\hfil\break
$\{\phi_{1,4''},\phi_{2,7''},\phi_{2,13''},\phi_{3,6},\phi_{3,10},\phi_{3,
2}\}$, $\{\phi_{2,9'},\phi_{2,15}\}$,
$\{\phi_{2,7'},\phi_{2,13'}\}$,
 $\{\phi_{2,11'},\phi_{2,5'}\}$, $\{\phi_{2,11''},\phi_{2,5''}\}$,\\ $\{\phi_{2,
9''},\phi_{2,3'}\}$, $\{\phi_{2,11'''},\phi_{2,5'''}\}$,
$\{\phi_{2,9'''}, \phi_{2,3''}\}$, $\{\phi_{2,7'''},\phi_{2,1}\}$,
$\{\phi_{3,4},\phi_{3,8},
\phi_{3,12}\}$\item[$a_0-a_1-b_1+b_2-c_0+c_2=0$]\hfil\break
$\{\phi_{1,16},\phi_{1,10''},\phi_{2,7'},\phi_{2,13'}\}$,
$\{\phi_{2,9'}, \phi_{2,15}\}$, $\{\phi_{2,11'},\phi_{2,5'}\}$,
$\{\phi_{2,7''},\phi_{2, 13''}\}$, $\{\phi_{2,11''},\phi_{2,5''}\}$,
$\{\phi_{2,9''},\phi_{2,3'}\}$,\\
 $\{\phi_{2,11'''},\phi_{2,5'''}\}$, $\{\phi_{2,9'''},\phi_{2,3''}\}$,
 $\{\phi_{2,7'''},\phi_{2,1}\}$, $\{\phi_{3,6},\phi_{3,10},\phi_{3,2}\}$,
$\{\phi_{3,4},\phi_{3,8},\phi_{3,12}\}$\item[$a_0-a_1-b_1+b_2+c_1-c_2=0$]\hfil\break
$\{\phi_{1,12''},\phi_{1,18'},\phi_{2,9'},\phi_{2,15}\}$,
$\{\phi_{2,7'}, \phi_{2,13'}\}$, $\{\phi_{2,11'},\phi_{2,5'}\}$,
$\{\phi_{2,7''},\phi_{2, 13''}\}$, $\{\phi_{2,11''},\phi_{2,5''}\}$,
$\{\phi_{2,9''},\phi_{2,3'}\}$,\\
 $\{\phi_{2,11'''},\phi_{2,5'''}\}$, $\{\phi_{2,9'''},\phi_{2,3''}\}$,
 $\{\phi_{2,7'''},\phi_{2,1}\}$, $\{\phi_{3,6},\phi_{3,10},\phi_{3,2}\}$,
$\{\phi_{3,4},\phi_{3,8},\phi_{3,12}\}$\item[$a_0-a_1-b_1+b_2+c_0-c_1=0$]\hfil\break
$\{\phi_{1,8'''},\phi_{1,14''},\phi_{2,11'},\phi_{2,5'}\}$,
$\{\phi_{2,9'}, \phi_{2,15}\}$, $\{\phi_{2,7'},\phi_{2,13'}\}$,
$\{\phi_{2,7''},\phi_{2, 13''}\}$, $\{\phi_{2,11''},\phi_{2,5''}\}$,
$\{\phi_{2,9''},\phi_{2,3'}\}$,\\
 $\{\phi_{2,11'''},\phi_{2,5'''}\}$, $\{\phi_{2,9'''},\phi_{2,3''}\}$,
 $\{\phi_{2,7'''},\phi_{2,1}\}$, $\{\phi_{3,6},\phi_{3,10},\phi_{3,2}\}$,
 $\{\phi_{3,4},\phi_{3,8},\phi_{3,12}\}$\item[$a_0-a_1=0$]\hfil\break
$\{\phi_{1,0},\phi_{1,6}\}$, $\{\phi_{1,4'},\phi_{1,10'}\}$,
$\{\phi_{1,8'}, \phi_{1,14'}\}$, $\{\phi_{1,4''},\phi_{1,10''}\}$,
$\{\phi_{1,8''},\phi_{1, 14''}\}$, $\{\phi_{1,12'},\phi_{1,18'}\}$,
$\{\phi_{1,8'''},\phi_{1,14'''}\}$,
 $\{\phi_{1,12''},\phi_{1,18''}\}$, $\{\phi_{1,16},\phi_{1,22}\}$, $\{\phi_{2,
9'},\phi_{2,15}\}$, $\{\phi_{2,7'},\phi_{2,13'}\}$,
$\{\phi_{2,11'},\phi_{2, 5'}\}$, $\{\phi_{2,7''},\phi_{2,13''}\}$,
$\{\phi_{2,11''},\phi_{2,5''}\}$,
 $\{\phi_{2,9''},\phi_{2,3'}\}$, $\{\phi_{2,11'''},\phi_{2,5'''}\}$,
 $\{\phi_{2,9'''},\phi_{2,3''}\}$, $\{\phi_{2,7'''},\phi_{2,1}\}$, $\{\phi_{3,
6},\phi_{3,4},\phi_{3,10},\phi_{3,8},\phi_{3,2},\phi_{3,12}\}$\item[$a_0-a_1+b_1-b_2-c_0+c_1=0$]\hfil\break
$\{\phi_{1,8''},\phi_{1,14'''},\phi_{2,11'},\phi_{2,5'}\}$,
$\{\phi_{2,9'}, \phi_{2,15}\}$, $\{\phi_{2,7'},\phi_{2,13'}\}$,
$\{\phi_{2,7''},\phi_{2, 13''}\}$, $\{\phi_{2,11''},\phi_{2,5''}\}$,
$\{\phi_{2,9''},\phi_{2,3'}\}$,\\
 $\{\phi_{2,11'''},\phi_{2,5'''}\}$, $\{\phi_{2,9'''},\phi_{2,3''}\}$,
 $\{\phi_{2,7'''},\phi_{2,1}\}$, $\{\phi_{3,6},\phi_{3,10},\phi_{3,2}\}$,
$\{\phi_{3,4},\phi_{3,8},\phi_{3,12}\}$\item[$a_0-a_1+b_1-b_2-c_1+c_2=0$]\hfil\break
$\{\phi_{1,12'},\phi_{1,18''},\phi_{2,9'},\phi_{2,15}\}$,
$\{\phi_{2,7'}, \phi_{2,13'}\}$, $\{\phi_{2,11'},\phi_{2,5'}\}$,
$\{\phi_{2,7''},\phi_{2, 13''}\}$, $\{\phi_{2,11''},\phi_{2,5''}\}$,
$\{\phi_{2,9''},\phi_{2,3'}\}$,\\
 $\{\phi_{2,11'''},\phi_{2,5'''}\}$, $\{\phi_{2,9'''},\phi_{2,3''}\}$,
 $\{\phi_{2,7'''},\phi_{2,1}\}$, $\{\phi_{3,6},\phi_{3,10},\phi_{3,2}\}$,
$\{\phi_{3,4},\phi_{3,8},\phi_{3,12}\}$\item[$a_0-a_1+b_1-b_2+c_0-c_2=0$]\hfil\break
$\{\phi_{1,4''},\phi_{1,22},\phi_{2,7'},\phi_{2,13'}\}$,
$\{\phi_{2,9'}, \phi_{2,15}\}$, $\{\phi_{2,11'},\phi_{2,5'}\}$,
$\{\phi_{2,7''},\phi_{2, 13''}\}$, $\{\phi_{2,11''},\phi_{2,5''}\}$,
$\{\phi_{2,9''},\phi_{2,3'}\}$,\\
 $\{\phi_{2,11'''},\phi_{2,5'''}\}$, $\{\phi_{2,9'''},\phi_{2,3''}\}$,
 $\{\phi_{2,7'''},\phi_{2,1}\}$, $\{\phi_{3,6},\phi_{3,10},\phi_{3,2}\}$,
$\{\phi_{3,4},\phi_{3,8},\phi_{3,12}\}$\item[$a_0-a_1+b_0-2b_1+b_2-2c_0+c_1+c_2=0$]\hfil\break
$\{\phi_{1,10''},\phi_{2,7''},\phi_{2,13''},\phi_{3,4},\phi_{3,8},\phi_{3,
12}\}$, $\{\phi_{2,9'},\phi_{2,15}\}$,
$\{\phi_{2,7'},\phi_{2,13'}\}$,
 $\{\phi_{2,11'},\phi_{2,5'}\}$, $\{\phi_{2,11''},\phi_{2,5''}\}$,\\ $\{\phi_{2,
9''},\phi_{2,3'}\}$, $\{\phi_{2,11'''},\phi_{2,5'''}\}$,
$\{\phi_{2,9'''}, \phi_{2,3''}\}$, $\{\phi_{2,7'''},\phi_{2,1}\}$,
$\{\phi_{3,6},\phi_{3,10},
\phi_{3,2}\}$\item[$a_0-a_1+b_0-2b_1+b_2+c_0-2c_1+c_2=0$]\hfil\break
$\{\phi_{1,14''},\phi_{2,11''},\phi_{2,5''},\phi_{3,4},\phi_{3,8},\phi_{3,
12}\}$, $\{\phi_{2,9'},\phi_{2,15}\}$,
$\{\phi_{2,7'},\phi_{2,13'}\}$,
 $\{\phi_{2,11'},\phi_{2,5'}\}$, $\{\phi_{2,7''},\phi_{2,13''}\}$,\\ $\{\phi_{2,
9''},\phi_{2,3'}\}$, $\{\phi_{2,11'''},\phi_{2,5'''}\}$,
$\{\phi_{2,9'''}, \phi_{2,3''}\}$, $\{\phi_{2,7'''},\phi_{2,1}\}$,
$\{\phi_{3,6},\phi_{3,10},
\phi_{3,2}\}$\item[$a_0-a_1+b_0-2b_1+b_2+c_0+c_1-2c_2=0$]\hfil\break
$\{\phi_{1,18'},\phi_{2,9''},\phi_{2,3'},\phi_{3,4},\phi_{3,8},\phi_{3,12}\}$,
 $\{\phi_{2,9'},\phi_{2,15}\}$, $\{\phi_{2,7'},\phi_{2,13'}\}$, $\{\phi_{2,
11'},\phi_{2,5'}\}$, $\{\phi_{2,7''},\phi_{2,13''}\}$,\\
$\{\phi_{2,11''}, \phi_{2,5''}\}$,
$\{\phi_{2,11'''},\phi_{2,5'''}\}$, $\{\phi_{2,9'''},\phi_{2,
3''}\}$, $\{\phi_{2,7'''},\phi_{2,1}\}$,
$\{\phi_{3,6},\phi_{3,10},\phi_{3,
2}\}$\item[$a_0-a_1+b_0-b_1-c_0+c_1=0$]\hfil\break
$\{\phi_{1,4'},\phi_{1,10''},\phi_{2,7'''},\phi_{2,1}\}$,
$\{\phi_{2,9'}, \phi_{2,15}\}$, $\{\phi_{2,7'},\phi_{2,13'}\}$,
$\{\phi_{2,11'},\phi_{2, 5'}\}$, $\{\phi_{2,7''},\phi_{2,13''}\}$,
$\{\phi_{2,11''},\phi_{2,5''}\}$,\\
 $\{\phi_{2,9''},\phi_{2,3'}\}$, $\{\phi_{2,11'''},\phi_{2,5'''}\}$,
 $\{\phi_{2,9'''},\phi_{2,3''}\}$, $\{\phi_{3,6},\phi_{3,10},\phi_{3,2}\}$,
$\{\phi_{3,4},\phi_{3,8},\phi_{3,12}\}$\item[$a_0-a_1+b_0-b_1-c_1+c_2=0$]\hfil\break
$\{\phi_{1,8'},\phi_{1,14''},\phi_{2,11'''},\phi_{2,5'''}\}$,
$\{\phi_{2,9'}, \phi_{2,15}\}$, $\{\phi_{2,7'},\phi_{2,13'}\}$,
$\{\phi_{2,11'},\phi_{2, 5'}\}$, $\{\phi_{2,7''},\phi_{2,13''}\}$,
$\{\phi_{2,11''},\phi_{2,5''}\}$,\\
 $\{\phi_{2,9''},\phi_{2,3'}\}$, $\{\phi_{2,9'''},\phi_{2,3''}\}$, $\{\phi_{2,
7'''},\phi_{2,1}\}$, $\{\phi_{3,6},\phi_{3,10},\phi_{3,2}\}$,
$\{\phi_{3,4},
\phi_{3,8},\phi_{3,12}\}$\item[$a_0-a_1+b_0-b_1+c_0-c_2=0$]\hfil\break
$\{\phi_{1,0},\phi_{1,18'},\phi_{2,9'''},\phi_{2,3''}\}$,
$\{\phi_{2,9'}, \phi_{2,15}\}$, $\{\phi_{2,7'},\phi_{2,13'}\}$,
$\{\phi_{2,11'},\phi_{2, 5'}\}$, $\{\phi_{2,7''},\phi_{2,13''}\}$,
$\{\phi_{2,11''},\phi_{2,5''}\}$,\\
 $\{\phi_{2,9''},\phi_{2,3'}\}$, $\{\phi_{2,11'''},\phi_{2,5'''}\}$,
 $\{\phi_{2,7'''},\phi_{2,1}\}$, $\{\phi_{3,6},\phi_{3,10},\phi_{3,2}\}$,
$\{\phi_{3,4},\phi_{3,8},\phi_{3,12}\}$\item[$a_0-a_1+b_0-b_2-c_0+c_2=0$]\hfil\break
$\{\phi_{1,8'},\phi_{1,14'''},\phi_{2,11''},\phi_{2,5''}\}$,
$\{\phi_{2,9'}, \phi_{2,15}\}$, $\{\phi_{2,7'},\phi_{2,13'}\}$,
$\{\phi_{2,11'},\phi_{2, 5'}\}$, $\{\phi_{2,7''},\phi_{2,13''}\}$,
$\{\phi_{2,9''},\phi_{2,3'}\}$,\\
 $\{\phi_{2,11'''},\phi_{2,5'''}\}$, $\{\phi_{2,9'''},\phi_{2,3''}\}$,
 $\{\phi_{2,7'''},\phi_{2,1}\}$, $\{\phi_{3,6},\phi_{3,10},\phi_{3,2}\}$,
$\{\phi_{3,4},\phi_{3,8},\phi_{3,12}\}$\item[$a_0-a_1+b_0-b_2+c_1-c_2=0$]\hfil\break
$\{\phi_{1,4'},\phi_{1,22},\phi_{2,7''},\phi_{2,13''}\}$,
$\{\phi_{2,9'}, \phi_{2,15}\}$, $\{\phi_{2,7'},\phi_{2,13'}\}$,
$\{\phi_{2,11'},\phi_{2, 5'}\}$, $\{\phi_{2,11''},\phi_{2,5''}\}$,
$\{\phi_{2,9''},\phi_{2,3'}\}$,\\
 $\{\phi_{2,11'''},\phi_{2,5'''}\}$, $\{\phi_{2,9'''},\phi_{2,3''}\}$,
 $\{\phi_{2,7'''},\phi_{2,1}\}$, $\{\phi_{3,6},\phi_{3,10},\phi_{3,2}\}$,
$\{\phi_{3,4},\phi_{3,8},\phi_{3,12}\}$\item[$a_0-a_1+b_0-b_2+c_0-c_1=0$]\hfil\break
$\{\phi_{1,0},\phi_{1,18''},\phi_{2,9''},\phi_{2,3'}\}$,
$\{\phi_{2,9'}, \phi_{2,15}\}$, $\{\phi_{2,7'},\phi_{2,13'}\}$,
$\{\phi_{2,11'},\phi_{2, 5'}\}$, $\{\phi_{2,7''},\phi_{2,13''}\}$,
$\{\phi_{2,11''},\phi_{2,5''}\}$,\\
 $\{\phi_{2,11'''},\phi_{2,5'''}\}$, $\{\phi_{2,9'''},\phi_{2,3''}\}$,
 $\{\phi_{2,7'''},\phi_{2,1}\}$, $\{\phi_{3,6},\phi_{3,10},\phi_{3,2}\}$,
$\{\phi_{3,4},\phi_{3,8},\phi_{3,12}\}$\item[$a_0-a_1+b_0+b_1-2b_2-2c_0+c_1+c_2=0$]\hfil\break
$\{\phi_{1,14'''},\phi_{2,11'''},\phi_{2,5'''},\phi_{3,4},\phi_{3,8},\phi_{3,
12}\}$, $\{\phi_{2,9'},\phi_{2,15}\}$,
$\{\phi_{2,7'},\phi_{2,13'}\}$,
 $\{\phi_{2,11'},\phi_{2,5'}\}$, $\{\phi_{2,7''},\phi_{2,13''}\}$,\\ $\{\phi_{2,
11''},\phi_{2,5''}\}$, $\{\phi_{2,9''},\phi_{2,3'}\}$,
$\{\phi_{2,9'''}, \phi_{2,3''}\}$, $\{\phi_{2,7'''},\phi_{2,1}\}$,
$\{\phi_{3,6},\phi_{3,10},
\phi_{3,2}\}$\item[$a_0-a_1+b_0+b_1-2b_2+c_0-2c_1+c_2=0$]\hfil\break
$\{\phi_{1,18''},\phi_{2,9'''},\phi_{2,3''},\phi_{3,4},\phi_{3,8},\phi_{3,
12}\}$, $\{\phi_{2,9'},\phi_{2,15}\}$,
$\{\phi_{2,7'},\phi_{2,13'}\}$,
 $\{\phi_{2,11'},\phi_{2,5'}\}$, $\{\phi_{2,7''},\phi_{2,13''}\}$,\\ $\{\phi_{2,
11''},\phi_{2,5''}\}$, $\{\phi_{2,9''},\phi_{2,3'}\}$,
$\{\phi_{2,11'''}, \phi_{2,5'''}\}$, $\{\phi_{2,7'''},\phi_{2,1}\}$,
$\{\phi_{3,6},\phi_{3,10},
\phi_{3,2}\}$\item[$a_0-a_1+b_0+b_1-2b_2+c_0+c_1-2c_2=0$]\hfil\break
$\{\phi_{1,22},\phi_{2,7'''},\phi_{2,1},\phi_{3,4},\phi_{3,8},\phi_{3,12}\}$,
 $\{\phi_{2,9'},\phi_{2,15}\}$, $\{\phi_{2,7'},\phi_{2,13'}\}$, $\{\phi_{2,
11'},\phi_{2,5'}\}$, $\{\phi_{2,7''},\phi_{2,13''}\}$,\\
$\{\phi_{2,11''}, \phi_{2,5''}\}$, $\{\phi_{2,9''},\phi_{2,3'}\}$,
$\{\phi_{2,11'''},\phi_{2, 5'''}\}$,
$\{\phi_{2,9'''},\phi_{2,3''}\}$, $\{\phi_{3,6},\phi_{3,10},\phi_{3,
2}\}$\item[$a_0-a_1+2b_0-b_1-b_2-c_0-c_1+2c_2=0$]\hfil\break
$\{\phi_{1,8'},\phi_{2,11'},\phi_{2,5'},\phi_{3,6},\phi_{3,10},\phi_{3,2}\}$,
 $\{\phi_{2,9'},\phi_{2,15}\}$, $\{\phi_{2,7'},\phi_{2,13'}\}$, $\{\phi_{2,
7''},\phi_{2,13''}\}$, $\{\phi_{2,11''},\phi_{2,5''}\}$,\\
$\{\phi_{2,9''}, \phi_{2,3'}\}$, $\{\phi_{2,11'''},\phi_{2,5'''}\}$,
$\{\phi_{2,9'''},\phi_{2, 3''}\}$, $\{\phi_{2,7'''},\phi_{2,1}\}$,
$\{\phi_{3,4},\phi_{3,8},\phi_{3,
12}\}$\item[$a_0-a_1+2b_0-b_1-b_2-c_0+2c_1-c_2=0$]\hfil\break
$\{\phi_{1,4'},\phi_{2,7'},\phi_{2,13'},\phi_{3,6},\phi_{3,10},\phi_{3,2}\}$,
 $\{\phi_{2,9'},\phi_{2,15}\}$, $\{\phi_{2,11'},\phi_{2,5'}\}$, $\{\phi_{2,
7''},\phi_{2,13''}\}$, $\{\phi_{2,11''},\phi_{2,5''}\}$,\\
$\{\phi_{2,9''}, \phi_{2,3'}\}$, $\{\phi_{2,11'''},\phi_{2,5'''}\}$,
$\{\phi_{2,9'''},\phi_{2, 3''}\}$, $\{\phi_{2,7'''},\phi_{2,1}\}$,
$\{\phi_{3,4},\phi_{3,8},\phi_{3,
12}\}$\item[$a_0-a_1+2b_0-b_1-b_2+2c_0-c_1-c_2=0$]\hfil\break
$\{\phi_{1,0},\phi_{2,9'},\phi_{2,15},\phi_{3,6},\phi_{3,10},\phi_{3,2}\}$,
 $\{\phi_{2,7'},\phi_{2,13'}\}$, $\{\phi_{2,11'},\phi_{2,5'}\}$, $\{\phi_{2,
7''},\phi_{2,13''}\}$, $\{\phi_{2,11''},\phi_{2,5''}\}$,\\
$\{\phi_{2,9''}, \phi_{2,3'}\}$, $\{\phi_{2,11'''},\phi_{2,5'''}\}$,
$\{\phi_{2,9'''},\phi_{2, 3''}\}$, $\{\phi_{2,7'''},\phi_{2,1}\}$,
$\{\phi_{3,4},\phi_{3,8},\phi_{3, 12}\}$
\end{description} 

\emph{\small Now, by Lemma $\ref{g7}$, 
the
generic Hecke algebra $\mathcal{H}(G_6)$ of $G_6$ is
 isomorphic to the subalgebra $\bar{H}:=<S,U>$
of the following specialization $H$ of $\mathcal{H}(G_7)$
$$\begin{array}{rccl}
    H & := & <S,T,U \,\,| &  STU=TUS=UST,\,\, T^3=1  \\
     &  &  & (S-x_0)(S-x_1)=0 \\
     &  &  & (U-z_0)(U-z_1)(U-z_2)=0>
  \end{array}$$
The algebra $H$ is the twisted symmetric algebra of the cyclic group
$C_3$ over the symmetric subalgebra $\bar{H}$ and this holds
 for all further
specializations of the parameters. If we denote by $\phi$ the
characters of $H$ and by $\psi$ the characters of $\bar{H}$, we have
$$\begin{array}{llllll}\scriptsize
    \mathrm{Ind}_{\bar{H}}^H(\psi_{1,0}) & = & \phi_{1,0}+\phi_{1,4''}+\phi_{1,8'''}   &  \mathrm{Ind}_{\bar{H}}^H(\psi_{1,4}) & = & \phi_{1,4'}+\phi_{1,8''}+\phi_{1,12''}\\
    \mathrm{Ind}_{\bar{H}}^H(\psi_{1,8}) & = & \phi_{1,8'}+\phi_{1,12'}+\phi_{1,16}   &  \mathrm{Ind}_{\bar{H}}^H(\psi_{1,6}) & = & \phi_{1,6}+\phi_{1,10''}+\phi_{1,14'''}\\
    \mathrm{Ind}_{\bar{H}}^H(\psi_{1,10}) & = & \phi_{1,10'}+\phi_{1,14''}+\phi_{1,18''}   &  \mathrm{Ind}_{\bar{H}}^H(\psi_{1,14}) & = & \phi_{1,14'}+\phi_{1,18'}+\phi_{1,22}\\
    \mathrm{Ind}_{\bar{H}}^H(\psi_{2,5''}) & = & \phi_{2,9'}+\phi_{2,13''}+\phi_{2,5'''}  &  \mathrm{Ind}_{\bar{H}}^H(\psi_{2,3''}) & = &\phi_{2,7'}+\phi_{2,11''}+\phi_{2,3''}\\
    \mathrm{Ind}_{\bar{H}}^H(\psi_{2,3'})& = & \phi_{2,11'}+\phi_{2,7'''}+\phi_{2,3'}  & \mathrm{Ind}_{\bar{H}}^H(\psi_{2,7}) & = &\phi_{2,7''}+\phi_{2,11'''}+\phi_{2,15}\\
    \mathrm{Ind}_{\bar{H}}^H(\psi_{2,1}) & = & \phi_{2,9''}+\phi_{2,5'}+\phi_{2,1}&   \mathrm{Ind}_{\bar{H}}^H(\psi_{2,5'}) & = & \phi_{2,9'''}+\phi_{2,13'}+\phi_{2,5''}\\
    \mathrm{Ind}_{\bar{H}}^H(\psi_{3,2}) & = & \phi_{3,6}+\phi_{3,10}+\phi_{3,2} &   \mathrm{Ind}_{\bar{H}}^H(\psi_{3,4}) & = & \phi_{3,4}+\phi_{3,8}+\phi_{3,12}
  \end{array}$$
 Let 
 $$\theta : \left\{ 
\begin{array}{ll} 
x_i \mapsto (\zeta_2)^i q^{a_i} &(0 \leq i <2),\\ 
z_k \mapsto (\zeta_3)^k q^{c_k} &(0 \leq k <3)
\end{array} \right. 
$$ 
be a cyclotomic specialization of $\mathcal{H}(G_6)$.
Let us consider the corresponding cyclotomic specialization
of $\mathcal{H}(G_7)$
 $$\vartheta : \left\{ 
\begin{array}{ll} 
x_i \mapsto (\zeta_2)^i q^{a_i} &(0 \leq i <2),\\ 
y_j \mapsto (\zeta_3)^j  &(0 \leq j <3),\\
z_k \mapsto (\zeta_3)^k q^{c_k} &(0 \leq k <3).
\end{array} \right. 
$$
Then $(\mathcal{H}(G_7))_\vartheta$ is the twisted symmetric algebra of the cyclic group $C_3$
over the symmetric subalgebra $(\mathcal{H}(G_6))_\theta$.
Therefore, the essential hyperplanes for $G_6$ are obtained from the essential hyperplanes
for $G_7$ by setting $b_0=b_1=b_2=0$.
If now, for example, $\theta$ is associated with no essential hyperplane for $G_6$, then
the Rouquier blocks of   $(\mathcal{H}(G_7))_\vartheta$ are:
\begin{center}
$\{\phi_{1,0},\phi_{1,4''},\phi_{1,8'''}\}$, 
$\{\phi_{1,4'},\phi_{1,8''},\phi_{1,12''}\}$,
 $\{\phi_{1,8'},\phi_{1,12'},\phi_{1,16}\}$,\\ $\{\phi_{1,6},\phi_{1,10''},\phi_{1,14'''}\}$,
 $\{\phi_{1,10'},\phi_{1,14''},\phi_{1,18''}\}$, $\{\phi_{1,14'}, \phi_{1,18'},\phi_{1,22}\}$,
$\{\phi_{2,9'},\phi_{2,13''},\phi_{2,5'''}, 
\phi_{2,7''}, \phi_{2,11'''}, \phi_{2,15}\}$,
$\{\phi_{2,7'},\phi_{2,11''},\phi_{2,3''}, \phi_{2,9'''},\phi_{2,13'},\phi_{2,5''}\}$,
$\{\phi_{2,11'},\phi_{2,7'''},\phi_{2,3'}, \phi_{2,9''},\phi_{2,5'},\phi_{2,1}\}$,
 $\{\phi_{3,6},
\phi_{3,10},\phi_{3,2}\}$,
$\{\phi_{3,4},\phi_{3,8},\phi_{3,12}\}$.\end{center}
By Clifford theory, the Rouquier blocks of $(\mathcal{H}(G_6))_\theta$, \ie the Rouquier blocks associated with no essential hyperplane for $G_6$ are:
\begin{center}
$\{\psi_{1,0}\}$, $\{\psi_{1,4}\}$, $\{\psi_{1,8}\}$, $\{\psi_{1,6}\}$, $\{\psi_{1,10}\}$, $\{\psi_{1,14}\}$,\\
$\{\psi_{2,5''},\psi_{2,7}\}$, $\{\psi_{2,3''},\psi_{2,5'}\}$,
$\{\psi_{2,3'}, \psi_{2,1}\}$, $\{\psi_{3,2}\}$, $\{\psi_{3,4}\}$.
\end{center}
In the same way, we obtain the Rouquier blocks associated with each essential hyperplane for $G_6$. Here we present only the non-trivial ones:
\begin{description}\scriptsize
\item[No essential hyperplane]\hfil\break
$\{\psi_{2,5''},\psi_{2,7}\}$, $\{\psi_{2,3''},\psi_{2,5'}\}$,
$\{\psi_{2,3'}, \psi_{2,1}\}$\item[$c_1-c_2=0$]\hfil\break
$\{\psi_{1,4},\psi_{1,8}\}$, $\{\psi_{1,10},\psi_{1,14}\}$,
$\{\psi_{2,5''}, \psi_{2,7}\}$,
$\{\psi_{2,3''},\psi_{2,3'},\psi_{2,1},\psi_{2,5'}\}$\item[$c_0-c_1=0$]\hfil\break
$\{\psi_{1,0},\psi_{1,4}\}$, $\{\psi_{1,6},\psi_{1,10}\}$,
$\{\psi_{2,5''}, \psi_{2,3''},\psi_{2,7},\psi_{2,5'}\}$,
$\{\psi_{2,3'},\psi_{2,1}\}$\item[$c_0-c_2=0$]\hfil\break
$\{\psi_{1,0},\psi_{1,8}\}$, $\{\psi_{1,6},\psi_{1,14}\}$,
$\{\psi_{2,5''}, \psi_{2,3'},\psi_{2,7},\psi_{2,1}\}$,
$\{\psi_{2,3''},\psi_{2,5'}\}$\item[$a_0-a_1-2c_0+c_1+c_2=0$]\hfil\break
$\{\psi_{1,6},\psi_{2,5''},\psi_{2,7},\psi_{3,4}\}$,
$\{\psi_{2,3''},\psi_{2, 5'}\}$,
$\{\psi_{2,3'},\psi_{2,1}\}$\item[$a_0-a_1+c_0-2c_1+c_2=0$]\hfil\break
$\{\psi_{1,10},\psi_{2,3''},\psi_{2,5'},\psi_{3,4}\}$,
$\{\psi_{2,5''}, \psi_{2,7}\}$,
$\{\psi_{2,3'},\psi_{2,1}\}$\item[$a_0-a_1+c_0+c_1-2c_2=0$]\hfil\break
$\{\psi_{1,14},\psi_{2,3'},\psi_{2,1},\psi_{3,4}\}$,
$\{\psi_{2,5''},\psi_{2, 7}\}$,
$\{\psi_{2,3''},\psi_{2,5'}\}$\item[$a_0-a_1-c_0-c_1+2c_2=0$]\hfil\break
$\{\psi_{1,8},\psi_{2,3'},\psi_{2,1},\psi_{3,2}\}$,
$\{\psi_{2,5''},\psi_{2, 7}\}$,
$\{\psi_{2,3''},\psi_{2,5'}\}$\item[$a_0-a_1-c_0+2c_1-c_2=0$]\hfil\break
$\{\psi_{1,4},\psi_{2,3''},\psi_{2,5'},\psi_{3,2}\}$,
$\{\psi_{2,5''},\psi_{2, 7}\}$,
$\{\psi_{2,3'},\psi_{2,1}\}$\item[$a_0-a_1+2c_0-c_1-c_2=0$]\hfil\break
$\{\psi_{1,0},\psi_{2,5''},\psi_{2,7},\psi_{3,2}\}$,
$\{\psi_{2,3''},\psi_{2, 5'}\}$,
$\{\psi_{2,3'},\psi_{2,1}\}$\item[$a_0-a_1-c_0+c_1=0$]\hfil\break
$\{\psi_{1,4},\psi_{1,6},\psi_{2,3'},\psi_{2,1}\}$,
$\{\psi_{2,5''},\psi_{2, 7}\}$,
$\{\psi_{2,3''},\psi_{2,5'}\}$\item[$a_0-a_1-c_1+c_2=0$]\hfil\break
$\{\psi_{1,8},\psi_{1,10},\psi_{2,5''},\psi_{2,7}\}$,
$\{\psi_{2,3''},\psi_{2, 5'}\}$,
$\{\psi_{2,3'},\psi_{2,1}\}$\item[$a_0-a_1+c_0-c_2=0$]\hfil\break
$\{\psi_{1,0},\psi_{1,14},\psi_{2,3''},\psi_{2,5'}\}$,
$\{\psi_{2,5''}, \psi_{2,7}\}$,
$\{\psi_{2,3'},\psi_{2,1}\}$\item[$a_0-a_1-c_0+c_2=0$]\hfil\break
$\{\psi_{1,8},\psi_{1,6},\psi_{2,3''},\psi_{2,5'}\}$,
$\{\psi_{2,5''},\psi_{2, 7}\}$,
$\{\psi_{2,3'},\psi_{2,1}\}$\item[$a_0-a_1+c_1-c_2=0$]\hfil\break
$\{\psi_{1,4},\psi_{1,14},\psi_{2,5''},\psi_{2,7}\}$,
$\{\psi_{2,3''},\psi_{2, 5'}\}$,
$\{\psi_{2,3'},\psi_{2,1}\}$\item[$a_0-a_1+c_0-c_1=0$]\hfil\break
$\{\psi_{1,0},\psi_{1,10},\psi_{2,3'},\psi_{2,1}\}$,
$\{\psi_{2,5''},\psi_{2, 7}\}$,
$\{\psi_{2,3''},\psi_{2,5'}\}$\item[$a_0-a_1=0$]\hfil\break
$\{\psi_{1,0},\psi_{1,6}\}$, $\{\psi_{1,4},\psi_{1,10}\}$,
$\{\psi_{1,8}, \psi_{1,14}\}$, $\{\psi_{2,5''},\psi_{2,7}\}$,
$\{\psi_{2,3''},\psi_{2,5'}\}$,
 $\{\psi_{2,3'},\psi_{2,1}\}$,
 $\{\psi_{3,2},\psi_{3,4}\}$
\end{description}}
\end{px}

Since it will take too many pages to give here the Rouquier blocks
associated with all essential hyperplanes for all exceptional irreducible complex
reflection groups, and in order to make it easier to work with them, we have stored these data in a computer file and
created two GAP functions which display them. These functions are
called \emph{AllBlocks} and \emph{DisplayAllBlocks} and they can be found
on my webpage, along with explanations for their use. Here is an example of the use of the second one
on the group $G_4$.

\begin{px}
\begin{verbatim}
gap> W:=ComplexReflectionGroup(4);
gap> DisplayAllBlocks(W);
No essential hyperplane
[["phi{1,0}"],["phi{1,4}"],["phi{1,8}"], ["phi{2,5}"],
["phi{2,3}"],["phi{2,1}"],["phi{3,2}"]]
c_1-c_2=0
[["phi{1,0}"],["phi{1,4}","phi{1,8}","phi{2,5}"],
["phi{2,3}","phi{2,1}"],["phi{3,2}"]]
c_0-c_1=0
[["phi{1,0}","phi{1,4}","phi{2,1}"],["phi{1,8}"],
["phi{2,5}","phi{2,3}"],["phi{3,2}"]]
c_0-c_2=0
[["phi{1,0}","phi{1,8}","phi{2,3}"],["phi{1,4}"],
["phi{2,5}","phi{2,1}"],["phi{3,2}"]]
2c_0-c_1-c_2=0
[["phi{1,0}","phi{2,5}","phi{3,2}"],
["phi{1,4}"],["phi{1,8}"], ["phi{2,3}"],["phi{2,1}"]]
c_0-2c_1+c_2=0
[["phi{1,0}"],["phi{1,4}","phi{2,3}","phi{3,2}"], ["phi{1,8}"],
["phi{2,5}"],["phi{2,1}"]]
c_0+c_1-2c_2=0
[["phi{1,0}"],["phi{1,4}"],["phi{1,8}","phi{2,1}","phi{3,2}"],
["phi{2,5}"],["phi{2,3}"]]
\end{verbatim}
\end{px}\

Let $W$ be any exceptional irreducible complex reflection group.
Now that we have the Rouquier blocks associated with no and each essential
hyperplane for $W$, we can determine the Rouquier blocks of any cyclotomic Hecke algebra
associated to $W$ with the use of Theorem $\ref{summary}$.
 We have also created the GAP functions \emph{RouquierBlocks} and
\emph{DisplayRouquierBlocks} (corresponding to \emph{AllBlocks} and \emph{DisplayAllBlocks}) which, given a cyclotomic specialization $\phi: u_{\mathcal{C},j} \mapsto \zeta_{e_\mathcal{C}}^j
q^{n_{\mathcal{C},j}}$, check to which essential hyperplanes the integers $n_{\mathcal{C},j}$ belong and, using the stored data, apply Theorem $\ref{summary}$ to return the Rouquier blocks of $\mathcal{H}_\phi$. 
We will give here an example of their use on $G_4$.

\begin{px}
\emph{\small The generic Hecke algebra of $G_4$ has a presantation of the form
$$\begin{array}{rccll}
    \mathcal{H}(G_4) & = & <S,T \,\,| &  STS=TST, & (S-u_0)(S-u_1)(S-u_2)=0 \\
     &  &  &  & (T-u_0)(T-u_1)(T-u_2)=0>
  \end{array}$$
If we want to calculate the Rouquier blocks of the cyclotomic Hecke
algebra
$$\begin{array}{rccll}
    \mathcal{H}_\phi & = & <S,T \,\,| &  STS=TST, & (S-1)(S-\zeta_3 q)(S-\zeta_3^2q^2)=0 \\
     &  &  &  & (T-1)(T-\zeta_3q)(T-\zeta_3^2q^2)=0>
  \end{array}$$
we use the following commands (the way to define a cyclotomic Hecke algebra in CHEVIE is explained in the GAP manual, cf., for example, \cite{jmichel}):}
\begin{verbatim}
gap> W:=ComplexReflectionGroup(4);
gap> H:=Hecke(W,[[1,E(3)*q,E(3)^2*q^2]]);
gap> DisplayRouquierBlocks(H);
[["phi{1,0}"],["phi{1,4}","phi{2,3}","phi{3,2}"],
["phi{1,8}"],["phi{2,5}"],[ "phi{2,1}"]]
\end{verbatim}
\end{px}

\section{The groups $G(d,1,r)$}

The group $G(d,1,r)$ is the group of all $r \times r$  monomial matrices with non-zero entries in $\mu_d$. It is isomorphic to the wreath product $\mu_d \wr \mathfrak{S}_r$ and its field of definition  is $\mathbb{Q}(\zeta_d)$. 

We will start by introducing some notations and results in combinatorics (cf.~\cite{BK}, \S 3A) which will be useful for the description of the Rouquier blocks of the cyclotomic Ariki-Koike algebras, \ie the cyclotomic Hecke algebras associated to the group $G(d,1,r)$.

\subsection{Combinatorics}

Let $\el=(\el_1,\el_2,\ldots,\el_h)$ be a \index{partition}\emph{partition}, \ie a finite decreasing sequence of positive integers
$$\el_1 \geq \el_2 \geq \ldots \geq \el_h \geq 1.$$
The integer
$$|\el|:=\el_1+\el_2+\ldots+\el_h$$
is called \index{size of a partition}\emph{the size of $\el$}. We also say that $\lambda$ \emph{is a partition of }
$|\el|$.
The integer $h$ is called \emph{the height of $\el$}\index{height of a partition} and we set $h_\el:=h$. To each partition $\el$ we associate its \index{$\beta$-number}\emph{$\beta$-number}, $\eb_\el=(\eb_1,\eb_2,\ldots,\eb_h)$, defined by
$$\eb_1:=h+\el_1-1,\eb_2:=h+\el_2-2,\ldots,\eb_h:=h+\el_h-h.$$
\begin{px}
\emph{\small If $\el=(4,2,2,1)$, then $\beta_\el=(7,4,3,1)$.}
\end{px}
Let $n$ be a non-negative integer. The \index{$\eb$-number shifted} \emph{$n$-shifted $\eb$-number} of $\el$ is the sequence of numbers defined by
$$\eb_\el[n]:=(\eb_1+n,\eb_2+n,\ldots,\eb_h+n,n-1,n-2,\ldots,1,0).$$
We have $\eb_\el[0]=\eb_\el$.
\begin{px}
\emph{\small If $\el=(4,2,2,1)$, then $\beta_\el[3]=(10,7,6,4,2,1,0)$.}
\end{px}

\subsubsection{Multipartitions}

Let $d$ be a positive integer and let $\el=(\el^{(0)},\el^{(1)},\ldots,\el^{(d-1)})$ be a \index{multipartition}$d$-partition, \ie a family of $d$ partitions indexed by the set $\{0,1,\ldots,d-1\}$. We set 
$$h^{(a)}:=h_{\el^{(a)}}, \,\,\, \eb^{(a)}:=\eb_{\el^{(a)}}$$
and we have
$$ \el^{(a)}=(\el_1^{(a)},\el_2^{(a)},\ldots,\el_{h^{(a)}}^{(a)}).$$
The integer
$$|\el|:=\sum_{a=0}^{d-1}|\el^{(a)}|$$
is called  \index{size of a multipartition}\emph{the size of $\el$}. We also say that $\lambda$ \emph{is a $d$-partition of}
$|\el|$.

\subsubsection{Ordinary symbols}

Let  $\el=(\el^{(0)},\el^{(1)},\ldots,\el^{(d-1)})$ be a $d$-partition.  We call \emph{$d$-height of $\el$} the family $(h^{(0)},h^{(1)},\ldots,h^{(d-1)})$ and we define the 
 \emph{height of $\el$}\index{height of a multipartition} to be the integer
$$h_\el:=\mathrm{max}\,\{h^{(a)} \,|\, 0 \leq a \leq d-1\}.$$

\begin{definition}\label{ordinary standard symbol}
The ordinary standard symbol \index{ordinary standard symbol} of $\el$ is the family of numbers defined by
$$B_\el=(B_\el^{(0)},B_\el^{(1)},\ldots,B_\el^{(d-1)}),$$
where, for all $a$ $(0 \leq a \leq d-1)$, we have
$$B_\el^{(a)}:=\eb^{(a)}[h_\el-h^{(a)}].$$
An ordinary symbol  \index{ordinary symbol} of $\el$ is a symbol obtained from the ordinary standard symbol by shifting all the rows by the same integer.
\end{definition}

The ordinary standard symbol of a $d$-partition $\el$ is of the form
$$
\begin{array}{cccccc}
 B_\el^{(0)} & = & b_1^{(0)} & b_2^{(0)} & \ldots & b_{h_\el}^{(0)} \\
 B_\el^{(1)} & = & b_1^{(1)} & b_2^{(1)} & \ldots & b_{h_\el}^{(1)} \\
  \vdots  &    & \vdots       &\vdots         &\vdots & \vdots  \\
 B_\el^{(d-1)} & = & b_1^{(d-1)} & b_2^{(d-1)} & \ldots & b_{h_\el}^{(d-1)} \\
\end{array}
$$

The \emph{ordinary content}  \index{ordinary content} of a $d$-partition of ordinary standard symbol $B_\el$ is the multiset
$$\mathrm{Cont}_\el := B_\el^{(0)} \cup B_\el^{(1)} \cup \ldots \cup B_\el^{(d-1)}$$
or (with the above notations) the polynomial defined by
$$\mathrm{Cont}_\el(x):= \sum_{0 \leq a < d}\sum_{1 \leq i \leq h_\el} x^{b_i^{(a)}}.$$

\begin{px}{\small \emph{If $d=2$ and $\el=((2,1),(3))$, then}
$$
B_\el=
\left(
\begin{array}{cc}
  3 & 1   \\
  4 & 0    
\end{array}
\right)
$$
\emph{\small We have
$\mathrm{Cont}_\el=\{0,1,3,4\}$ or $\mathrm{Cont}_\el(x)=1+x+x^3+x^4.$}}
\end{px}

\subsubsection{Charged symbols}

Let us suppose that we have a given \index{weight system}``weight system'', \ie a family of integers
$$m:=(m^{(0)},m^{(1)},\ldots,m^{(d-1)}).$$
Let  $\el=(\el^{(0)},\el^{(1)},\ldots,\el^{(d-1)})$ be a $d$-partition.  We call \emph{$(d,m)$-charged height of $\el$} the family $(hc^{(0)},hc^{(1)},\ldots,hc^{(d-1)})$, where
$$hc^{(0)}:=h^{(0)}-m^{(0)},hc^{(1)}:=h^{(1)}-m^{(1)},\ldots,hc^{(d-1)}:=h^{(d-1)}-m^{(d-1)}.$$
We define the \index{charged height of a multipartition}
 \emph{$m$-charged height of $\el$} to be the integer
$$hc_\el:=\mathrm{max}\,\{hc^{(a)} \,|\, 0 \leq a \leq d-1\}.$$

\begin{definition}\label{charged standard symbol}
The $m$-charged standard symbol\index{charged standard symbol} of $\el$ is the family of numbers defined by
$$Bc_\el=(Bc_\el^{(0)},Bc_\el^{(1)},\ldots,Bc_\el^{(d-1)}),$$
where, for all $a$ $(0 \leq a \leq d-1)$, we have
$$Bc_\el^{(a)}:=\eb^{(a)}[hc_\el-hc^{(a)}].$$
An $m$-charged symbol \index{charged symbol} of $\el$ is a symbol obtained from the $m$-charged standard symbol by shifting all the rows by the same integer.
\end{definition}
\begin{remark} \emph{The ordinary symbols correspond to the weight system }
\begin{center}
$m^{(0)}=m^{(1)}=\ldots=m^{(d-1)}=0.$
\end{center}
\end{remark}

The $m$-charged standard symbol of $\el$ is a tableau of numbers arranged into $d$ rows indexed by the set $\{0,1,\ldots,d-1\}$ such that the $a^{\mathrm{th}}$ row has length equal to $hc_\el+m^{(a)}$. For all $a$ $(0 \leq a \leq d-1)$, we set $l^{(a)}:=hc_\el+m^{(a)}$ and we denote by
$$
\begin{array}{cccccc}
 Bc_\el^{(a)} & = & bc_1^{(a)} & bc_2^{(a)} & \ldots & bc_{l^{(a)}}^{(a)} 
 \end{array}$$
the $a^{\mathrm{th}}$ row of the $m$-charged standard symbol.

The \emph{$m$-charged content} \index{charged content} of a $d$-partition of $m$-charged standard symbol $Bc_\el$ is the multiset
$$\mathrm{Contc}_\el := Bc_\el^{(0)} \cup Bc_\el^{(1)} \cup \ldots \cup Bc_\el^{(d-1)}$$
or (with the above notations) the polynomial defined by
$$\mathrm{Contc}_\el(x):= \sum_{0 \leq a < d}\sum_{1 \leq i \leq l^{(a)}}  x^{bc_i^{(a)}}.$$

\begin{px}{\small \emph{If $d=2$, $\el=((2,1),(3))$ and $m=(-1,2)$, then}
$$
Bc_\el=
\left(
\begin{array}{ccccc}
  3 & 1 &     &    & \\
  7 & 3 & 2 & 1 & 0    
\end{array}
\right)
$$
\emph{\small We have
$\mathrm{Contc}_\el=\{0,1,1,2,3,3,7\}$ or $\mathrm{Contc}_\el(x)=1+2x+x^2+2x^3+x^7.$}}
\end{px}
\subsection{Ariki-Koike algebras}

The \index{generic Ariki-Koike algebra}\emph{generic Ariki-Koike algebra} associated to $G(d,1,r)$ (cf.~\cite{ArKo}, \cite{BM}) is the algebra $\mathcal{H}_{d,r}$ generated over the Laurent polynomial ring in $d+1$ indeterminates  
$$\mathcal{L}_d:=\mathbb{Z}[u_0,u_0^{-1},u_1,u_1^{-1},\ldots,u_{d-1},u_{d-1}^{-1},x,x^{-1}]$$
by the elements $\mathrm{\textbf{s}},\mathrm{\textbf{t}}_1,\mathrm{\textbf{t}}_2,\ldots,\mathrm{\textbf{t}}_{r-1}$ satisfying the relations
\begin{itemize}
\item $\mathrm{\textbf{s}}\mathrm{\textbf{t}}_1\mathrm{\textbf{s}}\mathrm{\textbf{t}}_1=\mathrm{\textbf{t}}_1\mathrm{\textbf{s}}\mathrm{\textbf{t}}_1\mathrm{\textbf{s}}$, $\mathrm{\textbf{s}}\mathrm{\textbf{t}}_j=\mathrm{\textbf{t}}_j\mathrm{\textbf{s}} \textrm{ for } j\neq 1$,
\item $\mathrm{\textbf{t}}_j\mathrm{\textbf{t}}_{j+1}\mathrm{\textbf{t}}_j=\mathrm{\textbf{t}}_{j+1}\mathrm{\textbf{t}}_j\mathrm{\textbf{t}}_{j+1}$,  $ \mathrm{\textbf{t}}_i\mathrm{\textbf{t}}_j=\mathrm{\textbf{t}}_j\mathrm{\textbf{t}}_i \textrm{ for } |i-j|>1$,
\item $(\mathrm{\textbf{s}}-u_0)(\mathrm{\textbf{s}}-u_1)\ldots(\mathrm{\textbf{s}}-u_{d-1})=0$,
\item $(\mathrm{\textbf{t}}_j-x)(\mathrm{\textbf{t}}_j+1)=0$ for all $j=1,2,\ldots,r-1$.
\end{itemize}$ $\\
\begin{remark}\emph{ If the last relation in the above definition is replaced by
$$({\textbf{t}}_j-x)(\mathrm{\textbf{t}}_j-1)=0,$$
then we obtain a presentation of the generic Hecke algebra of $G(d,1,r)$. However, the second ``$-$'' becomes a ``$+$'', when we specialize via a cyclotomic specialization, so we might as well consider the generic Ariki-Koike algebra instead.}
\end{remark}\ $ $\\

For every $d$-partition $\el$ of $r$, we consider the free $\mathcal{L}_d$-module which has as basis the family of standard tableaux of $\el$. We can give to this module the structure of a $\mathcal{H}_{d,r}$-module (cf.~\cite{ArKo}, \cite{Ar1}, \cite{GrLe}) and hence obtain the \index{Specht module}\emph{Specht module} \textbf{Sp}$^\el$ associated to $\el$.  

Let $\mathcal{K}_d$ be the field of fractions of $\mathcal{L}_d$.
The  $\mathcal{K}_d\mathcal{H}_{d,r}$-module  $\mathcal{K}_d\textbf{Sp}^\el$, obtained by extension of scalars, is absolutely irreducible and every irreducible $\mathcal{K}_d\mathcal{H}_{d,r}$-module is isomorphic to a module of this type. Thus, $\mathcal{K}_d$ is a splitting field for $\mathcal{H}_{d,r}$. 
We denote by $\chi_\el$ the (absolutely) irreducible character of the $\mathcal{K}_d\mathcal{H}_{d,r}$-module  $\mathcal{K}_d\textbf{Sp}^\el$.

\subsection{Rouquier blocks, charged content and residues}

Let $q$ be an indeterminate and let
$$\phi : \left\{ 
\begin{array}{ll} 
u_a \mapsto \zeta_d^a q^{m_a}\,\,(0 \leq a <d),\\ 
x \mapsto q^n
\end{array} \right. 
$$
be a cyclotomic specialization of $\mathcal{H}_{d,r}$. Since the algebra $\mathcal{K}_d\mathcal{H}_{d,r}$ is split, we can deduce easily from Theorem $\ref{Semisimplicity Malle}$ and Proposition $\ref{cyclotomic Hecke is split ss}$ that the algebra $\mathbb{Q}(\zeta_d,q)(\mathcal{H}_{d,r})_\phi$ is split semisimple. Therefore, the Rouquier blocks of $(\mathcal{H}_{d,r})_\phi$ are the blocks of the algebra
$\mathcal{R}_{\mathbb{Q}(\zeta_d)}(q)(\mathcal{H}_{d,r})_\phi$, where
$$\mathcal{R}_{\mathbb{Q}(\zeta_d)}(q)=\mathbb{Z}[\zeta_d][q,q^{-1},(q^n-1)_{n \geq 1}^{-1}].$$

Theorem $3.13$ in \cite{BK} gives a description of the Rouquier blocks of $(\mathcal{H}_{d,r})_\phi$ when $n \neq 0$. However, in the proof it is supposed that $1-\zeta_d$ always belongs to a prime ideal of $\mathbb{Z}[\zeta_d]$. This is not correct, unless $d$ is a power of a prime number. Therefore, we will state here the part of the theorem that is correct and only for the case $n=1$.

\begin{theorem}\label{BK}
Let $\phi$ be a cyclotomic specialization of $\mathcal{H}_{d,r}$ such that $\phi(x) = q$. Let $\el$ and $\mu$ be two $d$-partitions of $r$.
If the irreducible characters $(\chi_\el)_\phi$ and $(\chi_\mu)_\phi$
are in the same Rouquier block of $(\mathcal{H}_{d,r})_\phi$, then $\mathrm{Contc}_\el=\mathrm{Contc}_\mu$ with respect to the weight system $m=(m_0,m_1,\ldots,m_{d-1})$. The converse holds when $d$ is a power of a prime number.
\end{theorem}
\begin{remark}\emph{ In \cite{Chlou2}, we have proved that the converse of Theorem $\ref{BK}$, and thus the description by Brou\'e and Kim, holds when $\phi$ is the ``spetsial'' cyclotomic specialization and $d$ is any positive integer.}
\end{remark}\ $ $\\

Set  $\mathcal{O}:=\mathcal{R}_{\mathbb{Q}(\zeta_d)}(q)$. Let $\mathfrak{p}$ be a prime ideal of  $\mathbb{Z}[\zeta_d]$ lying over a prime number $p$. By Proposition $\ref{Some properties of the Rouquier ring}$, the ring $\mathcal{O}$ is a Dedekind ring, whence $\mathcal{O}_{\mathfrak{p}\mathcal{O}}$ is a discrete valuation ring. Let us denote by $k_\mathfrak{p}$ the residue field of $\mathcal{O}_{\mathfrak{p}\mathcal{O}}$
and by $\pi_\mathfrak{p}$ the canonical surjection $\mathcal{O}_{\mathfrak{p}\mathcal{O}}\twoheadrightarrow k_\mathfrak{p}$. 
Following Corollary $\ref{residue field}$, the morphism $\pi_\mathfrak{p}$ induces a block bijection between
$\mathcal{O}_{\mathfrak{p}\mathcal{O}}(\mathcal{H}_{d,r})_\phi$ and 
$k_{\mathfrak{p}}(\mathcal{H}_{d,r})_\phi$.

\begin{definition}\label{diagram of a multipartition}
The diagram\index{diagram of a multipartition} of a $d$-partition $\lambda$ is the set
$$[\el]:=\{(i,j,a)\,|\,  (0 \leq a \leq d-1)(1\leq i \leq h^{(a)})(1 \leq j \leq \el_i^{(a)})\}.$$
A node\index{node} of $\el$ is any ordered triple $(i,j,a) \in [\el]$. 
The $\mathfrak{p}$-residue\index{residue} of the node $x=(i,j,a)$ with respect to  $\phi$ is
$$\mathrm{res}_{\mathfrak{p},\phi}(x):= \left\{ 
\begin{array}{ll} 
(\pi_\mathfrak{p}(j-i), \pi_\mathfrak{p}(\phi(u_a))), & \text{if } n=0 \text{ and } \pi_\mathfrak{p}(\phi(u_a)) \neq \pi_\mathfrak{p}(\phi(u_b))  \text{ for } b\neq a,\\  \\
\pi_\mathfrak{p}(\phi(u_ax^{j-i})), & \text{otherwise. }  
\end{array} \right. 
$$
\end{definition}

Let $\mathrm{Res}_{\mathfrak{p},\phi}:=\{\mathrm{res}_{\mathfrak{p},\phi}(x)\,|\,x \in [\el] \textrm{ for some $d$-partition $\el$ of $r$}\}$ be the set of all possible residues. For any $d$-partition $\el$ of $r$ and $f \in \mathrm{Res}_{\mathfrak{p},\phi}$, we set
$$C_f(\el):=| \{x \in [\el] \,|\, \mathrm{res}_{\mathfrak{p},\phi}(x)=f\}|.$$

\begin{definition}\label{residue equivalent}\index{residue equivalent}
Let $\el$ and $\mu$ be two $d$-partitions of $r$. We say that $\el$ and $\mu$ are $\mathfrak{p}$-residue equivalent  with respect to $\phi$ if $C_f(\el)=C_f(\mu)$ for all $f \in \mathrm{Res}_{\mathfrak{p},\phi}$.
\end{definition}

Then \cite{LyMa}, Theorem 2.13 implies that

\begin{theorem}\label{thm residue}
Let $\el$ and $\mu$ be two $d$-partitions of $r$.
The irreducible characters $(\chi_\el)_\phi$ and $(\chi_\mu)_\phi$ are in the same block of $\mathcal{O}_{\mathfrak{p}\mathcal{O}}
(\mathcal{H}_{d,r})_\phi$ if and only if $\el$ and $\mu$ are $\mathfrak{p}$-residue equivalent with respect to $\phi$.
\end{theorem}

\begin{corollary}\label{p1,p2}
Let $\mathfrak{p}_1$ and $\mathfrak{p}_2$ be two prime ideals of $\mathbb{Z}[\zeta_d]$ lying over the same prime number $p$. Then the blocks of $\mathcal{O}_{\mathfrak{p}_1\mathcal{O}}
(\mathcal{H}_{d,r})_\phi$ coincide with the blocks of $\mathcal{O}_{\mathfrak{p}_2\mathcal{O}}
(\mathcal{H}_{d,r})_\phi$.
\end{corollary}
\begin{apod}{Let $\mathfrak{p}$ be a prime ideal of  $\mathbb{Z}[\zeta_d]$ lying over  $p$ and let $a,b,c,d \in \mathbb{Z}$ such that $0 \leq a \leq b\leq d-1$. We have 
$\pi_\mathfrak{p}(\zeta_d^a q^c) = \pi_\mathfrak{p}(\zeta_d^b q^d)$ if and only if
$c=d$ and $\pi_\mathfrak{p}(\zeta_d^a)=\pi_\mathfrak{p}(\zeta_d^b)$.
If $\pi_\mathfrak{p}(\zeta_d^a)=\pi_\mathfrak{p}(\zeta_d^b)$, then the element
$\zeta_d^a-\zeta_d^b$ belongs to all the prime ideals lying over $p$.
Following the definition of $\mathfrak{p}$-residue, we deduce that
two $d$-partitions $\el$ and $\mu$ of $r$ are $\mathfrak{p}_1$-residue equivalent with respect to $\phi$ if and only if $\el$ and $\mu$ are $\mathfrak{p}_2$-residue equivalent with respect to $\phi$.}
\end{apod}

Theorem $\ref{thm residue}$, in combination with Proposition $\ref{semicontinuity  of Rouquier blocks}$, gives

\begin{proposition}\label{thm A}
Let $\el$ and $\mu$ be two $d$-partitions of $r$.
The irreducible characters $(\chi_\el)_\phi$ and $(\chi_\mu)_\phi$ are in the same Rouquier block of 
$(\mathcal{H}_{d,r})_\phi$
if and only if there exist a finite sequence
$\el_{(0)},\el_{(1)},\ldots,\el_{(m)}$ of $d$-partitions of $r$ and a finite
sequence $\mathfrak{p}_1,\ldots,\mathfrak{p}_m$ of $\phi$-bad prime
ideals for $G(d,1,r)$ such that
\begin{itemize}
  \item $\el_{(0)}=\el$ and $\el_{(m)}=\mu$,
  \item for all $i$ $(1\leq i \leq m)$,\,\,
         the $d$-partitions $\el_{(i-1)}$ and $\el_{(i)}$
         are $\mathfrak{p}_i$-residue equivalent with respect to $\phi$.
 \end{itemize}                   
\end{proposition}

\subsection{Essential hyperplanes}

The Schur elements of the algebra $\mathcal{K}_d\mathcal{H}_{d,r}$ have been independently  calculated by Geck, Iancu and Malle (\cite{GIM}) and by Mathas
(\cite{Mat}). Following their description by Theorem $\ref{schur elements of ArikiKoike}$, we deduce that
the essential hyperplanes for $G(d,1,r)$ are of the form
\begin{itemize}
\item $N=0$,
\item $kN+M_s-M_t=0$, where $0 \leq s<t<d$ and $-r<k<r$.
\end{itemize} 
The hyperplane $N=0$ is always essential for $G(d,1,r)$. 
Let $H$ be a hyperplane of the form $kN+M_s-M_t=0$, where $0 \leq s<t<d$ and $-r<k<r$.
The hyperplane
$H$
is essential for $G(d,1,r)$ if and only if there exists a prime ideal $\mathfrak{p}$ of $\mathbb{Z}[\zeta_d]$ such that $\zeta_d^s-\zeta_d^t \in \mathfrak{p}$. In this case, $H$ is $\mathfrak{p}$-essential for $G(d,1,r)$. In particular, if $\mathfrak{p}_1$ and $\mathfrak{p}_2$ are two prime ideals of $\mathbb{Z}[\zeta_d]$ lying over the same prime number $p$, then $H$ is $\mathfrak{p}_1$-essential if and only if it is  $\mathfrak{p}_2$-essential. 

\begin{px}\emph{\small The hyperplane $M_0=M_1$ is $(2)$-essential for $G(2,1,r)$, whereas it is not essential for $G(6,1,r)$, for all $r>0$.}
\end{px}

\subsection{Results}

Now we are going to determine the Rouquier blocks associated with no and each essential hyperplane for $G(d,1,r)$. All the results presented in this section have been first published in \cite{Chlou2}.

\begin{proposition}\label{no essential hyperplane}
The Rouquier blocks associated with no essential hyperplane for $G(d,1,r)$ are trivial.
\end{proposition}
\begin{apod}{Let $\phi$ be a cyclotomic specialization assoiciated with no essential hyperplane for $G(d,1,r)$. By Theorem  $\ref{schur elements of ArikiKoike}$, the coefficients of the Schur elements of 
$\mathcal{K}_d\mathcal{H}_{d,r}$ are units in $\mathbb{Z}[\zeta_d]$. We deduce that there are no $\phi$-bad prime ideals for $G(d,1,r)$, whence every irreducible character is a Rouquier block by itself.}
\end{apod}

\begin{proposition}\label{essential hyperplane N=0}
Let $\lambda, \mu$ be two $d$-partitions of $r$. The following two assertions are equivalent:
\begin{enumerate}[(i)]
\item The irreducible characters  $\chi_\lambda$ and  $\chi_\mu$ are in the same Rouquier block associated with the essential hyperplane $N=0$.
\item  We have $|\lambda^{(a)}|=|\mu^{(a)}|$ for all $a=0,1,\ldots,d-1$. 
\end{enumerate}
\end{proposition}
\begin{apod}{
Let $$\phi : \left\{ 
\begin{array}{ll} 
u_a \mapsto \zeta_d^a q^{m_a}\,\, (0 \leq a <d),\\ 
x \mapsto 1
\end{array} \right. 
$$
be a cyclotomic specialization associated with the essential hyperplane $N=0$. 
\begin{description}
\item[$(i) \Rightarrow (ii)$] Due to Proposition $\ref{thm A}$, it is enough to prove that
if two $d$-partitions  $\lambda, \mu$ of $r$ are $\mathfrak{p}$-residue equivalent with respect to $\phi$ for some prime ideal $\mathfrak{p}$ of $\mathbb{Z}[\zeta_d]$, then
$|\lambda^{(a)}|=|\mu^{(a)}|$ for all $a=0,1,\ldots,d-1$. Since
the integers $m_a\,(0 \leq a <d)$ do not belong to another essential hyperplane for $G(d,1,r)$, we have $\pi_\mathfrak{p}(\zeta_d^a q^{m_a}) \neq \pi_\mathfrak{p}(\zeta_d^b q^{m_b})$ for all $ 0 \leq a <b <d$. If $x=(i,j,a)$ is a node of $\el$ or $\mu$, then
$\mathrm{res}_{\mathfrak{p},\phi}(x)=(\pi_\mathfrak{p}(j-i), \pi_\mathfrak{p}(\zeta_d^a q^{m_a}))$. 
Since $\el$ and $\mu$  are $\mathfrak{p}$-residue equivalent, the number of nodes of $\el$ whose
$\mathfrak{p}$-residue's second entry is $ \pi_\mathfrak{p}(\zeta_d^a q^{m_a})$ must be equal to the number
of nodes of $\mu$ whose
$\mathfrak{p}$-residue's second entry is $ \pi_\mathfrak{p}(\zeta_d^a q^{m_a})$, for all $a=0,1,\ldots,d-1$.
We deduce that
$$\begin{array}{rcl}
|\lambda^{(a)}|=&|\{(i,j,a)\,|\,  (1\leq i \leq h_\el^{(a)})(1 \leq j \leq \el_i^{(a)})\}|&= \\
=& |\{(i,j,a)\,|\,  (1\leq i \leq h_\mu^{(a)})(1 \leq j \leq \mu_i^{(a)})\}|&=|\mu^{(a)}| 
\end{array}$$ 
for all $a=0,1,\ldots,d-1$.

\item[$(ii) \Rightarrow (i)$] Let  $a \in \{0,1,\ldots,d-1\}$.
It is enough to show that if $\el $ and $\mu$ are two distinct $d$-partitions of $r$ such that
\begin{center}
$|\el^{(a)}| =|\mu^{(a)}|$ and $\el^{(b)} = \mu^{(b)}$ for all $b \neq a$, 
\end{center}
then $(\chi_\el)_\phi$ and $(\chi_\mu)_\phi$ are in the same Rouquier block of $(\mathcal{H}_{d,r})_\phi$.
Set $l:=|\lambda^{(a)}|=|\mu^{(a)}|$. The partitions $\lambda^{(a)}$ and $\mu^{(a)}$ correspond to two distinct irreducible characters of the group $\mathfrak{S}_l$.
The cyclotomic Ariki-Koike algebra obtained from $\mathcal{H}_{1,l}$ via a cyclotomic specialization associated with the hyperplane $N=0$ is isomorphic to the group algebra $\mathbb{Z}[\mathfrak{S}_l]$.
For any finite group, it is known that  $1$ is the only block-idempotent of its  group algebra over $\mathbb{Z}$ (see also \cite{Rou}, \S3, Rem.1). Thus, all irreducible characters of  $\mathfrak{S}_l$ belong to the same Rouquier block of $\mathbb{Z}[\mathfrak{S}_l]$. 
Proposition $\ref{thm A}$ implies that there exist a finite sequence $\nu_{(0)},\nu_{(1)},\ldots,\nu_{(m)}$ of partitions of $l$ and a finite sequence $p_1, p_2, \ldots, p_m$ of prime numbers dividing the order of  $\mathfrak{S}_l$ such that
\begin{itemize}
\item $\nu_{(0)}=\lambda^{(a)}$ and $\nu_{(m)}=\mu^{(a)}$,
\item for all $i$ $(1\leq i \leq m)$,\,\,$\nu_{(i-1)}$ and $\nu_{(i)}$ are $(p_i)$-residue equivalent with respect to the cyclotomic specialization of $\mathcal{H}_{1,l}$ associated with the essential hyperplane $N=0$.
\end{itemize}
For  all $i$ $(1\leq i \leq m)$, we define $\nu_{d,i}$ to be the $d$-partition of $r$ such that
\begin{center}
$\nu_{d,i}^{(a)}:=\el_{(i)}$ and $\nu_{d,i}^{(b)} := \el^{(b)}$ for all $b \neq a$.
\end{center} 
Let $\mathfrak{p}_i$ be a prime ideal of $\mathbb{Z}[\zeta_d]$ lying over the prime number $p_i$. Then we have
 \begin{itemize}
\item $\nu_{d,0}=\lambda$ and $\nu_{d,m}=\mu$,
\item for all $i$ $(1\leq i \leq m)$,\,\,$\nu_{d,i-1}$ and $\nu_{d,i}$ are $\mathfrak{p}_i$-residue equivalent with respect to $\phi$.
\end{itemize}
By Proposition $\ref{thm A}$, the characters $(\chi_\el)_\phi$ and $(\chi_\mu)_\phi$ are in the same Rouquier block of $(\mathcal{H}_{d,r})_\phi$.}
\end{description}
\end{apod}

\begin{proposition}\label{essential hyperplane of type 1}
Let $\el, \mu$ be two $d$-partitions of $r$ and let $H$ be an essential hyperplane for $G(d,1,r)$ of the form $kN+M_s-M_t=0$, where $0 \leq s<t<d$ and $-r<k<r$.
The irreducible characters $\chi_\el$ and $\chi_\mu$ are in the same Rouquier block associated with the hyperplane $H$ if and only if the following conditions are satisfied:
\begin{enumerate}
\item We have $\el^{(a)}=\mu^{(a)}$ for all $a \notin \{s,t\}.$
\item If $\el^{st}:=(\el^{(s)},\el^{(t)})$ and $\mu^{st}:=(\mu^{(s)},\mu^{(t)})$, then
$\mathrm{Contc}_{\el^{st}}= \mathrm{Contc}_{\mu^{st}}$ with respect to the weight system $(0,k)$.
\end{enumerate}
\end{proposition}
\begin{apod}{Let $$\phi : \left\{ 
\begin{array}{ll} 
u_a \mapsto \zeta_d^a q^{m_a}\,\, (0 \leq j <d),\\ 
x \mapsto q^n
\end{array} \right. 
$$
 be a cyclotomic specialization associated with the essential hyperplane $H$. We can assume, without loss of generality, that $n=1$. We can also assume that $m_s=0$ and $m_t=k$.

Suppose that $(\chi_\el)_{\phi}$ and $(\chi_\mu)_{\phi}$ belong to  the same Rouquier block of $(\mathcal{H}_{d,r})_{\phi}$. By Theorem $\ref{BK}$, we have $\mathrm{Contc}_\el=\mathrm{Contc}_\mu$ with respect to the weight system $m=(m_0,m_1,\ldots,m_{d-1})$. Since the $m_a$, $a \notin \{s,t\}$, can take any value (as long as they do not belong to another essential hyperplane), the equality  $\mathrm{Contc}_\el=\mathrm{Contc}_\mu$ yields the first condition. Moreover, the $m$-charged standard symbols $Bc_\el$ and $Bc_\mu$ must have the same cardinality, whence $hc_\el=hc_\mu$. Therefore, we obtain
$$Bc_\el^{(a)}=\beta^{(a)}_\el[hc_\el-hc_\el^{(a)}]=\beta^{(a)}_\mu [hc_\mu-hc_\mu^{(a)}]=Bc_\mu^{(a)} \textrm{ for all } a \notin \{s,t\},$$
whence we deduce the following equality between multisets:
$$Bc_\el^{(s)} \cup Bc_\el^{(t)}=Bc_\mu^{(s)} \cup Bc_\mu^{(t)}.$$
We can assume that the $m_a$, $a \notin \{s,t\}$, are sufficiently large so that 
$$hc_\el \in \{hc_\el^{(s)},hc_\el^{(t)}\} \textrm{ and } hc_\mu \in \{hc_\mu^{(s)},hc_\mu^{(t)}\}.$$
If  $\el^{st}:=(\el^{(s)},\el^{(t)})$ and $\mu^{st}:=(\mu^{(s)},\mu^{(t)})$, then we have
$$Bc_{\el^{st}}^{(0)}=Bc_\el^{(s)},  Bc_{\el^{st}}^{(1)}=Bc_\el^{(t)}, 
Bc_{\mu^{st}}^{(0)}=Bc_\mu^{(s)}, Bc_{\mu^{st}}^{(1)}=Bc_\mu^{(t)}$$
with respect to the weight system $(0,k)$. By definition, we obtain
$\mathrm{Contc}_{\el^{st}}= \mathrm{Contc}_{\mu^{st}}$ with respect to the weight system $(0,k)$. 

Now let us suppose that the conditions $1$ and $2$ are satisfied. 
Since $H$ is an essential hyperplane for $G(d,1,r)$,  there exists a prime ideal $\mathfrak{p}$ of $\mathbb{Z}[\zeta_d]$ such that $\zeta_d^s-\zeta_d^t \in \mathfrak{p}$. We are going to show that the partitions $\el$ and $\mu$ are $\mathfrak{p}$-residue equivalent with respect to $\phi$. 
Thanks to the first condition, we only need to compare the $\mathfrak{p}$-residues of the nodes with third entry $s$ or $t$.

Set $l:=|\el^{st}|$. The first condition yields that $|\mu^{st}|=l$.   Let $\mathcal{H}_{2,l}$ be the generic Ariki-Koike algebra associated to the group $G(2,1,l)$. The algebra $\mathcal{H}_{2,l}$ is defined over the Laurent polynomial ring 
$$\mathbb{Z}[U_0,U_0^{-1},U_1,U_1^{-1},X,X^{-1}].$$ 
Let us consider the cyclotomic specialization 
$$\vartheta:U_0\mapsto 1, U_1 \mapsto -q^k,  X \mapsto q.$$
Due to Theorem $\ref{BK}$, the condition $2$ implies that the characters
 $(\chi_{\lambda^{st}})_{\vartheta}$ and  $(\chi_{\mu^{st}})_{\vartheta}$ belong to the same Rouquier block of $(\mathcal{H}_{2,l})_{\vartheta}$. We deduce that $kN+M_0-M_1=0$ is a $(2)$-essential hyperplane for $G(2,1,l)$ and that $\vartheta$ is associated with this hyperplane.
Following Proposition $\ref{thm A}$, 
 $\el^{st}$ and $\mu^{st}$ must be $(2)$-residue equivalent with respect to $\vartheta$. 
We have
\begin{itemize}
\item $(i,j,0) \in [\el^{st}]$ (resp.~$[\mu^{st}]$) if and only if $(i,j,s) \in [\el]$ (resp.~$[\mu]$).
Moreover, $\mathrm{res}_{(2),\vartheta}(i,j,0)=\pi_{(2)}(q^{j-i})$,
whereas $\mathrm{res}_{\mathfrak{p},\phi}(i,j,s)=\pi_{\mathfrak{p}}(\zeta_d^sq^{j-i})$.
\item $(i,j,1) \in [\el^{st}]$ (resp.~$[\mu^{st}]$) if and only if $(i,j,t) \in [\el]$ (resp.~$[\mu]$).
Moreover, $\mathrm{res}_{(2),\vartheta}(i,j,1)=\pi_{(2)}(-q^{k+j-i})$,
whereas $\mathrm{res}_{\mathfrak{p},\phi}(i,j,s)=\pi_{\mathfrak{p}}(\zeta_d^tq^{k+j-i})$.
\end{itemize}
Note that we have $\pi_{(2)}(1)=\pi_{(2)}(-1)$ and $\pi_{\mathfrak{p}}(\zeta_d^s)=\pi_{\mathfrak{p}}(\zeta_d^t)$. We deduce that $\el^{st}$ and $\mu^{st}$ are $(2)$-residue equivalent with respect to $\vartheta$ if and only if $\el$ and $\mu$ are $\mathfrak{p}$-residue equivalent with respect to $\phi$. }
\end{apod}

The following result is a corollary of the above proposition. However, we will show that it can also be obtained independently, with the use of the Morita equivalences established in \cite{DiMa}.

\begin{corollary}\label{second characterization}
Let $\el, \mu$ be two $d$-partitions of $r$ and let $H$ be an essential hyperplane for $G(d,1,r)$ of the form $kN+M_s-M_t=0$, where $0 \leq s<t<d$ and $-r<k<r$. Let $$\phi : \left\{ 
\begin{array}{ll} 
u_a \mapsto \zeta_d^a q^{m_a} \,\, (0 \leq a <d),\\ 
x \mapsto q^n
\end{array} \right. 
$$
 be a cyclotomic specialization associated with the essential hyperplane $H$.
The irreducible characters $(\chi_\el)_{\phi}$ and $(\chi_\mu)_{\phi}$ are in the same Rouquier block of $(\mathcal{H}_{d,r})_\phi$ if and only if the following conditions are satisfied:
\begin{enumerate}
\item We have $\el^{(a)}=\mu^{(a)}$ for all $a \notin \{s,t\}$.
\item  If $\el^{st}:=(\el^{(s)},\el^{(t)})$, 
$\mu^{st}:=(\mu^{(s)},\mu^{(t)})$ and $l:=|\el^{st}|=|\mu^{st}|$, 
then the characters $(\chi_{\el^{st}})_{\vartheta}$ and $(\chi_{\mu^{st}})_{\vartheta}$ belong to the same Rouquier block of the cyclotomic Ariki-Koike algebra of $G(2,1,l)$ obtained via the specialization  
$$\vartheta: U_0\mapsto q^{m_s}, U_1 \mapsto -q^{m_t}, X \mapsto q^n.$$
\end{enumerate}
\end{corollary}
\begin{apod}{Set $\mathcal{O}:=\mathcal{R}_{\mathbb{Q}(\zeta_d)}(q)$. Since $H$ is an essential hyperplane for $G(d,1,r)$,  there exists a prime ideal $\mathfrak{p}$ of $\mathbb{Z}[\zeta_d]$ such that $\zeta_d^s-\zeta_d^t \in \mathfrak{p}$. Due to Corollary $\ref{p1,p2}$, the Rouquier blocks of $(\mathcal{H}_{d,r})_\phi$ coincide with the blocks of $\mathcal{O}_{\mathfrak{p}\mathcal{O}}(\mathcal{H}_{d,r})_\phi$. 

From now on, all algebras are considered over the ring $\mathcal{O}_{\mathfrak{p}\mathcal{O}}$.
Following \cite{DiMa}, Theorem 1.1, we obtain that the algebra $(\mathcal{H}_{d,r})_\phi$ is Morita equivalent to the algebra
$$A:=\bigoplus_{\tiny
\begin{array}{c}
n_1,\ldots,n_{d-1} \geq 0\\
n_1+\ldots+n_{d-1}=r
\end{array}} (\mathcal{H}_{2,n_1})_{\phi'} \otimes \mathcal{H}(\mathfrak{S}_{n_2})_{\phi''}
\otimes\ldots\otimes  \mathcal{H}(\mathfrak{S}_{n_{d-1}})_{\phi''},
$$ where $\phi'$ is the restriction of $\phi$ to $\mathbb{Z}[u_s,u_s^{-1},u_t,u_t^{-1},x,x^{-1}]$
and $\phi''$ is the restriction of $\phi$ to $\mathbb{Z}[x,x^{-1}]$. Therefore, $(\mathcal{H}_{d,r})_\phi$ and $A$ have the same blocks.

Since $n \neq 0$, the blocks of $\mathcal{H}(\mathfrak{S}_{n_2})_{\phi''}$,$\ldots$, 
$\mathcal{H}(\mathfrak{S}_{n_2})_{\phi''}$ are trivial. Thus, we obtain that 
the irreducible characters $(\chi_\el)_{\phi}$ and $(\chi_\mu)_{\phi}$ are in the same (Rouquier) block of $(\mathcal{H}_{d,r})_{\phi}$ if and only if  the following conditions are satisfied:
\begin{enumerate}
\item We have $\el^{(a)}=\mu^{(a)}$ for all $a \notin \{s,t\}$.
\item If $\el^{st}:=(\el^{(s)},\el^{(t)})$, 
$\mu^{st}:=(\mu^{(s)},\mu^{(t)})$ and $l:=|\el^{st}|=|\mu^{st}|$, 
then the characters $(\chi_{\el^{st}})_{\phi'}$ and $(\chi_{\mu^{st}})_{\phi'}$ belong to the same block of
 $(\mathcal{H}_{2,l})_{\phi'}$. 
\end{enumerate} 
Theorem $\ref{thm residue}$ implies that the second condition holds if and only if  the $2$-partitions $\el^{st}$ and $\mu^{st}$ are
$\mathfrak{p}$-residue equivalent with respect to $\phi'$. Using the same argumentation as in the proof of Proposition $\ref{essential hyperplane of type 1}$, we obtain that
$\el^{st}$ and $\mu^{st}$ are
$\mathfrak{p}$-residue equivalent with respect to $\phi'$ if and only if they are
$(2)$-residue equivalent with respect to $\vartheta$,  \ie  if and only if the characters $(\chi_{\el^{st}})_{\vartheta}$ and 
$(\chi_{\mu^{st}})_{\vartheta}$ belong to the same Rouquier block of $(\mathcal{H}_{2,l})_{\vartheta}$.}
\end{apod}

\begin{px}\label{G(3,1,3)}
\emph{\small  Let $d:=3$ and $r:=3$. The irreducible characters of $G(3,1,3)$ are parametrized by the $3$-partitions of $3$. 
The generic Ariki-Koike algebra associated to $G(3,1,3)$ is the algebra $\mathcal{H}_{3,3}$ generated over the Laurent polynomial ring in $4$ indeterminates  
$$\mathbb{Z}[u_0,u_0^{-1},u_1,u_1^{-1},u_{2},u_{2}^{-1},x,x^{-1}]$$
by the elements $\mathrm{\textbf{s}},\mathrm{\textbf{t}}_1,\mathrm{\textbf{t}}_2$ satisfying the relations
\begin{itemize}
\item $\mathrm{\textbf{s}}\mathrm{\textbf{t}}_1\mathrm{\textbf{s}}\mathrm{\textbf{t}}_1=\mathrm{\textbf{t}}_1\mathrm{\textbf{s}}\mathrm{\textbf{t}}_1\mathrm{\textbf{s}}$, 
$\mathrm{\textbf{s}}\mathrm{\textbf{t}}_2=\mathrm{\textbf{t}}_2\mathrm{\textbf{s}}$,
 $\mathrm{\textbf{t}}_1\mathrm{\textbf{t}}_{2}\mathrm{\textbf{t}}_1=\mathrm{\textbf{t}}_{2}\mathrm{\textbf{t}}_1\mathrm{\textbf{t}}_{2}$,
\item $(\mathrm{\textbf{s}}-u_0)(\mathrm{\textbf{s}}-u_1)(\mathrm{\textbf{s}}-u_{2})=0$,
\item $(\mathrm{\textbf{t}}_1-x)(\mathrm{\textbf{t}}_1+1)=(\mathrm{\textbf{t}}_2-x)(\mathrm{\textbf{t}}_2+1)=0$.
\end{itemize}\
Let  $$\phi : \left\{ 
\begin{array}{ll} 
u_a \mapsto \zeta_3^a q^{m_a}\,\, (0 \leq a \leq 2),\\ 
x \mapsto q^n
\end{array} \right. 
$$ be a cyclotomic specialization of $\mathcal{H}_{3,3}$.
The essential hyperplanes for $G(3,1,3)$ are:
\begin{itemize}
\item  $N=0$.
\item  $kN+M_0-M_1=0$ for $k \in \{-2,-1,0,1,2\}$.
\item $kN+M_0-M_2=0$ for $k \in \{-2,-1,0,1,2\}$.
\item $kN+M_1-M_2=0$ for $k \in \{-2,-1,0,1,2\}$.
\end{itemize}
Let us suppose that $m_0=0$, $m_1=0$, $m_2=5$ and $n=1$.  These integers belong only to the essential hyperplane $M_0-M_1=0$.
Following Proposition $\ref{essential hyperplane of type 1}$, two irreducible characters $(\chi_\el)_{\phi}$, $(\chi_\mu)_{\phi}$ are in the same Rouquier block of $(\mathcal{H}_{3,3})_{\phi}$ if and only if
\begin{enumerate}
\item We have $\el^{(2)}=\mu^{(2)}$.
\item If $\el^{01}:=(\el^{(0)},\el^{(1)})$ and $\mu^{01}:=(\mu^{(0)},\mu^{(1)})$, then
$\mathrm{Contc}_{\el^{01}}= \mathrm{Contc}_{\mu^{01}}$ with respect to the weight system $(0,0)$, \ie
$\mathrm{Cont}_{\el^{01}}= \mathrm{Cont}_{\mu^{01}}$. 
\end{enumerate}
The first condition yields immediately that the irreducible characters corresponding to the $3$-partitions
$(\emptyset,\emptyset,(1,1,1))$, $(\emptyset,\emptyset,(2,1))$ and $(\emptyset,\emptyset,(3))$ are singletons.
Moreover, we have:}
\begin{center}
$B_{((1,1,1),\emptyset)}=\left(
\begin{array}{ccc}
  3 &2 &1     \\
  2 &1 &0 
\end{array}
\right),
$
$B_{(\emptyset,(1,1,1))}=\left(
\begin{array}{ccc}
  2 &1 &0     \\
  3 &2 &1 
\end{array}
\right),
$\\$ $\\$ $

$B_{((2,1),\emptyset)}=\left(
\begin{array}{cc}
  3 &1     \\
  1 &0 
\end{array}
\right),
$
$B_{(\emptyset,(2,1))}=\left(
\begin{array}{cc}
  1 &0     \\
  3  &1
\end{array}
\right),
$\\$ $\\$ $

$B_{((3),\emptyset)}=\left(
\begin{array}{c}
  3     \\
  0 
\end{array}
\right),
$
$B_{(\emptyset,(3))}=\left(
\begin{array}{c}
  0     \\
  3 
\end{array}
\right),
$\\$ $\\$ $

$B_{((1,1),(1))}=\left(
\begin{array}{cc}
  2  & 1   \\
  2  &0 
\end{array}
\right),
$
$B_{((1),(1,1))}=\left(
\begin{array}{cc}
  2 &0     \\
  2 &1
\end{array}
\right),
$\\$ $\\$ $

$B_{((2),(1))}=\left(
\begin{array}{c}
  2     \\
  1 
\end{array}
\right),
$
$B_{((1),(2))}=\left(
\begin{array}{c}
  1     \\
  2 
\end{array}
\right),
$\\$ $\\$ $

$B_{((1,1),\emptyset)}=\left(
\begin{array}{cc}
  2 & 1    \\
  1 & 0 
\end{array}
\right),
$
$B_{(\emptyset,(1,1))}=\left(
\begin{array}{cc}
  1 & 0     \\
  2 &1 
\end{array}
\right),
$\\$ $\\$ $

$B_{((2),\emptyset)}=\left(
\begin{array}{c}
  2     \\
  0 
\end{array}
\right),
$
$B_{(\emptyset,(2))}=\left(
\begin{array}{c}
  0     \\
  2 
\end{array}
\right),
$\\$ $\\$ $

$B_{((1),\emptyset)}=\left(
\begin{array}{c}
  1     \\
  0 
\end{array}
\right),
$
$B_{(\emptyset,(1))}=\left(
\begin{array}{c}
  0     \\
  1 
\end{array}
\right),
$\\$ $\\$ $

$B_{((1),(1))}=\left(
\begin{array}{c}
  1     \\
  1 
\end{array}
\right).
$\\$ $\\$ $

\end{center}
\emph{\small Hence, the Rouquier blocks of $(\mathcal{H}_{3,3})_\phi$ are:
\begin{enumerate}
\item $\{\chi_{((1),(1),(1))}\}$,
\item $\{\chi_{(\emptyset,\emptyset,(1,1,1))}\}$,
\item $\{\chi_{(\emptyset,\emptyset,(2,1))}\}$,
\item $\{\chi_{(\emptyset,\emptyset,(3))}\}$,
\item $\{\chi_{((1,1,1),\emptyset,\emptyset)},\chi_{(\emptyset,(1,1,1),\emptyset)}\}$,
\item $\{\chi_{((2,1),\emptyset,\emptyset)},\chi_{(\emptyset,(2,1),\emptyset)}\}$,
\item $\{\chi_{((3),\emptyset,\emptyset)},\chi_{(\emptyset,(3),\emptyset)}\}$,
\item $\{\chi_{((1,1),(1),\emptyset)},\chi_{((1),(1,1),\emptyset)}\}$,
\item $\{\chi_{((2),(1),\emptyset)},\chi_{((1),(2),\emptyset)}\}$,
\item $\{\chi_{((1,1),\emptyset,(1))},\chi_{(\emptyset,(1,1),(1))}\}$,
\item $\{\chi_{((2),\emptyset,(1))},\chi_{(\emptyset,(2),(1))}\}$,
\item $\{\chi_{((1),\emptyset,(1,1))},\chi_{(\emptyset,(1),(1,1))}\}$,
\item $\{\chi_{((1),\emptyset,(2))},\chi_{(\emptyset,(1),(2))}\}$.
\end{enumerate}
By definition, these are the Rouquier blocks associated with the $(1-\zeta_3)$-essential hyperplane $M_0-M_1=0$. }

\emph{\small If now we take $m_0=m_1=m_2=0$ and $n=1$, then the Rouquier blocks of
$(\mathcal{H}_{3,3})_\phi$ are unions of the Rouquier blocks associated with the essential hyperplanes
$M_0-M_1=0$, $M_0-M_2=0$ and $M_1-M_2=0$. Following Theorem $\ref{summary}$, the Rouquier
blocks of  $(\mathcal{H}_{3,3})_\phi$ are:
\begin{enumerate}
\item $\{\chi_{((1),(1),(1))}\}$,
\item $\{\chi_{((1,1,1),\emptyset,\emptyset)},\chi_{(\emptyset,(1,1,1),\emptyset)},\chi_{(\emptyset,\emptyset,(1,1,1))}\}$,
\item $\{\chi_{((2,1),\emptyset,\emptyset)},\chi_{(\emptyset,(2,1),\emptyset)},\chi_{(\emptyset,\emptyset,(2,1))}\}$,
\item $\{\chi_{((3),\emptyset,\emptyset)},\chi_{(\emptyset,(3),\emptyset)},\chi_{(\emptyset,\emptyset,(3))}\}$,
\item $\{\chi_{((1,1),(1),\emptyset)},\chi_{((1),(1,1),\emptyset)},\chi_{((1,1),\emptyset,(1))},\chi_{((1),\emptyset,(1,1))},\chi_{(\emptyset,(1,1),(1))},\chi_{(\emptyset,(1),(1,1))}\}$,
\item $\{\chi_{((2),(1),\emptyset)},\chi_{((1),(2),\emptyset)},\chi_{((2),\emptyset,(1))},\chi_{((1),\emptyset,(2))},\chi_{(\emptyset,(2),(1))},\chi_{(\emptyset,(1),(2))}\}$.
\end{enumerate}}
\end{px}

\section{The groups $G(2d,2,2)$}

Let $d\geq 1$. The group $G(2d,2,2)$ has $4d$ irreducible characters of degree $1$, 
$$\chi_{ijk} \,\,(0 \leq i,j\leq 1,\,\,0 \leq k < d),$$
and $d^2-d$ irreducible characters of degree $2$,
$$\chi_{kl}^{1},\, \chi_{kl}^{2}\,\,(0 \leq k \neq l < d),$$
where $\chi_{kl}^{1,2}=\chi_{lk}^{1,2}$.
The field of definition of $G(2d,2,2)$ is $\mathbb{Q}(\zeta_{2d})$.

The generic Hecke algebra of the group $G(2d,2,2)$ is the algebra $\mathcal{H}_{2d}$ generated over the Laurent  polynomial ring in $d+4$ indeterminates  
$$\mathbb{Z}[x_0,x_0^{-1},x_1,x_1^{-1},y_0,y_0^{-1},y_1,y_1^{-1},z_0,z_0^{-1},z_1,z_1^{-1},\ldots,z_{d-1},z_{d-1}^{-1}]$$
by the elements $\textbf{s},\textbf{t},\textbf{u}$ satisfying the relations
\begin{itemize}
\item $\bf stu=tus=ust$,
\item $(\textbf{s}-x_0)(\textbf{s}-x_1)=(\textbf{t}-y_0)(\textbf{t}-y_1)=(\textbf{u}-z_0)(\textbf{u}-z_1)\ldots(\textbf{u}-z_{d-1})=0$.
\end{itemize}

\subsection{Essential hyperplanes}

Let
$$\phi : \left\{ 
\begin{array}{ll} 
x_i \mapsto (-1)^iq^{a_i} & (0 \leq i <2)\\ 
y_j \mapsto (-1)^jq^{b_j} & (0 \leq j <2) \\
z_k \mapsto \zeta_d^kq^{c_k} & (0 \leq k <d)
\end{array} \right. 
$$
be a cyclotomic specialization of $\mathcal{H}_{2d}$.

The essential hyperplanes for $G(2d,2,2)$ are determined by the Schur elements of $\mathcal{H}_{2d}$.
The Schur elements of  $\mathcal{H}_{2d}$ have been calculated by Malle (\cite{Ma2}, Theorem $3.11$). Following their description (see subsection 6.7.3), the essential hyperplanes for $G(2d,2,2)$ are:
\begin{itemize}
\item $A_0-A_1=0$\,\,\,\,\,($2$-essential),
\item $B_0-B_1=0$\,\,\,\,\,($2$-essential),
\item $C_k-C_l=0$, where $0 \leq k<l<d$ and $\zeta_d^k-\zeta_d^l$ belongs to a prime ideal $\mathfrak{p}$ of 
 $\mathbb{Z}[\zeta_{2d}]$\,\,\,\,\,($\mathfrak{p}$-essential),
\item $A_i-A_{1-i}+B_j-B_{1-j}+C_k-C_l=0$, where $0 \leq i,j \leq 1$, $0 \leq k<l<d$ and $\zeta_d^k-\zeta_d^l$ belongs to a prime ideal $\mathfrak{p}$ of $\mathbb{Z}[\zeta_{2d}]$\,\,\,\,\,($\mathfrak{p}$-essential).
\end{itemize}
\begin{remark} \emph{When we say that a hyperplane is $2$-essential, we mean that it is $\mathfrak{I}$-essential for all prime ideals $\mathfrak{I}$ of $\mathbb{Z}[\zeta_{2d}]$ lying over $2$.}
\end{remark}

\subsection{Results}

In order to determine the Rouquier blocks associated with no and each essential hyperplane for $G(2d,2,2)$, we are going to use Proposition $\ref{aA}$. Following that result, if two irreducible characters
 $\chi_\phi$ and $\psi_\phi$ belong to the same Rouquier block of $(\mathcal{H}_{2d})_\phi$, then
 $$a_{\chi_\phi}+A_{\chi_\phi}=a_{\psi_\phi}+A_{\psi_\phi}.$$
Using the formulas for the Schur elements of $\mathcal{H}_{2d}$ given in the Appendix,
we can obtain the value of the sum $a_{\chi_\phi}+A_{\chi_\phi}$ for all $\chi \in \mathrm{Irr}(G(2d,2,2))$:

\begin{proposition}\label{values of aA}
Let $\chi \in \mathrm{Irr}(G(2d,2,2))$. If $\chi$ is a linear character $\chi_{ijk}$, then
$$a_{\chi_\phi}+A_{\chi_\phi}=d(a_i-a_{1-i}+b_j-b_{1-j}+2c_k)-2\sum_{l=0}^{d-1}c_l.$$
If $\chi$ is a character $\chi_{kl}^{1,2}$ of degree $2$, then
$$a_{\chi_\phi}+A_{\chi_\phi}=d(c_k+c_l)-2 \sum_{m=0}^{d-1}c_m.$$
\end{proposition}

Now we are ready to prove our main result (\cite{Chlou3}, Theorem 4.3):

\begin{theorem}\label{yes proof} For the group $G(2d,2,2)$, we have that:
\begin{enumerate}
\item The non-trivial Rouquier blocks associated with no essential hyperplane are 
$$\{\chi_{kl}^{1},\chi_{kl}^{2}\} \,\textrm{ for all } 0 \leq k<l<d.$$
\item The non-trivial Rouquier blocks associated with the $2$-essential hyperplane $A_0=A_1$ are
$$\{\chi_{0jk},\chi_{1jk}\}\, \textrm{ for all } 0 \leq j \leq 1 \textrm{ and } 0 \leq k<d,$$
$$\{\chi_{kl}^{1},\chi_{kl}^{2}\}\, \textrm{ for all } 0 \leq k<l<d.$$
\item The non-trivial Rouquier blocks associated with the $2$-essential hyperplane $B_0=B_1$ are
$$\{\chi_{i0k},\chi_{i1k}\} \,\textrm{ for all } 0 \leq i \leq 1 \textrm{ and } 0 \leq k<d,$$
$$\{\chi_{kl}^{1},\chi_{kl}^{2}\}\, \textrm{ for all } 0 \leq k<l<d.$$
\item The non-trivial Rouquier blocks associated with the $\mathfrak{p}$-essential hyperplane $C_k=C_l$ $(0 \leq k<l<d)$ are
$$ \{\chi_{ijk},\chi_{ijl}\}\, \textrm{ for all } 0 \leq i,j \leq 1,$$
$$\{\chi_{km}^{1},\chi_{km}^{2},\chi_{lm}^{1},\chi_{lm}^{2}\}\, \textrm{ for all } 0 \leq m <d \textrm{ with } m \notin \{k,l\},$$
$$\{\chi_{kl}^{1},\chi_{kl}^{2}\},$$
$$\{\chi_{rs}^{1},\chi_{rs}^{2}\}\, \textrm{ for all } 0 \leq r<s<d \textrm{ with } r,s \notin \{k,l\}.$$
\item The non-trivial Rouquier blocks associated with the $\mathfrak{p}$-essential hyperplane $A_i-A_{1-i}+B_j-B_{1-j}+C_k-C_l=0$ $(0 \leq i,j\leq 1,\,0 \leq k<l<d)$ are
$$\{\chi_{ijk},\chi_{1-i,1-j,l}, \chi_{kl}^{1},\chi_{kl}^{2} \}, $$
$$\{\chi_{rs}^{1},\chi_{rs}^{2}\}\, \textrm{ for all } 0 \leq r<s<d \textrm{ with } (r,s)  \neq (k,l).$$
\end{enumerate}
\end{theorem}
\begin{apod}{Let
$$\phi : \left\{ 
\begin{array}{lll} 
x_i \mapsto (-1)^iq^{a_i} & (0 \leq i <2)\\ 
y_j \mapsto (-1)^jq^{b_j} & (0 \leq j <2) \\
z_k \mapsto \zeta_d^kq^{c_k} & (0 \leq k <d)
\end{array} \right. 
$$
be a cyclotomic specialization of $\mathcal{H}_{2d}$.
\begin{enumerate}
\item If $\phi$ is a cyclotomic specialization associated with no essential hyperplane, then, by Proposition $\ref{Malle-Rouquier}$, each linear character is a Rouquier block by itself, whereas any character of degree $2$ is not. Due to the formulas of Proposition $\ref{values of aA}$,
Proposition $\ref{aA}$ yields that the character $\chi_{kl}^1$ $(0 \leq k<l<d)$ can  be in the same Rouquier block only with the character $\chi_{kl}^2$.
\item Suppose that $\phi$ is a cyclotomic specialization associated with the essential hyperplane
$A_0=A_1$. 
Since the
hyperplane $A_0=A_1$
 is not essential  for the characters of degree $2$,
Proposition $\ref{not essential for block}$ implies that $\{\chi_{kl}^{1},\chi_{kl}^{2}\}$ is a Rouquier block of $(\mathcal{H}_{2d})_\phi$ for all $0 \leq k<l<d$. Moreover, the hyperplane $A_0=A_1$ is $2$-essential for all characters of degree $1$ and thus, due to Proposition $\ref{Malle-Rouquier}$, there exist no linear character which is a block by itself. Due to the formulas of Proposition $\ref{values of aA}$,
Proposition $\ref{aA}$ yields that the character $\chi_{0jk}$ $(0 \leq j\leq 1,0 \leq k<d)$ can be in the same Rouquier block only with the character $\chi_{1jk}$.
\item If $\phi$ is a cyclotomic specialization associated with the essential hyperplane $B_0=B_1$, we proceed as in the previous case.
\item If $\phi$ is a cyclotomic specialization associated with the $\mathfrak{p}$-essential hyperplane
$C_k=C_l$, where $0 \leq k<l<d$, then the Rouquier blocks of $(\mathcal{H}_{2d})_\phi$ are unions of the Rouquier blocks associated with no essential hyperplane, due to Proposition $\ref{simple inclusion}$.
Hence, the characters 
 $\chi_{rs}^{1}$ and $\chi_{rs}^{2}$ are in the same Rouquier block of $(\mathcal{H}_{2d})_\phi$ for all $0 \leq r<s<d$.
Now, the hyperplane $C_k=C_l$ is $\mathfrak{p}$-essential for the following characters:
\begin{itemize}
\item $ \chi_{ijk}$, $\chi_{ijl}$, for all  $0 \leq i,j \leq 1,$
\item $\chi_{km}^{1,2}$, $\chi_{lm}^{1,2}$, for all $0 \leq m <d$ with $m \notin \{k,l\}.$
\end{itemize}
Due to the formulas of Proposition $\ref{values of aA}$,
Proposition $\ref{aA}$ yields  that 
\begin{itemize}
\item the character $\chi_{ijk}$ $(0 \leq i,j\leq 1)$ can be in the same Rouquier block only with the character $\chi_{ijl}$,
\item the character $\chi_{km}^{1}$ ($ 0 \leq m <d \textrm{ with } m \notin \{k,l\}$) can be in the same Rouquier block only with the characters $\chi_{km}^{2},\chi_{lm}^{1},\chi_{lm}^{2}$.
\end{itemize}
It remains to show that $\{\chi_{km}^{1},\chi_{km}^{2}\}$ ($ 0 \leq m <d \textrm{ with } m \notin \{k,l\}$) is not a Rouquier block of $(\mathcal{H}_{2d})_\phi$.
Following \cite{Ma2}, Table $3.10$, there exists an element $T_1$ of $\mathcal{H}_{2d}$ such that
$$\chi_{km}^1(T_1)=\chi_{km}^2(T_1)=x_0+x_1.$$
Suppose that $\{\chi_{km}^{1},\chi_{km}^{2}\}$ is a Rouquier block of $(\mathcal{H}_{2d})_\phi$ and set $y^{|\mu(\mathbb{Q}(\zeta_{2d}))|}:=q$. Then, 
 by Corollary $\ref{what we are searching}$, we must have
$$\frac{\phi(\chi_{km}^1(T_1))}{\phi(s_{\chi_{km}^1})}+\frac{\phi(\chi_{km}^2(T_1))}{\phi(s_{\chi_{km}^2})}
\in \mathcal{O},$$
where $ \mathcal{O}$ denotes the Rouquier ring of $\mathbb{Q}(\zeta_{2d})$. We have
$$\frac{\phi(\chi_{km}^1(T_1))}{\phi(s_{\chi_{km}^1})}+\frac{\phi(\chi_{km}^2(T_1))}{\phi(s_{\chi_{km}^2})}= \phi(x_0+x_1) \cdot (\frac{1}{\phi(s_{\chi_{km}^{1}})}+\frac{1}{\phi(s_{\chi_{km}^{2}})}),$$
where
$$\phi(x_0+x_1)=q^{a_0}-q^{a_1}=y^{a_0|\mu(\mathbb{Q}(\zeta_{2d}))|}-y^{a_1|\mu(\mathbb{Q}(\zeta_{2d}))|}.$$
Since $\phi$ is associated with the hyperplane $C_k=C_l$, we must have $a_0 \neq a_1$, whence
$\phi(x_0+x_1)^{-1} \in  \mathcal{O}$. We deduce that
$$\frac{1}{\phi(s_{\chi_{km}^{1}})}+\frac{1}{\phi(s_{\chi_{km}^{2}})} \in \mathcal{O}.$$
Using the formulas for the description of the Schur elements of $\chi_{km}^{1,2}$ given in the Appendix, we can easily calculate that the above element does not belong to the Rouquier ring.
\item Suppose that $\phi$ is a cyclotomic specialization associated with the $\mathfrak{p}$-essential hyperplane $A_i-A_{1-i}+B_j-B_{1-j}+C_k-C_l=0$, where $0 \leq i,j\leq 1$ and  $0 \leq k<l<d$. 
We have to distinguish two cases:
\begin{enumerate}
\item If $\mathfrak{p}$ is lying over an odd prime number, then
this hyperplane is 
$\mathfrak{p}$-essential for only three characters:
$\chi_{ijk}$, $\chi_{1-i,1-j,l}$ and either  $\chi_{kl}^{1}$ or $\chi_{kl}^{2}.$
If $\mathcal{O}$ is the Rouquier ring of $\mathbb{Q}(\zeta_{2d})$, then, by Proposition $\ref{Malle-Rouquier}$, these three characters belong to the same block of $\mathcal{O}_{\mathfrak{p}\mathcal{O}}(\mathcal{H}_{2d})_\phi$. All the remaining characters are blocks of $\mathcal{O}_{\mathfrak{p}\mathcal{O}}(\mathcal{H}_{2d})_\phi$ by themselves. Since the Rouquier blocks of $(\mathcal{H}_{2d})_\phi$ are unions of the Rouquier blocks associated with no essential hyperplane, we obtain the desired result.
\item If $\mathfrak{p}$ is lying over $2$, then the hyperplane $A_i-A_{1-i}+B_j-B_{1-j}+C_k-C_l=0$
is $\mathfrak{p}$-essential for the characters $\chi_{ijk}$, $\chi_{1-i,1-j,l}$, $\chi_{kl}^{1}$ and $\chi_{kl}^{2}.$ Using the same argumentation as in case $4$, we can show that the set
$\{\chi_{ijk}, \chi_{1-i,1-j,l}, \chi_{kl}^{1},\chi_{kl}^{2}\}$ is a Rouquier block of $(\mathcal{H}_{2d})_\phi$ (and not a union of two Rouquier blocks).
 Due to Proposition $\ref{not essential for block}$, the remaining Rouquier blocks associated with no essential hyperplane remain as they are.}
\end{enumerate}
\end{enumerate}
\end{apod}

\begin{px}\label{G(4,2,2)}\emph{\small
Let $d:=2$. The group $G(4,2,2)$ has $8$ irreducible characters of degree $1$, 
$\chi_{ijk} \,\,(0 \leq i,j,k\leq 1),$
and $2$ irreducible characters of degree $2$,
$\chi_{01}^{1,2}.$ 
The generic Hecke algebra of the group $G(4,2,2)$ is the algebra $\mathcal{H}_{4}$ generated over the Laurent  polynomial ring in $6$ indeterminates  
$$\mathbb{Z}[x_0,x_0^{-1},x_1,x_1^{-1},y_0,y_0^{-1},y_1,y_1^{-1},z_0,z_0^{-1},z_1,z_1^{-1}]$$
by the elements $\textbf{s},\textbf{t},\textbf{u}$ satisfying the relations
\begin{itemize}
\item $\bf stu=tus=ust$,
\item $(\textbf{s}-x_0)(\textbf{s}-x_1)=(\textbf{t}-y_0)(\textbf{t}-y_1)=(\textbf{u}-z_0)(\textbf{u}-z_1)=0$.
\end{itemize}
Let
$$\phi : \left\{ 
\begin{array}{ll} 
x_i \mapsto (-1)^iq^{a_i} & (0 \leq i <2)\\ 
y_j \mapsto (-1)^jq^{b_j}& (0 \leq j <2) \\
z_k \mapsto(-1)^kq^{c_k} & (0 \leq k <2)
\end{array} \right. 
$$
be a cyclotomic specialization of $\mathcal{H}_{4}$.
The essential hyperplanes for $G(4,2,2)$ are:
\begin{itemize}
\item  $H_1:\,A_0=A_1$,
\item  $H_2:\,B_0=B_1$,
\item  $H_3:\,C_0=C_1$,
\item $H_4:\,A_0-A_1+B_0-B_1+C_0-C_1=0$.
\item $H_5:\,A_0-A_1+B_1-B_0+C_0-C_1=0$.
\item $H_6:\,A_1-A_0+B_0-B_1+C_0-C_1=0$.
\item $H_7:\,A_1-A_0+B_1-B_0+C_0-C_1=0$.
\end{itemize}
The only non-trivial Rouquier block associated with no essential hyperplane is
$\{\chi_{01}^1,\chi_{01}^2\}$.
The Rouquier blocks associated with
\begin{itemize}
\item $H_1$ are: $\{\chi_{000},\chi_{100}\}$, $\{\chi_{001},\chi_{101}\}$,
$\{\chi_{010},\chi_{110}\}$, $\{\chi_{011},\chi_{111}\}$, $\{\chi_{01}^1,\chi_{01}^2\}$.
\item $H_2$ are: $\{\chi_{000},\chi_{010}\}$, $\{\chi_{001},\chi_{011}\}$,
$\{\chi_{100},\chi_{110}\}$, $\{\chi_{101},\chi_{111}\}$, $\{\chi_{01}^1,\chi_{01}^2\}$.
\item $H_3$ are: $\{\chi_{000},\chi_{001}\}$, $\{\chi_{010},\chi_{011}\}$,
$\{\chi_{100},\chi_{101}\}$, $\{\chi_{110},\chi_{111}\}$, $\{\chi_{01}^1,\chi_{01}^2\}$.
\item $H_4$ are: $\{\chi_{001}\}$, $\{\chi_{010}\}$, $\{\chi_{011}\}$,
$\{\chi_{100}\}$, $\{\chi_{101}\}$, $\{\chi_{110}\}$, $\{\chi_{000},\chi_{111},\chi_{01}^1,\chi_{01}^2\}$.
\item $H_5$ are: $\{\chi_{000}\}$, $\{\chi_{001}\}$, $\{\chi_{011}\}$,
$\{\chi_{100}\}$, $\{\chi_{110}\}$, $\{\chi_{111}\}$, $\{\chi_{010},\chi_{101},\chi_{01}^1,\chi_{01}^2\}$.
\item $H_6$ are: $\{\chi_{000}\}$, $\{\chi_{001}\}$, $\{\chi_{010}\}$,
$\{\chi_{101}\}$, $\{\chi_{110}\}$, $\{\chi_{111}\}$, $\{\chi_{100},\chi_{011},\chi_{01}^1,\chi_{01}^2\}$.
\item $H_7$ are: $\{\chi_{000}\}$, $\{\chi_{010}\}$, $\{\chi_{011}\}$,
$\{\chi_{100}\}$, $\{\chi_{101}\}$, $\{\chi_{111}\}$, $\{\chi_{110},\chi_{001},\chi_{01}^1,\chi_{01}^2\}$.
\end{itemize}
Let us take $a_0=2$, $a_1=4$, $b_0=3$, $b_1=1$ and $c_0=c_1=0$.
These integers belong to the essential hyperplanes $H_3$, $H_4$ and $H_7$.
By Theorem $\ref{summary}$, the Rouquier blocks of $(\mathcal{H}_{4})_\phi$ are
\begin{center}
$\{\chi_{000},\chi_{001},\chi_{110},\chi_{111},\chi_{01}^1,\chi_{01}^2\}$,
$\{\chi_{010},\chi_{011}\}$,
$\{\chi_{100},\chi_{101}\}$.
\end{center}}
\end{px}

\section{The groups $G(de,e,r)$}

All the results in this section have first appeared in \cite{Chlou3}.

\subsection {The groups $G(de,e,r)$, $r>2$}

We define the 
Hecke algebra of $G(de,e,r)$, $r >2$, to be the algebra $\mathcal{H}_{de,e,r}$ generated over the Laurent polynomial ring in $d+1$ indeterminates  
$$\mathbb{Z}[v_0,v_0^{-1},v_1,v_1^{-1},\ldots,v_{d-1},v_{d-1}^{-1},x,x^{-1}]$$
by the elements $a_0,a_1,\ldots,a_r$ satisfying the relations
\begin{itemize}
\item $(a_0-v_0)(a_0-v_1)\ldots(a_0-v_{d-1})=(a_j-x)(a_j+1)=0$ for $j=1,\ldots,r$,
\item $a_1a_3a_1=a_3a_1a_3$, $a_ja_{j+1}a_j=a_{j+1}a_ja_{j+1}$ for $j=2,\ldots,r-1$,
\item $a_1a_2a_3a_1a_2a_3=a_3a_1a_2a_3a_1a_2$,
\item $a_1a_j=a_ja_1$ for $j=4,\ldots,r$,
\item $a_i a_j=a_j a_i$  for $2 \leq i <j \leq r$ with $j-i>1$,
\item $a_0a_j=a_ja_0$ for $j=3,\ldots,r$,
\item $a_0a_1a_2=a_1a_2a_0$,
\item $\underbrace{a_2a_0a_1a_2a_1a_2a_1\ldots}_{e+1 \textrm{ factors}}=
\underbrace{a_0a_1a_2a_1a_2a_1a_2\ldots}_{e+1 \textrm{ factors}}$\,.
\end{itemize}
Let $$\phi : \left\{ 
\begin{array}{ll} 
v_j \mapsto \zeta_d^j q^{n_j}\,\, (0 \leq j <d),\\ 
x \mapsto q^n
\end{array} \right. 
$$
be a cyclotomic specialization of  $\mathcal{H}_{de,e,r}$. Following Theorem $\ref{summary}$, the Rouquier blocks of 
$(\mathcal{H}_{de,e,r})_\phi$ coincide with the Rouquier blocks of
$(\mathcal{H}_{de,e,r})_{\phi^e}$, where
$$\phi^e : \left\{ 
\begin{array}{ll} 
v_j \mapsto \zeta_d^j q^{en_j}\,\, (0 \leq j <d),\\ 
x \mapsto q^{en},
\end{array} \right. 
$$
since the integers $\{(n_j)_{0 \leq j <d},n\}$ and $\{(en_j)_{0 \leq j <d},en\}$ belong to the same essential hyperplanes for $G(de,e,r)$.

We now consider the generic Ariki-Koike algebra $\mathcal{H}_{de,r}$ generated over the ring
$$\mathbb{Z}[u_0,u_0^{-1},u_1,u_1^{-1},\ldots,u_{de-1},u_{de-1}^{-1},x,x^{-1}]$$
by the elements $\mathrm{\textbf{s}},\mathrm{\textbf{t}}_1,\mathrm{\textbf{t}}_2,\ldots,\mathrm{\textbf{t}}_{r-1}$ satisfying the relations described in subsection $5.3.2$. Let us consider the following cyclotomic specialization of $\mathcal{H}_{de,r}$:
$$\vartheta : \left\{ 
\begin{array}{ll} 
u_a \mapsto \zeta_{de}^{a} q^{m_a}\,\,(0 \leq a <de, m_a := n_{a \,\mathrm{mod}\,d}), \\ 
x \mapsto q^{en}.
\end{array} \right. 
$$
Following Lemma $\ref{gdeer}$, the algebra $(\mathcal{H}_{de,r})_\vartheta$ is the twisted symmetric algebra of the cyclic group $C_e$ over the symmetric subalgebra  $(\mathcal{H}_{de,e,r})_{\phi^e}$.

From now on, set $\mathcal{H}:=(\mathcal{H}_{de,r})_\vartheta$, $\bar{\mathcal{H}}:=
(\mathcal{H}_{de,e,r})_{\phi^e}$, $G:=C_e$, $K:=\mathbb{Q}(\zeta_{de})$ and let $\mathcal{R}_K(q)$ be the Rouquier ring of $K$.
Applying  Proposition $\ref{1.45}$ gives:

\begin{proposition}\label{first step}
The block-idempotents of $(Z\mathcal{R}_K(q)\bar{\mathcal{H}})^G$ coincide with the block-idempotents of $(Z\mathcal{R}_K(q)\mathcal{H})^{G^\vee}$.
\end{proposition}

The action of the cyclic group $G^\vee$ of order $e$ on $\mathrm{Irr}(K(q)\mathcal{H})$ corresponds to
the action generated by the cyclic permutation by $d$-packages on the $de$-partitions of $r$ (cf., for example, \cite{Ma4}, \S4.A):
$$\begin{array}{rl}
\tau_d: &(\el^{(0)},\ldots,\el^{(d-1)},\el^{(d)},\ldots,\el^{(2d-1)},\ldots,\el^{(de-d)},\ldots,\el^{(de-1)})\\ \mapsto &(\el^{(de-d)},\ldots,\el^{(de-1)},\el^{(0)},\ldots,\el^{(d-1)},\ldots, \el^{(de-2d)},\ldots,\el^{(de-d-1)}).
\end{array}$$
The $de$-partitions which are fixed by the action of 
$G^\vee$, \ie the $de$-partitions which are of the form
$$(\el^{(0)},\ldots,\el^{(d-1)},\el^{(0)},\ldots,\el^{(d-1)},\ldots,\el^{(0)},\ldots,\el^{(d-1)}),$$
 where the first $d$ partitions are repeated $e$ times, are called $d$-\emph{stuttering}.
 \index{d-stuttering partition}

\begin{proposition}\label{second step}
If $\el$ is a $de$-partition of $r$, then the characters $\chi_\el$ and $\chi_{\tau_d(\el)}$
 belong to the same Rouquier block of $\mathcal{H}$. In particular, the blocks of $\mathcal{R}_K(q)\mathcal{H}$ are stable under the action of $G^\vee$.
\end{proposition}
\begin{apod}{The symmetric group $\mathfrak{S}_{de}$ acts naturally on the set of $de$-partitions of $r$, and thus on
$\mathrm{Irr}(K(q)\mathcal{H})$:  If $\tau \in \mathfrak{S}_{de}$ and
$\el=(\el^{(0)},\el^{(1)},\ldots,\el^{(de-1)})$ is a  $de$-partition of $r$, then 
$\tau(\el):=(\el^{(\tau(0))},\el^{(\tau(1))},\ldots,\el^{(\tau(de-1))})$.
The action of  $G^\vee$  on $\mathrm{Irr}(K(q)\mathcal{H})$ corresponds to
the action of the cyclic subgroup  of order $e$ of  $\mathfrak{S}_{de}$ generated by the element
$$\tau_d=\prod_{j=0}^{d-1}\,\prod_{i=1}^{e-1}\sigma_{j,i}$$
where $\sigma_{j,i}$ denotes the transposition $(j,j+id)$. In order to prove that the characters
 $\chi_\el$ and $\chi_{\tau_d(\el)}$
 belong to the same Rouquier block of $\mathcal{H}$, it suffices to show that
 the characters $\chi_\el$ and $\chi_{\sigma_{j,i}(\el)}$
 are in the same Rouquier block of $\mathcal{H}$ for all  $j\,(0 \leq j < d)$ and  $i\,(0 \leq i < e)$.

Following Theorem $\ref{summary}$, the Rouquier blocks of $\mathcal{H}$ are unions of the Rouquier blocks associated with all the essential hyperplanes of the form
$$M_{s}=M_{t} \,\,(0 \leq s < t < de,\,s\equiv t\,\, \mathrm{mod}\,d).$$
Recall that the hyperplane
$M_{s}=M_{t}$ is actually essential for $G(de,1,r)$ if and only if
the element $\zeta_{de}^{s}-\zeta_{de}^{t}$ belongs to a prime ideal of
$\mathbb{Z}[\zeta_{de}]$.

Suppose that $e=p_1^{a_1}p_2^{a_2}\ldots p_m^{a_m}$, where the $p_k$ are distinct prime numbers. For $k \in \{1,2,\ldots,m\}$, we set
$c_k:= e/p_k^{a_k}.$ Then $\mathrm{gcd}(c_k)=1$ and by Bezout's theorem, there exist integers $(b_k)_{1 \leq k \leq m}$ such that $\sum_{k=1}^mb_kc_k=1$. We have $i=\sum_{k=1}^mi_k$, where $i_k:=ib_kc_k.$
The element $1-\zeta_e^{i_k}$ belongs to all the prime ideals of $\mathbb{Z}[\zeta_{de}]$ lying over the prime number $p_k$. Now set
 $$l_0:=0 \textrm{ and } l_k:=(l_{k-1}+i_k) \,\mathrm{ mod }\,e,\,\textrm{for all }k\,(1 \leq k \leq m).$$
We have that the element $\zeta_{de}^{j+l_{k-1}d}-\zeta_{de}^{j+l_{k}d}=\zeta_{de}^{j+l_{k-1}d}(1-\zeta_{e}^{i_{k}})$ belongs to all the prime ideals of $\mathbb{Z}[\zeta_{de}]$ lying over the prime number $p_k$.
Therefore, the hyperplane $M_{j+l_{k-1}d}=M_{j+l_kd}$ is essential for $G(de,1,r)$ for all $k$ $(1 \leq k \leq m)$.
Moreover, if we denote by $\tau_{j,i,k}$ the transposition $(j+l_{k-1}d,j+l_{k}d)$, then
\begin{itemize}
\item we have $\el^{(a)}=\tau_{j,i,k}(\el)^{(a)}$ for all $a \notin \{j+l_{k-1}d,j+l_{k}d\},$
\item the $2$-partitions $(\el^{(j+l_{k-1}d)},\el^{(j+l_{k}d)})$ and $(\tau_{j,i,k}(\el)^{(j+l_{k-1}d)},\tau_{j,i,k}(\el)^{(j+l_{k}d)})=(\el^{(j+l_{k}d)},\el^{(j+l_{k-1}d)})$ have the same ordinary content.
\end{itemize}
By Proposition $\ref{essential hyperplane of type 1}$, the characters
 $\chi_\el$ and $\chi_{\tau_{j,i,k}(\el)}$
 belong to the same Rouquier block associated with the essential hyperplane $M_{j+l_{k-1}d}=M_{j+l_kd}$ and thus, to the same Rouquier block of $\mathcal{H}$. We have
 \begin{center} $\sigma_{j,i}= \tau_{j,i,1} \circ \tau_{j,i,2} \circ \ldots \circ \tau_{j,i,m-1} \circ \tau_{j,i,m}\circ \tau_{j,i,m-1} \circ \ldots \circ \tau_{j,i,2} \circ \tau_{j,i,1}.$\end{center}
Consequently, the characters $\chi_\el$ and $\chi_{\sigma_{j,i}(\el)}$ belong to the same Rouquier block of $\mathcal{H}$ for all  $j\,(0 \leq j < d)$ and  $i\,(0 \leq i < e)$.}
\end{apod}

Thanks to the above result, Proposition $\ref{first step}$ now reads:

\begin{corollary}\label{third step}
The block-idempotents of $(Z\mathcal{R}_K(q)\bar{\mathcal{H}})^G$ coincide with the block-idempotents of $\mathcal{R}_K(q)\mathcal{H}$.
\end{corollary}

The following theorem demonstrates how we obtain the Rouquier blocks of $\bar{\mathcal{H}}$ from the Rouquier blocks of $\mathcal{H}$ (already determined in section 5.3).

\begin{theorem}\label{main cliff} Let $\el$ be a $de$-partition of $r$ and $\chi_\el$  the corresponding irreducible character of $G(de,1,r)$. We define $\mathrm{Irr}(K(q)\bar{\mathcal{H}})_\el$ to be the subset of $\mathrm{Irr}(K(q)\bar{\mathcal{H}})$ with the property:
$$\mathrm{Res}^{K(q)\mathcal{H}}_{K(q)\bar{\mathcal{H}}}(\chi_\el)=\sum_{\bar{\chi} \in \mathrm{Irr}(K(q)\bar{\mathcal{H}})_\el}\bar{\chi}.$$
Then
\begin{enumerate}
\item If $\el$ is $d$-stuttering and $\chi_\el$ is a block of $\mathcal{R}_K(q)\mathcal{H}$ by itself, then there are $e$ irreducible characters 
$(\bar{\chi})_{\bar{\chi} \in \mathrm{Irr}(K(q)\bar{\mathcal{H}})_\el}$. Each of these characters is a block of $\mathcal{R}_K(q)\bar{\mathcal{H}}$ by itself.
\item The other blocks of $\mathcal{R}_K(q)\mathcal{H}$ are in bijection with the remaining blocks of  $\mathcal{R}_K(q)\bar{\mathcal{H}}$ via the map of Proposition $\ref{1.45}$, i.e., the corresponding block-idempotents of $\mathcal{R}_K(q)\mathcal{H}$ coincide with the remaining block-idempotents of $\mathcal{R}_K(q)\bar{\mathcal{H}}$.
\end{enumerate}
\end{theorem}
\begin{apod}{
If $\el$ is a $d$-stuttering partition, then it is the only element in its orbit $\Omega$ under the action of $G^\vee$. Set $\bar{\Omega}:=\mathrm{Irr}(K(q)\bar{\mathcal{H}})_\el$.
By Proposition $\ref{1.42}$, we have
$|\Omega||\bar{\Omega}|=|G|=e$, whence 
$|\bar{\Omega}|=e$. Moreover, if $\bar{\chi} \in \bar{\Omega}$, then its Schur element $s_{\bar{\chi}}$ is equal to the Schur element $s_\el$ of $\chi_\el$. 
If $\chi_\el$ is a block of $\mathcal{R}_K(q)\mathcal{H}$ by itself, then, Propositions $\ref{semicontinuity of Rouquier blocks}$ and
 $\ref{Malle-Rouquier}$ imply that $s_\el=s_{\bar{\chi}}$ is invertible in $\mathcal{R}_K(q)$. Thus, $\bar{\chi}$ is a block of $\mathcal{R}_K(q)\bar{\mathcal{H}}$ by itself. 
 
 If $\el$ is $d$-stuttering and $\chi_\el$ is not a block of $\mathcal{R}_K(q)\mathcal{H}$ by itself, then, due to Theorem $\ref{summary}$, there exists a $de$-partition $\mu \neq \el$ such that $\chi_\el$ and $\chi_\mu$ belong to the same Rouquier block associated with an essential hyperplane $H$ for $G(de,1,r)$ such that the integers $\{(m_a)_{0 \leq a < de}, en\}$ belong to $H$. 
 If $H$ is $N=0,$ then, by Proposition $\ref{essential hyperplane N=0}$, we have $|\el^{(a)}|=|\mu^{(a)}|$ for all $a=0,1,\ldots,de-1$. Since 
$\el \neq \mu$, there exists $b \in \{0,1,\ldots,de-1\}$ such that $ \el^{(b)} \neq \mu^{(b)}$. If  $\nu$ is the $de$-partition of $r$ obtained from $\el$ by replacing $\el^{(b)}$ with $\mu^{(b)}$, then $\chi_\el$ and $\chi_{\nu}$ belong to the same block of $\mathcal{R}_K(q)\mathcal{H}$ and $\nu$ is not $d$-stuttering. 
If $H$ is of the form 
$kN+M_s-M_t=0, \textrm{ where } -r<k<r \textrm{ and } 0 \leq s<t<de$, then
$\el^{(a)}=\mu^{(a)}$ for all $a \notin \{s,t\}$. If $s \not\equiv t \,\mathrm{ mod }\,d$ or $e>2$, then
$\mu$ can not be $d$-stuttering. Suppose now that $s \equiv t \,\mathrm{ mod }\,d$ and $e=2$. 
 The description of $s_\el$ by Theorem 
 $\ref{schur elements of ArikiKoike}$ implies that 
 the hyperplane $M_s=M_t$ is not essential for $\chi_\el$. Due to  Proposition $\ref{not essential for block}$, we deduce that $k \neq 0$. 
Since the integers $\{(m_a)_{0 \leq a < de}, en\}$ belong to $H$ and $m_s=m_t$, we must have $n=0$.
If $\mu$ is $d$-stuttering, then $\mu^{(s)}=\mu^{(t)}$ and $|\mu^{(s)}|=|\mu^{(t)}|=|\el^{(t)}|=|\el^{(s)}|$.
Let $\nu$ be the $de$-partition obtained from $\el$ by replacing $\el^{(t)}$ with $\mu^{(t)}$. Then $\nu$ is not $d$-stuttering and the characters $\chi_\el$ and $\chi_{\nu}$ belong to the same Rouquier block associated with the essential hyperplane $N=0$. Since $n=0$, Theorem $\ref{summary}$ implies that 
 $\chi_\el$ and $\chi_{\nu}$ belong to the same block of $\mathcal{R}_K(q)\mathcal{H}$.
We will now show that the blocks of $\mathcal{R}_K(q)\mathcal{H}$ which contain at least one character corresponding to a not $d$-stuttering partition are in bijection with the remaining blocks of  $\mathcal{R}_K(q)\bar{\mathcal{H}}$ via the map of Proposition $\ref{1.45}$.

Suppose that $\el$ is not a $d$-stuttering partition and 
$b$ is the block containing  $\chi_\el$. 
Let $\bar{\chi} \in \mathrm{Irr}(K(q)\bar{\mathcal{H}})_\el$ and let $\bar{b}$ be the block of $\mathcal{R}_K(q)\bar{\mathcal{H}}$ which contains $\bar{\chi}$.
In order to establish the desired bijection, we have to show that $\bar{b}$ is stable under the action of $G$, \ie that $\bar{b}=\mathrm{Tr}(G,\bar{b}):=\sum_{g \in G/G_{\bar{b}}}g(\bar{b})$.  By Proposition $\ref{1.42}$, we have that $b=\mathrm{Tr}(G,\bar{b})$.

If $\el=(\el^{(0)},\ldots,\el^{(d-1)},\el^{(d)},\ldots,\el^{(2d-1)},\ldots,\el^{(ed-d)},\ldots,\el^{(ed-1)}),$ then, for $i=0,1,\ldots,e-1$, we define the $d$-partition $\el_{(i)}$ as follows:
$$\el_{(i)}:=(\el^{(id)},\el^{(id+1)},\ldots,\el^{(id+d-1)})$$
and we have
$$\el=(\el_{(0)},\el_{(1)},\ldots,\el_{(e-1)}).$$ 
Since $\el$ is not $d$-stuttering, there exists $m \in \{ 1,\ldots,e-1\}$ such that $\el_{(0)} \neq \el_{(m)}$. If $p$ is any prime divisor of $e$, we denote by $\el(p)$ the $de$-partition obtained from $\el$ by exchanging $\el_{(m)}$ and $\el_{(e/p)}$. Set
$$\sigma_p:=\prod_{j=0}^{d-1}\sigma_{j,m}\cdot\sigma_{j,e/p}\cdot\sigma_{j,m},$$
where $\sigma_{j,i}$ denotes the transposition $(j,j+id)$ for all $i$ $(0 \leq i <e)$.
Then $\el(p)=\sigma_p(\el)$. In the proof of Proposition $\ref{second step}$, we showed that the characters
$\chi_\el$ and $\chi_{\sigma_{j,i}(\el)}$
 are in the same Rouquier block of $\mathcal{H}$ for all  $j\,(0 \leq j < d)$ and  $i\,(0 \leq i < e)$.
Therefore, the characters $\chi_\el$ and
$\chi_{\el(p)}$ belong to the same block of $\mathcal{R}_K(q)\mathcal{H}$.
Moreover, by construction, the $de$-partition $\el(p)$ is not fixed by the generator $\tau_d^{e/p}$ of the unique subgroup of order $p$ of $G^\vee$. Thus, the order of  the stabilizer $G^\vee_{\chi_{\el(p)}}$ of
$\chi_{\el(p)}$  is prime to $p$.

By Proposition $\ref{1.42}$, we know that for each $\bar{\chi}_p \in
\mathrm{Irr}(K(q)\bar{\mathcal{H}})_{\el(p)}$, we have $|G^\vee_{\chi_{\el(p)}}||G_{\bar{\chi}_p}|=e$. Hence, $|G_{\bar{\chi}_p}|$ is divisible by the largest power of $p$ dividing $e$.  Since $b=\mathrm{Tr}(G,\bar{b})$,  the elements of $\mathrm{Irr}(K(q)\bar{\mathcal{H}})_{\el(p)}$ belong to blocks of $\mathcal{R}_K(q)\bar{\mathcal{H}}$ conjugate of $\bar{b}$ by $G$, whose stabilizer is $G_{\bar{b}}$. Following Lemma $\ref{1.43}$, we deduce that, for any prime number $p$,  $|G_{\bar{b}}|$ is divisible by the largest power of $p$ dividing $e$.
Thus, $G_{\bar{b}}=G$ and $\mathrm{Tr}(G,\bar{b})=\bar{b}$.}
\end{apod}

\begin{px}\label{G(3,3,3)}
\emph{\small Let $d:=1$, $e:=3$ and $r:=3$. The Hecke algebra of $G(3,3,3)$ is the algebra $\mathcal{H}_{3,3,3}$ generated over the Laurent polynomial ring $\mathbb{Z}[x,x^{-1}]$ by the elements $a_1,a_2,a_3$ satisfying the relations
\begin{itemize}
\item $a_1a_2a_1=a_2a_1a_2$,  $a_1a_3a_1=a_3a_1a_3$,  $a_2a_3a_2=a_3a_2a_3$,
\item $a_1a_2a_3a_1a_2a_3=a_3a_1a_2a_3a_1a_2$,
\item $(a_1-x)(a_1+1)=(a_2-x)(a_2+1)=(a_3-x)(a_3+1)=0$.
\end{itemize}
Let $\phi : x \mapsto q^n$ with $n\neq 0$ be a cyclotomic specialization of $\mathcal{H}_{3,3,3}$. We can apply Theorem $\ref{main cliff}$ and obtain the Rouquier blocks of $(\mathcal{H}_{3,3,3})_\phi$ from the
Rouquier blocks of $(\mathcal{H}_{3,3})_\vartheta$, where $\mathcal{H}_{3,3}$ is the generic Ariki-Koike algebra associated to $G(3,1,3)$ and 
$$\vartheta : \left\{ 
\begin{array}{ll} 
u_a \mapsto \zeta_{3}^{a} \,\, (0 \leq a \leq 2), \\ 
x \mapsto q^{n}.
\end{array} \right. 
$$
Since $n \neq 0$, the
Rouquier blocks of $(\mathcal{H}_{3,3})_\vartheta$ coincide with the
Rouquier blocks of $(\mathcal{H}_{3,3})_\theta$, where
$$\theta : \left\{ 
\begin{array}{ll} 
u_a \mapsto \zeta_{3}^{a}\,\,  (0 \leq a \leq 2), \\ 
x \mapsto q.
\end{array} \right. 
$$
The latter have been calculated in Example $\ref{G(3,1,3)}$ and are:
\begin{enumerate}
\item $\{\chi_{((1),(1),(1))}\}$,
\item $\{\chi_{((1,1,1),\emptyset,\emptyset)},\chi_{(\emptyset,(1,1,1),\emptyset)},\chi_{(\emptyset,\emptyset,(1,1,1))}\}$,
\item $\{\chi_{((2,1),\emptyset,\emptyset)},\chi_{(\emptyset,(2,1),\emptyset)},\chi_{(\emptyset,\emptyset,(2,1))}\}$,
\item $\{\chi_{((3),\emptyset,\emptyset)},\chi_{(\emptyset,(3),\emptyset)},\chi_{(\emptyset,\emptyset,(3))}\}$,
\item $\{\chi_{((1,1),(1),\emptyset)},\chi_{((1),(1,1),\emptyset)},\chi_{((1,1),\emptyset,(1))},\chi_{((1),\emptyset,(1,1))},\chi_{(\emptyset,(1,1),(1))},\chi_{(\emptyset,(1),(1,1))}\}$,
\item $\{\chi_{((2),(1),\emptyset)},\chi_{((1),(2),\emptyset)},\chi_{((2),\emptyset,(1))},\chi_{((1),\emptyset,(2))},\chi_{(\emptyset,(2),(1))},\chi_{(\emptyset,(1),(2))}\}$.
\end{enumerate}
Set $\mathcal{H}:=(\mathcal{H}_{3,3})_\vartheta$, $\bar{\mathcal{H}}:=
(\mathcal{H}_{3,3,3})_{\phi}$ and  $K:=\mathbb{Q}(\zeta_{3})$.
We have that
$$\mathrm{Irr}(K(q)\bar{\mathcal{H}})=\{\psi_1,\psi_2,\ldots,\psi_{10}\},$$
where
\begin{itemize}
\item $\psi_1=\mathrm{Res}^{K(q)\mathcal{H}}_{K(q)\bar{\mathcal{H}}}(\chi_{((1,1,1),\emptyset,\emptyset)})=\mathrm{Res}^{K(q)\mathcal{H}}_{K(q)\bar{\mathcal{H}}}(\chi_{(\emptyset,(1,1,1),\emptyset)})=\mathrm{Res}^{K(q)\mathcal{H}}_{K(q)\bar{\mathcal{H}}}(\chi_{(\emptyset,\emptyset,(1,1,1))})$,
\item $\psi_2=\mathrm{Res}^{K(q)\mathcal{H}}_{K(q)\bar{\mathcal{H}}}(\chi_{((2,1),\emptyset,\emptyset)})=\mathrm{Res}^{K(q)\mathcal{H}}_{K(q)\bar{\mathcal{H}}}(\chi_{(\emptyset,(2,1),\emptyset)})=\mathrm{Res}^{K(q)\mathcal{H}}_{K(q)\bar{\mathcal{H}}}(\chi_{(\emptyset,\emptyset,(2,1))})$,
\item $\psi_3=\mathrm{Res}^{K(q)\mathcal{H}}_{K(q)\bar{\mathcal{H}}}(\chi_{((3),\emptyset,\emptyset)})=\mathrm{Res}^{K(q)\mathcal{H}}_{K(q)\bar{\mathcal{H}}}(\chi_{(\emptyset,(3),\emptyset)})=\mathrm{Res}^{K(q)\mathcal{H}}_{K(q)\bar{\mathcal{H}}}(\chi_{(\emptyset,\emptyset,(3))})$,
\item $\psi_4=\mathrm{Res}^{K(q)\mathcal{H}}_{K(q)\bar{\mathcal{H}}}(\chi_{((1,1),(1),\emptyset)})=\mathrm{Res}^{K(q)\mathcal{H}}_{K(q)\bar{\mathcal{H}}}(\chi_{(\emptyset,(1,1),(1))})=\mathrm{Res}^{K(q)\mathcal{H}}_{K(q)\bar{\mathcal{H}}}(\chi_{((1),\emptyset,(1,1))})$,
\item $\psi_5=\mathrm{Res}^{K(q)\mathcal{H}}_{K(q)\bar{\mathcal{H}}}(\chi_{((1),(1,1),\emptyset)})=\mathrm{Res}^{K(q)\mathcal{H}}_{K(q)\bar{\mathcal{H}}}(\chi_{(\emptyset,(1),(1,1))})=\mathrm{Res}^{K(q)\mathcal{H}}_{K(q)\bar{\mathcal{H}}}(\chi_{((1,1),\emptyset,(1))})$,
\item $\psi_6=\mathrm{Res}^{K(q)\mathcal{H}}_{K(q)\bar{\mathcal{H}}}(\chi_{((2),(1),\emptyset)})=\mathrm{Res}^{K(q)\mathcal{H}}_{K(q)\bar{\mathcal{H}}}(\chi_{(\emptyset,(2),(1))})=\mathrm{Res}^{K(q)\mathcal{H}}_{K(q)\bar{\mathcal{H}}}(\chi_{((1),\emptyset,(2))})$,
\item $\psi_7=\mathrm{Res}^{K(q)\mathcal{H}}_{K(q)\bar{\mathcal{H}}}(\chi_{((1),(2),\emptyset)})=\mathrm{Res}^{K(q)\mathcal{H}}_{K(q)\bar{\mathcal{H}}}(\chi_{(\emptyset,(1),(2))})=\mathrm{Res}^{K(q)\mathcal{H}}_{K(q)\bar{\mathcal{H}}}(\chi_{((2),\emptyset,(1))})$,
\item $\psi_8+\psi_9+\psi_{10}=\mathrm{Res}^{K(q)\mathcal{H}}_{K(q)\bar{\mathcal{H}}}(\chi_{((1),(1),(1))})$.
\end{itemize}
Following Theorem $\ref{main cliff}$, the Rouquier blocks of $\bar{\mathcal{H}}$ are
\begin{center}
$\{\psi_1\}$, $\{\psi_2\}$, $\{\psi_3\}$, $\{\psi_4,\psi_5\}$, $\{\psi_6,\psi_7\}$, $\{\psi_8\}$, $\{\psi_9\}$, 
$\{\psi_{10}\}$.
\end{center}}
\end{px}

\subsection{The groups $G(de,e,2)$}

If the integer $e$ is odd, then everything that we said in the previous section applies to the case of $G(de,e,2)$.
Hence, we can obtain the Rouquier blocks of the cyclotomic Hecke algebras of $G(de,e,2)$ from those of $G(de,1,2)$.

If $e$ is even, then Clifford theory allows us to obtain the Rouquier blocks of the cyclotomic Hecke algebras of $G(de,e,2)$ from those of $G(de,2,2)$.

 Let $f,d \geq 1$. We denote by $\mathcal{H}_{2fd,2f,2}$ the generic
 Hecke algebra of $G(2fd,2f,2)$  generated over the Laurent polynomial ring in $d+4$ indeterminates  
$$\mathbb{Z}[x_0,x_0^{-1},x_1,x_1^{-1},y_0,y_0^{-1},y_1,y_1^{-1},u_0,u_0^{-1},u_1,u_1^{-1}\ldots,u_{d-1},u_{d-1}^{-1}],$$
by the elements $S,T,U$ satisfying the relations
\begin{itemize}
\item $STU=UST$, \,$TUS(TS)^{f-1}=U(ST)^{f}$,
\item $(S-x_0)(S-x_1)=(T-y_0)(T-y_1)=(U-u_0)(U-u_1)\ldots(U-u_{d-1})=0$. 
\end{itemize}
Let $$\phi : \left\{ 
\begin{array}{ll} 
x_i \mapsto (-1)^i q^{a_i} &(0 \leq i \leq 1),\\
y_j \mapsto (-1)^j q^{b_j}  &(0 \leq j \leq 1),\\
u_h \mapsto \zeta_{d}^{h} q^{e_h}  &(0 \leq h<d).  
\end{array} \right. 
$$
be a cyclotomic specialization of  $\mathcal{H}_{2fd,2f,2}$. Following Theorem $\ref{summary}$, the Rouquier blocks of 
$(\mathcal{H}_{2fd,2f,2})_\phi$ coincide with the Rouquier blocks of
$(\mathcal{H}_{2fd,2f,2})_{\phi^f}$, where
$$\phi^f : \left\{ 
\begin{array}{ll} 
x_i \mapsto (-1)^i q^{fa_i} &(0 \leq i \leq 1),\\
y_j \mapsto (-1)^j q^{fb_j} &(0 \leq j \leq 1),\\
u_h \mapsto \zeta_{d}^{h} q^{fe_h} & (0 \leq h <d), 
\end{array} \right. 
$$
since the integers $\{a_i,b_j,e_h\}$ and $\{fa_i,fb_j,fe_h\}$ belong to the same essential hyperplanes for $G(2fd,2f,2)$.

We now consider the generic Hecke algebra $\mathcal{H}_{2fd}$ of $G(2fd,2,2)$ generated over the ring
$$\mathbb{Z}[x_0,x_0^{-1},x_1,x_1^{-1},y_0,y_0^{-1},y_1,y_1^{-1},z_0,z_0^{-1},z_1,z_1^{-1}\ldots,z_{fd-1},z_{fd-1}^{-1}]$$
by the elements $\textbf{s},\textbf{t},\textbf{u}$ satisfying the relations
\begin{itemize}
\item $\bf stu=tus=ust$,
\item $(\textbf{s}-x_0)(\textbf{s}-x_1)=(\textbf{t}-y_0)(\textbf{t}-y_1)=(\textbf{u}-z_0)(\textbf{u}-z_1)\ldots(\textbf{u}-z_{fd-1})=0$.
\end{itemize}
Let us consider the following cyclotomic specialization of $\mathcal{H}_{2fd}$:
$$\vartheta : \left\{ 
\begin{array}{ll} 
x_i \mapsto (-1)^i q^{fa_i} &(0 \leq i \leq 1),\\
y_j \mapsto (-1)^j q^{fb_j} &(0 \leq j \leq 1),\\
z_k \mapsto \zeta_{fd}^{k} q^{c_k} &(0 \leq k <fd, c_k: = e_{k \,\mathrm{mod}\,d}).  
\end{array} \right. 
$$
Following Lemma $\ref{2fd}$, the algebra $(\mathcal{H}_{2fd})_\vartheta$ is the twisted symmetric algebra of the cyclic group $C_f$ over the symmetric subalgebra  $(\mathcal{H}_{2fd,2f,2})_{\phi^f}$.

From now on, set $\mathcal{H}:=(\mathcal{H}_{2fd})_\vartheta$, $\bar{\mathcal{H}}:=
(\mathcal{H}_{2fd,2f,2})_{\phi^f}$, $G:=C_f$, $K:=\mathbb{Q}(\zeta_{2fd})$, $y^{|\mu(K)|}:=q$ and let $\mathcal{R}_K(y)$ be the Rouquier ring of $K$.
Applying  Proposition $\ref{1.45}$ gives:

\begin{proposition}\label{1}
The block-idempotents of $(Z\mathcal{R}_K(y)\bar{\mathcal{H}})^G$ coincide with the block-idempotents of $(Z\mathcal{R}_K(y)\mathcal{H})^{G^\vee}$.
\end{proposition}

The action of the cyclic group $G^\vee$ of order $f$ on $\mathrm{Irr}(K(y)\mathcal{H})$ corresponds to the action 
$$\chi_{i,j,k} \mapsto \chi_{i,j,k+d }\,\, (0\leq i,j\leq 1, \,0 \leq k <fd),$$
$$\chi_{k,l}^{1,2} \mapsto \chi_{k+d ,l+d }^{1,2}\,\,
(0 \leq k <l<fd),$$
where all the indexes are considered  $\mathrm{mod}\,fd$. 

Let $\chi \in \mathrm{Irr}(K(y)\mathcal{H})$. If we denote by $\Omega$ the orbit of $\chi$ under the action of $G^\vee$, then $|\Omega|=f$. We define $\bar{\Omega}$ to be the subset of $\mathrm{Irr}(K(y)\bar{\mathcal{H}})$ with the property:
$$\mathrm{Res}^{K(y)\mathcal{H}}_{K(y)\bar{\mathcal{H}}}(\chi)=\sum_{\bar{\chi} \in \bar{\Omega}}\bar{\chi}.$$
By Proposition $\ref{1.42}$, we know that $|\Omega||\bar{\Omega}|=f$,
whence $|\bar{\Omega}|=1$. Since $\bar{\Omega}$ is also the orbit of $\bar{\chi}$ under the action of $G$, we deduce that the block-idempotents of $\mathcal{R}_K(y)\bar{\mathcal{H}}$ are fixed by the action of $G$.

With the help of the following lemma, we will show that the Rouquier blocks of 
$\mathcal{H}$ are also stable under the action of $G^\vee$. Here the results of Theorem $\ref{yes proof}$ are going to be used as definitions.

\begin{lemma}\label{three hyperplanes}
Let $k_1$, $k_2$ and $k_3$ be three distinct elements
of  $\{0,1,\ldots,fd-1\}$. If  the blocks of $\mathcal{R}_K(y)\mathcal{H}$ are unions of the Rouquier blocks associated with the (not necessarily essential) hyperplanes $C_{k_1}=C_{k_2}$ and $C_{k_2}=C_{k_3}$, then they are also unions of the Rouquier blocks associated with the (not necessarily essential) hyperplane
 $C_{k_1}=C_{k_3}$.
\end{lemma}
\begin{apod}{We only need to show that
\begin{enumerate}[(a)]
\item the characters $\chi_{i,j,k_1}$ and $\chi_{i,j,k_3}$ are in the same block of $\mathcal{R}_K(y)\mathcal{H}$ for all $0 \leq i,j \leq 1,$
\item the characters $\chi_{k_1,m}^{1,2}$ and $\chi_{k_3,m}^{1,2}$
are in the same block of $\mathcal{R}_K(y)\mathcal{H}$ for all
 $0 \leq m <fd$ with  $m \notin \{k_1,k_3\}.$
\end{enumerate}
Since the blocks of $\mathcal{R}_K(y)\mathcal{H}$ are unions of the Rouquier blocks associated with the  hyperplanes $C_{k_1}=C_{k_2}$ and $C_{k_2}=C_{k_3}$, Theorem $\ref{yes proof}$ implies that
\begin{enumerate}[(1)]
\item the characters $\chi_{i,j,k_1}$ and $\chi_{i,j,k_2}$ are in the same block of $\mathcal{R}_K(y)\mathcal{H}$ for all $0 \leq i,j \leq 1,$
\item the characters $\chi_{i,j,k_2}$ and $\chi_{i,j,k_3}$ are in the same block of $\mathcal{R}_K(y)\mathcal{H}$ for all $0 \leq i,j \leq 1,$
\item the characters $\chi_{k_1,m}^{1,2}$ and $\chi_{k_2,m}^{1,2}$
are in the same block of $\mathcal{R}_K(y)\mathcal{H}$ for all
 $0 \leq m <fd$ with  $m \notin \{k_1,k_2\},$
\item the characters $\chi_{k_2,m}^{1,2}$ and $\chi_{k_3,m}^{1,2}$
are in the same block of $\mathcal{R}_K(y)\mathcal{H}$ for all
 $0 \leq m <fd$ with  $m \notin \{k_2,k_3\}.$
\end{enumerate}
We immediately deduce (a) for all $0 \leq i,j \leq 1$ and (b) for all $0 \leq m <fd$ with $m \notin \{k_1,k_2,k_3\}.$
Finally, (3) implies that the characters $\chi_{k_1,k_3}^{1,2}$ and $\chi_{k_2,k_3}^{1,2}$
are in the same block of $\mathcal{R}_K(y)\mathcal{H}$, whereas by (4), $\chi_{k_1,k_2}^{1,2}$ and $\chi_{k_1,k_3}^{1,2}$ are also in the same block of $\mathcal{R}_K(y)\mathcal{H}$.
Thus, the characters $\chi_{k_1,k_2}^{1,2}$ and $\chi_{k_2,k_3}^{1,2}$ belong to the same Rouquier block of $\mathcal{H}$.}
\end{apod}

\begin{theorem}\label{2}
The Rouquier  blocks of 
$\mathcal{H}$ are stable under the action of $G^\vee$.
In particular, the block-idempotents of $\mathcal{R}_K(y)\bar{\mathcal{H}}$ coincide with the block-idempotents of $\mathcal{R}_K(y)\mathcal{H}$.
\end{theorem}
\begin{apod}{Following Theorem $\ref{summary}$, the Rouquier blocks of $\mathcal{H}$ are unions of the Rouquier blocks associated with all the essential hyperplanes of the form
$$C_{h+md}=C_{h+nd} \,\,(0 \leq h < d,\,0 \leq m < n < f).$$
Recall that the hyperplane
$C_{h+md}=C_{h+nd} $ is actually essential for $G(2fd,2,2)$ if and only if
the element $\zeta_{fd}^{h+md}-\zeta_{fd}^{h+nd}$ belongs to a prime ideal of
$\mathbb{Z}[\zeta_{2fd}]$.

Suppose that $f=p_1^{t_1}p_2^{t_2}\ldots p_r^{t_r}$, where the $p_i$ are distinct prime numbers. For $s \in \{1,2,\ldots,r\}$, we set
$h_s:= f/p_s^{t_s}.$ Then $\mathrm{gcd}(h_s)=1$ and by Bezout's theorem, there exist integers $(g_s)_{1 \leq s \leq r}$ such that $\sum_{s=1}^rg_sh_s=1$.
The element $1-\zeta_{f}^{g_sh_s}$ belongs to all the prime ideals of $\mathbb{Z}[\zeta_{2fd}]$ lying over the prime number $p_s$.  Let $h \in \{0,1,\ldots,d-1\}$ and $m \in \{0,1,\ldots,f-2\}$ and set
 $$l_0:=m \textrm{ and } l_s:=(l_{s-1}+g_sh_s) \,\,\mathrm{ mod }\,f, \textrm{ for all } s\,(1 \leq s \leq r).$$
We have that the element $\zeta_{fd}^{h+l_{s-1}d}-\zeta_{fd}^{h+l_{s}d}=\zeta_{fd}^{h+l_{s-1}d}(1-\zeta_{f}^{g_{s}h_s})$ belongs to all the prime ideals of $\mathbb{Z}[\zeta_{2fd}]$ lying over the prime number $p_s$.
Therefore, the hyperplane $C_{h+l_{s-1}d}=C_{h+l_sd}$ is essential for $G(2fd,2,2)$ for all $s$ $(1 \leq s \leq r)$.
Since $l_0=m$ and $l_r=m+1$, Lemma $\ref{three hyperplanes}$ implies that the Rouquier blocks of $\mathcal{H}$ are unions of the Rouquier blocks associated with the (not necessarily essential) hyperplane
$$C_{h+md}=C_{h+(m+1)d},$$
following their description by Theorem $\ref{yes proof}$. Since this holds for all
$m\,(0 \leq m \leq f-2)$, Lemma $\ref{three hyperplanes}$ again implies that 
the Rouquier blocks of $\mathcal{H}$ are unions of the Rouquier blocks associated with  all the hyperplanes of the form 
$$C_{h+md}=C_{h+nd} \,\,(0 \leq m < n < f),$$
for all $h\,(0 \leq h < d)$. Consequently, we obtain that
\begin{itemize}
\item the characters $(\chi_{i,j,h+md})_{ 0 \leq m < f}\,$ are in the same block of  $\mathcal{R}_K(y)\mathcal{H}$ for all $0\leq i,j\leq 1$ and $0\leq h <d,$
\item the characters $(\chi_{h+md,h+nd}^{1,2})_{ 0 \leq m<n< f}\,$ are in the same block of  $\mathcal{R}_K(y)\mathcal{H}$ for all  $0\leq h <d,$
\item the characters $(\chi_{h+md,h'+nd}^{1,2})_{ 0 \leq m,n< f}\,$ are in the same block of  $\mathcal{R}_K(y)\mathcal{H}$ for all  $0\leq h<h' <d$.
\end{itemize}
Thus, the blocks of $\mathcal{R}_K(y)\mathcal{H}$ are stable under the action of $G^\vee$.
Now, Proposition $\ref{1}$ implies that the block-idempotents of $\mathcal{R}_K(y)\bar{\mathcal{H}}$ coincide with the block-idempotents of $\mathcal{R}_K(y)\mathcal{H}$.
} 
\end{apod}

 Thanks to the above result, in order to determine the Rouquier blocks of $\bar{\mathcal{H}}$, it suffices to calculate the Rouquier blocks of $\mathcal{H}$: If $C$ is a Rouquier block of $\mathcal{H}$, then
 $\{\mathrm{Res}^{K(y)\mathcal{H}}_{K(y)\bar{\mathcal{H}}}(\chi)\,|\,\chi \in C\}$ is a Rouquier block of $\bar{\mathcal{H}}$.
 
 \begin{px}\label{G(4,4,2)} \emph{\small Let $f:=2$ and $d:=1$. The group $G(4,4,2)$ is isomorphic to the group $G(2,1,2)$. The generic Hecke algebra $\mathcal{H}_{4,2,2}$ of $G(4,4,2)$ is generated over the Laurent polynomial ring in $4$ indeterminates $$\mathbb{Z}[x_0,x_0^{-1},x_1,x_1^{-1},y_0,y_0^{-1},y_1,y_1^{-1}],$$
by the elements $S$ and $T$ satisfying the relations
\begin{itemize}
\item $(S-x_0)(S-x_1)=(T-y_0)(T-y_1)=0$.
\item $STST=TSTS.$ 
\end{itemize}
Let $$\phi : \left\{ 
\begin{array}{ll} 
x_i \mapsto (-1)^i q^{a_i} &(0 \leq i \leq 1),\\
y_j \mapsto (-1)^j q^{b_j} &(0 \leq j \leq 1) 
\end{array} \right. 
$$
be a cyclotomic specialization of  $\mathcal{H}_{4,2,2}$. Since $G(4,4,2) \simeq G(2,1,2)$, we can use the results on the Ariki-Koike algebras in order to determine the Rouquier blocks of  $(\mathcal{H}_{4,2,2})_\phi$. However, here we will demonstrate how we can apply Theorem $\ref{2}$ and obtain the Rouquier blocks of $(\mathcal{H}_{4,2,2})_\phi$ from the
Rouquier blocks of $(\mathcal{H}_{4})_\vartheta$, where $\mathcal{H}_{4}$ is the generic Hecke algebra associated to $G(4,2,2)$ and 
$$\vartheta : \left\{ 
\begin{array}{lll} 
x_i \mapsto (-1)^iq^{a_i} & (0 \leq i \leq 1),\\ 
y_j \mapsto (-1)^jq^{b_j} & (0 \leq j \leq 1), \\
z_k \mapsto(-1)^k & (0 \leq k \leq 1).
\end{array} \right. 
$$ 
Set $\mathcal{H}:=(\mathcal{H}_{4})_\vartheta$, $\bar{\mathcal{H}}:=
(\mathcal{H}_{4,4,2})_{\phi}$, $K:=\mathbb{Q}(i)$ and $y^{|\mu(K)|}:=q$.
We have that
$$\mathrm{Irr}(K(y)\bar{\mathcal{H}})=\{\chi_{((2),\emptyset)},\chi_{(\emptyset,(2))},\chi_{((1,1),\emptyset)},\chi_{(\emptyset,(1,1))},\chi_{((1),(1))}\},$$
where
\begin{itemize}
\item $\chi_{((2),\emptyset)}=\mathrm{Res}^{K(y)\mathcal{H}}_{K(y)\bar{\mathcal{H}}}(\chi_{000})=\mathrm{Res}^{K(y)\mathcal{H}}_{K(y)\bar{\mathcal{H}}}(\chi_{001})$,
\item $\chi_{(\emptyset,(2))}=\mathrm{Res}^{K(y)\mathcal{H}}_{K(y)\bar{\mathcal{H}}}(\chi_{010})=\mathrm{Res}^{K(y)\mathcal{H}}_{K(y)\bar{\mathcal{H}}}(\chi_{011})$,
\item $\chi_{((1,1),\emptyset)}=\mathrm{Res}^{K(y)\mathcal{H}}_{K(y)\bar{\mathcal{H}}}(\chi_{100})=\mathrm{Res}^{K(y)\mathcal{H}}_{K(y)\bar{\mathcal{H}}}(\chi_{101})$,
\item $\chi_{(\emptyset,(1,1))}=\mathrm{Res}^{K(y)\mathcal{H}}_{K(y)\bar{\mathcal{H}}}(\chi_{110})=\mathrm{Res}^{K(y)\mathcal{H}}_{K(y)\bar{\mathcal{H}}}(\chi_{111})$,
\item $\chi_{((1),(1))}=\mathrm{Res}^{K(y)\mathcal{H}}_{K(y)\bar{\mathcal{H}}}(\chi_{01}^1)=\mathrm{Res}^{K(y)\mathcal{H}}_{K(y)\bar{\mathcal{H}}}(\chi_{01}^2)$.
\end{itemize}
Following Theorem $\ref{2}$, the Rouquier blocks  associated with no essential hyperplane for $G(4,4,2)$ are trivial (as expected).
For $a_0=2$, $a_1=4$, $b_0=3$ and $b_1=1$, the Rouquier blocks of $\mathcal{H}$ have been calculated in Example $\ref{G(4,2,2)}$ and are:
\begin{center}
$\{\chi_{000},\chi_{001},\chi_{110},\chi_{111},\chi_{01}^1,\chi_{01}^2\}$,
$\{\chi_{010},\chi_{011}\}$,
$\{\chi_{100},\chi_{101}\}$.
\end{center}
Thanks to Theorem $\ref{2}$, we deduce that the Rouquier blocks of $\bar{\mathcal{H}}$ are:
\begin{center}
$\{\chi_{((2),\emptyset)}, \chi_{(\emptyset,(1,1))}, \chi_{((1),(1))}\}$,
$\{\chi_{(\emptyset,(2))}\}$,
$\{\chi_{((1,1),\emptyset)}\}$.
\end{center}
We can verify the above result with the use of Proposition $\ref{essential hyperplane of type 1}$, which yields that two irreducible characters $(\chi_\el)_\phi$ and $(\chi_\mu)_\phi$ are in the same Rouquier block of
$\bar{\mathcal{H}}$ if and only if $\mathrm{Contc}_\el=\mathrm{Contc}_\mu$ with respect to the weight system $(0,1)$.}
 \end{px}

\chapter{Appendix: Clifford theory and Schur elements for the Hecke algebras of complex reflection groups}

Let $W$ be a complex reflection group and let us denote by
$\mathcal{H}(W)$ its generic Hecke algebra. Suppose that the
assumptions $\ref{ypo}$ are satisfied. Let $W'$ be another complex
reflection group such that, for some specialization of the
parameters, $\mathcal{H}(W)$ becomes the twisted symmetric algebra of a
finite cyclic group $G$ over the symmetric subalgebra
$\mathcal{H}(W')$. Then, if we know the Schur elements and the blocks of $\mathcal{H}(W)$,
we can use Propositions $\ref{1.42}$ and $\ref{1.45}$ in order to
calculate the Schur elements and the blocks of $\mathcal{H}(W')$.

In particular, in all the cases of exceptional irreducible complex reflection groups that will be studied below, if we denote
by $\chi'$ the (irreducible) restriction to $\mathcal{H}(W')$ of an
irreducible character $\chi \in \mathrm{Irr}(\mathcal{H}(W))$, then
the corresponding Schur elements verify
$$s_\chi = |W:W'| s_{\chi'}.$$ 

Throughout the Appendix, we
denote by $\Phi_n$ the $n^{\mathrm{th}}$ $\mathbb{Q}$-cyclotomic polynomial, \ie the minimal polynomial  of $\zeta_n$ over $\mathbb{Q}$. The notations for the irreducible characters of the exceptional irreducible complex reflection groups are the ones used by the GAP package CHEVIE and are explained in subsection 5.2.3.
\section{The groups $G_4$,\ldots,$G_7$}

The following table gives the specializations of the parameters of
the generic Hecke algebra $\mathcal{H}(G_7)$,
$(x_0,x_1;y_0,y_1,y_2;z_0,z_1,z_2)$, which give the generic Hecke
algebras of the groups $G_4$, $G_5$ and $G_6$ (\cite{Ma2}, Table
4.6).

\begin{center}
\begin{tabular}{|c|c|c|c|c|}
  \hline
  Group & Index & S & T & U \\
  \hline
  $G_7$ & 1 & $x_0,x_1$ & $y_0,y_1,y_2$         & $z_0,z_1,z_2$\\
  $G_5$ & 2 & $1,-1$    & $y_0,y_1,y_2$         & $z_0,z_1,z_2$\\
  $G_6$ & 3 & $x_0,x_1$ & $1,\zeta_3,\zeta_3^2$ & $z_0,z_1,z_2$\\
  $G_4$ & 6 & $1,-1$    & $1,\zeta_3,\zeta_3^2$ & $z_0,z_1,z_2$\\
  \hline
\end{tabular}
\end{center}
$$\textrm{Specializations of the parameters for }\mathcal{H}(G_7)$$
\\

\begin{lemma}\label{g7}\
\begin{enumerate}
  \item The algebra $\mathcal{H}(G_7)$ specialized via
  $$(x_0,x_1;y_0,y_1,y_2;z_0,z_1,z_2) \mapsto
  (1,-1;y_0,y_1,y_2;z_0,z_1,z_2)$$
  is the twisted symmetric algebra of the
  cyclic group $C_2$ over the symmetric subalgebra $\mathcal{H}(G_5)$ with
  parameters $(y_0,y_1,y_2;z_0,z_1,z_2)$. 
  \item The algebra $\mathcal{H}(G_7)$ specialized via
  $$(x_0,x_1;y_0,y_1,y_2;z_0,z_1,z_2) \mapsto
  (x_0,x_1;1,\zeta_3,\zeta_3^2;z_0,z_1,z_2)$$
  is the twisted symmetric algebra of the
  cyclic group $C_3$ over the symmetric subalgebra $\mathcal{H}(G_6)$ with
  parameters $(x_0,x_1;z_0,z_1,z_2)$. 
  \item The algebra $\mathcal{H}(G_6)$ specialized via
  $$(x_0,x_1;z_0,z_1,z_2) \mapsto
  (1,-1;z_0,z_1,z_2)$$
  is the twisted symmetric algebra of the
  cyclic group $C_2$ over the symmetric subalgebra $\mathcal{H}(G_4)$ with
  parameters $(z_0,z_1,z_2)$. 
\end{enumerate}
\end{lemma}
\begin{apod}{\small We have
$$\begin{array}{rccl}
    \mathcal{H}(G_7) & = & <S,T,U\,\, | &  STU=TUS=UST,  \\
     &  &  & (S-x_0)(S-x_1)=0, \\
     &  &  & (T-y_0)(T-y_1)(T-y_2)=0, \\
     &  &  & (U-z_0)(U-z_1)(U-z_2)=0>.
  \end{array}$$
\begin{enumerate}
\item Let
  $$\begin{array}{rccl}
     A & := & <S,T,U \,\,| &  STU=TUS=UST, S^2=1,  \\
     &  &  & (T-y_0)(T-y_1)(T-y_2)=0, \\
     &  &  & (U-z_0)(U-z_1)(U-z_2)=0>
  \end{array}$$
   and $$\bar{A}:=<T,U>.$$
   Then $$A=\bar{A} \oplus S \bar{A} = \bar{A} \oplus  \bar{A}S\,\textrm{ and }\, \bar{A} \cong \mathcal{H}(G_5).$$
  \item Let
   $$\begin{array}{rccl}
     A & := & <S,T,U \,\,| &  STU=TUS=UST, T^3=1,  \\
     &  &  & (S-x_0)(S-x_1)=0, \\
     &  &  & (U-z_0)(U-z_1)(U-z_2)=0>
  \end{array}$$
  and $$\bar{A}:=<S,U>.$$
  Then $$A=\bigoplus_{i=0}^2 T^i \bar{A}=\bigoplus_{i=0}^2 \bar{A}T^i\,  \textrm{ and }\, \bar{A} \cong \mathcal{H}(G_6).$$
  \item Let
  $$\begin{array}{rccl}
     A & := & <S,U \,\,| &  SUSUSU=USUSUS, S^2=1,  \\
     &  &  & (U-z_0)(U-z_1)(U-z_2)=0>
  \end{array}$$
   and $$\bar{A}:=<U,SUS>.$$
   Then $$A=\bar{A} \oplus S \bar{A} =\bar{A} \oplus  \bar{A}S\,\textrm{ and } \,\bar{A} \cong \mathcal{H}(G_4).$$}
\end{enumerate}
\end{apod}

The Schur elements of all irreducible characters of
$\mathcal{H}(G_7)$ are calculated in \cite{Ma2}. They are
obtained via Galois transformations (permutation of indeterminates, permutation of roots of unity or combination of the two)
 from the following ones: {\footnotesize
\\
\\
$s_{\phi_{1,0}}=\Phi_{1}(x_0/x_1)\cdot\allowbreak\Phi_{1}(x_0y_0^2z_0^2/x_1y_1y_2z_1z_2)
\cdot\allowbreak\Phi_{1}(y_0/y_1)\cdot\allowbreak\Phi_{1}(y_0/y_2)\cdot
\allowbreak\Phi_{1}(z_0/z_1)\cdot\allowbreak\Phi_{1}(z_0/z_2)\cdot\allowbreak
\Phi_{1}(x_0y_0z_0/x_1y_1z_1)\cdot\allowbreak\Phi_{1}(x_0y_0z_0/x_1y_1z_2)\cdot
\allowbreak\Phi_{1}(x_0y_0z_0/x_1y_2z_1)\cdot\allowbreak
\Phi_{1}(x_0y_0z_0/x_1y_2z_2)$
\\
\\
$s_{\phi_{2,9'}}=2y_2/y_0\Phi_{1}(y_0/y_1)\cdot\allowbreak\Phi_{1}(y_2/y_0)\cdot\allowbreak
\Phi_{1}(z_1/z_0)\cdot\allowbreak\Phi_{1}(z_2/z_0)\cdot\allowbreak\
\Phi_{1}(r/x_0y_0z_0)\cdot\allowbreak\Phi_{1}(r/x_0y_2z_1)\cdot\allowbreak
\Phi_{1}(r/x_0y_2z_2)\cdot\allowbreak\Phi_{1}(r/x_1y_0z_0)\cdot\allowbreak
\Phi_{1}(r/x_1y_2z_1)\cdot\allowbreak\Phi_{1}(r/x_1y_2z_2)$\\
where $r=\root 2\of{x_0x_1y_1y_2z_1z_2}$
\\
\\
$s_{\phi_{3,6}}=3\Phi_{1}(x_1/x_0)\cdot\allowbreak\Phi_{1}(x_0y_0z_0/r)\cdot
\allowbreak\Phi_{1}(x_0y_0z_1/r)\cdot\allowbreak\Phi_{1}(x_0y_0z_2/r)\cdot
\allowbreak\Phi_{1}(x_0y_1z_0/r)\cdot\allowbreak\Phi_{1}(x_0y_1z_1/r)\cdot
\allowbreak\Phi_{1}(x_0y_1z_2/r)\cdot\allowbreak\Phi_{1}(x_0y_2z_0/r)\cdot
\allowbreak\Phi_{1}(x_0y_2z_1/r)\cdot\allowbreak\Phi_{1}(x_0y_2z_2/r)$\\
where $\ r=\root 3\of{x_0^2x_1y_0y_1y_2z_0z_1z_2}$}
\\
\\
Following Theorem $\ref{Semisimplicity Malle}$ and \cite{Ma4}, Table
8.1, if we set
$$\begin{array}{cl}
    X_i^{12}:=(\zeta_2)^{-i}x_i & (i=0,1), \\
    Y_j^{12}:=(\zeta_3)^{-j}y_j & (j=0,1,2), \\
    Z_k^{12}:=(\zeta_3)^{-k}z_k & (k=0,1,2),
  \end{array}$$
then $\mathbb{Q}(\zeta_{12})(X_0,X_1,Y_0,Y_1,Y_2,Z_0,Z_1,Z_2)$ is a
splitting field for $\mathcal{H}(G_7)$. Hence, the factorization of
the Schur elements over that field is as described by Theorem
$\ref{Schur element generic}$.

\section{The groups $G_8, \ldots, G_{15}$}

The following table gives the specializations of the parameters of
the generic Hecke algebra $\mathcal{H}(G_{11})$,
$(x_0,x_1;y_0,y_1,y_2;z_0,z_1,z_2,z_3)$, which give the generic
Hecke algebras of the groups $G_8,\ldots,G_{15}$ (\cite{Ma2}, Table
4.9).

\begin{center}
\begin{tabular}{|c|c|c|c|c|}
  \hline
  Group & Index & S & T & U \\
  \hline
  $G_{11}$ & 1  & $x_0,x_1$ & $y_0,y_1,y_2$         & $z_0,z_1,z_2,z_3$\\
  $G_{10}$ & 2  & $1,-1$    & $y_0,y_1,y_2$         & $z_1,z_1,z_2,z_3$\\
  $G_{15}$ & 2  & $x_0,x_1$ & $y_0,y_1,y_2$         & $\sqrt{u_0},\sqrt{u_1},-\sqrt{u_0},-\sqrt{u_1}$\\
  $G_9$    & 3  & $x_0,x_1$ & $1,\zeta_3,\zeta_3^2$ & $z_0,z_1,z_2,z_3$\\
  $G_{14}$ & 4  & $x_0,x_1$ & $y_0,y_1,y_2$         & $1,i,-1,-i$\\
  $G_8$    & 6  & $1,-1$    & $1,\zeta_3,\zeta_3^2$ & $z_0,z_1,z_2,z_3$\\
  $G_{13}$ & 6  & $x_0,x_1$ & $1,\zeta_3,\zeta_3^2$ & $\sqrt{u_0},\sqrt{u_1},-\sqrt{u_0},-\sqrt{u_1}$\\
  $G_{12}$ & 12 & $x_0,x_1$ & $1,\zeta_3,\zeta_3^2$ & $1,i,-1,-i$\\
  \hline
\end{tabular}
\end{center}

$$\textrm{Specializations of the parameters for }\mathcal{H}(G_{11})$$

\begin{lemma}\label{g11}\
\begin{enumerate}
 \item The algebra $\mathcal{H}(G_{11})$ specialized via
  $$(x_0,x_1;y_0,y_1,y_2;z_0,z_1,z_2,z_3) \mapsto
  (1,-1;y_0,y_1,y_2;z_0,z_1,z_2,z_3)$$
  is the twisted symmetric algebra of the
  cyclic group $C_2$ over the symmetric subalgebra $\mathcal{H}(G_{10})$ with
  parameters $(y_0,y_1,y_2;z_0,z_1,z_2,z_3)$. 
 \item The algebra $\mathcal{H}(G_{11})$ specialized via
  $$(x_0,x_1;y_0,y_1,y_2;z_0,z_1,z_2,z_3) \mapsto
  (x_0,x_1;1,\zeta_3,\zeta_3^2;z_0,z_1,z_2,z_3)$$
  is the twisted symmetric algebra of the
  cyclic group $C_3$ over the symmetric subalgebra $\mathcal{H}(G_9)$ with
  parameters $(x_0,x_1;z_0,z_1,z_2,z_3)$. 
 \item The algebra $\mathcal{H}(G_9)$ specialized via
  $$(x_0,x_1;z_0,z_1,z_2,z_3) \mapsto
  (1,-1;z_0,z_1,z_2,z_3)$$
  is the twisted symmetric algebra of the
  cyclic group $C_2$ over the symmetric subalgebra $\mathcal{H}(G_8)$ with
  parameters $(z_0,z_1,z_2,z_3)$. 
 \item The algebra $\mathcal{H}(G_{11})$ specialized via
  $$(x_0,x_1;y_0,y_1,y_2;z_0,z_1,z_2,z_3) \mapsto
  (x_0,x_1;y_0,y_1,y_2;1,i,-1,-i)$$
  is the twisted symmetric algebra of the
  cyclic group $C_4$ over the symmetric subalgebra $\mathcal{H}(G_{14})$ with
  parameters $(x_0,x_1;y_0,y_1,y_2)$. 
 \item The algebra $\mathcal{H}(G_{14})$ specialized via
  $$(x_0,x_1;y_0,y_1,y_2) \mapsto
  (x_0,x_1;1,\zeta_3,\zeta_3^2)$$
  is the twisted symmetric algebra of the
  cyclic group $C_3$ over the symmetric subalgebra $\mathcal{H}(G_{12})$ with
  parameters $(x_0,x_1)$. 
 \item The algebra $\mathcal{H}(G_{11})$ specialized via
  $$(x_0,x_1;y_0,y_1,y_2;z_0,z_1,z_2,z_3) \mapsto
  (x_0,x_1;y_0,y_1,y_2;\sqrt{u_0},\sqrt{u_1},-\sqrt{u_0},-\sqrt{u_1})$$
  is the twisted symmetric algebra of the
  cyclic group $C_2$ over the symmetric subalgebra $\mathcal{H}(G_{15})$ with
  parameters $(x_0,x_1;y_0,y_1,y_2;u_0,u_1)$. 
 \item The algebra $\mathcal{H}(G_{15})$ specialized via
  $$(x_0,x_1;y_0,y_1,y_2;u_0,u_1) \mapsto
  (x_0,x_1;1,\zeta_3,\zeta_3^2;u_0,u_1)$$
  is the twisted symmetric algebra of the
  cyclic group $C_3$ over the symmetric subalgebra $\mathcal{H}(G_{13})$ with
  parameters $(x_0,x_1;u_0,u_1)$. 
\end{enumerate}
\end{lemma}
\begin{apod}{\small We have $$\begin{array}{rccl}
    \mathcal{H}(G_{11}) & = & <S,T,U \,\,| &  STU=TUS=UST,  \\
     &  &  & (S-x_0)(S-x_1)=0, \\
     &  &  & (T-y_0)(T-y_1)(T-y_2)=0, \\
     &  &  & (U-z_0)(U-z_1)(U-z_2)(U-z_3)=0>.
  \end{array}$$
\begin{enumerate}
\item Let
  $$\begin{array}{rccl}
     A & := & <S,T,U \,\,| &  STU=TUS=UST, S^2=1,  \\
     &  &  & (T-y_0)(T-y_1)(T-y_2)=0, \\
     &  &  & (U-z_0)(U-z_1)(U-z_2)(U-z_3)=0>
  \end{array}$$
   and $$\bar{A}:=<T,U>.$$
   Then $$ A=\bar{A} \oplus S \bar{A}=\bar{A} \oplus  \bar{A}S\, \textrm{ and }\, \bar{A} \cong \mathcal{H}(G_{10}).$$
\item Let
   $$\begin{array}{rccl}
     A & := & <S,T,U \,\,| &  STU=TUS=UST, T^3=1,  \\
     &  &  & (S-x_0)(S-x_1)=0, \\
     &  &  & (U-z_0)(U-z_1)(U-z_2)(U-z_3)=0>
  \end{array}$$
  and $$\bar{A}:=<S,U>.$$ Then
  $$A=\bigoplus_{i=0}^2 T^i \bar{A}=\bigoplus_{i=0}^2  \bar{A}T^i \, \textrm{ and }\, \bar{A}  \cong \mathcal{H}(G_9).$$
\item Let
  $$\begin{array}{rccl}
     A & := & <S,U \,\,| &  SUSUSU=USUSUS, S^2=1,  \\
     &  &  & (U-z_0)(U-z_1)(U-z_2)(U-z_3)=0>
  \end{array}$$
   and $$\bar{A}:=<U,SUS>.$$
   Then $$A=\bar{A} \oplus S \bar{A}=\bar{A} \oplus  \bar{A}S\, \textrm{ and }\, \bar{A}  \cong \mathcal{H}(G_8).$$
\item Let
   $$\begin{array}{rccl}
    A& := & <S,T,U \,\,| &  STU=TUS=UST, U^4=1,  \\
     &  &  & (S-x_0)(S-x_1)=0,\\
     &  &  & (T-y_0)(T-y_1)(T-y_2)=0>
  \end{array}$$
  and $$\bar{A}:=<S,T>.$$ Then
  $$A=\bigoplus_{i=0}^3 U^i \bar{A}=\bigoplus_{i=0}^3  \bar{A}U^i \, \textrm{ and }\, \bar{A}\cong \mathcal{H}(G_{14}).$$
\item Let
   $$\begin{array}{rccl}
    A& := & <S,T \,\,| &  STSTSTST=TSTSTSTS, T^3=1,  \\
     &  &  & (S-x_0)(S-x_1)=0 >
  \end{array}$$
  and $$\bar{A}:=<S,TST^2,T^2ST>.$$
  Then
  $$A=\bigoplus_{i=0}^2 T^i \bar{A}=\bigoplus_{i=0}^2  \bar{A}T^i \,\textrm{ and }\, \bar{A}  \cong \mathcal{H}(G_{12}).$$
 \item Let
  $$\begin{array}{rccl}
     A & := & <S,T,U \,\,| &  STU=TUS=UST,  \\
     &  &  & (S-x_0)(S-x_1)=0, \\
     &  &  & (T-y_0)(T-y_1)(T-y_2)=0, \\
     &  &  & (U^2-u_0)(U^2-u_1)=0>
  \end{array}$$
   and $$\bar{A}:=<S,T,U^2>.$$
   Then $$ A=\bar{A} \oplus U \bar{A}=\bar{A} \oplus  \bar{A}U\, \textrm{ and }\, \bar{A} \cong \mathcal{H}(G_{15}).$$
 \item Let
   $$\begin{array}{rccl}
     A & := & <U^2,S,T \,\,| & STU^2=U^2ST, U^2STST=TU^2STS, T^3=1,  \\
     &  &  & (S-x_0)(S-x_1)=0, \\
     &  &  & (U^2-u_0)(U^2-u_1)=0>
  \end{array}$$
  and $$\bar{A}:=<U^2,S,T^2ST>.$$ Then
  $$ A=\bigoplus_{i=0}^2 T^i \bar{A}=\bigoplus_{i=0}^2  \bar{A}T^i \,\textrm{ and }\, \bar{A} \cong \mathcal{H}(G_{13}).$$}
\end{enumerate}
\end{apod}

The Schur elements of all irreducible characters of
$\mathcal{H}(G_{11})$ are calculated in \cite{Ma2}. They are
obtained via Galois transformations from the following ones:
\\
\\
{\footnotesize
$s_{\phi_{1,0}}=\Phi_{1}(x_0/x_1)\cdot\allowbreak\Phi_{1}(y_0/y_1)\cdot\allowbreak
\Phi_{1}(y_0/y_2)\cdot\allowbreak\Phi_{1}(z_0/z_1)\cdot\allowbreak\Phi_{1}(
z_0/z_2)\cdot\allowbreak\Phi_{1}(z_0/z_3)\cdot\allowbreak\Phi_{1}(x_0y_0
z_0/x_1y_1z_1)\cdot\allowbreak\Phi_{1}(x_0y_0z_0/x_1y_1z_2)\cdot\allowbreak
\Phi_{1}(x_0y_0z_0/x_1y_1z_3)\cdot\allowbreak\Phi_{1}(x_0y_0z_0/x_1y_2z_
1)\cdot\allowbreak \Phi_{1}(x_0y_0z_0/x_1y_2z_2)\cdot\allowbreak
\\ \Phi_{1}(x_0y_0z_0/x_1y_2z_3)\cdot\allowbreak\Phi_{1}(x_0y_0^2z_0^2/x_1y_1y_2z_1z_2)
\cdot\allowbreak
\Phi_{1}(x_0y_0^2z_0^2/x_1y_1y_2z_1z_3)\cdot\allowbreak
\Phi_{1}(x_0y_0^2z_0^2/x_1y_1y_2z_2z_3)\cdot\allowbreak\Phi_{1}(x_0^2y_0^2z_0^3/x_1^2y_1y_2z_1z_2z_3)$
\\
\\
$s_{\phi_{2,1}}=-2z_1/z_0\Phi_{1}(y_0/y_2)\cdot\allowbreak\Phi_{1}(y_1/y_2)\cdot
\allowbreak\Phi_{1}(z_0/z_2)\cdot\allowbreak\Phi_{1}(z_0/z_3)\cdot\allowbreak
\Phi_{1}(z_1/z_2)\cdot\allowbreak\Phi_{1}(z_1/z_3)\cdot\allowbreak\Phi_{1}
(y_0z_0z_1/y_2z_2z_3)\cdot\allowbreak\Phi_{1}(y_1z_0z_1/y_2z_2z_3)\cdot
\allowbreak\Phi_{1}(r/x_0y_2z_2)\cdot\allowbreak\Phi_{1}(r/x_0y_2z_3)\cdot
\allowbreak\Phi_{1}(r/x_1y_2z_2)\cdot\allowbreak\\
\Phi_{1}(r/x_1y_2z_3)\cdot
\allowbreak\Phi_{1}(r/x_0y_0z_1)\cdot\allowbreak\Phi_{1}(r/x_0y_1z_1)\
\cdot\allowbreak\Phi_{1}(r/x_1y_0z_1)\cdot\allowbreak\Phi_{1}(r/x_1y_1z_1)$\\
where $r=\root 2\of{x_0x_1y_0y_1z_0z_1}$
\\
\\
$s_{\phi_{3,2}}=3\Phi_{1}(x_1/x_0)\cdot\allowbreak\Phi_{1}(z_1/z_3)\cdot\allowbreak
\Phi_{1}(z_2/z_3)\cdot\allowbreak\Phi_{1}(z_0/z_3)\cdot\allowbreak\Phi_{1}
(r/x_1y_0z_3)\cdot\allowbreak\Phi_{1}(r/x_1y_1z_3)\cdot\allowbreak\Phi_1
(r/x_1y_2z_3)\cdot\allowbreak\Phi_{1}(x_0y_0z_0/r)\cdot\allowbreak\Phi
_{1}(x_0y_0z_1/r)\cdot\allowbreak\Phi_{1}(x_0y_0z_2/r)\cdot\allowbreak\Phi_{1}
(x_0y_1z_0/r)\cdot\allowbreak
\\ \Phi_{1}(x_0y_1z_1/r)\cdot\allowbreak\
\Phi_{1}(x_0y_1z_2/r)\cdot\allowbreak\Phi_{1}(x_0y_2z_0/r)\cdot\allowbreak
\Phi_{1}(x_0y_2z_1/r)\cdot\allowbreak\Phi_{1}(x_0y_2z_2/r)$\\
where $r=\root 3\of{x_0^2x_1y_0y_1y_2z_0z_1z_2}$
\\
\\
$s_{\phi_{4,21}}=-4\Phi_{1}(y_0/y_1)\cdot\allowbreak\Phi_{1}(y_0/y_2)\cdot\allowbreak\
\Phi_{1}(r/x_0y_0z_0)\cdot\allowbreak\Phi_{1}(r/x_1y_0z_0)\cdot\allowbreak
\Phi_{1}(x_0y_0z_1/r)\cdot\allowbreak\Phi_{1}(x_0y_0z_2/r)\cdot\allowbreak
\Phi_{1}(x_0y_0z_3/r)\cdot\allowbreak\Phi_{1}(x_1y_0z_1/r)\cdot\allowbreak
\Phi_{1}(x_1y_0z_2/r)\cdot\allowbreak\Phi_{1}(x_1y_0z_3/r)\cdot\allowbreak
\\
\Phi_{1}(x_0x_1y_0y_1z_0z_1/r^2)\cdot\allowbreak\Phi_{1}(x_0x_1y_0y_1
z_0z_2/r^2)\cdot\allowbreak\Phi_{1}(x_0x_1y_0y_1z_0z_3/r^2)\cdot\allowbreak
\\ \Phi_{1}(x_0x_1y_0y_2z_0z_1/r^2)\cdot\allowbreak\Phi_{1}(x_0x_1y_0y_2z_0z_2/r^2)
\cdot\allowbreak\Phi_{1}(x_0x_1y_0y_2z_0z_3/r^2)$\\
where $r=\root 4\of{x_0^2x_1^2y_0^2y_1y_2z_0z_1z_2z_3}$}
\\
\\
Following Theorem $\ref{Semisimplicity Malle}$ and \cite{Ma4}, Table
8.1, if we set
$$\begin{array}{cl}
    X_i^{24}:=(\zeta_2)^{-i}x_i & (i=0,1), \\
    Y_j^{24}:=(\zeta_3)^{-j}y_j & (j=0,1,2), \\
    Z_k^{24}:=(\zeta_4)^{-k}z_k & (k=0,1,2,3),
  \end{array}$$
then $\mathbb{Q}(\zeta_{24})(X_0,X_1,Y_0,Y_1,Y_2,Z_0,Z_1,Z_2,Z_3)$
is a splitting field for $\mathcal{H}(G_{11})$. Hence, the
factorization of the Schur elements over that field is as described
by Theorem $\ref{Schur element generic}$.

\section{The groups $G_{16}, \ldots, G_{22}$}

The following table gives the specializations of the parameters of
the generic Hecke algebra $\mathcal{H}(G_{19})$,
$(x_0,x_1;y_0,y_1,y_2;z_0,z_1,z_2,z_3,z_4)$, which give the generic
Hecke algebras of the groups $G_{16},\ldots,G_{22}$ (\cite{Ma2},
Table 4.12).

\begin{center}
\begin{tabular}{|c|c|c|c|c|}
  \hline
  Group & Index & S & T & U \\
  \hline
  $G_{19}$ & 1  & $x_0,x_1$ & $y_0,y_1,y_2$         & $z_0,z_1,z_2,z_3,z_4$\\
  $G_{18}$ & 2  & $1,-1$    & $y_0,y_1,y_2$         & $z_0,z_1,z_2,z_3,z_4$\\
  $G_{17}$ & 3  & $x_0,x_1$ & $1,\zeta_3,\zeta_3^2$ & $z_0,z_1,z_2,z_3,z_4$\\
  $G_{21}$ & 5  & $x_0,x_1$ & $y_0,y_1,y_2$         & $1,\zeta_5,\zeta_5^2,\zeta_5^3,\zeta_5^4$\\
  $G_{16}$ & 6  & $1,-1$    & $1,\zeta_3,\zeta_3^2$ & $z_0,z_1,z_2,z_3,z_4$\\
  $G_{20}$ & 10 & $1,-1$    & $y_0,y_1,y_2$         & $1,\zeta_5,\zeta_5^2,\zeta_5^3,\zeta_5^4$\\
  $G_{22}$ & 15 & $x_0,x_1$ & $1,\zeta_3,\zeta_3^2$ & $1,\zeta_5,\zeta_5^2,\zeta_5^3,\zeta_5^4$\\
  \hline
\end{tabular}
\end{center}

$$\textrm{Specializations of the parameters for }\mathcal{H}(G_{19})$$

\begin{lemma}\label{g19}\
\begin{enumerate}
 \item The algebra $\mathcal{H}(G_{19})$ specialized via
  $$(x_0,x_1;y_0,y_1,y_2;z_0,z_1,z_2,z_3,z_4) \mapsto
  (1,-1;y_0,y_1,y_2;z_0,z_1,z_2,z_3,z_4)$$
  is the twisted symmetric algebra of the
  cyclic group $C_2$ over the symmetric subalgebra $\mathcal{H}(G_{18})$ with
  parameters $(y_0,y_1,y_2;z_0,z_1,z_2,z_3,z_4)$. 
 \item The algebra $\mathcal{H}(G_{19})$ specialized via
  $$(x_0,x_1;y_0,y_1,y_2;z_0,z_1,z_2,z_3,z_4) \mapsto
  (x_0,x_1;1,\zeta_3,\zeta_3^2;z_0,z_1,z_2,z_3,z_4)$$
  is the twisted symmetric algebra of the
  cyclic group $C_3$ over the symmetric subalgebra $\mathcal{H}(G_{17})$ with
  parameters $(x_0,x_1;z_0,z_1,z_2,z_3,z_4)$. 
\item The algebra $\mathcal{H}(G_{17})$ specialized via
  $$(x_0,x_1;z_0,z_1,z_2,z_3,z_4) \mapsto
  (1,-1;z_0,z_1,z_2,z_3,z_4)$$
  is the twisted symmetric algebra of the
  cyclic group $C_2$ over the symmetric subalgebra $\mathcal{H}(G_{16})$ with
  parameters $(z_0,z_1,z_2,z_3,z_4)$. 
\item The algebra $\mathcal{H}(G_{19})$ specialized via
  $$(x_0,x_1;y_0,y_1,y_2;z_0,z_1,z_2,z_3,z_4) \mapsto
  (x_0,x_1;y_0,y_1,y_2;1,\zeta_5,\zeta_5^2,\zeta_5^3,\zeta_5^4)$$
  is the twisted symmetric algebra of the
  cyclic group $C_5$ over the symmetric subalgebra $\mathcal{H}(G_{21})$ with
  parameters $(x_0,x_1;y_0,y_1,y_2)$. 
\item The algebra $\mathcal{H}(G_{21})$ specialized via
  $$(x_0,x_1;y_0,y_1,y_2) \mapsto
  (1,-1;y_0,y_1,y_2)$$
  is the twisted symmetric algebra of the
  cyclic group $C_2$ over the symmetric subalgebra $\mathcal{H}(G_{20})$ with
  parameters $(y_0,y_1,y_2)$. 
\item The algebra $\mathcal{H}(G_{21})$ specialized via
  $$(x_0,x_1;y_0,y_1,y_2) \mapsto
  (x_0,x_1;1,\zeta_3,\zeta_3^2)$$
  is the twisted symmetric algebra of the
  cyclic group $C_3$ over the symmetric subalgebra $\mathcal{H}(G_{22})$ with
  parameters $(x_0,x_1)$. 
\end{enumerate}
\end{lemma}
\begin{apod}{\small We have $$\begin{array}{rccl}
    \mathcal{H}(G_{19}) & = & <S,T,U \,\,| &  STU=TUS=UST,  \\
     &  &  & (S-x_0)(S-x_1)=0, \\
     &  &  & (T-y_0)(T-y_1)(T-y_2)=0, \\
     &  &  & (U-z_0)(U-z_1)(U-z_2)(U-z_3)(U-z_4)=0>.
  \end{array}$$
\begin{enumerate}
\item Let
  $$\begin{array}{rccl}
     A & := & <S,T,U \,\,| &  STU=TUS=UST, S^2=1,  \\
     &  &  & (T-y_0)(T-y_1)(T-y_2)=0, \\
     &  &  & (U-z_0)(U-z_1)(U-z_2)(U-z_3)(U-z_4)=0>
  \end{array}$$
   and $$\bar{A}:=<T,U>.$$ Then
   $$A=\bar{A} \oplus S \bar{A}=\bar{A} \oplus  \bar{A}S\, \textrm{ and }\, \bar{A} \cong \mathcal{H}(G_{18}).$$
  \item Let
   $$\begin{array}{rccl}
     A & := & <S,T,U \,\,| &  STU=TUS=UST, T^3=1,  \\
     &  &  & (S-x_0)(S-x_1)=0, \\
     &  &  & (U-z_0)(U-z_1)(U-z_2)(U-z_3)(U-z_4)=0>
  \end{array}$$
  and $$\bar{A}:=<S,U>.$$
  Then $$A=\bigoplus_{i=0}^2 T^i \bar{A} =\bigoplus_{i=0}^2  \bar{A}T^i \,\textrm{ and } \,\bar{A} \cong \mathcal{H}(G_{17}).$$
  \item Let
  $$\begin{array}{rccl}
     A & := & <S,U \,\,| &  SUSUSU=USUSUS, S^2=1,  \\
     &  &  & (U-z_0)(U-z_1)(U-z_2)(U-z_3)(U-z_4)=0>
  \end{array}$$
   and $$\bar{A}:=<U,SUS>.$$ Then
   $$A=\bar{A} \oplus S \bar{A}=\bar{A} \oplus  \bar{A}S\, \textrm{ and }\, \bar{A} \cong
   \mathcal{H}(G_{16}).$$
   \item Let
   $$\begin{array}{rccl}
    A& := & <S,T,U \,\,| &  STU=TUS=UST, U^5=1,  \\
     &  &  & (S-x_0)(S-x_1)=0, \\
     &  &  & (T-y_0)(T-y_1)(T-y_2)=0>
  \end{array}$$
  and $$\bar{A}:=<S,T>.$$ Then
  $$A=\bigoplus_{i=0}^4 U^i \bar{A}=\bigoplus_{i=0}^4  \bar{A}U^i \,  \textrm{ and } \bar{A}\, \cong \mathcal{H}(G_{21}).$$
\item Let
   $$\begin{array}{rccl}
    A& := & <S,T \,\,| &  STSTSTSTST=TSTSTSTSTS, S^2=1,  \\
     &  &  & (T-y_0)(T-y_1)(T-y_2)=0>
  \end{array}$$
  and $$\bar{A}:=<T,STS>.$$ Then
  $$ A=\bar{A} \oplus S\bar{A} =\bar{A} \oplus  \bar{A}S\,\textrm{ and } \,\bar{A} \cong \mathcal{H}(G_{20}).$$
\item Let
   $$\begin{array}{rccl}
    A& := & <S,T \,\,| &  STSTSTSTST=TSTSTSTSTS, T^3=1,  \\
     &  &  & (S-x_0)(S-x_1)=0 >
  \end{array}$$
  and $$\bar{A}:=<S,TST^2,T^2ST>.$$  Then
  $$A=\bigoplus_{i=0}^2 T^i \bar{A}=\bigoplus_{i=0}^2  \bar{A}T^i \, \textrm{ and } \,\bar{A}\cong
  \mathcal{H}(G_{22}).$$}
\end{enumerate}
\end{apod}

The Schur elements of all irreducible characters of
$\mathcal{H}(G_{19})$ are calculated in \cite{Ma2}.  They are
obtained via Galois transformations from the following ones:
\\
\\
{\footnotesize
$s_{\phi_{1,0}}=\Phi_{1}(x_0/x_1)\cdot\allowbreak\Phi_{1}(y_0/y_1)\cdot\allowbreak
\Phi_{1}(y_0/y_2)\cdot\allowbreak\Phi_{1}(z_0/z_1)\cdot\allowbreak\Phi_{1}(
z_0/z_2)\cdot\allowbreak\Phi_{1}(z_0/z_3)\cdot\allowbreak\Phi_{1}(z_0/z_4)
\cdot\allowbreak\Phi_{1}(x_0y_0z_0/x_1y_1z_1)\cdot\allowbreak\Phi_{1}(
x_0y_0z_0/x_1y_1z_2)\cdot\allowbreak\Phi_{1}(x_0y_0z_0/x_1y_1z_3)\cdot
\allowbreak\Phi_{1}(x_0y_0z_0/x_1y_1z_4)\cdot\allowbreak\\
\Phi_{1}(x_0y_0z_0/x_1y_2z_1)
\cdot\allowbreak\Phi_{1}(x_0y_0z_0/x_1y_2z_2)\cdot\allowbreak
\Phi_{1}(x_0y_0z_0/x_1y_2z_3)\cdot\allowbreak\Phi_{1}(x_0y_0z_0/x_1y_2z_4)\cdot
\allowbreak\\
\Phi_{1}(x_0y_0^2z_0^2/x_1y_1y_2z_1z_2)\cdot\allowbreak\Phi_1
(x_0y_0^2z_0^2/x_1y_1y_2z_1z_3)\cdot\allowbreak
\Phi_{1}(x_0y_0^2z_0^2/x_1y_1y_2z_1z_4) \cdot\allowbreak\\
\Phi_{1}(x_0y_0^2z_0^2/x_1y_1y_2z_2z_3)\cdot
\allowbreak\Phi_{1}(x_0y_0^2z_0^2/x_1y_1y_2z_2z_4)\cdot\allowbreak\Phi_{1}
(x_0y_0^2z_0^2/x_1y_1y_2z_3z_4)\cdot\allowbreak\\
\Phi_{1}(x_0^2y_0^2z_0^3/x_1^2y_1y_2z_1z_2z_3)\cdot\allowbreak\Phi_{1}(x_0^2y_0^2z_0^3/x_1^2y_1y_2z_1z_2z_4)
\cdot\allowbreak\Phi_{1}(x_0^2y_0^2z_0^3/x_1^2y_1y_2z_1z_3z_4)\cdot\allowbreak\\
\Phi_{1}(x_0^2y_0^2z_0^3/x_1^2y_1y_2z_2z_3z_4)\cdot\allowbreak
\Phi_{1}(x_0^2y_0^3z_0^4/x_1^2y_1^2y_2z_1z_2z_3z_4)\cdot\allowbreak\Phi_{1}
(x_0^2y_0^3z_0^4/x_1^2y_1y_2^2z_1z_2z_3z_4)\cdot\allowbreak\\
\Phi_{1}(x_0^3y_0^4z_0^4/x_1^3y_1^2y_2^2z_1z_2z_3z_4)$
\\
\\
$s_{\phi_{2,31'}}=-2\Phi_{1}(y_0/y_2)\cdot\allowbreak\Phi_{1}(y_1/y_2)\cdot\allowbreak
\Phi_{1}(z_0/z_2)\cdot\allowbreak\Phi_{1}(z_0/z_3)\cdot\allowbreak\Phi_{1}
(z_0/z_4)\cdot\allowbreak\Phi_{1}(z_1/z_2)\cdot\allowbreak\Phi_{1}(z_1/
z_3)\cdot\allowbreak\Phi_{1}(z_1/z_4)\cdot\allowbreak\Phi_{1}(y_0z_0z_1/
y_2z_2z_3)\cdot\allowbreak\Phi_{1}(y_0z_0z_1/y_2z_2z_4)\cdot\allowbreak
\Phi_{1}(y_0z_0z_1/y_2z_3z_4)\cdot\allowbreak\Phi_{1}(y_1z_0z_1/y_2z_2z_3)
\cdot\allowbreak\Phi_{1}(y_1z_0z_1/y_2z_2z_4)\cdot\allowbreak\Phi_{1}(y_1z_0z_1/y_2z_3z_4)
\cdot\allowbreak\Phi_{1}(y_0y_1z_0z_1^2/y_2^2z_2z_3z_4)\cdot
\allowbreak\Phi_{1}(y_0y_1z_0^2z_1/y_2^2z_2z_3z_4)\cdot\allowbreak\Phi_{1}
(r/x_0y_0z_0)\cdot\allowbreak\Phi_{1}(x_0y_0z_1/r)\cdot\allowbreak\Phi_{1}
(x_1y_0z_0/r)\cdot\allowbreak\Phi_{1}(r/x_1y_0z_1)\cdot\allowbreak\
\Phi_{1}(r/x_1y_2z_2)\cdot\allowbreak\Phi_{1}(r/x_1y_2z_3)\cdot\allowbreak
\Phi_{1}(r/x_1y_2z_4)\cdot\allowbreak\Phi_{1}(r/x_0y_2z_2)\cdot\allowbreak
\Phi_{1}(r/x_0y_2z_3)\cdot\allowbreak\Phi_{1}(r/x_0y_2z_4)\cdot\allowbreak
\Phi_{1}(rz_0z_1/x_0y_2z_2z_3z_4)\cdot\allowbreak\\ \Phi_{1}(rz_0z_1/x_1y_2z_2z_3z_4)$\\
where $r=\root 2\of{x_0x_1y_0y_1z_0z_1}$
\\
\\
$s_{\phi_{3,22'}}=3\Phi_{1}(x_1/x_0)\cdot\allowbreak\Phi_{1}(z_0/z_3)\cdot\allowbreak
\Phi_{1}(z_0/z_4)\cdot\allowbreak\Phi_{1}(z_1/z_3)\cdot\allowbreak\Phi_{1}
(z_1/z_4)\cdot\allowbreak\Phi_{1}(z_2/z_3)\cdot\allowbreak\Phi_{1}(z_2/z_4)
\cdot\allowbreak\Phi_{1}(x_0z_0z_1/x_1z_3z_4)\cdot\allowbreak\Phi_{1}
(x_0z_0z_2/x_1z_3z_4)\cdot\allowbreak\Phi_{1}(x_0z_1z_2/x_1z_3z_4)\cdot
\allowbreak\Phi_{1}(r/x_1y_0z_3)\cdot\allowbreak\Phi_{1}(r/x_1y_0z_4)\cdot
\allowbreak\Phi_{1}(r/x_1y_1z_3)\cdot\allowbreak\Phi_{1}(r/x_1y_1z_4)\cdot
\allowbreak\Phi_{1}(r/x_1y_2z_3)\cdot\allowbreak\Phi_{1}(r/x_1y_2z_4)\cdot
\allowbreak\Phi_{1}(x_0y_0z_0/r)\cdot\allowbreak\Phi_{1}(x_0y_0z_1/r)\
\cdot\allowbreak\Phi_{1}(x_0y_0z_2/r)\cdot\allowbreak\Phi_{1}(x_0y_1z_0/r)
\cdot\allowbreak\Phi_{1}(x_0y_1z_1/r)\cdot\allowbreak\Phi_{1}(x_0y_1z_2/r)
\cdot\allowbreak\Phi_{1}(x_0y_2z_0/r)\cdot\allowbreak\Phi_{1}(x_0y_2z_1/r)
\cdot\allowbreak\Phi_{1}(x_0y_2z_2/r)\cdot\allowbreak\Phi_{1}(r^2/x_0x_1y_0y_1z_3z_4)
\cdot\allowbreak\Phi_{1}(r^2/x_0x_1y_0y_2z_3z_4)\cdot\allowbreak
\Phi_{1}(r^2/x_0x_1y_1y_2z_3z_4)$\\
where $r=\root 3\of{x_0^2x_1y_0y\ _1y_2z_0z_1z_2}$
\\
\\
$s_{\phi_{4,18}}=-4\Phi_{1}(y_1/y_0)\cdot\allowbreak\Phi_{1}(y_0/y_2)\cdot
\allowbreak\Phi_{1}(z_0/z_4)\cdot\allowbreak\Phi_{1}(z_1/z_4)\cdot\allowbreak
\Phi_{1}(z_2/z_4)\cdot\allowbreak\Phi_{1}(z_3/z_4)\cdot\allowbreak\Phi_{1}
(x_0y_0z_0/r)\cdot\allowbreak\Phi_{1}(x_0y_0z_1/r)\cdot\allowbreak\Phi_{1}
(x_0y_0z_2/r)\cdot\allowbreak\Phi_{1}(x_0y_0z_3/r)\cdot\allowbreak\Phi_{1}(
x_1y_0z_0/r)\cdot\allowbreak\Phi_{1}(x_1y_0z_1/r)\cdot\allowbreak
\Phi_{1}(x_1y_0z_2/r)\cdot\allowbreak\Phi_{1}(x_1y_0z_3/r)\cdot\allowbreak
\Phi_{1}(r/x_0y_1z_4)\cdot\allowbreak\Phi_{1}(r/x_1y_1z_4)\cdot\allowbreak
\Phi_{1}(r/x_0y_2z_4)\cdot\allowbreak\Phi_{1}(r/x_1y_2z_4)\cdot\allowbreak\\
\Phi_{1}(r^2/x_0x_1y_0y_1z_0z_1)\cdot\allowbreak\Phi_{1}(r^2/x_0x_1y_0y_1z_0z_2)
\cdot\allowbreak\Phi_{1}(x_0x_1y_0y_1z_0z_3/r^2)\cdot\allowbreak\Phi_{1}
(x_0x_1y_0y_1z_1z_2/r^2)\cdot\allowbreak\Phi_{1}(r^2/x_0x_1y_0y_1z_1z_3
)\cdot\allowbreak\Phi_{1}(r^2/x_0x_1y_0y_1z_2z_3)\cdot\allowbreak
\Phi_{1}(r^2/x_0x_1y_1y_2z_0z_4)\cdot\allowbreak\Phi_{1}(r^2/x_0x_1y_1y_2z_1z_4)
\cdot\allowbreak\Phi_{1}(r^2/x_0x_1y_1y_2z_2z_4)\cdot\allowbreak\Phi_{1}
(r ^2/x_0x_1y_1y_2z_3z_4)$\\
where $r=\root 4\of{x_0^2x_1^2y_0^2y_1y_2z_0z_1z_2z_3}$
\\
\\
$s_{\phi_{5,16}}=5\Phi_{1}(x_0/x_1)\cdot\allowbreak\Phi_{1}(y_2/y_0)\cdot\allowbreak
\Phi_{1}(y_2/y_1)\cdot\allowbreak\Phi_{1}(x_0y_0z_0/r)\cdot\allowbreak
\Phi_{1}(x_0y_0z_1/r)\cdot\allowbreak\Phi_{1}(x_0y_0z_2/r)\cdot\allowbreak
\Phi_{1}(x_0y_0z_3/r)\cdot\allowbreak\Phi_{1}(x_0y_0z_4/r)\cdot\allowbreak
\ \Phi_{1}(x_0y_1z_0/r)\cdot\allowbreak\Phi_{1}(x_0y_1z_1/r)\cdot
\allowbreak\Phi_{1}(x_0y_1z_2/r)\cdot\allowbreak\Phi_{1}(x_0y_1z_3/r)\cdot
\allowbreak\Phi_{1}(x_0y_1z_4/r)\cdot\allowbreak\Phi_{1}(r/x_1y_2z_0)\cdot
\allowbreak\Phi_{1}(r/x_1y_2z_1)\cdot\allowbreak\Phi_{1}(r/x_1y_2z_2)\cdot
\allowbreak\Phi_{1}(r/x_1y_2z_3)\cdot\allowbreak\Phi_{1}(r/x_1y_2z_4)\cdot
\allowbreak\Phi_{1}(x_0x_1y_0y_1z_0z_1/r^2)\cdot\allowbreak\Phi_{1}
(x_0x_1y_0y_1z_0z_2/r^2)\cdot\allowbreak\\
\Phi_{1}(x_0x_1y_0y_1z_0z_3/r^2)\cdot
\allowbreak\Phi_{1}(x_0x_1y_0y_1z_0z_4/r^2)\cdot\allowbreak\Phi_{1}
(x_0x_1y_0y_1z_1z_2/r^2)\cdot\allowbreak\Phi_{1}(x_0x_1y_0y_1z_1z_3/r^2)
\cdot\allowbreak\Phi_{1}(x_0x_1y_0y_1z_1z_4/r^2)\cdot\allowbreak\Phi_{1}
(x_0x_1y_0y_1z_2z_3/r^2)\cdot\allowbreak\Phi_{1}(x_0x_1y_0y_1z_2z_4/r^2)\cdot\allowbreak
\Phi_{1}(x_0x_1y_0y_1z_3z_4/r^2)$\\
where $r=\root 5\of{x_0^3x_1^2y_0^2y_1^2y _2z_0z_1z_2z_3z_4}$
\\
\\
$s_{\phi_{6,15}}=-6\Phi_{1}(z_0/z_1)\cdot\allowbreak\Phi_{1}(z_0/z_2)\cdot\allowbreak
\Phi_{1}(z_0/z_3)\cdot\allowbreak\Phi_{1}(z_0/z_4)\cdot\allowbreak\Phi_{1}
(r/x_0y_0z_0)\cdot\allowbreak\Phi_{1}(r/x_0y_1z_0)\cdot\allowbreak\Phi_
{1}(r/x_0y_2z_0)\cdot\allowbreak\Phi_{1}(x_1y_0z_0/r)\cdot\allowbreak\Phi_{1}
(x_1y_1z_0/r)\cdot\allowbreak\Phi_{1}(x_1y_2z_0/r)\cdot\allowbreak
\Phi_{1}(x_0x_1y_0y_1z_0z_1/r^2)\cdot\allowbreak\\
\Phi_{1}(x_0x_1y_0y_1z_0z_2/r^2)\cdot\allowbreak\Phi_{1}(x_0x_1y_0y_1z_0z_3/r^2)\cdot\allowbreak
\Phi_{1}(x_0x_1y_0y_1z_0z_4/r^2)\cdot\allowbreak\Phi_{1}(x_0x_1y_0y_2z_0z_1/r^2)
\cdot\allowbreak\Phi_{1}(x_0x_1y_0y_2z_0z_2/r^2)\cdot\allowbreak\Phi_{1}
(x_0x_1y_0y_2z_0z_3/r^2)\cdot\allowbreak\Phi_{1}(x_0x_1y_0y_2z_0z_4/r^2)
\cdot\allowbreak\Phi_{1}(x_0x_1y_1y_2z_0z_1/r^2)\cdot\allowbreak\Phi_{1}
(x_0x_1y_1y_2z_0z_2/r^2)\cdot\allowbreak\Phi_{1}(x_0x_1y_1y_2z_0z_3/r^2)
\cdot\allowbreak\Phi_{1}(x_0x_1y_1y_2z_0z_4/r^2)\cdot\allowbreak\\
\Phi_{1} (x_0^2x_1y_0y_1y_2z_0z_1z_2/r^3)\cdot\allowbreak
\Phi_{1}(x_0^2x_1y_0y_1y_2z_0z_1
z_3/r^3)\cdot\allowbreak\Phi_{1}(x_0^2x_1y_0y_1y_2z_0z_1z_4/r^3)\cdot
\allowbreak\\
\Phi_{1}(x_0^2x_1y_0y_1y_2z_0z_2z_3/r^3)\cdot\allowbreak
\Phi_{1} (x_0^2x_1y_0y_1y_2z_0z_2z_4/r^3)\cdot\allowbreak
\Phi_{1}(x_0^2x_1y_0y_1y_2z_0z _3z_4/r^3)$\\ where
$r=\root6\of{x_0^3x_1^3y_0^2y_1^2y_2^2z_0^2z_1z_2z_3z_4}$}
\\
\\
Following Theorem $\ref{Semisimplicity Malle}$ and \cite{Ma4}, Table
8.1, if we set
$$\begin{array}{cl}
    X_i^{60}:=(\zeta_2)^{-i}x_i & (i=0,1), \\
    Y_j^{60}:=(\zeta_3)^{-j}y_j & (j=0,1,2), \\
    Z_k^{60}:=(\zeta_5)^{-k}z_k & (k=0,1,2,3,4),
  \end{array}$$
then
$\mathbb{Q}(\zeta_{60})(X_0,X_1,Y_0,Y_1,Y_2,Z_0,Z_1,Z_2,Z_3,Z_4)$ is
a splitting field for $\mathcal{H}(G_{19})$. Hence, the factorization
of the Schur elements over that field is as described by Theorem
$\ref{Schur element generic}$.

\section{The groups $G_{25}$, $G_{26}$}

The following table gives the specialization of the parameters of
the generic Hecke algebra $\mathcal{H}(G_{26})$,
$(x_0,x_1;y_0,y_1,y_2)$, which gives the generic Hecke algebra of the
group $G_{25}$ (\cite{Ma5}, Theorem 6.3).

\begin{center}
\begin{tabular}{|c|c|c|c|}
  \hline
  Group & Index & S & T \\
  \hline
  $G_{26}$ & 1  & $x_0,x_1$ & $y_0,y_1,y_2$\\
  $G_{25}$ & 2  & $1,-1$    & $y_0,y_1,y_2$\\
  \hline
\end{tabular}
\end{center}

$$\textrm{Specialization of the parameters for }\mathcal{H}(G_{26})$$

\begin{lemma}\label{g26}\
The algebra $\mathcal{H}(G_{26})$ specialized via
$$(x_0,x_1;y_0,y_1,y_2) \mapsto
(1,-1;y_0,y_1,y_2)$$ is the twisted symmetric algebra of the cyclic
group $C_2$ over the symmetric subalgebra $\mathcal{H}(G_{25})$ with
parameters $(y_0,y_1,y_2)$. 
\end{lemma}
\begin{apod}
{\small We have $$\begin{array}{rccl}
    \mathcal{H}(G_{26}) & = & <S,T,U \,\,| &  STST=TSTS, UTU=TUT, SU=US,  \\
     &  &  & (S-x_0)(S-x_1)=0, \\
     &  &  & (T-y_0)(T-y_1)(T-y_2)=0, \\
     &  &  & (U-y_0)(U-y_1)(U-y_2)=0>.
  \end{array}$$
  Let
$$\begin{array}{rccl}
     A & := & <S,T,U \,\,| & STST=TSTS, UTU=TUT, SU=US, S^2=1,  \\
     &  &  & (T-y_0)(T-y_1)(T-y_2)=0, \\
     &  &  & (U-y_0)(U-y_1)(U-y_2)=0>
  \end{array}$$
   and $$\bar{A}:=<SUS,T,U>.$$ Then
   $$A=\bar{A} \oplus S \bar{A}=\bar{A} \oplus  \bar{A}S\, \textrm{ and } \,\bar{A} \cong \mathcal{H}(G_{25}).$$}
\end{apod}

The Schur elements of all irreducible characters of
$\mathcal{H}(G_{26})$ are calculated in \cite{Ma5}. They are
obtained via Galois transformations from the following ones:
\\
\\
{\footnotesize
$s_{\phi_{1,0}}=-\Phi_{1}(x_0/x_1)\cdot\allowbreak\Phi_{1}(y_0/y_1)\cdot\allowbreak
\Phi_{1}(y_0/y_2)\cdot\allowbreak\Phi_{2}(x_0y_0/x_1y_1)\cdot\allowbreak
\Phi_{2}(x_0y_0/x_1y_2)\cdot\allowbreak\Phi_{1}(x_0y_0^2/x_1y_1^2)\cdot
\allowbreak\Phi_{1}(x_0y_0^2/x_1y_2^2)\cdot\allowbreak\Phi_{2}(x_0y_0^3/x_1y_1^2y_2)
\cdot\allowbreak\Phi_{2}(x_0y_0^3/x_1y_1y_2^2)\cdot\allowbreak
\Phi_{6}(x_0y_0^2/x_1y_1y_2)\cdot\allowbreak\Phi_{2}(y_0^2/y_1y_2)\cdot
\allowbreak\Phi_{6}(y_0/y_1)\cdot\allowbreak\Phi_{6}(y_0/y_2)$
\\
\\
$s_{\phi_{2,3}}=y_1/y_0\Phi_{1}(x_0/x_1)\cdot\allowbreak\Phi_{1}(y_0/y_2)\cdot
\allowbreak\Phi_{1}(y_1/y_2)\cdot\allowbreak\Phi_{1}(x_0y_0/x_1y_2)\cdot
\allowbreak\Phi_{1}(x_0y_1/x_1y_2)\cdot\allowbreak\Phi_{2}(x_0y_0/x_1y_2)\cdot
\allowbreak\Phi_{2}(x_0y_1/x_1y_2)\cdot\allowbreak\Phi_{2}(x_0y_0/x_1y_1)
\cdot\allowbreak\Phi_{2}(x_0y_1/x_1y_0)\cdot\allowbreak\Phi_{6}
(x_0y_0y_1/x_1y_2^2)\cdot\allowbreak\Phi_{2}(y_0y_1/y_2^2)\cdot\allowbreak\Phi_{6}
(y_0/y_1)$
\\
\\
$s_{\phi_{3,6}}=-\Phi_{1}(x_0/x_1)\cdot\allowbreak\Phi_{3}(x_0/x_1)\cdot
\allowbreak\Phi_{2}(x_0y_0/x_1y_1)\cdot\allowbreak\Phi_{2}(x_0y_0/x_1y_2)\cdot
\allowbreak\Phi_{2}(x_0y_1/x_1y_0)\cdot\allowbreak\Phi_{2}(x_0y_1/x_1y_2)\cdot
\allowbreak\Phi_{2}(x_0y_2/x_1y_0)\cdot\allowbreak\Phi_{2}(x_0y_2/x_1y_1)
\cdot\allowbreak\Phi_{2}(y_0y_1/y_2^2)\cdot\allowbreak\Phi_{2}(y_0y_2/y_1^2
)\cdot\allowbreak\Phi_{2}(y_1y_2/y_0^2)$
\\
\\
$s_{\phi_{3,1}}=-\Phi_{1}(x_1/x_0)\cdot\allowbreak\Phi_{1}(y_0/y_1)\cdot
\allowbreak\Phi_{1}(y_0/y_2)\cdot\allowbreak\Phi_{2}(y_0/y_2)\cdot\allowbreak\Phi_{1}
(y_1/y_2)\cdot\allowbreak\Phi_{2}(y_0y_1/y_2^2)\cdot\allowbreak\Phi_{2}(
y_0^2/y_1y_2)\cdot\allowbreak\Phi_{6}(y_0/y_2)\cdot\allowbreak\Phi_{2}
(x_0y_0/x_1y_2)\cdot\allowbreak\Phi_{2}(x_0y_1/x_1y_0)\cdot\allowbreak
\Phi_{1}(x_0y_0^2/x_1y_1^2)\cdot\allowbreak\\
\Phi_{2}(x_0y_0^2y_1/x_1y_2^3)$
\\
\\
$s_{\phi_{6,2}}=\Phi_{1}(x_0/x_1)\cdot\allowbreak\Phi_{1}(y_1/y_0)\cdot\allowbreak
\Phi_{1}(y_0/y_2)\cdot\allowbreak\Phi_{1}(y_1/y_2)\cdot\allowbreak\Phi_{2}(
y_2/y_0)\cdot\allowbreak\Phi_{6}(y_0/y_2)\cdot\allowbreak\Phi_{2}(y_0y_2
/y_1^2)\cdot\allowbreak\Phi_{1}(x_0y_1/x_1y_2)\cdot\allowbreak\Phi_{2}
(x_1y_0/x_0y_2)\cdot\allowbreak\Phi_{2}(x_0y_1/x_1y_2)\cdot\allowbreak
\Phi_{2}(x_0y_0^3/x_1y_1^2y_2)$
\\
\\
$s_{\phi_{8,3}}=2\Phi_{1}(y_0/y_1)\cdot\allowbreak\Phi_{1}(y_0/y_2)\cdot\allowbreak
\Phi_{2}(x_1y_2/x_0y_1)\cdot\allowbreak\Phi_{2}(x_1y_1/x_0y_2)\cdot
\allowbreak\Phi_{2}(ry_0/x_1y_2^2)\cdot\allowbreak\Phi_{2}(ry_0/x_1y_1^2)\cdot
\allowbreak\Phi_{1}(ry_2/x_1y_0y_1)\cdot\allowbreak\Phi_{1}(ry_1/x_1y_0y_2)
\cdot\allowbreak\Phi_{3}(ry_0/x_1y_1y_2)\cdot\allowbreak\Phi_{3}(ry_0/x_0
y_1y_2)$\\
where $r=\root 2\of{-x_0x_1y_1y_2}$
\\
\\
$s_{\phi_{9,7}}=\Phi_{1}(\zeta_3^2)\cdot\allowbreak\Phi_{6}(y_0/y_1)\cdot\allowbreak
\ \Phi_{6}(y_2/y_0)\cdot\allowbreak\Phi_{6}(y_1/y_2)\cdot\allowbreak
\Phi_{2}(\zeta_3x_0y_1y_2/x_1y_0^2)\cdot\allowbreak
\Phi_{2}(\zeta_3x_0y_0y_2/x_1y_1^2)\cdot\allowbreak\Phi_{2}(\zeta_3x_0y_0y_1/x_1y_2^2)\cdot
\allowbreak\Phi_{1}(x_1/x_0)\cdot\allowbreak\Phi_{1}(\zeta_3x_0/x_1)$}
\\
\\
Following Theorem $\ref{Semisimplicity Malle}$ and \cite{Ma4}, Table
8.2, if we set
$$\begin{array}{cl}
    X_i^6:=(\zeta_2)^{-i}x_i & (i=0,1), \\
    Y_j^6:=(\zeta_3)^{-j}y_j & (j=0,1,2),
  \end{array}$$
then $\mathbb{Q}(\zeta_3)(X_1,X_2,Y_1,Y_2,Y_2)$ is a splitting field
for $\mathcal{H}(G_{26})$. Hence, the factorization of the Schur
elements over that field is as described by Theorem
$\ref{Schur element generic}$.

\section{The group $G_{28}$ (``$F_4$'')}

Let $\mathcal{H}(G_{28})$ be the generic Hecke algebra of the real
reflection group $G_{28}$ over the ring
$\mathbb{Z}[x_0^\pm,x_1^\pm,y_0^\pm,y_1^\pm]$. We have
$$\begin{array}{rccl}
     \mathcal{H}(G_{28}) & = & <S_1,S_2,T_1,T_2 \,\,| &S_1S_2S_1=S_2S_1S_2,\,\, T_1T_2T_1=T_2T_1T_2, \\
     &  &  &  S_1T_1=T_1S_1,\,\, S_1T_2=T_2S_1,\,\, S_2T_2=T_2S_2, \\
     &  &  &  S_2T_1S_2T_1=T_1S_2T_1S_2,\\
     &  &  & (S_i-x_0)(S_i-x_1)=(T_i-y_0)(T_i-y_1)=0>.
  \end{array}$$
The Schur elements of all irreducible characters of
$\mathcal{H}(G_{28})$ have been calculated in \cite{Lu79b}. They
are obtained via Galois transformations from the following
ones:
\\
\\
{\footnotesize
$s_{\phi_{1,0}}=\Phi_{1}(y_0/y_1)\cdot\allowbreak\Phi_{6}(y_0/y_1)\cdot\allowbreak
\Phi_{1}(x_0/x_1)\cdot\allowbreak\Phi_{6}(x_0/x_1)\cdot\allowbreak\Phi_{1}(
x_0y_0^2/x_1y_1^2)\cdot\allowbreak
\Phi_{6}(x_0y_0/x_1y_1)\cdot\allowbreak
\Phi_{1}(x_0^2y_0/x_1^2y_1)\cdot\allowbreak\Phi_{4}(x_0y_0/x_1y_1)\cdot
\allowbreak\Phi_{2}(x_0y_0/x_1y_1)\cdot\allowbreak\Phi_{2}(x_0y_0/x_1y_1)$
\\
\\
$s_{\phi_{2,4''}}=-y_1/y_0\Phi_{6}(y_0/y_1)\cdot\allowbreak\Phi_{3}(x_0/x_1)\cdot
\allowbreak\Phi_{6}(x_0/x_1)\cdot\allowbreak\Phi_{1}(x_0/x_1)\cdot\allowbreak\
\Phi_{1}(x_0/x_1)\cdot\allowbreak\Phi_{1}(x_0^2y_0/x_1^2y_1)\cdot
\allowbreak\Phi_{2}(x_0y_0/x_1y_1)\cdot\allowbreak\Phi_{2}(x_0y_1/x_1y_0)\cdot
\allowbreak\Phi_{1}(x_0^2y_1/x_1^2y_0)$
\\
\\
$s_{\phi_{4,8}}=2\Phi_{6}(y_0/y_1)\cdot\allowbreak\Phi_{6}(x_1/x_0)\cdot\allowbreak
\Phi_{2}(x_0y_1/x_1y_0)\cdot\allowbreak\Phi_{2}(x_0y_1/x_1y_0)\cdot
\allowbreak\Phi_{2}(x_1y_1/x_0y_0)\cdot\allowbreak\Phi_{2}(x_0y_0/x_1y_1)$
\\
\\
$s_{\phi_{4,1}}=\Phi_{1}(y_0/y_1)\cdot\allowbreak\Phi_{6}(y_0/y_1)\cdot\allowbreak
\Phi_{1}(x_1/x_0)\cdot\allowbreak\Phi_{6}(x_0/x_1)\cdot\allowbreak\Phi_{2}(
x_0y_1/x_1y_0)\cdot\allowbreak
\Phi_{6}(x_0y_0/x_1y_1)\cdot\allowbreak
\Phi_{2}(x_0y_0/x_1y_1)\cdot\allowbreak\Phi_{2}(x_0y_0/x_1y_1)$
\\
\\
$s_{\phi_{6,6''}}=3\Phi_{1}(y_1/y_0)\cdot\allowbreak\Phi_{1}(y_1/y_0)\cdot\allowbreak
\Phi_{1}(x_1/x_0)\cdot\allowbreak\Phi_{1}(x_1/x_0)\cdot\allowbreak\Phi_{6}
(x_0y_0/x_1y_1)\cdot\allowbreak
\Phi_{2}(x_0y_1/x_1y_0)\cdot\allowbreak \Phi_{2}(x_1y_0/x_0y_1)$
\\
\\
$s_{\phi_{8,3''}}=-y_1/y_0\Phi_{6}(y_0/y_1)\cdot\allowbreak\Phi_{6}(x_0/x_1)\cdot
\allowbreak\Phi_{1}(x_0/x_1)\cdot\allowbreak\Phi_{1}(x_1/x_0)\cdot\allowbreak
\Phi_{3}(x_0/x_1)\cdot\allowbreak\Phi_{1}(x_0y_1^2/x_1y_0^2)\cdot
\allowbreak\Phi_{1}(x_0y_0^2/x_1y_1^2)$
\\
\\
$s_{\phi_{9,2}}=\Phi_{1}(y_0/y_1)\cdot\allowbreak\Phi_{1}(x_0/x_1)\cdot\allowbreak
\Phi_{1}(x_0y_1^2/x_1y_0^2)\cdot\allowbreak\Phi_{4}(x_0y_0/x_1y_1)\cdot
\allowbreak\Phi_{1}(x_1^2y_0/x_0^2y_1)\cdot\allowbreak\Phi_{2}(x_0y_0/x_1y_1)
\cdot\allowbreak\Phi_{2}(x_0y_0/x_1y_1)$
\\
\\
$s_{\phi_{12,4}}=6\Phi_{3}(y_0/y_1)\cdot\allowbreak\Phi_{3}(x_1/x_0)\cdot\allowbreak
\Phi_{2}(x_0y_1/x_1y_0)\cdot\allowbreak\Phi_{2}(x_0y_1/x_1y_0)\cdot
\allowbreak\Phi_{2}(x_0y_0/x_1y_1)\cdot\allowbreak\Phi_{2}(x_1y_1/x_0y_0)$
\\
\\
$s_{\phi_{16,5}}=2x_1y_1/x_0y_0\Phi_{6}(y_0/y_1)\cdot\allowbreak\Phi_{6}(x_1/x_0)\cdot
\allowbreak\Phi_{4}(x_0y_1/x_1y_0)\cdot\allowbreak\Phi_{4}(x_0y_0/x_1y_1)$}
\\
\\
Following Theorem $\ref{Semisimplicity Malle}$, if we set
$$\begin{array}{cl}
    X_i^2:=(\zeta_2)^{-i}x_i & (i=0,1), \\
    Y_j^2:=(\zeta_2)^{-j}y_j & (j=0,1),
  \end{array}$$
then $\mathbb{Q}(X_0,X_1,Y_0,Y_1)$ is a splitting field for
$\mathcal{H}(G_{28})$. Hence, the factorization of the Schur elements
over that field is as described by Theorem $\ref{Schur element
generic}$.

\section{The group $G_{32}$}

Let $\mathcal{H}(G_{32})$ be the generic Hecke algebra of the
complex reflection group $G_{32}$ over the ring
$\mathbb{Z}[x_0^\pm,x_1^\pm,x_2^\pm]$. We have
$$\begin{array}{rccl}
     \mathcal{H}(G_{32}) & = & <S_1,S_2,S_3,S_4 \,\,| &S_iS_{i+1}S_i=S_{i+1}S_iS_{i+1},  \\
     &  &  &  S_iS_j=S_jS_i \textrm{ when }|i-j|>1, \\
     &  &  & (S_i-x_0)(S_i-x_1)(S_i-x_2)=0>.
  \end{array}$$
The Schur elements of all irreducible characters of
$\mathcal{H}(G_{32})$ have been calculated in \cite{Ma5}. They
are obtained via Galois transformations from the following
ones:
\\
\\
{\footnotesize
$s_{\phi_{1,0}}=\Phi_{1}(x_0/x_2)\cdot\allowbreak\Phi_{1}(x_0/x_2)\cdot
\allowbreak\Phi_{1}(x_0/x_1)\cdot\allowbreak\Phi_{1}(x_0/x_1)\cdot\allowbreak
\Phi_{1}(x_0^3/x_1x_2^2)\cdot\allowbreak\Phi_{1}(x_0^3/x_1^2x_2)
\cdot\allowbreak\Phi_{1}(x_0^5/x_1^3x_2^2)\cdot\allowbreak
\Phi_{1}(x_0^5/x_1^2x_2^3)
\cdot\allowbreak\Phi_{2}(x_0^4/x_1x_2^3)\cdot\allowbreak\Phi_{2}(x_0^4/x_1^3x_2)\cdot\allowbreak\Phi_{2}(x_0^2/x_1x_2)\cdot\allowbreak
\Phi_{2}(x_0^2/x_1x_2)\cdot\allowbreak\Phi_{6}(x_0/x_2)\cdot\allowbreak
\Phi_{6}(x_0/x_1)\cdot\allowbreak\Phi_{6}(x_0^3/x_1x_2^2)\cdot\allowbreak
\Phi_{6}(x_0^3/x_1^2x_2)\cdot\allowbreak\Phi_{6}(x_0^2/x_1x_2)\cdot\allowbreak
\Phi_{4}(x_0^2/x_1x_2)\cdot\allowbreak\Phi_{4}(x_0/x_2)\cdot\allowbreak
\Phi_{4}(x_0/x_1)\cdot\allowbreak\Phi_{3}(x_0^2/x_1x_2)\cdot\allowbreak
\Phi_{10}(x_0/x_2)\cdot\allowbreak\Phi_{10}(x_0/x_1)\cdot\allowbreak
\Phi_{5}(x_0^2/x_1x_2)$
\\
\\
$s_{\phi_{4,1}}=\Phi_{1}(x_0^4/x_1x_2^3)\cdot\allowbreak\Phi_{1}(x_0^3/x_1x_2^2)
\cdot\allowbreak\Phi_{1}(x_0^3/x_1^2x_2)\cdot\allowbreak\Phi_{1}(x_0^2x_1
/x_2^3)\cdot\allowbreak\Phi_{1}(x_1/x_0)\cdot\allowbreak\Phi_{1}(x_1/x_2
)\cdot\allowbreak\Phi_{1}(x_0/x_2)\cdot\allowbreak\Phi_{1}(x_0/x_2)\cdot
\allowbreak\Phi_{2}(x_0^5/x_1x_2^4)\cdot\allowbreak\Phi_{2}(x_0^3x_1/x_2^4)
\cdot\allowbreak\Phi_{2}(x_0^3/x_1x_2^2)\cdot\allowbreak\Phi_{2}(x_0^2
/x_1x_2)\cdot\allowbreak\Phi_{2}(x_0/x_2)\cdot\allowbreak\Phi_{2}(x_0x_1
/x_2^2)\cdot\allowbreak\Phi_{6}(x_0/x_2)\cdot\allowbreak\Phi_{6}(x_0/x_1)
\cdot\allowbreak\Phi_{4}(x_0/x_2)\cdot\allowbreak\Phi_{3}(x_0^2/x_1x_2
)\cdot\allowbreak\Phi_{10}(x_0/x_1)\cdot\allowbreak\Phi_{15}(x_0/x_2)$
\\
\\
$s_{\phi_{5,4}}=\Phi_{1}(x_0^3x_1^2/x_2^5)\cdot\allowbreak\Phi_{1}(x_0^2x_1/x_2^3
)\cdot\allowbreak\Phi_{1}(x_0/x_2)\cdot\allowbreak\Phi_{1}(x_0/x_2)\cdot
\allowbreak\Phi_{1}(x_1/x_0)\cdot\allowbreak\Phi_{1}(x_1/x_0)\cdot\allowbreak
\Phi_{1}(x_1/x_2)\cdot\allowbreak\Phi_{1}(x_1/x_2)\cdot\allowbreak
\Phi_{2}(x_0^3/x_1x_2^2)\cdot\allowbreak\Phi_{2}(x_0x_1^2/x_2^3)\cdot\allowbreak\Phi_{2}(x_1/x_2)\cdot\allowbreak\Phi_{2}(x_0^2/x_1x_2)\cdot
\allowbreak\Phi_{2}(x_0x_1/x_2^2)\cdot\allowbreak\Phi_{2}(x_0^4x_1/x_2^5)\cdot\allowbreak\Phi_{2}(x_0/x_2)\cdot\allowbreak\Phi_{2}(x_0/x_2)\cdot
\allowbreak\Phi_{6}(x_0/x_2)\cdot\allowbreak\Phi_{6}(x_0/x_2)\cdot\allowbreak\Phi_{6}(x_0/x_1)\cdot\allowbreak\Phi_{4}(x_0/x_1)\cdot\allowbreak\Phi_{3}(x_0x_1/x_2^2)\cdot\allowbreak\Phi_{12}(x_0/x_2)$
\\
\\
$s_{\phi_{6,8}}=x_1^2/x_0^2\Phi_{1}(x_0/x_1)\cdot\allowbreak\Phi_{1}(x_1/x_0)\cdot\allowbreak\Phi_{1}(x_0/x_2)\cdot\allowbreak
\Phi_{1}(x_1/x_2)\cdot\allowbreak\Phi_{1}(x_0/x_2)\cdot\allowbreak\Phi_{1}(x_1/x_2)\cdot\allowbreak\Phi_{1}(x_0x_1^2/x_2^3)\cdot\allowbreak
\Phi_{1}(x_0^2x_1/x_2^3)\cdot\allowbreak\Phi_{2}(x_0x_1^2/x_2^3)\cdot\allowbreak\Phi_{2}(x_0^2x_1/x_2^3)\cdot\allowbreak\Phi_{2}(x_0^2/x_1x_2)\cdot
\allowbreak\Phi_{2}(x_1^2/x_0x_2)\cdot\allowbreak\\
\Phi_{2}(x_0x_1/x_2^2)\cdot\allowbreak\Phi_{2}(x_0x_1/x_2^2)\cdot\allowbreak\Phi_{2}(x_1/x_2)\cdot
\allowbreak\Phi_{2}(x_0/x_2)\cdot\allowbreak\Phi_{6}(x_0x_1/x_2^2)\cdot
\allowbreak\Phi_{6}(x_0/x_2)\cdot\allowbreak\Phi_{6}(x_1/x_2)\cdot\allowbreak\Phi_{10}(x_0/x_1)\cdot\allowbreak\Phi_{5}(x_0x_1/x_2^2)$
\\
\\
$s_{\phi_{10,2}}=\Phi_{1}(x_0^2x_1/x_2^3)\cdot\allowbreak\Phi_{1}
(x_0^3/x_1x_2^2)\cdot\allowbreak\Phi_{1}(x_1/x_2)\cdot\allowbreak\Phi_{1}
(x_1/x_2)\cdot\allowbreak\Phi_{1}(x_1/x_0)\cdot\allowbreak\Phi_{1}(x_2/x_0)\cdot\allowbreak\Phi_{1}(x_0/x_2)\cdot\allowbreak\Phi_{1}(x_0/x_2)
\cdot\allowbreak\Phi_{2}(x_0^4/x_1^3x_2)\cdot\allowbreak\Phi_{2}(x_0^3/x_1x_2^2)\cdot\allowbreak\Phi_{2}(x_0x_1/x_2^2)\cdot\allowbreak\Phi_{2}
(x_0x_2/x_1^2)\cdot\allowbreak\Phi_{2}(x_0/x_2)\cdot\allowbreak\Phi_{2}(x_0/x_2)\cdot\allowbreak\Phi_{2}(x_1/x_2)\cdot\allowbreak\Phi_{2}(x_0^2/
x_1x_2)\cdot\allowbreak\Phi_{6}(x_0^2x_1/x_2^3)\cdot\allowbreak\Phi_{6}(x_0/x_2)\cdot\allowbreak\Phi_{6}(x_0/x_1)\cdot\allowbreak\Phi_{4}(x_0/x_2)
\cdot\allowbreak\Phi_{3}(x_0^2/x_1x_2)$
\\
\\
$s_{\phi_{15,6}}=\Phi_{1}(x_0^3/x_1x_2^2)\cdot\allowbreak\Phi_{1}(x_0^3/x_1^2x_2)\cdot\allowbreak\Phi_{1}(x_0/x_1)\cdot\allowbreak\Phi_{1}(x_0/x_1)
\cdot\allowbreak\Phi_{1}(x_2/x_1)\cdot\allowbreak\Phi_{1}(x_2/x_1)\cdot\allowbreak\Phi_{1}(x_2/x_0)\cdot\allowbreak\Phi_{1}(x_0/x_2)\cdot
\allowbreak\Phi_{2}(x_0^3x_2/x_1^4)\cdot\allowbreak\Phi_{2}(x_0^3x_1/x_2^4)\cdot\allowbreak\Phi_{2}(x_1^2/x_0x_2)\cdot\allowbreak\Phi_{2}(x_0x_1/x_2^2)
\cdot\allowbreak\Phi_{2}(x_0/x_2)\cdot\allowbreak\Phi_{2}(x_0/x_1)\cdot\allowbreak\Phi_{2}(x_0^2/x_1x_2)\cdot\allowbreak\Phi_{2}(x_0^2/x_1x_2)
\cdot\allowbreak\Phi_{6}(x_0^2/x_1x_2)\cdot\allowbreak\Phi_{6}(x_0/x_1)\cdot\allowbreak\Phi_{6}(x_0/x_2)\cdot\allowbreak\Phi_{4}(x_0^2/x_1x_2)$
\\
\\
$s_{\phi_{15,8}}=\Phi_{1}(x_1^2x_2/x_0^3)\cdot\allowbreak\Phi_{1}(x_0^2x_2/x_1^3)\cdot\allowbreak\Phi_{1}(x_0/x_2)\cdot\allowbreak\Phi_{1}(x_0/
x_2)\cdot\allowbreak\Phi_{1}(x_1/x_2)\cdot\allowbreak\Phi_{1}(x_1/x_2)\cdot\allowbreak\Phi_{1}(x_1/x_0)\cdot\allowbreak\Phi_{1}(x_0/x_1)\cdot
\allowbreak\Phi_{2}(x_1x_2/x_0^2)\cdot\allowbreak\Phi_{2}(x_0x_2/x_1^2)\cdot\allowbreak\Phi_{2}(x_0x_1/x_2^2)\cdot\allowbreak\Phi_{2}(x_0x_1/x_2^2)
\cdot\allowbreak\Phi_{2}(x_1/x_2)\cdot\allowbreak\Phi_{2}(x_1/x_2)\cdot\allowbreak\Phi_{2}(x_0/x_2)\cdot\allowbreak\Phi_{2}(x_0/x_2)\cdot
\allowbreak\Phi_{6}(x_0/x_2)\cdot\allowbreak\Phi_{6}(x_1/x_2)\cdot\allowbreak
\Phi_{4}(x_0x_1/x_2^2)$
\\
\\
$s_{\phi_{20,3}}=\Phi_{1}(x_0^2x_2/x_1^3)\cdot\allowbreak\Phi_{1}(x_0/x_1)\cdot\allowbreak\Phi_{1}(x_0/x_1)\cdot
\allowbreak\Phi_{1}(x_2/x_0)\cdot\allowbreak\Phi_{1}(x_0/x_2)\cdot\allowbreak\Phi_{1}(x_0^4/x_1x_2^3)\cdot\allowbreak\Phi_{1}(x_0^3/x_1x_2^2)\cdot
\allowbreak\Phi_{1}(x_1/x_2)\cdot\allowbreak\Phi_{2}(x_0x_1^2/x_2^3)\cdot\allowbreak\Phi_{2}(x_1^2/x_0x_2)\cdot\allowbreak\Phi_{2}(x_0^2/x_1x_2)
\cdot\allowbreak\Phi_{2}(x_0x_1/x_2^2)\cdot\allowbreak\Phi_{2}(x_0/x_2)\cdot\allowbreak\Phi_{2}(x_2/x_0)\cdot\allowbreak\Phi_{6}(x_0^3/x_1^2x_2)
\cdot\allowbreak\Phi_{6}(x_0/x_2)\cdot\allowbreak\Phi_{4}(x_0/x_2)\cdot\allowbreak\Phi_{3}(x_0x_1/x_2^2)$
\\
\\
$s_{\phi_{20,5}}=-\Phi_{1}(x_1^3/x_0^2x_2)\cdot\allowbreak\Phi_{1}(x_0^3/x_1^2x_2)\cdot\allowbreak\Phi_{1}(x_0x_1^2/x_2^3)\cdot
\allowbreak\Phi_{1}(x_0^2x_1/x_2^3)\cdot\allowbreak\Phi_{1}(x_2/x_1)\cdot\allowbreak\Phi_{1}(x_2/x_1)\cdot\allowbreak\Phi_{1}(x_0/x_2)\cdot\allowbreak
\Phi_{1}(x_0/x_2)\cdot\allowbreak\Phi_{2}(x_1^3/x_0x_2^2)\cdot\allowbreak\Phi_{2}(x_0^3/x_1x_2^2)\cdot\allowbreak\Phi_{2}(x_0x_1/x_2^2)\cdot
\allowbreak\Phi_{2}(x_0x_1/x_2^2)\cdot\allowbreak\Phi_{6}(x_0x_1/x_2^2)\cdot\allowbreak\Phi_{6}(x_1/x_0)\cdot\allowbreak\Phi_{6}(x_1/x_2)\cdot
\allowbreak\Phi_{6}(x_0/x_2)\cdot\allowbreak\Phi_{3}(x_0x_1/x_2^2)$
\\
\\
$s_{\phi_{20,7}}=\Phi_{1}(x_0^3x_1/x_2^4)\cdot\allowbreak\Phi_{1}(x_0x_1^2/x_2^3)\cdot\allowbreak\Phi_{1}(x_1/x_0)
\cdot\allowbreak\Phi_{1}(x_1/x_0)\cdot\allowbreak\Phi_{1}(x_2/x_0)\cdot\allowbreak\Phi_{1}(x_0/x_2)\cdot\allowbreak\Phi_{1}(x_1/x_2)\cdot
\allowbreak\Phi_{1}(x_1/x_2)\cdot\allowbreak\Phi_{2}(x_0/x_2)\cdot\allowbreak\Phi_{2}(x_2/x_0)\cdot\allowbreak\Phi_{2}(x_0^3x_2/x_1^4)\cdot
\allowbreak\Phi_{2}(x_0x_1^2/x_2^3)\cdot\allowbreak\Phi_{2}(x_0x_1/x_2^2)\cdot\allowbreak\Phi_{2}(x_1^2/x_0x_2)\cdot\allowbreak\Phi_{6}(x_0/x_2)\cdot
\allowbreak\Phi_{6}(x_0/x_2)\cdot\allowbreak\Phi_{6}(x_1/x_2)\cdot\allowbreak\Phi_{3}(x_0^2/x_1x_2)$
\\
\\
$s_{\phi_{20,12}}=2\Phi_{1}(x_2/x_1)\cdot\allowbreak\Phi_{1}(x_1/x_2)\cdot\allowbreak\Phi_{1}(x_2/x_0)\cdot\allowbreak\Phi_{1}(x_1/x_0)\cdot\allowbreak\Phi_{1}(x_2/x_0)\
\cdot\allowbreak\Phi_{1}(x_1/x_0)\cdot\allowbreak\Phi_{1}(x_0/x_2)\cdot\allowbreak\Phi_{1}(x_0/x_1)\cdot\allowbreak\Phi_{2}(x_0x_1^2/x_2^3)\cdot
\allowbreak\Phi_{2}(x_0x_2^2/x_1^3)\cdot\allowbreak\Phi_{2}(x_0^2/x_1x_2)\cdot\allowbreak\Phi_{2}(x_1x_2/x_0^2)\cdot\allowbreak\Phi_{2}(x_0/x_1)\
\cdot\allowbreak\Phi_{2}(x_0/x_2)\cdot\allowbreak\Phi_{2}(x_0/x_1)\cdot\allowbreak\Phi_{2}(x_0/x_2)\cdot\allowbreak\Phi_{6}(x_2/x_1)\cdot\allowbreak\Phi_{6}(x_1/x_2)\cdot\allowbreak\Phi_{3}(x_0^2/x_1x_2)$
\\
\\
$s_{\phi_{24,6}}=\Phi_{1}(x_1^3/x_0^2x_2)
\cdot\allowbreak\Phi_{1}(x_2^3/x_0^2x_1)\cdot\allowbreak\Phi_{1}(x_0/x_1)\cdot\allowbreak\Phi_{1}(x_0/x_2)\cdot\allowbreak\Phi_{1}(x_0/x_1)\cdot\
\allowbreak\Phi_{1}(x_0/x_2)\cdot\allowbreak\Phi_{1}(x_2/x_1)\cdot\allowbreak\Phi_{1}(x_1/x_2)\cdot\allowbreak\Phi_{2}(x_0/x_1)\cdot\allowbreak\
\Phi_{2}(x_0/x_2)\cdot\allowbreak\Phi_{2}(x_0x_2/x_1^2)\cdot\allowbreak\Phi_{2}(x_0x_1/x_2^2)\cdot\allowbreak\Phi_{6}(x_0/x_1)\cdot\allowbreak
\Phi_{6}(x_0/x_2)\cdot\allowbreak\Phi_{4}(x_0/x_1)\cdot\allowbreak\Phi_{4}(x_0/x_2)\cdot\allowbreak\Phi_{5}(x_0^2/x_1x_2)$
\\
\\
$s_{\phi_{30,4}}=\Phi_{1}(x_0^5/x_1^3x_2^2)\cdot\allowbreak\Phi_{1}(x_0/x_2)\cdot\allowbreak\Phi_{1}(x_1/x_0)\cdot\allowbreak\Phi_{1}(x_1/x_0)
\cdot\allowbreak\Phi_{1}(x_1/x_2)\cdot\allowbreak\Phi_{1}(x_1/x_2)\cdot\allowbreak\Phi_{1}(x_0/x_2)\cdot\allowbreak\Phi_{1}(x_0/x_2)\cdot
\allowbreak\Phi_{2}(x_1/x_0)\cdot\allowbreak\Phi_{2}(x_1/x_2)\cdot\allowbreak\Phi_{2}(x_2/x_0)\cdot\allowbreak\Phi_{2}(x_0^5/x_1x_2^4)\cdot\allowbreak
\Phi_{2}(x_0x_2^2/x_1^3)\cdot\allowbreak\Phi_{2}(x_0/x_2)\cdot\allowbreak\Phi_{2}(x_0x_2/x_1^2)\cdot\allowbreak\Phi_{2}(x_0^2/x_1x_2)\cdot
\allowbreak\Phi_{6}(x_0/x_1)\cdot\allowbreak\Phi_{6}(x_1/x_2)\cdot\allowbreak\Phi_{6}(x_0/x_2)\cdot\allowbreak\Phi_{4}(x_0/x_2)$
\\
\\
$s_{\phi_{30,12'}}=\Phi_{1}(x_1^5/x_0^3x_2^2)\cdot\allowbreak\Phi_{1}(x_1/x_2)\cdot
\allowbreak\Phi_{1}(x_0/x_1)\cdot\allowbreak\Phi_{1}(x_0/x_1)\cdot\allowbreak\Phi_{1}(x_0/x_2)\cdot\allowbreak\Phi_{1}(x_0/x_2)\cdot\allowbreak\Phi_{1}(x_1/x_2)\cdot\allowbreak\Phi_{1}(x_1/x_2)\cdot
\allowbreak\Phi_{2}(x_0/x_1)\cdot\allowbreak\Phi_{2}(x_0/x_2)\cdot\allowbreak\Phi_{2}(x_2/x_1)\cdot\allowbreak\Phi_{2}(x_1^5/x_0x_2^4)\cdot
\allowbreak\Phi_{2}(x_1x_2^2/x_0^3)\cdot\allowbreak\Phi_{2}(x_1/x_2)\cdot\allowbreak\Phi_{2}(x_1x_2/x_0^2)\cdot\allowbreak\Phi_{2}(x_1^2/x_0x_2)\cdot\allowbreak\Phi_{6}(x_1/x_0)\cdot\allowbreak\Phi_{6}(x_0/x_2)\cdot\allowbreak\Phi_{6}(x_1/x_2)\cdot\allowbreak\Phi_{4}(x_1/x_2)$
\\
\\
$s_{\phi_{36,5}}=\Phi_{1}(1/\zeta_3)\cdot\allowbreak\Phi_{1}(x_0/x_2)\cdot\allowbreak\Phi_{1}(x_1/x_0)\cdot\allowbreak\Phi_{1}(\zeta_3
x_0^2/x_1x_2)\cdot\allowbreak\Phi_{1}(\zeta_3^2x_0x_2/x_1^2)\cdot\allowbreak\Phi_{1}(x_2^2/\zeta_3^2x_0x_1)\cdot\allowbreak\Phi_{2}(x_1x_2/x_0^2)
\cdot\allowbreak\Phi_{2}(x_0^2x_2/\zeta_3x_1^3)\cdot\allowbreak\Phi_{2}(\zeta_3^2x_0^2x_1/x_2^3)\cdot\allowbreak\Phi_{6}(x_0^2/x_1x_2)\cdot
\allowbreak\Phi_{6}(x_0/x_1)\cdot\allowbreak\Phi_{6}(x_0/x_2)\cdot\allowbreak\Phi_{6}(x_1/x_2)\cdot\allowbreak\Phi_{5}(\zeta_3x_0/x_2)\cdot\allowbreak
\Phi_{5}(\zeta_3x_0/x_1)$
\\
\\
$s_{\phi_{40,8}}=\Phi_{1}(x_0^3x_1^2/x_2^5)\cdot\allowbreak\Phi_{1}(x_1/x_0)\cdot\allowbreak\Phi_{1}(x_0/x_1)\cdot\allowbreak\Phi_{1}(x_0/x_2)\cdot
\allowbreak\Phi_{1}(x_2/x_0)\cdot\allowbreak\Phi_{1}(x_0/x_2)\cdot\allowbreak\Phi_{1}(x_1^2x_2/x_0^3)\cdot\allowbreak\Phi_{1}(x_2/x_1)\cdot\allowbreak
\Phi_{2}(x_0/x_1)\cdot\allowbreak\Phi_{2}(x_0x_1/x_2^2)\cdot\allowbreak\Phi_{2}(x_0^2/x_1x_2)\cdot\allowbreak\Phi_{2}(x_0/x_2)\cdot\allowbreak\
\Phi_{6}(x_0/x_1)\cdot\allowbreak\Phi_{6}(x_1/x_2)\cdot\allowbreak\Phi_{4}(x_0/x_1)\cdot\allowbreak\Phi_{4}(x_0/x_2)\cdot\allowbreak\Phi_{3}(x_0x_2/x_1^2)$
\\
\\
$s_{\phi_{45,6}}=\Phi_{1}(\zeta_3)\cdot\allowbreak\Phi_{1}(\zeta_3^2x_0^2/x_1x_2)\cdot\allowbreak\Phi_{1}(\zeta_3x_0x_2/x_1^2)
\cdot\allowbreak\Phi_{1}(\zeta_3x_0x_1/x_2^2)\cdot\allowbreak\Phi_{1}(x_2/x_0)\cdot\allowbreak\Phi_{1}(x_1/x_2)\cdot\allowbreak\Phi_{2}(\zeta_3^2x_0^2x_2/x_1^3)\cdot\allowbreak\Phi_{2}(\zeta_3^2x_1^2x_2/x_0^3)
\cdot\allowbreak\Phi_{2}(x_0/\zeta_3^2x_2)\cdot\allowbreak\Phi_{2}(\zeta_3x_1/
x_2)\cdot\allowbreak\Phi_{2}(x_0x_1/x_2^2)\cdot\allowbreak\Phi_{6}(x_1/x_0)\cdot\allowbreak\Phi_{6}(x_1/x_2)\cdot\allowbreak\Phi_{6}(x_0/x_2)
\cdot\allowbreak\Phi_{6}(x_0x_1/x_2^2)\cdot\allowbreak\Phi_{4}(\zeta_3x_0/x_2)\cdot\allowbreak\Phi_{4}(\zeta_3x_1/x_2)$
\\
\\
$s_{\phi_{60,7}}=\Phi_{1}(x_0x_1^2/x_2^3)\cdot\allowbreak\Phi_{1}(x_0/x_1)\cdot\allowbreak\Phi_{1}(x_0/x_1)\cdot\allowbreak\Phi_{1}(x_1/x_2)\cdot\allowbreak\Phi_{1}(x_1/x_2)
\cdot\allowbreak\Phi_{1}(x_0/x_2)\cdot\allowbreak\Phi_{1}(x_0/x_2)\cdot\allowbreak\Phi_{1}(x_1x_2^2/x_0^3)\cdot\allowbreak\Phi_{2}(x_0/x_1)\cdot\allowbreak\Phi_{2}(x_0^4x_1/x_2^5)\cdot\allowbreak\Phi_{2}(x_1^2/x_0x_2)\cdot\allowbreak\Phi_{2}(x_1x_2/x_0^2)\cdot\allowbreak\Phi_{2}(x_1/x_2)\cdot\allowbreak
\Phi_{2}(x_1/x_2)\cdot\allowbreak\Phi_{6}(x_0/x_1)\cdot\allowbreak\Phi_{6}(x_2/x_1)\cdot\allowbreak\Phi_{4}(x_0/x_1)$
\\
\\
$s_{\phi_{60,11''}}=\Phi_{1}(x_1x_2^3/x_0^4)\cdot\allowbreak\Phi_{1}(x_0/x_1)\cdot\allowbreak\Phi_{1}(x_0/x_1)\cdot\allowbreak\Phi_{1}(x_1/x_2)\cdot\allowbreak\Phi_{1}(x_1/x_2)\cdot\allowbreak\Phi_{1}(x_2/x_0)\cdot\allowbreak\Phi_{1}(x_2/x_0)\cdot\allowbreak\Phi_{1}(x_0x_2^3/x_1^4)\cdot\allowbreak\Phi_{2}
(x_1^2/x_0x_2)\cdot\allowbreak\Phi_{2}(x_0x_1/x_2^2)\cdot\allowbreak\Phi_{2}(x_0x_1/x_2^2)\cdot\allowbreak\Phi_{2}(x_0^2/x_1x_2)\cdot\allowbreak\Phi_{2}(x_0/x_1)\cdot\allowbreak\Phi_{2}(x_0/x_1)\cdot\allowbreak\Phi_{6}(x_0x_1/x_2^2)\cdot\allowbreak\Phi_{6}(x_1/x_0)$
\\
\\
$s_{\phi_{60,12}}=2\Phi_{1}(x_0x_2^2/x_1^3)\cdot\allowbreak\Phi_{1}(x_0x_1^2/x_2^3)\cdot\allowbreak\Phi_{1}(x_1/x_0)\cdot\allowbreak\Phi_{1}(x_1/x_0)\cdot\allowbreak\Phi_{1}(x_1/x_0)\cdot\allowbreak\Phi_{1}(x_2/x_0)\cdot\allowbreak\Phi_{1}(x_2/x_0)\cdot\allowbreak\Phi_{1}(x_2/x_0)\cdot\allowbreak\Phi_{2}(x_2/x_1)\cdot\allowbreak\Phi_{2}
(x_1/x_2)\cdot\allowbreak\Phi_{2}(x_0^2/x_1x_2)\cdot\allowbreak\Phi_{2}
(x_0^2/x_1x_2)\cdot\allowbreak\Phi_{6}(x_0^2/x_1x_2)\cdot\allowbreak\Phi_{6}(x_1/x_2)\cdot\allowbreak\Phi_{6}(x_2/x_1)\cdot\allowbreak\Phi_{4}
(x_0/x_1)\cdot\allowbreak\Phi_{4}(x_0/x_2)$
\\
\\
$s_{\phi_{64,8}}=2\Phi_{1}(rx_1/x_2^2)\cdot
\allowbreak\Phi_{1}(x_2^2/rx_0)\cdot\allowbreak\Phi_{1}(x_0/x_2)\cdot\allowbreak\Phi_{1}(x_2/x_1)\cdot\allowbreak\Phi_{1}(x_0/x_1)\cdot\allowbreak
\Phi_{1}(x_1/x_0)\cdot\allowbreak\Phi_{1}(x_0^3/x_1x_2^2)\cdot\allowbreak\Phi_{1}(x_0x_2^2/x_1^3)\cdot\allowbreak\Phi_{2}(rx_0^2/x_1^2x_2)\cdot
\allowbreak\Phi_{2}(rx_1^2/x_0^2x_2)\cdot\allowbreak\Phi_{3}(x_0x_1/x_2^2)\cdot\allowbreak\Phi_{10}(r/x_2)\cdot\allowbreak\Phi_{15}(r/x_0)$\\
where $r=\root 2\of{x_0x_1}$
\\
\\
$s_{\phi_{80,9}}=2\Phi_{1}(x_0x_2^2/x_1^3)\cdot\allowbreak
\Phi_{1}(x_0x_1^2/x_2^3)\cdot\allowbreak\Phi_{1}(x_0/x_1)
\cdot\allowbreak\Phi_{1}(x_0/x_1)\cdot\allowbreak\Phi_{1}(x_2/x_0)\cdot\allowbreak\Phi_{1}(x_2/x_0)\cdot\allowbreak\Phi_{1}(x_2/x_1)\cdot
\allowbreak\Phi_{1}(x_2/x_1)\cdot\allowbreak\Phi_{2}(x_0/x_2)\cdot\allowbreak\Phi_{2}(x_0/x_1)\cdot\allowbreak\Phi_{4}(x_1x_2/x_0^2)\cdot\allowbreak
\Phi_{4}(x_0/x_1)\cdot\allowbreak\Phi_{4}(x_0/x_2)\cdot\allowbreak\Phi_{3}(x_0^2/x_1x_2)\cdot\allowbreak\Phi_{12}(x_1/x_2)$
\\
\\
$s_{\phi_{81,10}}=3\Phi_{2}(rx_2/x_0^2)\cdot\allowbreak\Phi_{2}(rx_2/x_1^2)\cdot\allowbreak\Phi_{2}(rx_0/x_2^2)\cdot\allowbreak\Phi_{2}
(rx_0/x_1^2)\cdot\allowbreak\Phi_{2}(rx_1/x_0^2)\cdot\allowbreak\Phi_{2}(rx_1/x_2^2)\cdot\allowbreak\Phi_{2}(x_0x_1/x_2^2)\cdot\allowbreak\Phi_{2}(x_0x_2/
x_1^2)\cdot\allowbreak\Phi_{2}(x_1x_2/x_0^2)\cdot\allowbreak\Phi_{2}(r/x_2)
\cdot\allowbreak\Phi_{2}(r/x_0)\cdot\allowbreak\Phi_{2}(r/x_1)\cdot
\allowbreak\Phi_{4}(r^2/x_0x_1)\cdot\allowbreak\Phi_{4}(r^2/x_0x_2)\cdot
\allowbreak\Phi_{4}(r^2/x_1x_2)\cdot\allowbreak\Phi_{5}(r/x_0)\cdot\allowbreak
\Phi_{5}(r/x_2)\cdot\allowbreak\Phi_{5}(r/x_1)$\\
where $r=\root 3\of{x_0x_1x_2}$}
\\
\\
Following Theorem $\ref{Semisimplicity Malle}$ and \cite{Ma4}, Table
8.2, if we set
$$X_i^6:=(\zeta_3)^{-i}x_i  \,\,\,(i=0,1,2),$$
then $\mathbb{Q}(\zeta_3)(X_0,X_1,X_2)$ is a splitting field for
$\mathcal{H}(G_{32})$. Hence, the factorization of the Schur elements
over that field is as described by Theorem $\ref{Schur element
generic}$.

\section{The groups $G(de,e,r)$}

The generic Hecke algebras of the groups $G(de,e,r)$ are presented in Chapter $5$. Here we will only give some applications of Clifford theory and a description of their Schur elements.

\subsection{The groups $G(de,e,r)$, $r>2$}

Proposition 1.6 of \cite{Ariki} yields the specialization of the parameters of 
the generic Hecke algebra $\mathcal{H}(G(de,1,r))$,
$(x_0,x_1;u_0,u_1,\ldots,u_{de-1})$, which gives the generic
Hecke algebra of the group $G(de,e,r)$.

\begin{lemma}\label{gdeer}
The algebra $\mathcal{H}(G(de,1,r))$ specialized via
 $$\left\{ 
\begin{array}{ll} 
x_i \mapsto x_i &(0 \leq i \leq 1),\\
u_k \mapsto \zeta_{e}^{[k/d]} v_{k \,\mathrm{mod}\,d}^{1/e} &(0 \leq k \leq de-1)  
\end{array} \right. 
$$
 is the twisted symmetric algebra of the
 cyclic group $C_e$ over the symmetric subalgebra $\mathcal{H}(G(de,e,r))$ with
 parameters $(x_0,x_1;v_0,v_1,\ldots,v_{d-1})$. 
\end{lemma}
\begin{apod}{\small 
The algebra $\mathcal{H}(G(de,1,r))$ is generated 
by the elements $\mathrm{\textbf{s}},\mathrm{\textbf{t}}_1,\mathrm{\textbf{t}}_2,\ldots,\mathrm{\textbf{t}}_{r-1}$ satisfying the relations
\begin{itemize}
\item $\mathrm{\textbf{s}}\mathrm{\textbf{t}}_1\mathrm{\textbf{s}}\mathrm{\textbf{t}}_1=\mathrm{\textbf{t}}_1\mathrm{\textbf{s}}\mathrm{\textbf{t}}_1\mathrm{\textbf{s}}$, $\mathrm{\textbf{s}}\mathrm{\textbf{t}}_j=\mathrm{\textbf{t}}_j\mathrm{\textbf{s}} \textrm{ for } j\neq 1$,
\item $\mathrm{\textbf{t}}_j\mathrm{\textbf{t}}_{j+1}\mathrm{\textbf{t}}_j=\mathrm{\textbf{t}}_{j+1}\mathrm{\textbf{t}}_j\mathrm{\textbf{t}}_{j+1}$,  $ \mathrm{\textbf{t}}_i\mathrm{\textbf{t}}_j=\mathrm{\textbf{t}}_j\mathrm{\textbf{t}}_i \textrm{ for } |i-j|>1$,
\item $(\mathrm{\textbf{s}}-u_0)(\mathrm{\textbf{s}}-u_1)\ldots(\mathrm{\textbf{s}}-u_{de-1})=(\mathrm{\textbf{t}}_j-x_0)(\mathrm{\textbf{t}}_j-x_1)=0$.
\end{itemize}
Let $A$ be the algebra obtained from $\mathcal{H}(G(de,1,r))$ via the given specialization, \ie the algebra generated by the elements $\mathrm{\textbf{s}},\mathrm{\textbf{t}}_1,\mathrm{\textbf{t}}_2,\ldots,\mathrm{\textbf{t}}_{r-1}$ satisfying the same braid relations as above, as well as:
\begin{center}
$(\mathrm{\textbf{s}}^e-v_0)(\mathrm{\textbf{s}}^e-v_1)\ldots(\mathrm{\textbf{s}}^e-v_{d-1})=(\mathrm{\textbf{t}}_j-x_0)(\mathrm{\textbf{t}}_j-x_1)=0$.
\end{center}
If  $\bar{A}:=<\mathrm{\textbf{s}}^e,\tilde{\textbf{t}}_1:=\textbf{s}^{-1}\textbf{t}_1\textbf{s},\mathrm{\textbf{t}}_1,\mathrm{\textbf{t}}_2,\ldots,\mathrm{\textbf{t}}_{r-1}>$, then
   $$A=\bigoplus_{i=0}^{e-1} \mathrm{\textbf{s}}^i \bar{A}=\bigoplus_{i=0}^{e-1} \bar{A}\mathrm{\textbf{s}}^i\,  \textrm{ and }\, \bar{A} \cong \mathcal{H}(G(de,e,r)).$$}
\end{apod}

For $x_1=-1$, the algebra $\mathcal{H}(G(de,1,r))$ becomes the generic Ariki-Koike algebra $\mathcal{H}_{de,r}$ associated to $G(de,1,r)$. Set $n:=de$ and $x:=x_0$. The following result, which has been obtained independently by Geck, Iancu and Malle (\cite{GIM}) and by Mathas
(\cite{Mat}),
 gives a description of the Schur elements of $\mathcal{H}_{n,r}$. Recall that the irreducible characters of $G(n,1,r)$ are parametrized by the $n$-partitions of $r$.

\begin{theorem}\label{schur elements of ArikiKoike}
Let $\el$ be an $n$-partition of $r$ with ordinary standard symbol
$B_\el=(B_\el^{(0)}, B_\el^{(1)}, \ldots, B_\el^{(n-1)})$. Fix $L \geq h_\el$, where $h_\el$ is the height of $\el$.
We set $B_{\el,L}:=(B_\el^{(0)}[L-h_\el], B_\el^{(1)}[L-h_\el], \ldots, B_\el^{(n-1)}[L-h_\el])=(B_{\el,L}^{(0)}, B^{(1)}_{\el,L}, \ldots, B^{(n-1)}_{\el,L})$ and
$B_{\el,L}^{(s)}=(b_1^{(s)},b_2^{(s)},\ldots,
b_{L}^{(s)})$.  Let $a_{L}:=r(n-1)+\binom{ n}{ 2}\binom{ L}{ 2}$ and $b_{L}:=nL(L-1)(2nL-n-3)/12$. Then the Schur element of the irreducible character $\chi_\el$ is given by the formulae $s_\el=(-1)^{a_{L}} x^{b_{L}}(x-1)^{-r}(u_0u_1\ldots u_{n-1})^{-r}\nu_\el/ \delta_\el$, where
$$
\nu_\el=
\prod_{0\leq s<t<n}(u_s-u_t)^L\prod_{0 \leq s,t <n}\prod_{b_s \in B_{\el,L}^{(s)}}\prod_{1 \leq k \leq b_s} 
(x^ku_s-u_t)$$
and
$$\delta_\el=\prod_{0\leq s< t <n}\prod_{(b_s,b_t) \in B_{\el,L}^{(s)}\times B_{\el,L}^{(t)}}(x^{b_s}u_s-x^{b_t}u_t) \prod_{0 \leq s <n} \prod_{1 \leq i < j \leq L}(x^{b_i^{(s)}}u_s-x^{b_j^{(s)}}u_s).
$$
\end{theorem}
Following \cite{Brou}, Table
1, the field of definition of $G(n,1,r)$ is $K:=\mathbb{Q}(\zeta_{n})$. By Theorem $\ref{Semisimplicity Malle}$, if we set
$$\begin{array}{ll}
    X^{|\mu(K)|}:=x, &  \\
    U_k^{|\mu(K)|}:=(\zeta_n)^{-k}u_k & (k=0,1,\ldots,n-1),
  \end{array}$$
then the algebra $K(X,U_0,U_1,\ldots,U_{n-1})\mathcal{H}_{n,r}$ is split semisimple. We easily deduce that the factorization of the Schur elements
of this algebra is as described by Theorem $\ref{Schur element
generic}$.

\subsection{The groups $G(de,e,2)$, $e$ odd}

Lemma $\ref{gdeer}$ holds when $r=2$ and $e$ is odd.

\subsection{The groups $G(de,e,2)$, $e$ even}

Suppose that $e=2f$ for some $f\geq1$. Proposition 1.6 of \cite{Ariki} yields  the specialization of the parameters of 
the generic Hecke algebra $\mathcal{H}(G(2fd,2,2))$,
$(x_0,x_1;y_0,y_1;z_0,z_1,\ldots,z_{fd-1})$, which gives the generic
Hecke algebra of the group $G(2fd,2f,2)$.

\begin{lemma}\label{2fd}
The algebra $\mathcal{H}(G(2fd,2,2))$ specialized via
 $$\left\{ 
\begin{array}{ll} 
x_i \mapsto x_i &(0 \leq i \leq 1),\\
y_j \mapsto y_j  &(0 \leq j \leq 1),\\
z_k \mapsto \zeta_{f}^{[k/d]} u_{k \,\mathrm{mod}\,d}^{1/f} &(0 \leq k \leq fd-1)
\end{array} \right. 
$$
 is the twisted symmetric algebra of the
 cyclic group $C_f$ over the symmetric subalgebra $\mathcal{H}(G(2fd,2f,2))$ with
 parameters $(x_0,x_1;y_0,y_1;u_0,u_1,\ldots,u_{d-1})$. 
\end{lemma}
\begin{apod}{\small We have $$\begin{array}{rccl}
    \mathcal{H}(G(2fd,2,2)) & = & <S,T,U \,\,| &  STU=TUS=UST,  \\
     &  &  & (S-x_0)(S-x_1)=0, \\
     &  &  & (T-y_0)(T-y_1)=0, \\
     &  &  & (U-z_0)(U-z_1)\ldots(U-z_{fd-1})=0>.
  \end{array}$$
  Let
$$\begin{array}{rccl}
     A & := & <S,T,U \,\,| & STU=TUS=UST,  \\
     &  &  & (S-x_0)(S-x_1)=0, \\
     &  &  & (T-y_0)(T-y_1)=0, \\
     &  &  & (U^f-u_0)(U^f-u_1)\ldots(U^f-u_{d-1})=0>
  \end{array}$$
   and $$\bar{A}:=<S,T,U^f>.$$ Then
   $$A=\bigoplus_{i=0}^{f-1} U^i \bar{A}=\bigoplus_{i=0}^{f-1} \bar{A}U^i\,  \textrm{ and }\, \bar{A} \cong \mathcal{H}(G(2fd,2f,2)).$$}
\end{apod}

Set $n:=fd=de/2$. The group $G(2n,2,2)$ has $4n$ irreducible characters of degree $1$, 
$$\chi_{ijk} \,\,(0 \leq i,j\leq 1)\,(0 \leq k < n),$$
and $n^2-n$ irreducible characters of degree $2$,
$$\chi_{kl}^{1},\, \chi_{kl}^{2}\,\,(0 \leq k < l < n).$$

Following \cite{Ma2}, Theorem 3.11, the Schur elements of the irreducible characters of $\mathcal{H}(G(2n,2,2))$ are:
\\
\\
{\footnotesize $s_{\chi_{ijk}}=\Phi_1(x_ix_{1-i}^{-1}) \cdot \Phi_1(y_jy_{1-j}^{-1})\cdot \prod_{l=0,\,l\neq k}^{n-1} (\Phi_1(z_kz_l^{-1})
\cdot \Phi_1(x_ix_{1-i}^{-1} y_jy_{1-j}^{-1}z_kz_l^{-1}))$
\\
\\
$s_{\chi_{kl}^{1,2}}=-2 \prod_{m=0,\,m\neq k,l}^{n-1}( \Phi_1(z_kz_m^{-1})\cdot \Phi_1(z_lz_m^{-1})) \cdot\allowbreak
\prod_{i=0}^1 (\Phi_1(X_iX_{1-i}^{-1}Y_iY_{1-i}^{-1}Z_kZ_l^{-1})\cdot
\Phi_1(X_iX_{1-i}^{-1}Y_{1-i}Y_{i}^{-1}Z_lZ_k^{-1}))$\\
where $X_i^2:=x_i$, $Y_j^2:=y_j$, $Z_k^2:=z_k$.} 
\\
\\
Following \cite{Brou}, Table
1, the field of definition of $G(2n,2,2)$ is $K:=\mathbb{Q}(\zeta_{2n})$. By Theorem $\ref{Semisimplicity Malle}$, if we set
$$\begin{array}{cl}
    \mathcal{X}_i^{|\mu(K)|}:=(\zeta_2)^{-i}x_i & (i=0,1), \\
    \mathcal{Y}_j^{|\mu(K)|}:=(\zeta_2)^{-j}y_j & (j=0,1), \\
    \mathcal{Z}_k^{|\mu(K)|}:=(\zeta_n)^{-k}z_k & (k=0,1,\ldots,n-1),
  \end{array}$$
then the algebra $K(\mathcal{X}_0,\mathcal{X}_1,\mathcal{Y}_0,\mathcal{Y}_1,\mathcal{Z}_0,\mathcal{Z}_1,\ldots,\mathcal{Z}_{n-1})\mathcal{H}(G(2n,2,2))$ is split semisimple. Hence, the factorization of the Schur elements
of this algebra is as described by Theorem $\ref{Schur element
generic}$.

\printindex

\end{document}